\documentclass[10pt, reqno]{amsart}

\usepackage{amssymb}
\usepackage{mathrsfs,bbold}
\usepackage{amsmath, amsthm}
\usepackage{enumerate}
\usepackage{dsfont}
\usepackage{color}
\usepackage[a4paper]{geometry}
\usepackage[utf8]{inputenc}   
\usepackage{stmaryrd}
\usepackage{titlesec}
\usepackage{mathabx}
\usepackage{appendix}
\usepackage{todonotes, fancyhdr}
\usepackage{chngcntr}
\usepackage{cancel}

\usepackage{tikz-cd}
\usetikzlibrary{calc}

\usepackage{setspace}
\renewcommand{\baselinestretch}{0.99}

\usepackage[all]{xy}
\numberwithin{subsection}{section}
\numberwithin{subsubsection}{subsection}
\numberwithin{equation}{section} 

\newenvironment{Dem}[1][\unskip]{%
    \begin{list}{\hspace{1.15cm}{\sf \textbf{{\small Proof #1 --}}}}{%
        \setlength{\topsep}{0pt}%
        \setlength{\leftmargin}{0pt}%
        \setlength{\rightmargin}{0pt}%
        \setlength{\listparindent}{0pt}%
        \setlength{\itemindent}{0pt}%
        \setlength{\parsep}{0pt}%
        \addtolength{\leftmargin}{0pt}%
        \addtolength{\rightmargin}{0pt}%
    } \item }{\hfill $\rhd$\end{list}\smallskip}

\newenvironment{Dem*}[1][\unskip]{%
    \begin{list}{\hspace{0cm}{\sf \textbf{{\small Proof #1 --}}}}{%
        \setlength{\topsep}{0pt}%
        \setlength{\leftmargin}{0pt}%
        \setlength{\rightmargin}{0pt}%
        \setlength{\listparindent}{0pt}%
        \setlength{\itemindent}{0pt}%
        \setlength{\parsep}{0pt}%
        \addtolength{\leftmargin}{20pt}%
        \addtolength{\rightmargin}{0pt}%
    } \item }{\hfill $\rhd$\end{list}\smallskip}

\newenvironment{Rem}[1][\unskip]{%
    \begin{list}{\hspace{1.15cm}{\textsf{\textbf{{\small \textsl{Remark}} 
--}}}}{%
        \setlength{\topsep}{0pt}%
        \setlength{\leftmargin}{0pt}%
        \setlength{\rightmargin}{0pt}%
        \setlength{\listparindent}{0pt}%
        \setlength{\itemindent}{0pt}%
        \setlength{\parsep}{0pt}%
        \addtolength{\leftmargin}{0pt}%
        \addtolength{\rightmargin}{0pt}%
    } \item }{\hfill \end{list}\smallskip}

\newenvironment{Examples}[1][\unskip]{%
    \begin{list}{\hspace{1.15cm}{\textsf{\textbf{{\small \textsl{Examples}} }}}}{%
        \setlength{\topsep}{0pt}%
        \setlength{\leftmargin}{0pt}%
        \setlength{\rightmargin}{0pt}%
        \setlength{\listparindent}{0pt}%
        \setlength{\itemindent}{0pt}%
        \setlength{\parsep}{0pt}%
        \addtolength{\leftmargin}{20pt}%
        \addtolength{\rightmargin}{0pt}%
    } \item }{\hfill \end{list}\smallskip}

\renewcommand\thesection       {\arabic{section}}
\renewcommand\thesubsection    {\thesection{\boldmath $.$}\arabic{subsection}}
\renewcommand\thesubsubsection    {\thesection{\boldmath $.$}\arabic{subsection}{\boldmath $.$}\arabic{subsubsection}}

\titleformat{\section}[block]
{\filcenter\normalfont\sffamily\bfseries\Large}
{{\hspace{-0.7cm}}\thesection \hspace{0.2em} --\vspace{0.3cm}}{0.5em}{}

\titleformat{\subsection}[block]
{\filcenter\normalfont\sffamily\bfseries\large}  						  
{\hspace{-0.7cm}\thesubsection \hspace{0.5em} \vspace{0.3cm}}{.5em}{}  
\titlespacing{\subsection}{-0pc}{1.5ex plus .1ex minus .2ex}{0pc}

\titleformat{\subsubsection}[block]
{\normalfont\sffamily\bfseries}					  
{\thesubsubsection \vspace{0.3cm}}{.5em}{}  
\titlespacing{\subsection}{-0pc}{1.5ex plus .1ex minus .2ex}{0pc}

\newtheoremstyle{mystyle}
{3pt}               
{3pt}               
{\it }                      
{}                      
{\sffamily\bfseries}             
{}                      
{0.5em}                 
{#1 #2{\Large$.$}  }

\theoremstyle{mystyle}

\newtheorem{thm}{Theorem}
\newtheorem*{thm*}{Theorem}

\newtheorem{cor}[thm]{\hspace{-0.15cm}  {Corollary} }
\newtheorem{lem}[thm]{\hspace{-0.14cm}  {Lemma} }
\newtheorem{prop}[thm]{\hspace{-0.13cm} {Proposition}}
\newtheorem{defn}[thm]{ \hspace{-0.31cm} {Definition}}
\newtheorem*{rem*}{\hspace{-0.15cm} {Remark}}

\newtheoremstyle{mystyle2}
{3pt}               
{3pt}               
{\it }                      
{}                      
{\sffamily\bfseries}             
{}                      
{0.5em}                 
{\llap{#2 }#1{\hspace{0.2cm}--}}

\theoremstyle{mystyle2}

\newtheorem*{definition*}{Definition}
\newtheorem*{theorem*}{Theorem}
\newtheorem*{Remark*}{Remark}
\newtheorem*{lem*} {Lemma}
\newtheorem*{defn*} {Definition}
\newtheorem*{prop*} {Proposition}
\newtheorem*{cor*} {Corollary}

\newtheoremstyle{mystyle3}
{3pt}               
{3pt}               
{\it }                      
{}                      
{\sffamily\bfseries}             
{}                      
{0.5em}                 
{#1#2){\hspace{0.2cm}--}}

\theoremstyle{mystyle3}

\newtheorem{assumA}{Assumption (A}
\newtheorem{assumB}{Assumption (B}
\newtheorem{assumC}{Assumption (C}
\newtheorem{assumD}{Assumption (D}

\newcommand{\refA}[1]{\textbf{\textsf{(A\ref{#1})}}}
\newcommand{\REFA}{\textbf{\textsf{(A)}}}
\newcommand{\refB}[1]{\textbf{\textsf{(B\ref{#1})}}}
\newcommand{\REFB}{\textbf{\textsf{(B)}}}
\newcommand{\refC}[1]{\textbf{\textsf{(C\ref{#1})}}}

\newcommand{\refD}[1]{\textbf{\textsf{(D\ref{#1})}}}
\newcommand{\REFD}{\textbf{\textsf{(D)}}}

\usepackage{tikz-cd}
\usetikzlibrary{calc}
\usetikzlibrary{shapes.misc}
\usetikzlibrary{shapes.symbols}
\usetikzlibrary{snakes}
\usetikzlibrary{decorations}
\usetikzlibrary{decorations.markings}

\makeatletter
\pgfdeclareshape{crosscircle}
{
  \inheritsavedanchors[from=circle] 
  \inheritanchorborder[from=circle]
  \inheritanchor[from=circle]{north}
  \inheritanchor[from=circle]{north west}
  \inheritanchor[from=circle]{north east}
  \inheritanchor[from=circle]{center}
  \inheritanchor[from=circle]{west}
  \inheritanchor[from=circle]{east}
  \inheritanchor[from=circle]{mid}
  \inheritanchor[from=circle]{mid west}
  \inheritanchor[from=circle]{mid east}
  \inheritanchor[from=circle]{base}
  \inheritanchor[from=circle]{base west}
  \inheritanchor[from=circle]{base east}
  \inheritanchor[from=circle]{south}
  \inheritanchor[from=circle]{south west}
  \inheritanchor[from=circle]{south east}
  \inheritbackgroundpath[from=circle]
  \foregroundpath{
    \centerpoint%
    \pgf@xc=\pgf@x%
    \pgf@yc=\pgf@y%
    \pgfutil@tempdima=\radius%
    \pgfmathsetlength{\pgf@xb}{\pgfkeysvalueof{/pgf/outer xsep}}%
    \pgfmathsetlength{\pgf@yb}{\pgfkeysvalueof{/pgf/outer ysep}}%
    \ifdim\pgf@xb<\pgf@yb%
      \advance\pgfutil@tempdima by-\pgf@yb%
    \else%
      \advance\pgfutil@tempdima by-\pgf@xb%
    \fi%
    \pgfpathmoveto{\pgfpointadd{\pgfqpoint{\pgf@xc}{\pgf@yc}}{\pgfqpoint{-0.707107\pgfutil@tempdima}{-0.707107\pgfutil@tempdima}}}
    \pgfpathlineto{\pgfpointadd{\pgfqpoint{\pgf@xc}{\pgf@yc}}{\pgfqpoint{0.707107\pgfutil@tempdima}{0.707107\pgfutil@tempdima}}}
    \pgfpathmoveto{\pgfpointadd{\pgfqpoint{\pgf@xc}{\pgf@yc}}{\pgfqpoint{-0.707107\pgfutil@tempdima}{0.707107\pgfutil@tempdima}}}
    \pgfpathlineto{\pgfpointadd{\pgfqpoint{\pgf@xc}{\pgf@yc}}{\pgfqpoint{0.707107\pgfutil@tempdima}{-0.707107\pgfutil@tempdima}}}
  }
}
\makeatother

\colorlet{symbols}{black}    
\colorlet{testcolor}{green!60!black}
\colorlet{supcolor}{red!60!black}

\tikzset{
	root/.style={circle, fill=testcolor!70, draw=testcolor, inner sep=1pt, minimum size=0.5mm},
	dot/.style={circle, draw=black, fill=black, inner sep=0pt, minimum size=0.2mm},
	noise/.style={circle, draw=black, fill=white, inner sep=0pt, minimum size=1.4mm},
	noiseblue/.style={circle, fill=blue!40, draw=blue, inner sep=0pt, minimum size=1mm},
	noisegray/.style={circle, fill=gray!40, draw=gray, inner sep=0pt, minimum size=1mm},
	blackdot/.style={circle, draw=black, fill=black, inner sep=0pt, minimum size=1.2mm},
	K/.style= {semithick, shorten >=0pt,shorten <=0pt,-},
	DK/.style={thick, densely dotted, shorten >=0pt,shorten <=0pt},   
	}

\makeatletter
\DeclareRobustCommand\bigop[1]{%
  \mathop{\vphantom{\sum}\mathpalette\bigop@{#1}}\slimits@
}
\newcommand{\bigop@}[2]{%
  \vcenter{%
    \sbox\z@{$#1\sum$}%
    \hbox{\resizebox{\ifx#1\displaystyle.9\fi\dimexpr\ht\z@+\dp\z@}{!}{$\m@th#2$}}%
  }%
}
\makeatother          
             
\newcommand{\Bigstar}{\DOTSB\bigop{\bigstar}}

\definecolor{orange}{rgb}{1.0, 0.5, 0.0}

\newcommand{\iden}{\text{\rm Id}}

\newcommand{\unit}{\mathbf{1}}

\newcommand{\admi}{\mathcal{J}}
\newcommand{\remind}{\mathcal{N}}

\newcommand{\copro}{\Delta^+}
\newcommand{\comul}{\Delta}

\newcommand{\trino}[1]{|\!|\!|#1|\!|\!|}
\newcommand{\brarb}[1]{[\hspace{-0.5mm}]#1[\hspace{-0.5mm}]}

\newcommand{\bsf}{\boldsymbol{f}}
\newcommand{\bsg}{\boldsymbol{g}}
\newcommand{\bsh}{\boldsymbol{h}}
\newcommand{\mcF}{\boldsymbol{p}}
\newcommand{\bsu}{\boldsymbol{u}}
\newcommand{\bsv}{\boldsymbol{v}}
\newcommand{\bsw}{\boldsymbol{w}}

\newcommand{\bstau}{\boldsymbol{\tau}}
\newcommand{\bssigma}{\boldsymbol{\sigma}}

\newcommand{\bseta}{\boldsymbol{\eta}}
\newcommand{\bsvarphi}{\boldsymbol{\varphi}}
\newcommand{\mcFsi}{\boldsymbol{\psi}}

\newcommand{\bbF}{\mathbb{F}}
\newcommand{\bbD}{\mathbb{D}}

\newcommand{\sfT}{\mathsf{T}}
\newcommand{\sfU}{\mathsf{U}}
\newcommand{\sfD}{\mathsf{D}}

\newcommand{\wangle}[1]{\langle\!\langle#1\rangle\!\rangle}

\newcommand{\graftingatv}{\overset{\sf e}{\curvearrowright}_{(v)}}
\newcommand{\grafting}{\overset{\sf e}{\curvearrowright}}
\newcommand{\flatgrafting}{\overset{\sf e}{\curvearrowright}_{\flat}}
\newcommand{\rootgrafting}{\overset{\sf e}{\curvearrowright}_1}
\newcommand{\nonrootgrafting}{\overset{\sf e}{\curvearrowright}_2}

\newcommand{\mfs}{\mathfrak{s}}

\newcommand{\ssk}{\smallskip}

\renewcommand{\epsilon}{\varepsilon}

\newcommand\bbC{\mathbb{C}}
\newcommand\bbE{\mathbb{E}}
\newcommand{\bbG}{\mathbb{G}}

\newcommand\bbI{\mathbb{I}}
\newcommand\bbN{\mathbb{N}}
\newcommand\bbR{\mathbb{R}}
\newcommand{\bbT}{\mathbb{T}}
\newcommand{\bbX}{\mathbb{X}}
\newcommand\bbZ{\mathbb{Z}}

\newcommand{\mcB}{\mathcal{B}} 
\newcommand{\mcC}{\mathcal{C}} 
\newcommand{\mcD}{\mathcal{D}}

\renewcommand{\mcF}{\mathcal{F}}
\newcommand{\mcG}{\mathcal{G}}

\newcommand{\mcI}{\mathcal{I}}
\newcommand{\mcJ}{\mathcal{J}}
\newcommand{\mcK}{\mathcal{K}}
\newcommand{\mcL}{\mathcal{L}}
\newcommand{\mcM}{\mathcal{M}}
\newcommand{\mcN}{\mathcal{N}}

\newcommand{\mcP}{\mathcal{P}}
\newcommand{\mcQ}{\mathcal{Q}}
\newcommand\mcS{\mathcal{S}}

\newcommand\mcU{\mathcal U}

\newcommand\mcV{\mathcal V}

\newcommand{\bfA}{\mathbf{A}}

\newcommand{\bfG}{\mathbf{G}}

\newcommand{\bfK}{\mathbf{K}}

\newcommand{\bfT}{\mathbf{T}}

\newcommand{\Xplus}{X_+}

\makeatletter
\newcommand*{\defeq}{\mathrel{\rlap{%
                     \raisebox{0.3ex}{$\m@th\cdot$}}%
                     \raisebox{-0.3ex}{$\m@th\cdot$}}%
                     =}
\makeatother

\makeatletter
\newcommand*{\eqdef}{=\mathrel{\rlap{%
                     \raisebox{0.3ex}{$\m@th\cdot$}}%
                     \raisebox{-0.3ex}{$\m@th\cdot$}}%
                     }
\makeatother

\begin{document}

\begin{center}
{\Huge\sffamily{A tourist's guide to regularity structures   \vspace{0.15cm}   \\
and singular stochastic PDEs   \vspace{0.5cm}}}
\end{center}

\begin{center}
{\sf I. BAILLEUL}\footnote{I.B. was partially supported by the ANR via the ANR-16-CE40-0020-01 grant.} \& {\sf M. HOSHINO}\footnote{M. H. was partially supported by the JSPS KAKENHI Grant Number 19K14556.}
\end{center}

\vspace{1cm}

\begin{center}
\begin{minipage}{0.8\textwidth}
\renewcommand\baselinestretch{0.7} \scriptsize \textbf{\textsf{\noindent Abstract.}} {\sf We give an essentially self-contained treatment of the fundamental analytic and algebraic features of regularity structures and its applications to the study of singular stochastic PDEs.   }
\end{minipage}
\end{center}

\vspace{0.6cm}

{\sf 
\begin{center}
\begin{minipage}[t]{11cm}
\baselineskip =0.35cm
{\scriptsize 

\center{\textbf{Contents}}

\vspace{0.1cm}

\textbf{ 1.~Introduction\dotfill 
\pageref{SectionIntro}}

\textbf{ 2.~Basics on regularity structures\dotfill 
\pageref{SectionBasicsRS}}

\textbf{ 3.~Regularity structures built from integration operators\dotfill
\pageref{SectionIntegration}}

\textbf{ 4.~Solving singular PDEs within regularity structures\dotfill 
\pageref{SectionSolvingPDEs}}

\textbf{ 5.~Renormalization structures\dotfill 
\pageref{SectionConcreteRenormStructure}}

\textbf{ 6.~Multi-pre-Lie structure and renormalized equations\dotfill 
\pageref{SectionMultiAndrenormalizedEquations}}

\textbf{ 7.~The BHZ character\dotfill 
\pageref{SectionBHZCharacter}}

\textbf{ 8.~The manifold of solutions\dotfill 
\pageref{SectionManifold}}

\textbf{ 9.~Building regularity and renormalization structures\dotfill 
\pageref{SectionBuildingRS}}

\textbf{ A.~Summary of notations\dotfill
\pageref{SectionAppendixSummary}}

\textbf{ B.~Basics from algebra\dotfill
\pageref{SectionAppendixAlgebra}}

\textbf{ C.~Technical proofs\dotfill
\pageref{SectionAppendixProofs}}

\textbf{ D.~Comments\dotfill
\pageref{SectionAppendixComments}}

}\end{minipage}
\end{center}
}   \vspace{1cm}

\vspace{0.6cm}

\section{Introduction}
\label{SectionIntro}

The class of singular stochastic partial differential equations (PDEs) is 
characterized by the appearance in their formulation of ill-defined products due to the presence in the equation of distributions with low regularity, typically realizations of random distributions. Here are three typical examples.
\begin{itemize}
   \item[-] The $2$ or $3$-dimensional parabolic Anderson model equation (PAM)
   \begin{equation} \label{EqPAM}
   (\partial_t - \Delta_x)u = u\xi,
   \end{equation}
   with $\xi$ a space white noise. It represents the evolution of a Brownian particle in a $2$ or $3$-dimensional white noise environment in the torus. (The operator $\Delta_x$ stands here for the $2$ or $3$-dimensional 
Laplacian.)  \vspace{0.1cm}
      
   \item[-] The scalar $\Phi^4_3$ equation from quantum field theory
   \begin{equation} \label{EqPhi43}
   (\partial_t - \Delta_x)u = -u^3 + \zeta,
   \end{equation}
   with $\zeta$ a $3$-dimensional spacetime white noise and $\Delta_x$ the $3$-dimensional Laplacian in the torus or the Euclidean space. Its invariant measure is the scalar $\Phi^4_3$ measure from quantum field theory. 
  \vspace{0.1cm}
   
   \item[-] The generalized (KPZ) equation
   \begin{equation} \label{EqGKPZ}
   (\partial_t-\partial_x^2)u = f(u)\zeta + g(u)\vert\partial_x u\vert^2,
   \end{equation}
   with $\zeta$ a $1$-dimensional spacetime white noise. In a more sophisticated form, it provides amongst others a description of the random motion of a rubber on a Riemannian manifold under a random perturbation of the mean curvature flow motion.    \vspace{0.1cm}
\end{itemize}

A $d$-dimensional space white noise has H\"older regularity $-d/2-\kappa$, and a $d$-dimensional spacetime white noise has H\"older regularity $-d/2-1-\kappa$ under the parabolic scaling, almost surely for every positive $\kappa$. Whereas one expects from the heat operator that its inverse regularizes a distribution by $2$, this is not sufficient to make sense of 
any of the products $u\xi, u^3, f(u)\zeta, \vert\partial_x u\vert^2, g(u)\vert\partial_x u\vert^2$ above, as the product of two H\"older distributions is well-defined if and only if the sum of their regularity exponents 
is positive. Why then bother about such equations? It happens that they appear as scaling limits of a number of microscopic nonlinear random dynamics where the strength of the nonlinearity and the randomness balance each other. Many microscopic random systems exhibit this feature as you will 
see from reading Corwin \& Shen's nice review \cite{CorwinShen} on singular stochastic PDEs.

\ssk

A typical statement in the theory of regularity structures \cite{Hai14, BHZ, ChandraHairer, BCCH18} about a singular stochastic PDE takes the following informal form, stated here in restricted generality. Consider a subcritical singular stochastic PDE
\begin{equation} \label{EqGenericPDE}
(\partial_t-\Delta_x)u = f(u,\partial u)\zeta + g(u,\partial_x u) \eqdef F(u,\partial u; \zeta)
\end{equation}
driven by a possibly multi-dimensional irregular random noise $\zeta$ that is almost surely of spacetime regularity $\alpha-2$, for a deterministic constant $\alpha\in\bbR$. (The notion of `subcriticality' will be properly defined later in the text.) We talk of the sufficiently regular function $F$ as a `nonlinearity' -- even though a particular $F$ could depend linearly or affinely of one or all of its arguments. Denote by $\frak{F}$ 
the space of nonlinearities that are affine functions of the noise argument. For each $\varepsilon\in(0,1]$, denote also by $\zeta_\epsilon$ a regularized version of the noise which converges to $\zeta$ as $\varepsilon$ goes to $0$, obtained for instance by convolution $\zeta_\epsilon\defeq \zeta*\varrho_\varepsilon$ with a deterministic smooth mollifier $\varrho_\varepsilon$. Write 
$$
u_\varepsilon=\textsc{Sol}(\zeta_\epsilon ; F)
$$ 
for the solution to the well-posed parabolic equation 
$$
(\partial_t-\Delta_x)u_\varepsilon = F(u_\varepsilon,\partial_x u_\varepsilon; \zeta_\epsilon)
$$
started at time $0$ from a given (regular enough) fixed initial condition.

\medskip

\noindent \textbf{\textsf{Meta-theorem 1. (What it means to be a solution) --}} \emph{The following three points hold true.}\vspace{0.1cm}

\begin{itemize}\setlength{\itemsep}{0.2cm}
   \item  \emph{One can associate to each subcritical singular stochastic 
PDE a finite dimensional unbounded Lie group called the renormalization group. Denote by $k$ its generic elements.}   
   
   \item \emph{This group acts explicitly on the right on the nonlinearity space $\frak{F}$}
\begin{equation} \label{EqGroupAction}
\big(k,F\big)\mapsto F^{(k)}  \in \frak{F}.
\end{equation}

   \item \emph{There exists some deterministic (typically diverging) elements $(k_\epsilon)_{0<\epsilon\leq 1}$ of the renormalization group such that, for any element $k$ of the renormalization group, the solutions 
$$
\overline{u}_\epsilon^{(k)} \defeq \textsc{Sol}\Big(\zeta_\epsilon; (F^{(k_\epsilon)})^{(k)}\Big)
$$ 
to the well-posed stochastic PDE
$$
(\partial_t-\Delta_x)\overline{u}_\epsilon^{(k)} = \big(F^{(k_\epsilon)}\big)^{(k)} \big(\overline{u}_\epsilon^{(k)},\partial_x \overline{u}_\epsilon^{(k)}; \zeta_\epsilon\big)
$$ 
with given initial condition, converge in probability in an appropriate function/distribution space to some $\overline{u}^{(k)}$ as $\epsilon$ goes to $0$.}
\end{itemize} 
\vspace{0.1cm}
\emph{A solution to a singular stochastic PDE is not a single function or distribution, but rather the \emph{family} $(\overline{u}^{(k)})_k$ of functions/distributions indexed by the renormalization group.  } 

\medskip

The $k_\epsilon$ typically diverge as $\epsilon>0$ goes to $0$ but there are situations where they remain bounded, or are even constant, like in \cite{BrunedGabrielHairerZambotti}. We stick to the tradition and talk about any of the above limit functions/distributions as a solution to equation \eqref{EqGenericPDE}. We talk of 
the family of solutions. To have a picture in mind, consider the family of maps 
\begin{equation} \label{EqDefnSEpsilon}
S_\epsilon(x) = \big(x-1/\epsilon\big)^2,
\end{equation}
on $\bbR$. It explodes in every fixed interval as $\epsilon$ goes to $0$, 
but remains finite, and converges, in a moving window $S_\epsilon(x+1/\epsilon)$, where it is equal to $x^2$. It also converges in the other moving window $S_\epsilon(x+1/\epsilon+1)$, where it is equal to $(x+1)^2$. No 
given moving window is a priori better than another. In this parallel, the function $\textsc{Sol}(\zeta_\epsilon ; \cdot)$ plays the role of $S_\epsilon$, with the infinite dimensional nonlinearity space $\frak{F}$ in the role of the state space $\bbR$. The role of the translations $x\mapsto 
x+1/\epsilon$, is played by the group action \eqref{EqGroupAction} of $k_\epsilon$ on the space of nonlinearities $\frak{F}$. The explicit action of $k_\epsilon$ on the space of nonlinearities gives some formulas of the form 
$$
F^{(k_\epsilon)}(u,\partial_x u ; \zeta) = f(u,\partial u)\zeta + g(u,\partial_x u) + C_\epsilon(u,\partial u),
$$
for some functions $C_\epsilon$ built from $f,g$ and their derivatives, and from $k_\varepsilon$, so $\overline{u}_\epsilon\defeq\textsc{Sol}\big(\zeta_\epsilon; F^{(k_\epsilon)}\big)$ is the solution to the equation 
\begin{equation} \label{EqrenormalizedPDE}
(\partial_t-\Delta_x)\overline{u}_\epsilon = f(\overline{u}_\epsilon,\partial_x \overline{u}_\epsilon)\zeta_\epsilon + g(\overline{u}_\epsilon,\partial_x \overline{u}_\epsilon) + C_\epsilon(\overline{u}_\epsilon,\partial_x \overline{u}_\epsilon),
\end{equation}
with given initial condition. We talk of the function $C_\epsilon$ as a counterterm.

\ssk

In a robust solution theory for differential equations a solution to a differential equation ends up being a continuous function of the parameters in the equation. In the case of equation \eqref{EqGenericPDE} the parameters are the functions $f,g$, the noise $\zeta$ and the initial condition of the equation. While it is unreasonable and wrong to expect that the solutions from Meta-theorem 1 are continuous functions of each realization of the noise, they happen to be continuous functions of a measurable functional of the latter build by probabilistic means. We talk about that functional of the noise as an {\it enhanced noise}.   

\ssk

\noindent \textbf{\textsf{Meta-theorem 2 (Continuity of a solution with respect to the enhanced noise) --}} \emph{For any subcritical singular stochastic PDE \eqref{EqGenericPDE}, and for any $\zeta$ in a class of random noises including space or spacetime white noises, there is a measurable functional $\sf \Pi$ of the noise taking values in a metric space and such that any individual solution of equation \eqref{EqGenericPDE} is a continuous function of $\sf \Pi$.   }

\ssk

This is a fundamental point to be compared with the fact that in the stochastic calculus approach to stochastic (partial) differential equations, the solutions to the equations are only measurable functionals of the noise. In a setting where both approaches can be used and coincide the regularity structures point of view provides a factorization of the measurable solution map of stochastic calculus under the form of the composition of a measurable function $\sf \Pi$ of the noise with a continuous function of $\sf \Pi$. A number of probabilistic statements about $\sf \Pi$ are then automatically transferred to the solution of the equation by continuity. The support of the law of the random variable $\sf \Pi$ determines for instance the support of the law of the solution to the equation. A large deviation result for the laws of a family of random $\sf \Pi$ is also automatically transported by continuity into a large deviation result for the laws of the corresponding family of solutions of the equation. 

The functional $\sf \Pi$ is built as a limit {\it in probability} of some elementary functionals of the noise. This is in the end the reason why the convergence result in Meta-theorem 1 holds in probability. (The convergence is almost sure along an appropriate sequence.)

\ssk

How is it that one can prove such statements? The starting point is that solutions of singular stochastic PDEs are not expected to be any kind of H\"older functions or distributions. Rather, under an assumption on the equation captured by the notion of subcriticality, we expect any possible solution
\begin{equation} \label{EqDescriptionLocal}
u(\cdot) \simeq \sum_\tau u_\tau(x) ({\sf \Pi}_x\tau)(\cdot)
\end{equation}
to be described locally in terms of a {\it finite} number of equation-dependent reference functions or distributions ${\sf \Pi}_x\tau$ that are polynomial functionals of the noise. One of these symbols is denoted by $\bf 1$ and ${\sf \Pi}_x{\bf 1}$ is the constant function equal to $1$. Unlike $u$, which may be a distribution in some situations, the $u_\tau$ are always some functions. Making the parallel with a classical situation, one could talk of the family $(u_\tau)_\tau$ as a {\it jet} for $u$. A fundamental result in regularity structures gives some mild ($({\sf \Pi}_x\tau)_{x,\tau}$-dependent) coherence conditions on the $(u_\tau)_\tau$ under which a quantified version of the relation \eqref{EqDescriptionLocal} determines a unique function/distribution $u$ satisfying it. By trading $u$ for $(u_\tau)_\tau$ the theory of regularity structures then provides a complete description of the local structure of the possible solutions to a given singular stochastic PDE, in terms of their local expansion coefficients with respect to some equation-dependent polynomial functionals ${\sf \Pi}_x\tau$ of the noise. The theory actually turns the problem upside down by reformulating any singular stochastic PDEs as an equation with unknown some tuple $(u_\tau)_\tau$ of local coefficients satisfying a priori the above mentioned coherence condition. Since \eqref{EqDescriptionLocal} gives a local description of $u$ near any state space point $x$, and $u\simeq u_{\bf 1}(x)$ at first order near $x$, an equation of the form
$$
\mathcal{L}u = F(u\,;\,\zeta)
$$
will be rewritten near $x$ as

\begin{equation} \label{EqLocalRewriting}
\sum_\tau u_\tau(x) \, \mcL({\sf \Pi}_x\tau) \simeq \sum_k\frac{(\partial_u^kF)(u_{\bf 1}(x))}{k!} \, \sum_{\tau_1,\dots,\tau_k\neq {\bf 1}}u_{\tau_1}(x)\cdots u_{\tau_k}(x)\,({\sf \Pi}_x\tau_1)\cdots ({\sf \Pi}_x\tau_k)\zeta.
\end{equation}
Identifying the terms on both sides provides a triangular system for the ${\sf \Pi}_x\tau$ and a family of equations for the $u_\tau(x)$. In the example of the parabolic Anderson model equation one has for instance $F(u;\zeta)=u\zeta$, the expansion is indexed by some trees $\tau\in\{\begin{tikzpicture}[scale=0.3,baseline=0.08cm] \node at (0,0)  [dot] (1) {}; \node at (0,0.8)  [noise] (2) {}; \draw[K] (1) to (2); \end{tikzpicture}, 
 \begin{tikzpicture}[scale=0.3,baseline=0.1cm] \node at (0,0)  [dot] (1) {}; \node at (0,0.6) [noise] (2) {}; \node at (0,1.4) [noise] (3) {}; \draw[K] (1) to (2); \draw[K] (2) to (3); \end{tikzpicture}, \dots\}$, and the system reads
\begin{equation*} \begin{split}
\mcL({\sf \Pi}_x\begin{tikzpicture}[scale=0.3,baseline=0.08cm] \node at (0,0)  [dot] (1) {}; \node at (0,0.8)  [noise] (2) {}; \draw[K] (1) to (2); \end{tikzpicture}) &= \zeta,   \\
\mcL({\sf \Pi}_x \begin{tikzpicture}[scale=0.3,baseline=0.1cm] \node at (0,0)  [dot] (1) {}; \node at (0,0.6) [noise] (2) {}; \node at (0,1.4) [noise] (3) {}; \draw[K] (1) to (2); \draw[K] (2) to (3); \end{tikzpicture}) &= ({\sf \Pi}_x\begin{tikzpicture}[scale=0.3,baseline=0.08cm] \node at (0,0)  [dot] (1) {}; \node at (0,0.8)  [noise] (2) {}; \draw[K] (1) to (2); \end{tikzpicture}) \zeta, \quad \textrm{etc.}   \\
u_\tau(x) &= u_{\bf 1}(x), \; \textrm{for all }\tau,
\end{split} \end{equation*} 
since the nonlinearity is linear here. The base point $x$ in ${\sf \Pi}_x\tau$ means that rather than simply defining ${\sf \Pi}_x\tau$ as $\mcL^{-1}$ of the right hand side of its defining equation, we only keep the Taylor remainder at $x$ of this quantity, at some $\tau$-dependent order. Of course the singular feature of the equation has not disappeared as the products $({\sf \Pi}_x\tau_1)\cdots ({\sf \Pi}_x\tau_k)\zeta$ in \eqref{EqLocalRewriting} are still problematic. However, proceeding this way we have isolated the product problem in the problem of making sense of the reference functions/distributions ${\sf \Pi}_x\tau$, a task which has little to do with the actual task of solving the equation. This is the very point where the fact that the noise is random plays a crucial role. It allows indeed to build the reference functions/distributions $({\sf \Pi}_x\tau)(\omega)$ not as some functions of a realization $\zeta(\omega)$ of the noise but rather as  some random variables jointly defined with the noise on a common probability space. This realizes a wonderful decoupling of probability and analysis. {\it To the former} the task of building the enhanced noise: Functions or distribution-valued random variables ${\sf \Pi}_x\tau$ that involve the noise only. {\it To the latter} the task of solving uniquely an equation in  the side space  of local coefficients built from the enhanced noise, regardless of any multiplication problem. The construction of the ${\sf \Pi}_x\tau$ is done by a limiting procedure called \textit{renormalization}, after similar procedures used in quantum field theory to tackle similar problems.

\medskip

The fundamentals of the theory of regularity structures were built gradually by M. Hairer and his co-authors in four groundbreaking works \cite{Hai14, BHZ, ChandraHairer, BCCH18}. In paper \cite{Hai14} M. Hairer sets the analytic framework of regularity structures and provides an ad hoc study of the renormalization problem for the parabolic Anderson model equation \eqref{EqPAM} and scalar $\Phi^4_3$ equation \eqref{EqPhi43}. The algebra involved in the renormalization process of a large class of singular stochastic PDEs was unveiled in Bruned, Hairer and Zambotti's work \cite{BHZ}. The proof that the renormalization algorithm provided in \cite{BHZ} converges was given by Chandra \& Hairer in \cite{ChandraHairer}. Last, the fact that the renormalization can be `implemented' at the level of the equation was proved in Bruned, Chandra, Chevyrev and Hairer's work \cite{BCCH18}, giving a wonderful analogue of the equivalence of the ``subtraction scheme'' versus ``counterterms'' approaches to renormalization problems in quantum field theory. Altogether these four works provide a black box for the local well-posedness theory of subcritical singular stochastic PDEs. This work gives an essentially self-contained short treatment of the fundamental analytic and algebraic features of regularity structures and its applications to the study of singular stochastic PDEs that contains the essential points of the works \cite{Hai14, BHZ, BCCH18}. It is intended for readers who already have an idea of the subject and who wish to understand in depth the mechanics at work. We hope nonetheless that even a newcomer to the field may grasp the matter by following the road taken here. Regularity structures and the fundamental tools are developed in generality within a highly abstract setting. No trees are in particular involved in the analysis before we actually construct an example of regularity structure adapted to the study of the generalized (KPZ) equation \eqref{EqGKPZ} in Section \ref{SectionBuildingRS}. When it comes to applying these tools to singular stochastic PDEs we trade generality for the concrete example of the generalized (KPZ) equation \eqref{EqGKPZ}, that involves all the difficulties of the most general case. We do not treat Chandra \& Hairer's work \cite{ChandraHairer} constructing the measurable functional $\sf \Pi$ of the noise involved in Meta-theorem 2 using Bruned, Hairer and Zambotti's renormalization process \cite{BHZ}.

\medskip

We stress here that Hairer's approach to singular stochastic PDEs is somewhat orthogonal to the purely probabilistic approaches of stochastic PDEs pioneered by Pardoux, Walsh or da Prato \& Zabczyk using martingale technics, described in \cite{DPZ} for instance. No knowledge of these approaches is needed to understand what follows.

\medskip

It is our aim here to give a concise self-contained version of what seems to us to be the most important features of the 433(=236+118+79) pages of the works \cite{Hai14, BHZ, BCCH18}. A number of comments about different statements, concepts, other works, are deferred to Appendix {\sf \ref{SectionAppendixComments}} so as to keep focused in the main body of the text. The reader is invited to read this section at any point along her/his reading. We expect that the reader will see from the present work the simplicity that governs the architecture of the theory. The climb may be hard but the view after the walk is stunning. 

The probabilistic/renormalization side of the analysis of singular stochastic PDEs was not mature yet when this work was first written. This is why the deep work \cite{ChandraHairer} of Chandra \& Hairer was left aside in this tourist's guide. The situation has changed after the seminal work \cite{LOTT} of Linares, Otto, Tempelmayr \& Tsatsoulis and some subsequent works by Hairer \& Steele \cite{HS} and Bailleul \& Hoshino \cite{BH23}. We refer the reader to \cite{BHReview} for a review of the subject and further references.

\ssk

Besides the original articles \cite{Hai14, BHZ, ChandraHairer, BCCH18}, Hairer's lectures notes \cite{HairerBrazil, HairerTakagi}, the book \cite{FrizHairerBook} by Friz and Hairer, Chandra \& Weber's article \cite{ChandraWeber} and Berglund's book \cite{Berglund}, provide other accessible accounts of part of the material presented here. The work \cite{CorwinShen} of Corwin and Shen provides a nice non-technical overview of the context in which singular stochastic PDEs arise.

\medskip

We will introduce the different pieces of the puzzle one after the other to arrive at a clear understanding of the mathematical form of the above meta-theorems. This will be done along the following lines. 

\ssk

\noindent {\it 1. Concrete regularity structures, models and modelled distributions.} We will first set the scene to talk of the local behavior of functions/distributions
\begin{equation} \label{EqLocalPicture}
f(\cdot) \sim \sum_\tau f_\tau(x) ({\sf \Pi}^{\sf g}_x\tau)(\cdot),
\end{equation}
near each spacetime point $x$, giving a generalization of the notion of jet, in terms of reference functions/distributions $({\sf \Pi}^{\sf g}_x\tau)(\cdot)$ parametrized by $\tau$ in a finite set and by the spacetime points $x$. This will involve the setting of {\it (concrete) regularity structures} $\mathscr{T}=\big((T^+,\Delta^+),(T,\Delta)\big)$, and {\it models} $\sf M=(\Pi, g)$, from which the reference functions/distributions $({\sf \Pi}^{\sf g}_x\tau)(\cdot)$ are built. See Sections \ref{SectionConcreteRS} and \ref{SubsectionModelsAndCo} for a detailed descriptions of these objects. The notation $({\sf \Pi}^{\sf g}_x\tau)(\cdot)$ will be defined in Definition \ref{DefnModel}. In the same way as a family of functions $\{f_k\}_{k\in\bbN^d}$ on $\bbR^d$ needs to satisfy a quantitative consistency condition for a function $f$ satisfying
$$
f(\cdot) \sim \sum_{k\in\bbN^d} f_k(x) (\cdot-x)^k, \qquad \textrm{near all $x$,}
$$
to exist, a collection of functions $(f_\tau)_\tau$ needs to satisfy a quantitative consistency condition for a distribution $f$ satisfying \eqref{EqLocalPicture} to exist. This condition will involve the notion of {\it modelled distribution} and {\it reconstruction operator} $\textbf{\textsf{R}}^{\sf M}$, with a notion of consistency that will depend on the component $\sf g$ of the model $\sf M$. One can think of a modelled distribution as the consistent jet of an function or a distribution. At that stage, given a regularity structure and a model on it, we will have a convenient way of representing a class of functions/distributions on the state space -- not all of them. Given a (system of) singular stochastic PDE(s), a good choice of concrete regularity structure will allow to represent the set of functions/distributions that appear in a naive analysis of the equation via a Picard iteration by some elements of our class. Unlike what happens in the study of controlled ordinary differential equations driven by an $\ell$-dimension control, there is no universal concrete regularity structure for the set of all singular stochastic PDEs. One associates to each (system of) subcritical singular stochastic PDE(s) a specific regularity structure. 

\ssk

\noindent {\it 2. Lifting the equation as an equation in the space of consistent jets.} The regularity structure associated with equation \eqref{EqGenericPDE} is built from a noise symbol and some operators $(\mcI_n)_{n\in\bbN^{d+1}}$ that play the role in the regularity structure of the operators $\partial^n(\partial_t-\Delta_x)^{-1}$, involved in the Picard fixed point formulation of the equation. (The letter $\partial=(\partial_t,\partial_x)$ stands here for the time/space derivative operator.) One proceeds then by formulating the equation as a fixed point problem in the space of consistent jets of functions/distributions $(f_\tau)_\tau$, encoded in the notion of modelled distribution. This will require to introduce a tweaked version $\mathcal{K}^{\sf M}$ of the operator $\mcI$, as the latter does not produce consistent jets from consistent jets -- the notion of consistency depends on $\sf g$ whereas the operator $\mcI$ does not. The equation on the jet space will happen then to have a unique solution in small time under some mild conditions. This solution will be a continuous function of all the parameters in the equations, the model in particular. Along the way we will turn the initial analytical multiplication problem into the problem of defining some models enjoying some appropriate properties -- the so called {\it admissible models}. We will see that it is straightforward to construct what is called the canonical lift of a regularized version $\zeta_\epsilon$ of the noise $\zeta$ as an admissible model ${\sf M}^{\epsilon}$; this can be done for any smooth noise. In those terms, the solution $u_\epsilon$ to a well-posed (system of) singular stochastic PDE(s) driven by some smooth noise $\zeta_\epsilon$ can be written as the reconstruction
\begin{equation} \label{EqSolProjection}
u_\epsilon = \textsf{\textbf{R}}^{{\sf M}^\epsilon}(\bsu_\epsilon)
\end{equation}
of a consistent jet $\bsu_\epsilon$ obtained as the fixed point 
\begin{equation} \label{EqFixedPoint}
\bsu_\epsilon = \Phi\big({\sf M}^\epsilon, \bsu_\epsilon\big)
\end{equation}
of a map $\Phi$ that depends continuously on its model argument. So does the reconstruction map $\textsf{\textbf{R}}^{{\sf M}^\epsilon}$. 

\ssk

\noindent {\it 3. Renormalized models and renormalized equation.} However the model ${\sf M}^\epsilon$ does not converge in the appropriate space as the positive regularization parameter $\epsilon$ goes to $0$, so a solution to the (system of) singular stochastic PDE(s) under study cannot be defined as the limit of the $u_\epsilon$ as $\epsilon$ goes to $0$. The situation is similar to what happens to the function $S_\epsilon$ from \eqref{EqDefnSEpsilon}. One has to look at ${\sf M}^{\epsilon}$ in a moving window to obtain a finite limit. The {\it renormalization group} will provide us precisely with this possibility, and will provide us in particular with a family of renormalized canonical models $({}^{k_\epsilon}{\sf M}^\epsilon)_{0<\epsilon\leq 1}$. To make the final step from here to the meta-theorems, we will see that this action of the renormalization group on the set of models has a dual action on the space $\frak{F}$ of nonlinearities. The $^{k_\epsilon}{\sf M}^\epsilon$-reconstruction $\overline{u}_\epsilon$ of the unique $^{k_\epsilon}{\sf M}^\epsilon$-dependent fixed point equation in the space of jets will happen to solve a `renormalized' version of the singular equation \eqref{EqGenericPDE}, with some additional $\epsilon$-dependent terms diverging most of the time as the regularization parameter $\epsilon$ tends to $0$, as in \eqref{EqrenormalizedPDE}. The continuity of both the solution of the fixed point equation \eqref{EqFixedPoint} and the reconstruction map, as functions of the underlying model, will ensure the convergence of $\overline{u}_\epsilon$ to some limit $\overline{u}$ for some converging renormalized models $^{k_\epsilon}{\sf M}^\epsilon$ with limit $\overline{\sf M}$, say. The limit function $\overline{u}$ will satisfy a system
$$
\overline{u} = \textsf{\textbf{R}}^{\overline{\sf M}}(\overline\bsu),\quad \overline\bsu = \Phi\big(\overline{\sf M}, \overline\bsu\big),
$$
similar to the system of equations \eqref{EqSolProjection} and \eqref{EqFixedPoint} satisfied by $u_\epsilon$. It is in this sense that $\overline{u}$ will deserve to be called a solution of the singular stochastic PDE under study. Think of $\overline{u}$ as a function/distribution defined from its `Taylor' jet $\overline\bsu$, with the latter solution of a fixed point problem.

\medskip

\noindent \textbf{\textsf{A word about algebra.}} It is one of the features of the theory of regularity structures that algebra plays an important role, unlike what one usually encounters in the analytic study of PDEs. This is partly due to the choice of description of the objects involved in the analysis, in terms of ``jets-like'' quantities. Elementary consistency requirements directly bring algebra into play, under the form of Hopf algebras and actions of the latter on some vector spaces. This is what concrete regularity structures are. The appearance of algebra in the study of singular stochastic PDEs is also due to the fact that the renormalization algorithm used to define the random variables that play the role of a number of ill-defined polynomial functionals of the noise is conveniently encoded in an algebraic structure that we call \textit{renormalization structure}; it differs from a concrete regularity structure. These two points involve Hopf algebras. A last piece of algebra is also needed under the form of some pre-Lie algebras. This is an algebraic structure that behaves as the differentiation operation $(f,g)\mapsto g'f$, derivative of $g$ in the direction of $f$. Using an algebraic language sheds a gentle light on the meaning of the renormalization process at the level of the equation. Pre-Lie algebras are the ingredient that we use to understand how to build the counterterm $C_\epsilon$ in Meta-theorem {\sf 1} and equation \eqref{EqrenormalizedPDE}. 

The analysis or probability-oriented reader should not be frightened by the perspective of working with some algebraic tools; we will hardly need anything more than a few definitions and elementary facts that are direct consequences of the latter; everything else is proved. We refer the reader to Manchon's lecture notes \cite{Manchon}, or the first four chapters of Sweedler's book \cite{Sweedler}, for some accessible references on Hopf algebras, and to Foissy's work \cite{FoissyTypedTrees} for the basics on pre-Lie algebras; all we need is elementary and recalled below. Appendix {\sf \ref{SectionAppendixAlgebra}} contains in any case all the results from algebra that we use without proving them, with precise pointers to the literature.

\medskip

\textit{\textbf{Organisation of the article.}} This Tourist's Guide has been organized as follows. Basics on regularity structures are introduced in Section \ref{SectionBasicsRS}, under the form of concrete regularity structures. The reconstruction theorem, that ensures that a consistent jet describes a distribution in the state space is proved there, in Theorem \ref{ThmReconstructionRS}. This allows to formulate in Section \ref{SectionSolvingPDEs} a singular PDE as an equation in a space of modelled distributions over a regularity structure associated with the singular PDE. A fixed point argument is used in Section \ref{SectionSolvingPDEs} to prove a local in time well-posedness result in a space of modelled distributions. Despite their possible differences, the regularity structures built for the study of different subcritical elliptic or parabolic 
singular PDEs all involve the construction of the counterpart of a (or several) regularizing convolution operator(s) and the proof of its (/their) 
continuity properties in spaces of modelled distributions. This is done in Section \ref{SectionIntegration}. Section \ref{SectionConcreteRenormStructure} sets the scene of renormalization structures. They encode the renormalization algorithm used to build the random variables whose realizations play the role of a finite number of reference functions/distributions. The renormalization algorithm is described in Section \ref{SectionBHZCharacter}. The dual action of the renormalization operation on the genuine singular PDE is clarified by the introduction of some pre-Lie structures; this is done in Section \ref{SectionMultiAndrenormalizedEquations}. Nothing so far requires a deep understanding of how one builds the regularity or the renormalization structure associated with a given singular stochastic PDE. It suffices to assume that they satisfy a small number of simple assumptions to run the analysis. A summary of the assumptions can be found at the begining of Section \ref{SectionBuildingRS}. Section \ref{SectionBuildingRS} is dedicated to constructing explicitly such structures in the example of the generalized (KPZ) equation. A summary of our notations is given in Appendix {\sf \ref{SectionAppendixSummary}}. Appendix {\sf \ref{SectionAppendixAlgebra}} contains a number of elementary facts from algebra that we use without proof in the text. Precise pointers to the proofs of these facts are given. Appendix {\sf \ref{SectionAppendixProofs}} contains the proof of technical results that were not given in the body of the text to keep concentrated on the essential features of the method. A number of comments about the notions, the statements or the literature are collected in Appendix {\sf \ref{SectionAppendixComments}}. The reader is invited to read them at any time.

\medskip

\noindent \textbf{\textsf{Notations --}} {\it We use a number of greek letters with different meanings. As a rule, $\alpha, \beta,\gamma$ stand for real numbers, while $\tau, \sigma, \mu, \nu, \eta, \varphi,\psi$ stand for elements of regularity or renormalization structures. 
\begin{itemize}
   \item Given two statements $\frak{a}$ and $\frak{a}^+$, we agree to write $\frak{a}^{(+)}$ to mean both the statement $\frak{a}$ and the statement $\frak{a}^+$.
   \item Denote by $e_i=(0,\dots,0,\overset{i}{1},0,\dots,0)\in\bbR^d$ the $i$-th basis vector of $\bbR^d$.
   \item Denote by $\bbN$ the set of non-negative integers. For each $k=(k_i)_{i=1}^d\in\bbN^d$ and $x=(x_i)_{i=1}^d\in\bbR^d$, we use the notations
   $$
   k!\defeq k_1!\cdots k_d!,\qquad
   x^k\defeq x_1^{k_1}\cdots x_d^{k_d}.
   $$
   For any $k=(k_i)_{i=1}^d,\ell=(\ell_i)_{i=1}^d$, we write $\ell\le k$ if $\ell_i\le k_i$ for any $i$, and then define
   $$
   \binom{k}{\ell}\defeq\frac{k!}{\ell!(k-\ell)!}=\binom{k_1}{\ell_1}\cdots\binom{k_d}{\ell_d}.
   $$
\end{itemize}
All the notations introduced along the way are gathered in Appendix {\sf\ref{SectionAppendixSummary}}, with pointers to the section where they are introduced.}

\medskip

\noindent \textbf{\textsf{Assumptions --}} {\it We emphasize along the way a number of `assumptions' on some regularity structures; they are summarized at the begining of Section \ref{SectionBuildingRS}. All these properties are satisfied by the regularity structures used to study some singular stochastic PDEs; so, strictly speaking, they are not assumptions. We first present the theory of regularity structures independently of its applications to singular stochastic PDEs and then gradually introduce some more specific features of the particular (tree-indexed) structures that are used in that setting. We find it convenient to state in the form of some `assumption' some special features that a regularity structure can have to emphasize their role in the proofs of some particular results.}

\vfill \pagebreak

\section{Basics on regularity structures}
\label{SectionBasicsRS}

Hairer's theory of regularity structures builds on Gubinelli's approach \cite{GubinelliControlled} to T. Lyons' theory of rough paths and rough differential equations \cite{Lyons98}. This is a theory of controlled ordinary differential equations 
$$
dz_t = V(z_t) d{\bf X}_t,\quad z_t\in\bbR^k,
$$
driven by irregular controls $\bf X$. Gubinelli's notion of `path controlled by a rough path $\bf X$' gives a Taylor-like description of a path around each time $s$ in terms of some `monomials' given by the different components of the increments of the rough path $\bf X$ between the running time $t$ and the fixed time $s$. This notion of controlled path turns out to be stable by (regular enough) nonlinear maps and by the operator $d^{-1}$ defining the integral against the reference rough path $\bf X$. These facts allow to formulate controlled ordinary differential equations driven by a rough path as an integral equation in a space of controlled paths and to prove local well-posedness of the equation under some mild regularity assumptions on the vector fields involved in the dynamics by some fixed point arguments. There is no need to know anything about rough path in this work.) Hairer chooses a similar angle to build his theory of singular stochastic PDEs, with an important add on. In the context of the above controlled differential equation, the theory of regularity structures provides a setting in which one has not only a pointwise description of a potential solution path $z$ but also a `local' description of the $\bbR^k$-valued distribution $V(z_t)d{\bf X}_t$ on some time interval around every fixed time $s$. This way, all the terms that are involved in the analysis of the controlled differential equation are described in terms of some local expansion.

We start this section by explaining in simple terms why the strategy of local expansion devices automatically brings algebra into play, independently of any problem of dynamical nature. This is the content of Section \ref{SectionLocalExpansionDevices}. The backbone of these expansion devices is encoded in the notion of concrete regularity structure introduced in Section \ref{SectionConcreteRS}. The reference functions/distributions ${\sf \Pi}_x\tau$ that we informally used in the local expansions \eqref{EqDescriptionLocal} and \eqref{EqLocalDescription} play the role classically devoted to the Taylor monomials $(\cdot-x)^n$. The ${\sf \Pi}_x\tau$ are built from two primary objects $\sf \Pi$ and $\sf g$ that will play in the sequel a crucial role. They define jointly what is called a model. In this setting, the quantification of the approximation \eqref{EqDescriptionLocal} leads to the notion of modelled distribution. This condition on the collection $(u_\tau(x))_{x,\tau}$ of `coefficients' turns out to garantee the existence of a function/distribution which is indeed well approximated by the finite sum $\sum_\tau u_\tau(x) {\sf \Pi}_x\tau$ near any state space point $x$; a condition for uniqueness of the approximated function/distribution is also known. This fact is the content of the reconstruction theorem. Models, modelled distributions and their reconstruction are presented in Section \ref{SubsectionModelsAndCo}. Section \ref{section products and derivatives} describes how one can make sense of some nonlinear operations on modelled distributions and gives the properties of a derivative operator acting on the space of modelled distributions. 

No PDE-related matters are involved in this section. Its main purpose is to give a general analysis of some arbitrary expansion devices. It is only in Section \ref{SectionIntegration} that we will introduce some particular features of such devices that are involved in the application of this general machinery to the study of singular stochastic PDEs. Together with the introduction of the main objects of regularity structures the main result of this section is the reconstruction theorem, Theorem \ref{ThmReconstructionRS}.

\bigskip

\subsection{Algebra as the mechanics of local expansion devices}
\label{SectionLocalExpansionDevices}

Regularity structures are the backbone of expansion devices for the local description of functions and distributions in (an open set of) a Euclidean space, say $\bbR^d$. (The isotropic nature of the Euclidean space plays no role here, and as a matter of fact everything works in a non-isotropic setting. We stick to a Euclidean setting here for simplicity.)

\ssk

{\it 1. The Taylor expansion device.} The usual notion of local description of a function near a point $x\in\bbR^d$ involves the Taylor expansion operation and amounts to comparing a function to a polynomial centered at $x$
\begin{equation}   \label{EqLocalDescription}
f(\cdot) \simeq \sum_n f_n(x)\,(\cdot-x)^n, \quad \textrm{near }x.
\end{equation}
The sum over $n$ is finite, the approximation quantified and we end up describing the class of H\"older functions with real positive regularity exponents. By the binomial expansion, one gets a local description of $f$ near any other point $y$ writing
\begin{equation}   \label{EqDescriptionf}   \begin{split}
f(\cdot) \simeq \sum_{\ell\leq n} f_n(x)\,{n\choose{\ell}} (\cdot-y)^\ell (y-x)^{n-\ell} \simeq \sum_\ell \left(\sum_{n;\ell\leq n} f_n(x)\,{n\choose{\ell}} (y-x)^{n-\ell}\right)(\cdot-y)^\ell.
\end{split}\end{equation}
This expansion brings the important insight that the coefficient $f_\ell(y)$ should be compared with the polynomial centered at $x$
\begin{equation}   \label{EqLocalDescription2}
f_\ell(y) \simeq \sum_{n\geq\ell} f_n(x)\,{n\choose{\ell}}(y-x)^{n-\ell}, \quad \textrm{near }x
\end{equation}
where the coefficients $f_n(x)$ are associated with some multiindices $n\ge\ell$. This relation between the coefficients in the expansion at different points is a {\it consistency condition}. Conversely, this consistency condition for the coefficients $(f_n)_{n\in\bbN^d}$ is required to ensure the existence of a function $f$ satisfying \eqref{EqLocalDescription}.

\ssk

{\it 2. The reference objects for general expansion devices.} A more general local description device involves an $\bbR^d$-indexed collection of functions or distributions $(\Pi_x\tau)(\cdot)$, with labels $\tau$ in a given finite set $\mcB$. We will consider the functions/distributions locally described as
\begin{equation}   \label{EqLocalDescription3}
f(\cdot) \simeq \sum_\tau f_\tau(x)(\Pi_x\tau)(\cdot), \quad\textrm{near each }x\in\bbR^d,
\end{equation} 
for some coefficients $f_\tau(x)$. The above expression implicitly assumes that the coefficients $f_\tau(x)$ are function of $x$. One has $\mcB=\bbN^d$ and $(\Pi_xk)(\cdot) = (\cdot - x)^k$, in the usual, Taylor, polynomial setting. So, what is an alternative of the binomial expansion? Like in the former setting, it seems meaningful to impose that $\Pi_x\tau$ is a linear combination of $\Pi_y\sigma$ with several labels $\sigma$ near ay nother point $y$. Denoting by $T$ the real vector space spanned by $\mcB$, we express the situation by the identity
\begin{equation}  \label{EqTransitionMapsRelation}
(\Pi_x\tau)(\cdot) = \big(\Pi_y(\Gamma_{yx}\tau)\big)(\cdot),
\end{equation}
with a linear map $\Gamma_{yx} : T\rightarrow T$. Moreover, since the roles of $x$ and $y$ are exchangeable, the linear maps $\Gamma_{yx}$ are invertible, with $\Gamma_{yx}^{-1}=\Gamma_{xy}$, and one has a group action of an $\bbR^d\times\bbR^d$-indexed group on the local description structure $T$.   

\ssk

While one uses the same polynomial-type local description \eqref{EqLocalDescription2} for the $f_n$ as we do for $f$ in \eqref{EqLocalDescription}, in the usual H\"older setting, there is no reason in a more general local description device to use the same reference objects for $f$ and for its local coefficients. This is in particular the case if the $(\Pi_x\tau)(\cdot)$ are meant to describe distributions, among others, while it makes sense to some use functions only as reference objects to describe the functions $f_\tau$. For this reason, we introduce another set $\mcB^+$ and associate to each $\mu\in\mcB^+$ with a reference function ${\sf g}_{yx}(\mu)$ playing the role of $(y-x)^{n-\ell}$ in \eqref{EqLocalDescription2}. More precisely, we assume that each $f_\sigma$ can be compared with a finite linear combination
\begin{equation} \label{EqTransitionfTau}
f_\sigma(y) \simeq \sum_{\tau\in\mcB,\, \mu\in\mcB^+} c_\mu^{\tau,\sigma}f_\tau(x) {\sf g}_{yx}(\mu), \quad\textrm{near }x
\end{equation}
with some constant $c_\mu^{\tau,\sigma}$ depending on $\tau\in\mcB,\sigma\in\mcB^+,\mu\in\mcB^+$. In the Taylor polynomial setting, one has $\mcB^+=\mcB=\bbN^d$, ${\sf g}_{yx}(k)=(y-x)^k$ and $c_k^{n,\ell}={\bf1}_{k=n-\ell}{n\choose\ell}$. In this model situation, all the $\sigma$ appearing in the above formula have a `size' greater than or equal to the size of $\tau$. Such a hierarchy will be encoded later as a graded structure on $T$. Unlike the model Taylor polynomial setting, the elements of $\mcB^+$ are associated with some functions only, while the elements of $\mcB$ may be associated with some distributions. Note that we take care not to write $(\Pi_x\tau)(y)$ as $\Pi_x\tau$ may be a distribution. It would be consistent to write ${\sf g}_x(\mu)(y)$ but we stick to the established and convenient notation ${\sf g}_{yx}(\mu)$. 

\ssk

{\it 3. Consistency relations.} With the above notation, one has
$$
f(\cdot) \simeq \sum_{\sigma\in\mcB} f_\sigma(y)(\Pi_y\sigma)(\cdot) \simeq \sum_{\tau,\sigma\in\mcB,\,\mu\in\mcB^+} c_\mu^{\tau,\sigma}f_\tau(x) \, {\sf g}_{yx}(\mu)(\Pi_y\sigma)(\cdot).
$$
Comparing this expansion to the original expansion \eqref{EqLocalDescription3} at the point $x$, we would have the explicit representation
$$
\Gamma_{yx}\tau = \sum_{\sigma\in\mcB,\,\mu\in\mcB^+}c_\mu^{\tau,\sigma}{\sf g}_{yx}(\mu) \sigma.
$$
To lighten the notation it will be convenient to introduce a notation under the form of the formal ratio 
$$
\tau/\sigma=\sum_{\mu\in\mcB^+}c_\mu^{\tau,\sigma}\mu.
$$
It is an element of the real vector space $T^+$ spanned by $\mcB^+$. Then the above consistency formulas read
$$
f_\sigma(y) \simeq \sum_{\tau\in\mcB}f_\tau(x) \, {\sf g}_{yx}(\tau/\sigma),\qquad
\Gamma_{yx}\tau = \sum_{\sigma\in\mcB}{\sf g}_{yx}(\tau/\sigma)\sigma.
$$
We can also derive a transitive relation similar to \eqref{EqTransitionMapsRelation} for $\Gamma_{zy}$ and $\Gamma_{yx}$. Indeed we can develop the coefficient $f_\tau(x)$ in \eqref{EqTransitionfTau}, and re-indexing the labels, one gets
\begin{equation} \label{EqDeltaDeltaPlus} \begin{split}
f_\eta(z) &\simeq \sum_{\sigma\in\mcB} f_\sigma(y)\,{\sf g}_{zy}(\sigma/\eta)
\simeq \sum_{\sigma,\eta\in\mcB} f_\tau(x)\,{\sf g}_{zy}(\sigma/\eta) \, {\sf g}_{yx}(\tau/\sigma).
\end{split} \end{equation}
In the comparison with \eqref{EqTransitionfTau}, it would be natural to impose the identity
\begin{equation}  \label{EqTransitionMapsRelation2}
\sum_{\sigma\in\mcB} {\sf g}_{zy}(\sigma/\eta) \, {\sf g}_{yx}(\tau/\sigma) 
= 
{\sf g}_{zx}(\tau/\eta).
\end{equation}
Together with \eqref{EqTransitionMapsRelation}, this is a fundamental transition relation analogous to the binomial expansion.

\ssk

{\it 4. Algebra encodes consistency relations.} There is a way to describe the skeleton of the relations \eqref{EqTransitionMapsRelation} and \eqref{EqTransitionMapsRelation2} with no mention of any explicit reference functions or distributions. We introduce the splitting maps 
\begin{align*}
\Delta &: T\rightarrow T\otimes T^+, \quad \Delta\tau = \sum_{\sigma} \sigma\otimes (\tau/\sigma),\\
\Delta^+ &: T^+\rightarrow T^+\otimes T^+,\quad \Delta^+(\tau/\eta) = \sum_{\sigma} (\sigma/\eta)\otimes (\tau/\sigma).
\end{align*}
(Since there is no unique way to express an element of $T^+$ in the form $\tau/\sigma$, the `definition' of $\Delta^+(\tau/\eta)$ may seem problematic. However, this formula actually holds for the ``concrete regularity structure" defined below.) Using this notation the transition map reads
$$
\Gamma_{yx}=(\iden\otimes{\sf g}_{yx})\Delta,
$$
and \eqref{EqTransitionMapsRelation} and \eqref{EqTransitionMapsRelation2} take the form
$$
\Pi_x=(\Pi_y\otimes{\sf g}_{yx})\Delta,\qquad
{\sf g}_{zx}=({\sf g}_{zy}\otimes{\sf g}_{yx})\Delta^+.
$$
An important property of these maps is given by the following relations.
\begin{equation}\label{EqRightComodule}   
\begin{split}
(\Delta\otimes \textrm{Id})\Delta &= (\textrm{Id}\otimes \Delta^+)\Delta,\\
(\textrm{Id}\otimes \Delta^+)\Delta^+ &= (\Delta^+\otimes \textrm{Id})\Delta^+.
\end{split}
\end{equation}
For instance, the transitive relation $\Gamma_{zx}=\Gamma_{zy}\Gamma_{yx}$ is encoded by the first relation as follows.
\begin{align*}
\Gamma_{zy}\Gamma_{yx}&=(\iden\otimes{\sf g}_{zy}\otimes{\sf g}_{yx})(\Delta\otimes\iden)\Delta\\
&=(\iden\otimes{\sf g}_{zy}\otimes{\sf g}_{yx})(\iden\otimes\Delta^+)\Delta=(\iden\otimes{\sf g}_{zx})\Delta
=\Gamma_{zx}.
\end{align*}   
Based on the model case of the Taylor expansion device, we ask the family of reference functions ${\sf g}_{yx}(\mu), \mu\in\mcB^+$ to be sufficiently rich to describe locally an \textit{algebra} of functions. The cheapest way to ensure that property is to assume that the linear span $T^+$ of $\mcB^+$ has an algebra structure and that the maps ${\sf g}_{yx}$ on $T^+$ are some characters of this algebra, that is they are multiplicative maps. We do not impose an algebra structure on $T$ as the elements of $\mcB$ are meant to be associated with some distributions. In the end, we find that it is natural to introduce a graded vector space $T$ and a graded algebra $T^+$ with some splitting maps $\Delta$ and $\Delta^{+}$. Using the assumed invertibility of the transition maps $\Gamma_{yx}$, an elementary fact from algebra then leads directly to the Hopf algebra structure that appears below in the definition of a concrete regularity structure. (The curious reader can see Proposition \ref{PropBialgebraHopf} in Appendix {\sf \ref{SectionAppendixAlgebra}}. We do not need to understand the details of its simple proof now.)

Note that the dimension $d$ of the state space, or the fact that it is Euclidean, play no role in this discussion.

\ssk

We choose to record the essential features of this discussion in the definition of a `concrete' regularity structure given below; this is a special form of the more general notion of regularity structure from Hairer' seminal work \cite{Hai14}. The reader should keep in mind that the entire algebraic setting can be understood at a basic level from the above consistency requirements on a given local description device. We invite the reader to look at Appendix {\sf \ref{SectionAppendixAlgebra}} and read the definitions and basic properties of bialgebras, Hopf algebras, and comodules. It is better to read this appendix in the light of the preceding discussion.

\bigskip

\subsection{Regularity structures}
\label{SectionConcreteRS}

We define in this section a particular form of regularity structure that turns out to be sufficient for the study of (systems of) singular stochastic PDE(s). A number of notations are fixed here. A general regularity structure is defined in \cite[Section 2]{Hai14} and recalled at the end of this section. The following particular form can be essentially found in \cite[Section 8]{Hai14} and \cite{BHZ}. The following definition is to be read in the light of the discussion of Section \ref{SectionLocalExpansionDevices} and will be best understood by recalling first from Appendix {\sf \ref{SectionAppendixAlgebra}} the definition of a Hopf algebra and the meaning of connectedness in this setting.

\medskip

\begin{defn*}
A \textsf{\textbf{concrete regularity structure}} $\mathscr{T}=(T^+,T)$ 
is the pair of graded vector spaces   \vspace{-0.1cm}
$$
T^+ = \bigoplus_{\alpha\in A^+} T_\alpha^+, \qquad T = \bigoplus_{\beta\in A} T_\beta   \vspace{-0.1cm}
$$
such that the following holds.   \vspace{0.1cm}

\begin{itemize}
\setlength{\itemsep}{0.1cm}
   \item The vector spaces $T_\alpha^+$ and $T_\beta$ are finite dimensional.   \vspace{0.1cm}
   
   \item The vector space $T^+$ is a connected graded bialgebra with unit 
${\bf1}_+$, counit ${\bf1}_+'$, coproduct $\Delta^+:T^+\to T^+\otimes T^+$, and grading $A^+\subset[0,\infty)$.   \vspace{0.1cm}

   \item The index set $A$ for $T$ is a locally finite subset of $\bbR$ bounded below. The vector space $T$ is a right comodule over $T^+$, that is $T$ is equipped with a splitting map $\Delta : T \rightarrow T\otimes T^+$ which satisfies 
\begin{align} \label{eq T is comodule}
(\Delta\otimes \textrm{\emph{Id}})\Delta = (\textrm{\emph{Id}}\otimes \Delta^+)\Delta,\quad\textrm{and}\quad (\iden\otimes{\bf1}_+')\Delta=\iden. 
\end{align}
Moreover, for any $\beta\in A$ 
\begin{align}\label{eq T is graded}
\Delta T_\beta \subset \bigoplus_{\alpha\geq 0}T_{\beta-\alpha}\otimes T_\alpha^+.
\end{align}
\end{itemize}
We denote by 
$$
\mathscr{T} \defeq  \big((T^+,\Delta^+), (T,\Delta)\big)
$$
a concrete regularity structure.
\end{defn*}

\medskip

Recall that $\bigoplus$ denotes an algebraic sum, so each element of $T^{(+)}$ is a linear combination of elements of $T_\alpha^{(+)}$ for finitely many $\alpha$, even if $A^{(+)}$ is an infinite set. The different elements that appear in the definition of a concrete regularity structure will acquire later a concrete meaning. The elements of $T$ and $T^+$ will index some expansion devices for the study of a given (system of) singular stochastic PDE(s) -- remember that each equation will have its own regularity structure. We saw in Section \ref{SectionLocalExpansionDevices} the meaning of the splitting maps $\Delta$ and $\Delta^+$ and their intertwining/coherence relations \eqref{EqRightComodule} in terms of some expansion device. Recall from the definition of the graded bialgebra given in Appendix {\sf \ref{SectionAppendixAlgebra}} that $T_0^+=\langle{\bf1}_+\rangle$ and $T_\alpha^+T_\beta^+\subset T_{\alpha+\beta}^+$, for any $\alpha,\beta\in A^+$. By Proposition \ref{PropBialgebraHopf}, the bialgebra $T^+$ is indeed a Hopf algebra; we denote by $S_+$ its antipode. Denoting by $\mcM : T^+\otimes T^+\rightarrow T^+$ the multiplication operator $\mcM(a\otimes b)\defeq ab$, and by ${\bf1}_+'$ the counit of $T^+$ -- think of it as a dual vector to the vector $\textbf{\textsf{1}}_+$, the antipode $S_+$ is characterized by the identity
$$
\mcM(\textrm{Id}\otimes S_+)\Delta^+\tau = \mcM(S_+\otimes\textrm{Id})\Delta^+\tau = {\bf1}_+'(\tau)\textbf{\textsf{1}}_+.
$$
Moreover, the coproduct $\Delta^+$ satisfies $\Delta^+{\bf1}_+={\bf1}_+\otimes{\bf1}_+$, and 
   \begin{equation}   \label{EqDefnDeltaPlus}
   \Delta^+\tau \in \left\{\tau\otimes {\bf 1}_+ + {\bf 1}_+\otimes \tau + \sum_{0<\alpha'<\alpha} T^+_{\alpha'}\otimes T^+_{\alpha-\alpha'}\right\},
   \end{equation}
for any $\tau\in T_\alpha^+$ with $\alpha>0$. Similarly, it is straightforward from \eqref{eq T is comodule} and \eqref{eq T is graded} to check that
\begin{equation}\label{EqDelta}  
\Delta \tau \in \left\{\tau\otimes {\bf 1}_+ + \sum_{\beta'<\beta} T_{\beta'}\otimes T^+_{\beta-\beta'}\right\}
\end{equation}
for any $\tau\in T_\beta$. This identity will later imply for a set $(\Pi_x\tau)_{x,\tau}$ of reference functions/distributions that $\Pi_{x'}\tau\simeq\Pi_x\tau$, up to some terms of smaller `homogeneity' $\beta'<\beta$.

For an arbitrary element $\tau$ in $T$, set
\begin{equation*}
\tau = \sum_{\beta\in A} \tau_\beta \in \bigoplus_{\beta\in A} T_\beta.
\end{equation*}
We use a similar notation for elements of $T^+$. An element $\tau$ of $T_\alpha^{(+)}$ is said to be \emph{homogeneous} and is assigned \textsf{\textbf{homogeneity}} $|\tau| \defeq  \alpha$. The homogeneous spaces $T_{\beta}$ and $T^{+}_\alpha$ being finite dimensional, all norms on them are equivalent; we use a generic notation $\|\cdot\|_\beta$ or $\|\cdot\|_\alpha$ for norms on these spaces. For simplicity, we write
\begin{align}
\label{EqNotationDirectSum}
\|\tau\|_\alpha \defeq \|\tau_\alpha\|_\alpha.
\end{align}
Note that we do not assume any relation between the linear spaces $T_\alpha^+$ and $T_\beta$ at that stage. Note also that the homogeneity function $\vert\cdot\vert$ takes values in $\bbR$, and that the parameter $\beta$ in \eqref{EqDelta} can be non-positive, unlike in \eqref{EqDefnDeltaPlus}.

\medskip

\noindent \textbf{\textsf{Notations.}} $\bullet$ {\it Let $\mcB_\alpha^+$ 
 and $\mcB_\beta$ be bases of $T_\alpha^+$ and $T_\beta$, respectively. Set 
$$
\mcB^+ \defeq  \bigcup_{\alpha\in A^+} \mcB_\alpha^+, \quad \mcB \defeq  \bigcup_{\beta\in A} \mcB_\beta.
$$   }\vspace{0.15cm}

$\bullet$ \textit{Recall from the end of Section \ref{SectionIntro} our convention about statements of the form $\frak{s}^{(+)}$. Given $\sigma, \tau \in \mcB^{(+)}$, we use the notation $\sigma\leq^{(+)}\tau$ to mean that 
$\sigma=\tau$ or $|\sigma|<|\tau|$; we write $\tau/^{(+)}\sigma$ for 
the element of $T^+$ defined by the expansion
\begin{equation}\label{EqCoprodQuotient}
\Delta^{(+)}\tau = \sum_{\sigma\in\mcB^{(+)},\, \sigma\le^{(+)}\tau } \sigma\otimes (\tau/^{(+)}\sigma).
\end{equation}
Write $\sigma<^{(+)}\tau$ to mean further that $\sigma$ is different from 
$\tau$. The notations $\tau/^{(+)}\sigma$ and $\sigma<^{(+)}\tau$ are only used for $\tau$ and $\sigma$ in $\mcB^{(+)}$.}

\bigskip

It should be noted that, for any regularity structure $\mathscr{T}=(T^+,T)$, the pair
$$
\mathscr{T}^+ \defeq  \big((T^+,\Delta^+),(T^+,\Delta^+)\big)
$$ 
also define a regularity structure. The polynomial regularity structure defined later takes such a form. Also, an algebraic structure arising from branched rough paths considered in \cite{GubinelliBranched} essentially takes the above form. In general, we consider two distinct spaces $T^+$ and $T$ to encode distributions by the vector space $T$ which is not defined as an algebra.

\medskip

Interpreting the splitting maps $\Delta$ and $\Delta^+$ as chopping elements into pieces, keep in mind that $\tau/^{(+)}\sigma$ can be a sum of elements of $T^+$, in case $\sigma$ appears `at different places' as a part of $\tau$. This will be particularly clear in formula \eqref{eq coproduct on polynomial} below for the polynomial regularity structure, where 
the binomial coefficient $\binom{n}{\ell}$ will account for the number of 
$X^\ell$ inside $X^n$, for $0\leq \ell\leq n$. Note that for $\sigma<^{(+)}\tau$ in $\mcB^{(+)}$, we have
\begin{equation}
\label{EqCoprodQuotient+}
\begin{split}
\Delta^+(\tau/^{(+)}\sigma) &= \sum_{\sigma \leq^{(+)} \eta \leq^{(+)} \tau}(\eta/^{(+)}\sigma)\otimes(\tau/^{(+)}\eta)   \\
                                     &= (\tau/^{(+)}\sigma)\otimes\unit_+ + \unit_+\otimes(\tau/^{(+)}\sigma) + \sum_{\sigma <^{(+)} \eta <^{(+)} \tau}(\eta/^{(+)}\sigma)\otimes(\tau/^{(+)}\eta).
\end{split}
\end{equation}
These two identities are direct consequences of the \textit{co-associativity properties} 
$$
(\Delta^{(+)}\otimes \textrm{Id})\Delta^{(+)} = ( \textrm{Id}\otimes \Delta^+)\Delta^{(+)},
$$
of the coproduct $\Delta^{(+)}$, obtained by identifying the corresponding terms in the left and right hand sides. In the setting of singular stochastic PDEs where the elements of $T$ are (decorated) trees, $\tau/\sigma$ will be a product of trees, and each of these trees will eventually be involved in the action of re-centering the corresponding analytic objects 
to a given running point, while leaving the trunk tree $\sigma$ untouched. The definition of a model given in Section \ref{SubsectionModelsAndCo} illustrates exactly this picture. 

\medskip

Here are two examples of regularity structures.
\begin{itemize}
\item[--] Let symbols $X_1,\dots, X_d$ be given. For $n\in\bbN^d$, set $X^n \defeq  X_1^{n_1}\cdots X_d^{n_d}$; this is an element of the free commutative algebra with unit ${\bf1}(\defeq X^0)$ generated by the $X_i$. We can see that $T_X\defeq \textrm{span}\big\{X^n;n\in\bbN^d\big\}$ is a bialgebra with the coproduct
\begin{align}\label{eq coproduct on polynomial}
\Delta_\textrm{pol} X^n \defeq  \sum_{\ell\leq n} \binom{n}{\ell} X^\ell\otimes X^{n-\ell}.   
\end{align}
Let $\frak{s}=(\frak{s}_i)\in(\bbN\backslash\{0\})^d$ be an integer-valued fixed vector, called a \textsf{\textbf{scaling}}. This vector accounts for the natural scaling properties in the different directions of $\bbR^d$ for the problem at hand. If for instance $\bbR^d$ does not stand for the isotropic 
Euclidean space but rather for a non-isotropic space with topology $\bbR^d$, as a Lie group, different directions will naturally have different homogeneities, depending on the geometry of the space. We define the scaled 
degree of $n\in\bbN^d$ by
$$
|n|_{\frak{s}}=\sum_{i=1}^d\frak{s}_in_i.
$$
Then the definition $T_\alpha=\textrm{span}\{X^n;\vert n\vert_{\frak{s}}=\alpha\}$ gives a grading for the bialgebra $T_X$. Since $T_0=\textrm{span}\{\bf1\}$, the space $T_X$ is a connected graded bialgebra. Thus it is indeed a Hopf algbera; the antipode is actually given by $S_+X^n=(-X)^n$. The \textsf{\textbf{polynomial regularity structure}} is given by
$$
\mathscr{T}_X\defeq \big((T_X,\Delta_\textrm{pol}),(T_X,\Delta_\textrm{pol})\big).
$$
   
   \item[--] To have another picture in mind, think of $T$ and $T^+$ as sets of possibly labelled rooted trees, with $T^+$ consisting only of trees with positive tree homogeneities -- a homogeneity is assigned to each labelled tree. This notion of homogeneity induces the decomposition \eqref{EqDelta} of $T$ into linear spaces spanned by trees with equal homogeneities; a similar decomposition holds for $T^+$. The coproduct $\Delta^{(+)}\tau$ is typically a sum over subtrees $\sigma$ of $\tau$ with the same root as $\tau$, and $\tau/\sigma$ is the quotient tree obtained from $\tau$ by identifying $\sigma$ with the root; this quotient tree is better seen as a product of trees. See Section \ref{SectionBuildingRS} for constructions of regularity structures of this sort associated with singular stochastic PDEs. For such regularity structures, the minimum regularity of the elements of $T$ is given by the minimum regularity of the noises in the equation. One can leave aside trees by the time we arrive at Section \ref{SectionBuildingRS} and work in the abstract setting of this section throughout.
\end{itemize}

\medskip

A group of linear transforms on the space $T^+$ has an important role in \cite{Hai14}. 

\medskip

\begin{defn*}
A \textbf{\textsf{character}} $g$ on the Hopf algebra $T^+$ is a linear map $g:T^+\to\bbR$, such that $g({\bf1}_+)=1$ and $g(\tau_1\tau_2)=g(\tau_1)g(\tau_2)$, for any $\tau_1,\tau_2\in T^+$. 
The set $G^+$ of all characters of the algebra $T^+$ turns into a group with the convolution product $*$ defined by 
$$
(g_1*g_2)\tau \defeq  (g_1\otimes g_2)\Delta^+\tau, \quad \tau\in T^+,
$$
where we identify the tensor product of two real numbers with their product.
The unit of $G^+$ is the counit ${\bf1}_+'$ of $T^+$, and the inverse of $g\in G^+$ is given by $g^{-1}\defeq g\circ S_+$, where $S_+$ is the antipode of $T^+$.
The group $G^+$ is called \textbf{\textsf{structure group}}.
\end{defn*}

\medskip

Think of the usual convolution product $(g_1*g_2)(x) = \int g_1(y)g_2(x-y)dy$, where one first splits $x$ into $y$ and $x-y$, then apply $f$ and $g$ to each piece, before taking the product and summing over all possible splittings. The group $G^+$ acts on $T$ from left. One associates to a character $g$ of $T^+$ the map 
$$
\widehat g \defeq  (\textrm{Id}\otimes g)\Delta : T\to T,
$$ 
from $T$ to itself.  We have indeed
\begin{equation}\label{Section2:GactionAssociative}
\widehat{g_1*g_2} = \widehat{g_1} \circ \widehat{g_2}
\end{equation}
for any $g_1,g_2\in G^+$, as a direct consequence of the comodule property \eqref{eq T is comodule}. Indeed,
\begin{align*}
\widehat{g_1*g_2} &= (\iden\otimes(g_1*g_2))\Delta
=(\iden\otimes g_1\otimes g_2)(\iden\otimes\Delta^+)\Delta\\
&=(\iden\otimes g_1\otimes g_2)(\Delta\otimes\iden)\Delta=
(\widehat{g_1}\otimes g_2)\Delta\\
&=\widehat{g_1} \circ \widehat{g_2}.
\end{align*}
Also, for any $\tau\in T_\beta$,
$$
\Big(\widehat{g}(\tau) - \tau\Big) \in \bigoplus_{\beta'<\beta}T_{\beta'},
$$
as a consequence of the structural identity \eqref{EqDelta}. Similarly, one defines the action of $G^+$ on $T^+$ by
$$
\widehat{g}^+\defeq (\textrm{Id}\otimes g)\Delta^+:T^+\to T^+.
$$
for $g\in G^+$. This operator also satisfies the properties similar to $\widehat{g}$.

\medskip

In the end of this section, we recall from \cite{Hai14} the original definition of regularity structures.

\medskip

\begin{defn*}
A \textbf{\textsf{regularity structure}} $\mathscr{T} = \big( \bfA,\bfT,\bfG \big)$ consists of the following.
\begin{itemize} \setlength{\itemsep}{0.1cm}
\item[(1)]
$\bfA$: a subset of $\bbR$ such that the set $\{\alpha\in \bfA\, ;\, \alpha<\gamma\}$ is finite for every $\gamma\in\bbR$.
\item[(2)]
$\bfT=\bigoplus_{\alpha\in \bfA}\bfT_\alpha$: an algebraic sum of Banach spaces $(\bfT_\alpha,\|\cdot\|_\alpha)$.
\item[(3)]
$\bfG$: a group of continuous linear operators on $\bfT$ such that, for any $\Gamma\in \bfG$ and $\alpha\in\bfA$,
$$
(\Gamma-\iden)\bfT_\alpha\subset \bigoplus_{\beta\in \bfA,\, \beta<\alpha}\bfT_\beta.
$$
\end{itemize}
\end{defn*}

\medskip

A concrete regularity structure $\big((T^+,\Delta^+), (T,\Delta)\big)$ turns into a regularity structure in the above sense, by setting $\bfA=A$, $\bfT=T$, and $\bfG=\{\widehat{g}\}_{g\in G^+}$.
The map $\widehat{g}$ is denoted by $\Gamma_g$ in Hairer's work \cite{Hai14}, where $G=\{\Gamma_g\ ;\, g\in G^+\}$ is defined as a (particular form of) structure group. In this article, we prefer the former Fourier-like notation, which is consistent with the fact that the `hat' map defines a linear representation of $G^+$ into $L(T)$.

\bigskip

\subsection{Models and modelled distributions}
\label{SubsectionModelsAndCo}

The preceding section contains the algebraic backbone of regularity structures. Its analytic flesh is introduced in this section on models and modelled distributions. This analytic setting depends on which (system of) singular stochastic PDE(s) one studies. We will not use the same function spaces to analyze a class of equations involving only the heat operator $(\partial_t-\Delta_x)$ and a system of two equations involving $(\partial_t-\Delta_x)$ for one and an operator with a different scaling for the other, like $(\partial_t + (-\Delta_x)^a)$, or simply $(-\Delta_x)^a$, with 
$0<a\neq1$, for the other. {\it We choose to concentrate in the present work on parabolic equations involving the heat operator only}. We will thus work throughout with the parabolic space $\bbR\times\bbR^d$, with generic point $x=(x_0,x')=(x_i)_{i=0}^d$, equipped with the distance function
$$
d(x,y) = d\big((x_0,x'),(y_0,y')\big) = \sqrt{\vert x_0-y_0\vert} + \vert x'-y'\vert.
$$
The H\"older spaces introduced in the next paragraph of this section will play a prominent role. They are used in the second paragraph to define models over a given regularity structure. Models give us the reference functions/distributions $(\Pi_x\tau)(\cdot)$ and ${\sf g}_{yx}(\mu)$ that we will use in our expansion devices to describe potential solutions of given singular stochastic PDEs. Expansion devices associate to each spacetime point a distribution meant to give the local description of a globally 
defined distribution. There is however no reason that such a globally defined distribution exists if no condition on its local `jets' is imposed. The appropriate consistency condition is encoded in the definition of a modelled distribution. Under this consistency condition, it is a fundamental fact that all these local descriptions can be patched together to define a unique globally defined distribution locally that is close to its local description, everywhere. This is what the reconstruction theorem does 
for us. We end this section with a paragraph on the special properties of 
modelled distributions representing functions.

\vfill \pagebreak

\hfill \textcolor{black}{\textsf{\textbf{ {\Large \S} Function spaces}}}

\medskip

Set
$$
\frak{s} \defeq  (2,1,\dots,1)\in\bbN\times\bbN^d,
$$
and define, for any multi-index $n=(n_0,n_1,\dots,n_d)\in\bbN\times\bbN^d$, the scaled degree of $n$ by 
$$
|n|_{\frak{s}} \defeq  2n_0+n_1+\cdots+n_d.
$$ 
Throughout this article, we define some analytic tools by using the specific kernel approach as introduced in Otto \& Weber \cite{OttoWeber}, instead of the local presentation as in \cite{Hai14}. We define a non-positive elliptic operator on $\bbR\times\bbR^d$
$$
\mcG \defeq  \partial_{x_0}^2 - \Delta_{x'}^2
$$
and denote by 
$$
P_t \defeq  e^{t\mcG}
$$ 
its semigroup, and by $p_t(x,y)$ its kernel with respect to Lebesgue 
measure. It is a symmetric function of $(x,y)$ that satisfies the scaling 
property
$$
p_t(x,y) = t^{-(d+2)/4}p\left(\Big(t^{-\frak{s}_j/4}(x_j-y_j)\Big)_{0\leq j\leq d}\right)
$$
for a Schwartz function $p\in\mcS(\bbR\times\bbR^d)$. The estimate
\begin{equation} \label{EqMomentEstimate} 
\int \big\vert \partial^n_x p_t(x,y) \big\vert \, d^a(x,y)\,dy \lesssim t^{\frac{a-\vert n\vert_{\frak{s}}}{4}},
\end{equation}
holds as a consequence for any multiindex $n\in\bbN\times\bbN^d$ and any positive exponent $a$. For a fixed positive integer $N\in\bbN$, we define operators $Q_t^{(N)}$ and $P_t^{(N)}$ setting
$$
Q_t^{(N)} \defeq  (-t\mcG)^N e^{t\mcG},\quad 
P_t^{(N)}\defeq  \int_t^\infty Q_s^{(N)} \frac{ds}s.
$$ 
This implies that $P_t^{(N)} = {\sf P}(-t\mcG) e^{t\mcG}$, for a monic polynomial ${\sf P}$ of degree $N$. One has in particular $P_t^{(1)}=P_t$, and
\begin{equation} \label{EqContinuousLPDecomposition}
P_t^{(N)} = \int_t^1 Q_s^{(N)}\,\frac{ds}{s} + P_1^{(N)}.
\end{equation}
(Those who know a little about Littlewood-Paley decomposition will recognize in $Q_t^{(N)}$ the counterpart of the Littlewood-Paley projectors $\Delta_i$ and in the integral with respect to the measure $ds/s$ the counterpart of the uniform measure on the integers; the integral operator associated with $P_1^{(N)}$ plays the role of $\Delta_{-1}$; this is an infinitely smoothing operator.) 

\ssk

\begin{defn*}
Fix $N\geq 1$ and pick a real number $\alpha < 4N$. We define the \textbf{\textsf{$\alpha$-H\"older space}} $\mcC^\alpha(\bbR\times\bbR^d)$ as the 
set of tempered distributions on $\bbR\times\bbR^d$ with finite $\mcC^\alpha$-norm 
defined by 
\begin{equation} \label{EqDefnHolderNorm}
\|\Lambda\|_{\mcC^\alpha} \defeq  \big\|P_1^{(N)}(\Lambda)\big\|_{L^\infty(\bbR\times\bbR^d)} + \underset{0<t\leq 1}{\sup}\,t^{-\frac{\alpha}{4}}\big\|Q_t^{(N)}(\Lambda)\big\|_{L^\infty(\bbR\times\bbR^d)}.
\end{equation}
\end{defn*}

\ssk

We work with the elements of $\mcS'$ (the space of tempered distributions) satisfying global H\"older estimates.
In contrast, Hairer \cite{Hai14} considered the elements of $\mcD'$ (the dual of compactly supported smooth functions) satisfying local H\"older estimates.
Such a technical difference is not important here, but when considering distributions diverging at infinity, the former definition should be appropriately modified by incorporating weight functions. See \cite{BH23} for instance.

The constraint on $\alpha$ of $\mcC^\alpha$ comes from the fact that all polynomials of scaled degree no greater than $4N$ are in the kernel of the operator $Q_t^{(N)}$. The above definition of the H\"older spaces depends on $N$, for the range of regularity exponents considered; write momentarily $\mcC^\alpha_N(\bbR\times\bbR^d)$. We remark that if $\alpha/4<N<N'$ are given then one can prove that the $N$ and $N'$-dependent norms are equivalent on $\mcC^\alpha_N(\bbR\times\bbR^d)$ -- this is a classical fact, worked out e.g. in Appendix {\sf A} of \cite{BB1}. In the sequel, the exponent $N$ is fixed once and for all to a large enough value depending on the problem at hand, so we do not record it in the notations for the H\"older spaces. More generally, one can define H\"older spaces using other elliptic operators than $\mcG$ with the same `scaling properties' as $\mcG$; the spaces will be identical and the different norms equivalent. We will use this remark only in the proof of Proposition \ref{PropSchauderEstimates} on the classical Schauder estimates. One can also show that for a positive non-integer regularity exponent $a$ the space $\mcC^a(\bbR\times\bbR^d)$ coincides with the usual space of $a$-H\"older functions, for the parabolic distance $d$, with equivalent norms. See e.g. the proof of Proposition 2.5 in \cite{BB1}.

\ssk

Note that, if $\alpha<0$, then the equivalence
\begin{equation} \label{EqEquivalentNegativeHolderNorm}
\|\Lambda\|_{\mcC^\alpha} \simeq \sup_{0<t\le1}t^{-\frac\alpha{4}} \big\|P^{(N)}_t(\Lambda)\big\|_{L^\infty(\bbR\times\bbR^d)} \simeq \sup_{0<t\le1}t^{-\frac\alpha{4}} \|P_t(\Lambda)\|_{L^\infty(\bbR\times\bbR^d)}
\end{equation}
holds. The middle term is bounded by the right one, because $P_1^{(1)}=P_1$ and $Q_t^{(1)}=\varphi(t\mcG)P_{t/2}$, with a uniformly bounded 
operator $\varphi(t\mcG)$. The other direction follows from identity \eqref{EqContinuousLPDecomposition} relating the operators $P$ and $Q^{(1)}$. 

\bigskip

\begin{Remark*}
Otto \& Weber \cite{OttoWeber} were the first to use the semigroup generated by $\mcG$ in a singular stochastic PDE setting.
\end{Remark*}

\bigskip

\hfill \textsf{\textbf{ {\Large \S} Models}}

\medskip

Recall from the introduction of Section \ref{SectionBasicsRS} the intuitive motivation for introducing regularity structures. Whereas the algebra involved in the use of local description devices is captured by the notion of regularity structure, the actual family of functions and distributions involved in these local descriptions is captured by the notion of model over a regularity structure. 

\smallskip

\begin{defn} \label{DefnModel}
A \textbf{\textsf{model over a regularity structure}} $\mathscr{T}$ is a pair $({\sf g}, {\sf \Pi})$ of maps
$$
{\sf g} : \bbR\times\bbR^d \rightarrow G^+, \qquad {\sf \Pi} : T\rightarrow \mcS'(\bbR\times\bbR^d)
$$
with the following properties.   \vspace{0.1cm}
\begin{itemize}
\setlength{\itemsep}{0.1cm}
   \item Set ${\sf g}_{yx}\defeq {\sf g}_y*{\sf g}_x^{-1}$, for each $x,y\in\bbR\times\bbR^d$. For each exponent $\gamma\in\bbR$, one has
\begin{align}\label{EqEstimateGammayx}
\|{\sf g}\|_\gamma\defeq  \sup_{\tau\in\mcB^+,|\tau|<\gamma}\sup_{x,y\in\bbR\times\bbR^d}\frac{|{\sf g}_{yx}(\tau)|}{d(y,x)^{|\tau|}} < \infty.
\end{align}

   \item The map ${\sf\Pi}$ is linear. Set 
   $$
   {\sf\Pi}_x^{\sf g}\defeq ({\sf\Pi}\otimes{\sf g}_x^{-1})\Delta
   $$ 
   for each $x\in\bbR\times\bbR^d$. For each exponent $\gamma\in\bbR$, one has 
\begin{align}\label{EqestimatePix}
\|{\sf\Pi}^{\sf g}\|_\gamma \defeq  \sup_{\sigma\in\mcB, |\sigma|<\gamma}\sup_{x\in\bbR\times\bbR^d,\, 0<t\le 1} t^{-\frac{|\sigma|}{4}} \big|\big\langle {\sf\Pi}_x^{\sf g}\sigma , p_t(x,\cdot) \big\rangle\big| <  \infty.
\end{align}
\end{itemize}
We also define a pseudo-distance on the space of models over a given regularity structure setting for each $\gamma\in\bbR$
\begin{equation} \label{EqPseudoDistanceModels}
{\sf d}_\gamma{\sf (M,M')} \defeq  
\sup_{\substack{\tau\in\mcB^+,\, |\tau|<\gamma \\ x,y\in\bbR\times\bbR^d}}\frac{\big|{\sf g}_{yx}(\tau)-{\sf g}'_{yx}(\tau)\big|}{d(y,x)^{|\tau|}} 
+ \sup_{\substack{\sigma\in\mcB,\,|\sigma|<\gamma \\ x\in\bbR\times\bbR^d \\ 0<t\leq 1}}\,t^{-\frac{|\sigma|}{4}}\,\Big|\big\langle {\sf\Pi}_x^{\sf g}\sigma - {\sf\Pi'}_x^{\sf g'}\sigma , p_t(x,\cdot) \big\rangle\Big|. 
\end{equation}
\end{defn}

\medskip

\noindent By the analytic properties of $\sf(g,\Pi)$, we have
\begin{align}\label{EqEstimatehatgyx}
\|\widehat{\sf g}\|_\gamma\defeq\sup_{\sigma,\tau\in\mcB,\,|\sigma|\le|\tau|<\gamma}\sup_{x,y\in\bbR\times\bbR^d}\frac{\big|\big(\widehat{{\sf g}_{yx}}(\tau)\big)_\sigma\big|}{d(y,x)^{|\tau|-|\sigma|}}<\infty,
\end{align}
where $(\cdot)_\sigma$ denotes the $\sigma$-components of elements of $T$.
In Hairer's original work \cite{Hai14}, a pair $(\Pi=\{\Pi_x\}_x,\Gamma=\{\Gamma_{yx}\}_{x,y})$ of the family of linear operators 
$$
\Pi_x:T\to\mcS'(\bbR\times\bbR^d),\qquad \Gamma_{yx}:T\to T
$$
is called a model if it satisfies analytic conditions \eqref{EqestimatePix} and \eqref{EqEstimatehatgyx} and algebraic conditions
$$
\Pi_x=\Pi_y\Gamma_{yx},\qquad
\Gamma_{xx}=\iden,\qquad
\Gamma_{zx}=\Gamma_{zy}\Gamma_{yx}.
$$
In our setting, the choice of operators $\Pi_x\defeq{\sf\Pi}_x^{\sf g}$ and $\Gamma_{yx}=\widehat{{\sf g}_{yx}}$ provides a model in the original meaning. The algebraic conditions satisfied by $\Pi_x$ and $\Gamma_{yx}$ are encoded by the algebraic properties of $\Delta$ and $\Delta^+$.
Indeed, since ${\sf g}_y^{-1}*{\sf g}_{yx}={\sf g}_x^{-1}$,
\begin{align}\label{EqTransitionRelation}
\begin{aligned}
{\sf\Pi}_y^{\sf g}\circ\widehat{{\sf g}_{yx}} 
&=({\sf\Pi}\otimes {\sf g}_y^{-1}\otimes {\sf g}_{yx})(\Delta\otimes\iden)\Delta\\
&=({\sf\Pi}\otimes {\sf g}_y^{-1}\otimes {\sf g}_{yx})(\iden\otimes\Delta^+)\Delta\\
&=\big({\sf\Pi}\otimes ({\sf g}_y^{-1}*{\sf g}_{yx})\big)\Delta={\sf\Pi}_x^{\sf g}.
\end{aligned}
\end{align}
So the above choice is more specific, but such a specific choice can be found in \cite[Section 8]{Hai14} and \cite{BHZ, ChandraHairer}.
We also remark that for the analytic estimates in original definition \cite{Hai14} the supremum over $x,y$ are local, and a family compactly supported test functions are considered in the condition \eqref{EqestimatePix}, instead of a single function $p_t(x,\cdot)$. For simplicity, we use the global estimates \eqref{EqEstimateGammayx} and \eqref{EqestimatePix} in this article. 

Emphasize that ${\sf g}$ acts on $T^+$, while $\sf \Pi$ acts on $T$, and \textit{note that ${\sf g}$ plays on $T^+$ the same role as ${\sf \Pi}$ on $T$}: For $\tau\in T^+$ and $\sigma\in T$, one has
\begin{equation}
\label{EqParallel}
{\sf g}_{yx}(\tau) = \Big({\sf g}_y(\cdot)\otimes {\sf g}_x^{-1}(\cdot)\Big)\Delta^+\tau, \quad
({\sf \Pi}_x^{\sf g}\sigma)(y) = \Big({\sf \Pi}(\cdot)(y)\otimes {\sf g}_x^{-1}(\cdot)\Big)\Delta\sigma,
\end{equation}
in a distributional sense for the latter. Therefore the maps
$$
\big({\sf \Pi}^{(\sf g)}\tau\big)(x) \defeq  {\sf g}_x(\tau), \qquad x\in\bbR\times\bbR^d,\;\tau\in T^+
$$
define a model $\sf \big(g,\Pi^{(g)}\big)$ on $\mathscr{T}^+=(T^+,T^+)$.

\ssk

In the class of problems we consider, it is sufficient in each problem to fix $\gamma\in\bbR$ to a large enough value; we omit as a consequence this parameter from the notations, unless necessary. We emphasize the dependence of $\Pi_x$ on $\sf g$ using our notation. We stress that ${\sf\Pi}\tau$ is only an element of $\mcS'(\bbR\times\bbR^d)$. Think of ${\sf\Pi}$ as an interpretation operator for the symbols $\tau$, with $\tau$ encoding the structure of the analytic object ${\sf\Pi}\tau$. One can think of ${\sf\Pi}^{\sf g}_x\tau=({\sf\Pi}\otimes{\sf g}_x^{-1})\Delta\tau$, as ${\sf \Pi}\tau$ `fully recentered' at $x$, to give it a concrete meaning. The splitting map $\Delta$ identifies the different sets of internal pieces of $\tau$ that can be `recentered' to the point $x$ by the action of the map ${\sf g}_x^{-1}$, with the full recentering operation on ${\sf \Pi}\tau$ being the result of all these recentering operations. Condition \eqref{EqestimatePix} conveys the idea that ${\sf\Pi}^{\sf g}_x\tau$ behaves \textit{at point $x$} like an element of $\mcC^{|\tau|}(\bbR\times\bbR^d)$, as a result of this full recentering operation. We will see in Section \ref{SectionBuildingRS} concrete examples of recentering operations that can be understood as replacing a function by its Taylor remainder of a certain degree.    

The following immediate consequence of the bound \eqref{EqestimatePix} will be useful in the next section.

\medskip

\begin{prop} \label{PropModelsOnKernels}
One has
$$
\sup_{x\in\bbR\times\bbR^d,\, 0<t\le 1} t^{-\frac{|\tau|-|n|_{\frak{s}}}{4}}
\big|\big\langle{\sf\Pi}_x^{\sf g}\tau , \partial_x^n p_t(x,\cdot) \big\rangle\big| < \infty,
$$
for any model $({\sf\Pi},{\sf g})$ on $\mathscr{T}$, $\tau\in\mcB$, and $n\in\bbN\times\bbN^d$.
\end{prop}

\medskip

\begin{Dem}
By the semigroup property,
$$
\partial_x^np_t(x,y)=\int\partial_x^np_{\frac{t}2}(x,z)p_{\frac{t}2}(z,y)dz.
$$
We need to apply the distribution ${\sf\Pi}_x^{\sf g}\tau$ to the kernel $p_{\frac{t}2}(z,\cdot)$. Using the expansion \eqref{EqCoprodQuotient} of $\Delta\tau$ and the relation \eqref{EqTransitionRelation} to write ${\sf \Pi}_x^{\sf g}\tau$ in terms of ${\sf \Pi}_z^{\sf g}\tau$, one has
\begin{align*}
\big\vert \big\langle {\sf\Pi}_x^{\sf g}\tau, p_{\frac{t}2}(z,\cdot) \big\rangle \big\vert = \Big\vert\sum_{\sigma\le\tau} {\sf g}_{zx}(\tau/\sigma) \big\langle {\sf\Pi}_z^{\sf g}\sigma, p_{\frac{t}2}(z,\cdot) \big\rangle\Big\vert \lesssim \sum_{\sigma\le \tau} d(z,x)^{|\tau|-|\sigma|}\, t^{\frac{|\sigma|}{4}}
\end{align*}
from the bound \eqref{EqestimatePix} in the definition of a model. Using the bound \eqref{EqMomentEstimate} on the moments of the heat kernel, we then have
\begin{align*}
\Big|\big\langle{\sf\Pi}_x^{\sf g}\tau , \partial_x^n p_t(x,\cdot) \big\rangle\Big|
&\lesssim \sum_{\sigma\le \tau} t^{\frac{|\sigma|}{4}} \int \big\vert\partial_x^np_{\frac{t}2}(x,z)\big\vert \, d(z,x)^{|\tau|-|\sigma|} dz \\
&\lesssim \sum_{\sigma\le \tau} \,t^{\frac{|\sigma|}{4}} t^{\frac{|\tau|-|\sigma|-|n|_{\frak{s}}}{4}}
\lesssim t^{\frac{|\tau|-|n|_{\frak{s}}}{4}}.
\end{align*}
\end{Dem}

\medskip

We close this paragraph by the remark about the situation where all the ${\sf \Pi}\tau$ are some continuous functions. Then it follows from the bound on $\big\langle{\sf\Pi}_x^{\sf g}\tau , p_t(x,\cdot) \big\rangle$, and the fact that $p_t(x,\cdot)$ is converging to a Dirac mass at $x$, that the function ${\sf \Pi}^{\sf g}_x\tau$ satisfies $({\sf \Pi}^{\sf g}_x\tau)(x) = 0$, for all $\tau\in T$ such that $\vert\tau\vert>0$. This will be the case of the smooth (possibly renormalized) 
models from Section \ref{SectionMultiAndrenormalizedEquations}.

\bigskip

\hfill \textsf{\textbf{ {\Large \S} Modelled distributions and their reconstruction}}

\medskip

Think of a $T$-valued function $\bsf$ on $\bbR\times\bbR^d$ as the data needed to associate with each spacetime point $x\in\bbR\times\bbR^d$ the local description ${\sf \Pi}^{\sf g}_x\bsf(x)$ of a possibly globally defined distribution close to ${\sf \Pi}^{\sf g}_x\bsf(x)$ near each $x$.  There is no reason that such a globally defined object exists if one does not impose relations between the different components of $\bsf$. This is what the next definition does. For a real number $\gamma\geq\beta_0$ set
\vspace{-0.1cm}$$
T_{<\gamma} \defeq  \bigoplus_{\beta<\gamma} T_\beta, \qquad T_{<\gamma}^+ \defeq  \bigoplus_{\alpha<\gamma} T_\alpha^+.
$$
Recall from \eqref{EqNotationDirectSum} the meaning of the notation $\|\tau\|_\alpha$, for $\alpha\in A$ and $\tau\in T$.

\ssk

\begin{defn} \label{DefnModelledDistribution}
Let ${\sf g} : \bbR\times\bbR^d\rightarrow G^+$ be a function satisfying \eqref{EqEstimateGammayx}. Fix a regularity exponent $\gamma\in\bbR$. One defines the space $\mcD^\gamma(T, {\sf g})$ of \textbf{\textsf{distributions modelled on the regularity structure $\mathscr{T}$, with transition ${\sf g}$}}, as the space of functions $\bsf : \bbR\times\bbR^d\rightarrow T_{<\gamma}$ such that
\begin{equation*}
\begin{split}
&\brarb{\bsf}_{\mcD^\gamma} \defeq  \max_{\beta<\gamma}\, \sup_{x\in \bbR\times\bbR^d}\, \big\| \bsf(x) \big\|_{\beta} < \infty,   \\
&\|\bsf\|_{\mcD^\gamma} \defeq  \max_{\beta<\gamma}\sup_{x,y\in \bbR\times\bbR^d}\frac{\big\| \bsf(y) - \widehat{{\sf g}_{yx}}\bsf(x) \big\|_{\beta}} 
{d(y,x)^{\gamma-\beta}} < \infty.
\end{split}
\end{equation*}
Set $\trino{\bsf}_{\mcD^\gamma} \defeq  \brarb{\bsf}_{\mcD^\gamma} + \|\bsf\|_{\mcD^\gamma}$. We also define the pseudo-distance between two modelled distributions $\bsf\in\mcD^\gamma(T,\sf g)$ and $\bsf'\in\mcD^\gamma(T,\sf g')$ defined for two distinct models models $\sf M=(g,\Pi)$ and $\sf M'=(g',\Pi')$, by setting
$$
d\big(\bsf,\bsf'\big) \defeq  \underset{x,y\in\bbR\times\bbR^d}{\sup}\,\underset{\beta<\gamma}{\max}\,\left\{ \big\|\bsf(x) - \bsf'(x)\big\|_\beta + \frac{\Big\| \Big\{\bsf(y)-\widehat{{\sf g}_{yx}}\bsf(x)\Big\} -  \Big\{\bsf'(y)-\widehat{{\sf g}'_{yx}}\bsf'(x)\Big\}\Big\|_\beta}{d(y,x)^{\gamma-\beta}}  \right\}.
$$
\end{defn}

\ssk

(As in the definition of models, we choose a global bound to define modelled distributions -- the original definition in \cite[Section 2]{Hai14} relies on the local bounds. Since we consider the equations on the domain $(x_0,x')\in [0,T]\times \bbT$, there is no difference between global and local bounds on the spatial variable $x'$. In Section \ref{SectionSolvingPDEs} below, we will consider a weighted norm with respect to the time variable $x_0$.)

For a basis element $\sigma\in\mcB\subset T$, and an arbitrary element $h$ in $T$, denote by $h_\sigma$ its component on $\sigma$ in the basis $\mcB$. For a modelled distribution $\bsf(\cdot) = \sum_{\sigma\in\mcB} f_\sigma(\cdot) \, \sigma$ in $\mcD^\gamma(T, {\sf g})$, and $\sigma_0\in\mcB$, we have
\begin{equation} \label{EqIncrementMdF} \begin{split}
\big(\bsf(y) - \widehat{{\sf g}_{yx}}\bsf(x)\big)_{\sigma_0} &= f_{\sigma_0}(y) - \sum_{\tau \geq \sigma_0} {\sf g}_{yx}(\tau/\sigma_0) \, f_\tau(x)  \\
&= f_{\sigma_0}(y) - f_{\sigma_0}(x) - \sum_{\tau > \sigma_0} {\sf g}_{yx}(\tau/\sigma_0) \, f_\tau(x).
\end{split} \end{equation}

\medskip

{\small {\sl \begin{Examples} $\bullet$ The archetype of a modelled distribution is given by the lift 
\begin{equation} \label{EqLiftHolder}
\bsf(x) \defeq  \sum_{\vert n\vert_{\frak{s}} < \gamma} \frac{f^{(n)}(x)}{n!} X^n,
\end{equation}
in the polynomial regularity structure of a $\gamma$-H\"older real valued 
function $f$ on $\bbR\times\bbR^d$ with a positive regularity exponent $\gamma$. The identities \eqref{EqIncrementMdF} become in that case the Taylor expansions
\begin{equation} \label{EqConditionJetHolder}
f^{(n)}(y) - f^{(n)}(x) - \sum_{\vert\ell\vert_{\frak{s}}<\gamma-\vert n\vert_{\frak{s}}} \frac1{\ell!} f^{(n+\ell)}(x) (y-x)^\ell = O\Big(d(y,x)^{\gamma-\vert n\vert_\frak{s}}\Big)
\end{equation}
satisfied by each $f^{(n)}$. Note here that the function $f$ on $\bbR\times\bbR^d$ is $\gamma$-H\"older iff there exists a family $(f_n)_{n\in\bbN\times\bbN^d},\vert n\vert_\frak{s}<\gamma$ of functions on $\bbR\times\bbR^d$ satisfying $f_0=f$ and the condition \eqref{EqConditionJetHolder} with $f_n$ in the role of $f^{(n)}$. The ``if" part holds because one can get $f_n=f^{(n)}$ from \eqref{EqConditionJetHolder} inductively. The ``only if" part is the classical and elementary fact for the isotropic case, and can be found in Appendix {\sf A} of \cite{Hai14} for anisotropic cases. So the notion of modelled distribution with values in the polynomial regularity structure captures exactly the classical notion of regularity.

\ssk

$\bullet$ Given a basis element $\tau\in \mcB$, set
\begin{equation}
\label{EqDefnMdTau}
\bsh^\tau(x) \defeq  
\widehat{{\sf g}_x}\tau-\tau =
\sum_{\sigma<\tau}{\sf g}_x(\tau/\sigma)\sigma.
\end{equation}
It follows from identity \eqref{eq T is graded} in the definition of a concrete regularity structure that $\bsh^\tau$ takes values in $T_{<\vert\tau\vert}$. 
Since $\widehat{{\sf g}_{yx}}\circ\widehat{{\sf g}_x}=\widehat{{\sf g}_y}$ by \eqref{Section2:GactionAssociative}, it follows that
\begin{align*}
\bsh^\tau(y)-\widehat{{\sf g}_{yx}}\big(\bsh^\tau(x)\big)
&=(\widehat{{\sf g}_y}\tau-\tau)-\widehat{{\sf g}_{yx}}(\widehat{{\sf g}_x}\tau-\tau)\\
&=(\widehat{{\sf g}_y}\tau-\tau)-(\widehat{{\sf g}_y}\tau-\widehat{{\sf g}_{yx}}\tau)\\
&=\widehat{{\sf g}_{yx}}\tau-\tau=
\sum_{\sigma<\tau}{\sf g}_{yx}(\tau/\sigma)\sigma.
\end{align*}
The size estimate $\big|{\sf g}_{yx}(\tau/\sigma)\big| \lesssim d(y,x)^{|\tau|-|\sigma|}$ required from the $\sf g$-component of a model, then shows that $\bsh^\tau$ is a modelled distribution in $\mcD^{|\tau|}(T_{<|\tau|}, 
{\sf g})$.

\ssk

$\bullet$ If $\bsf(\cdot) = \sum_{\sigma\in\mcB}f_\sigma(\cdot)\,\sigma$, is an element of $\mcD^\gamma(T,{\sf g})$, then, for each $\tau\in\mcB$, the $T^+$-valued function 
$$
\bsf/\tau (\cdot)\defeq  \sum_{\sigma\geq\tau}f_\sigma(\cdot)\,\sigma/\tau.
$$
is an element of $\mcD^{\gamma-\vert\tau\vert}(T^+,{\sf g})$, where we denote by $\mcD^\gamma(T^+,{\sf g})$ the space of $T^+$-valued modelled distributions with transition ${\sf g}$. Recall that $\mathscr{T}^+=(T^+,T^+)$ is also a regularity structure.
\end{Examples}}}

\ssk

The next statement says that the consistency condition encoded in the notion of modelled distribution $\bsf$ ensures the existence of a globally defined object close to ${\sf \Pi}^{\sf g}_x\bsf(x)$ near each $x\in\bbR\times\bbR^d$, and gives condition for uniqueness. Recall $A$ stands for the index set in the grading of $T$ and set
$$
\beta_0 \defeq  \min A.
$$ 

\medskip

\begin{thm} \label{ThmReconstructionRS} \textbf{\textsf{(Reconstruction theorem)}}
Let $\mathscr{T}$ be a concrete regularity structure and ${\sf M}=(\sf g,\Pi)$ be a model over $\mathscr{T}$. Fix a regularity exponent $\gamma\in\bbR\backslash\{0\}$. There exists a linear continuous operator 
$$
\textsf{\textbf{R}}^{\sf M} : \mcD^\gamma(T, {\sf g}) \rightarrow \mcC^{\beta_0\wedge0}(\bbR\times\bbR^d)
$$
satisfying the property
\begin{equation}
\label{EqReconstructionCondition}
\Big|\big\langle \textsf{\textbf{R}}^{\sf M}\bsf - {\sf\Pi}_x^{\sf g}\bsf(x) , p_t(x,\cdot)\big\rangle\Big|  \lesssim  \|{\sf \Pi}^{\sf g}\|\, \big\|\bsf\big\|_{\mcD^\gamma} \, t^{\frac{\gamma}{4}},
\end{equation}
uniformly in $\bsf\in\mcD^\gamma(T, {\sf g}), x\in\bbR\times\bbR^d$ and $0<t\leq 1$. Such an operator is unique if the exponent $\gamma$ is positive.
\end{thm}

\medskip

A distribution $\textbf{\textsf{R}}^{\sf M}\bsf$ satisfying identity \eqref{EqReconstructionCondition} is called a \textbf{\textsf{reconstruction of the modelled distribution}} $\bsf$. When $\gamma=0$, the existence of a reconstruction is not ensured by \eqref{EqReconstructionCondition} in 
general. See Example 5.5 in \cite{CaravennaZambotti}. We will see as a particular case of Corollary \ref{Cor reconstruction of smooth model} that the lift \eqref{EqLiftHolder} of a $\gamma$-H\"older function $f$ in the polynomial regularity structure has indeed $f$ as a reconstruction.

Notice from the definition of ${\sf\Pi}^{\sf g}_x$, we have the relation
\begin{align*}
({\sf\Pi}_x^{\sf g}\otimes{\sf g}_x)\Delta
&=({\sf\Pi}\otimes{\sf g}_x^{-1}\otimes{\sf g}_x)(\Delta\otimes\iden)\Delta\\
&=({\sf\Pi}\otimes{\sf g}_x^{-1}\otimes{\sf g}_x)(\iden\otimes\Delta^+)\Delta\\
&=({\sf\Pi}\otimes{\bf1}_+')\Delta={\sf\Pi},
\end{align*}
and thus
$$
{\sf\Pi}\tau={\sf\Pi}_x^{\sf g}(\widehat{{\sf g}_x}\tau)={\sf\Pi}_x^{\sf g}\tau+{\sf\Pi}_x^{\sf g}\bsh^\tau.
$$
Therefore the constraint $\big|\big\langle{\sf\Pi}^{\sf g}_x\tau , p_t(x,\cdot) \big\rangle\big|  \lesssim t^{|\tau|/4}$, that needs to be satisfied by a model, is equivalent to the estimate
\begin{equation} \label{EqEquivalentPiTau}
\Big|\big\langle {\sf \Pi}\tau - {\sf\Pi}^{\sf g}_x \bsh^\tau(x) , p_t(x,\cdot)\big\rangle\Big| = \Big|\big\langle {\sf \Pi}\tau - \sum_{\sigma<\tau} {\sf g}_x(\tau/\sigma){\sf\Pi}^{\sf g}_x\sigma , p_t(x,\cdot)\big\rangle\Big|  \lesssim t^{|\tau|/4},
\end{equation}
which says that ${\sf \Pi}\tau$ is a/the reconstruction of the modelled distribution $\bsh^\tau$ from \eqref{EqDefnMdTau}, depending on whether $\vert\tau\vert\leq 0$ or $\vert\tau\vert>0$. Since, for $\vert\tau\vert<0$, the difference $(\ast)$ of two reconstructions of $\bsh^\tau$ satisfies
$$
\big\vert\big\langle(\ast),p_t(x,\cdot)\big\rangle\big\vert \lesssim t^{\vert\tau\vert/4}
$$
for all $x\in\bbR\times\bbR^d$, this difference is a $\mcC^{\vert\tau\vert}(\bbR\times\bbR^d)$ distribution from identity \eqref{EqEquivalentNegativeHolderNorm}. So the estimate \eqref{EqEquivalentPiTau} shows in particular that we could require from scratch that the $\sf \Pi$ map of a model 
of $\mathscr{T}$ takes values in $\mcC^{\beta_0\wedge 0}(\bbR\times\bbR^d)$ rather than $\mcS'(\bbR\times\bbR^d)$. The case $\vert\tau\vert=0$ does not cause any problem as we assume that the only element of $T$ of null homogeneity is $\textbf{\textsf{1}}$. 

\ssk

We will only work with $\mcD^\gamma(T,\sf g)$-spaces with positive regularity exponents $\gamma$ in our study of singular stochastic PDEs. We only 
give a proof of the reconstruction theorem in that setting, following Otto \& Weber's nice approach \cite{OttoWeber}. An extension to the inhomogeneous integral kernels can be found in \cite{Hos23}. See Friz and Hairer's lecture notes \cite{FrizHairerBook} for another treatment along these lines. See Hairer's original work \cite{Hai14} or the references given in Appendix {\sf \ref{SectionAppendixComments}} for a proof of Theorem \ref{ThmReconstructionRS} in the case $\gamma\leq 0$.

\medskip

\begin{Dem}
{\it Existence --} We construct explicitly a reconstruction operator. Note first that since
\begin{equation*} \begin{split}
\Big({\sf \Pi}^{\sf g}_y\bsf(y) - {\sf \Pi}^{\sf g}_x\bsf(x)\Big)(\cdot) &= \Big({\sf \Pi}^{\sf g}_y\big(\bsf(y) - \widehat{{\sf g}_{yx}}\bsf(x)\big)\Big)(\cdot)   \\
&=\sum_{\tau\in\mcB} 
\big(\bsf(y) - \widehat{{\sf g}_{yx}}\bsf(x)\big)_\tau
\,\big({\sf \Pi}^{\sf g}_y\tau\big)(\cdot)
\end{split} \end{equation*}
one has 
\begin{equation*} \begin{split}
\left\vert\big\langle {\sf \Pi}^{\sf g}_y\bsf(y) - {\sf \Pi}^{\sf g}_x\bsf(x) , p_t(y,\cdot) \big\rangle\right\vert \lesssim \sum_{\tau\in\mcB, \vert\tau\vert<\gamma} d(y,x)^{\gamma-\vert\tau\vert} t^{\frac{\vert\tau\vert}{4}},
\end{split} \end{equation*}
from the bounds on models and modelled distributions. For $0<s\leq t\leq 1$ and $x\in\bbR\times\bbR^d$, set 
$$
\textbf{\textsf{I}}_s^t(x) \defeq  \int p_{t-s}(x,y) \, \big\langle {\sf \Pi}^{\sf g}_y\bsf(y) , p_s(y,\cdot) \big\rangle\,dy.
$$
We will obtain the distribution ${\sf R}^{\sf M}\bsf$ from $\textbf{\textsf{I}}_s^t$ under the form $\lim_{t\downarrow 0}\lim_{s\downarrow 0}\textbf{\textsf{I}}_s^t$, with the limits taken in that order, with $s$ sent to $0$ first and then $t$ sent to $0$. First, from the bounds on modelled 
distributions, we have
$$
\textbf{\textsf{I}}_t^t(x) = \big\langle {\sf \Pi}^{\sf g}_x\bsf(x) , p_t(x,\cdot) \big\rangle
\quad \textrm{and} \quad
\left\vert \textbf{\textsf{I}}_t^t(x) \right\vert
\le
\sum_{\tau\in\mcB} \left\vert f_\tau(x) \right\vert \left\vert \big\langle {\sf \Pi}^{\sf g}_x\tau , p_t(x,\cdot) \big\rangle \right\vert 
\lesssim
t^{\frac{\beta_0}{4}},
$$
and moreover, for $0<s'<s<t\leq 1$, we have from the semigroup property of the kernel $p$ the $x$-uniform estimate
$$
\left\vert \textbf{\textsf{I}}_{s'}^t(x) - \textbf{\textsf{I}}_s^t(x) \right\vert 
= \left\vert \int p_{t-s}(x,z)p_{s-s'}(z,y) \big\langle {\sf \Pi}^{\sf g}_{y}\bsf(y) - {\sf \Pi}^{\sf g}_z\bsf(z) , p_{s'}(y,\cdot)\big\rangle\,dzdy\right\vert
$$
\begin{equation*} \begin{split}
&\leq \sum_{\tau\in\mcB, \vert\tau\vert<\gamma} \int p_{t-s}(x,z)p_{s-s'}(z,y)\,d(y,z)^{\gamma-\vert\tau\vert}(s')^{\frac{\vert\tau\vert}{4}}\,dzdy   \\
&\lesssim \sum_{\tau\in\mcB, \vert\tau\vert<\gamma} (s-s')^{\frac{\gamma-\vert\tau\vert}{4}} (s')^{\frac{\vert\tau\vert}{4}}.
\end{split}
\end{equation*}
For $s'\in[s/2,s)$, this implies
\begin{equation} \label{*EqReconstructionCondition}
\left\vert \textbf{\textsf{I}}_{s'}^t(x) - \textbf{\textsf{I}}_s^t(x) \right\vert 
\lesssim s^{\frac{\gamma}{4}}.
\end{equation}
For $s'\in(0,s/2)$, by taking $n\in\bbN$ such that $s'\in\big[s/2^{n+1},s/2^n\big)$, we have
\begin{align*}
\left\vert \textbf{\textsf{I}}_{s'}^t(x) - \textbf{\textsf{I}}_s^t(x) \right\vert 
&\le\sum_{m=0}^{n-1} \left\vert \textbf{\textsf{I}}_{s/2^m}^t(x) - \textbf{\textsf{I}}_{s/2^{m+1}}^t(x) \right\vert 
+\left\vert \textbf{\textsf{I}}_{s/2^n}^t(x) - \textbf{\textsf{I}}_{s'}^t(x) \right\vert \\
&\lesssim \sum_{m=0}^{n-1}(s/2^m)^{\frac{\gamma}{4}}+(s/2^n)^{\frac{\gamma}{4}}
\lesssim s^{\frac{\gamma}{4}}.
\end{align*}
Thus the bound \eqref{*EqReconstructionCondition} holds uniformly over $0<s'<s$. Hence the (locally in $t$) uniform limit
$$
\textbf{\textsf{I}}_0^t(x) \defeq 
\underset{s\rightarrow 0}{\lim}\, \textbf{\textsf{I}}_s^t(x)
$$
exists, since $\gamma$ is positive. As the identity $P_{t'}\textbf{\textsf{I}}_0^t=\textbf{\textsf{I}}_0^{t+t'}$ follows from the semigroup property, we see from \eqref{EqEquivalentNegativeHolderNorm} that $\{\textbf{\textsf{I}}_0^t\}_{0<t\le1}$ is bounded in the space $\mcC^{\beta_0}(\bbR\times\bbR^d)$. (Note that all of the above estimates on $\textbf{\textsf{I}}_s^t$ holds over $0<s\le t\le2$, since the bounds on ${\sf\Pi}_x^{\sf g}\tau$ can be extended to $0<t\le2$ by a similar argument to Proposition \ref{PropModelsOnKernels}.) Therefore by noting the continuity of $t\mapsto P_t\Lambda$
$$
\|(P_t-\iden)\Lambda\|_{\mcC^{\beta_0-\varepsilon}}\lesssim t^{\frac\varepsilon4}\|\Lambda\|_{\mcC^{\beta_0}}
$$
for any $\varepsilon>0$ and $\Lambda\in\mcC^{\beta_0}$ (see e.g., Lemma 2.15 of \cite{Hos23}), for $0<s<t\le1$ we have
$$
\|\textbf{\textsf{I}}_0^t-\textbf{\textsf{I}}_0^s\|_{\mcC^{\beta_0-\varepsilon}}
=\|(P_{t-s}-\iden)\textbf{\textsf{I}}_0^s\|_{\mcC^{\beta_0-\varepsilon}}
\lesssim (t-s)^{\frac{\varepsilon}4}\|\textbf{\textsf{I}}_0^s\|_{\mcC^{\beta_0}}\lesssim (t-s)^{\frac\varepsilon4}.
$$
Hence $\{\textbf{\textsf{I}}_0^t\}_{0<t\le1}$ converges in $\mcC^{\beta_0-\epsilon}(\bbR\times\bbR^d)$ as $t$ goes to $0$, for any $\epsilon>0$. Denote its limit by $\textbf{\textsf{R}}^{\sf M}\bsf$. Since
$$
\big\langle \textbf{\textsf{R}}^{\sf M}\bsf, p_t(x,\cdot) \big\rangle = 
\lim_{s\to0} \big\langle \textbf{\textsf{I}}_0^{s} , p_t(x,\cdot) \big\rangle = \lim_{s\to0} \textbf{\textsf{I}}_0^{t+s}(x) = \textbf{\textsf{I}}_0^{t}(x),
$$
we have actually $\textbf{\textsf{R}}^{\sf M}\bsf\in\mcC^{\beta_0}(\bbR\times\bbR^d)$ from the $x$-uniform bound $\vert \textbf{\textsf{I}}_0^{t}(x) \vert \lesssim t^{\beta_0/4}$. Letting $s=t$ and sending $s'$ to $0$ in \eqref{*EqReconstructionCondition} we can check that $\textbf{\textsf{R}}^{\sf M}\bsf$ satisfies the bound \eqref{EqReconstructionCondition}. 

\smallskip

{\it Uniqueness -- } To prove uniqueness of the reconstruction operator on $\mcD^\gamma(T,\sf g)$ when the regularity exponent $\gamma$ is positive, we start from the identity
$$
\big\vert \big\langle\textbf{\textsf{R}}^{\sf M}\bsf - (\textbf{\textsf{R}}^{\sf M})'\bsf , p_t(x,\cdot)\big\rangle\big\vert \lesssim t^{\frac\gamma4},
$$
satisfied uniformly in $x\in\bbR\times\bbR^d$ by any other reconstruction 
operator $(\textbf{\textsf{R}}^{\sf M})'$. As for any Schwartz function $\varphi\in\mcS(\bbR\times\bbR^d)$ the convolutions $\int \varphi(x)p_t(x,z)dx$, converge to $\varphi$ in the smooth topology, one has from the symmetry of the kernels $p_t$ and the fact that $\gamma$ is positive
$$
\big\langle \textbf{\textsf{R}}^{\sf M}\bsf - \big(\textbf{\textsf{R}}^{\sf M}\big)'\bsf , \varphi \big\rangle = \underset{t\rightarrow 0}{\lim} 
\int \Big\langle  \textbf{\textsf{R}}^{\sf M}\bsf - (\textbf{\textsf{R}}^{\sf M})'\bsf , p_t(x,\cdot) \Big\rangle\,\varphi(x)dx = \underset{t\rightarrow 0}{\lim} \, O\big(t^{\frac\gamma4}\big) = 0.
$$
\end{Dem}

\medskip

One can use Proposition \ref{PropModelsOnKernels} to improve estimate \eqref{EqReconstructionCondition} under the form
\begin{equation} \label{EqImprovedReconstructionEstimate}
\Big|\big\langle \textsf{\textbf{R}}^{\sf M}\bsf - {\sf\Pi}_x^{\sf g}\bsf(x) , \partial_x^n p_t(x,\cdot)\big\rangle\Big|  \lesssim t^{\frac{\gamma-\vert n\vert_\frak{s}}{4}},
\end{equation}
uniformly in $x\in\bbR\times\bbR^d$, for each $n\in\bbN\times\bbN^d$. 

We now state two standard properties of the reconstruction operator. The following fact implies that $\big\langle\textbf{\textsf{R}}^{\sf M}\bsf , \varphi\big\rangle$ depends only on the restriction of $\bsf$ to the support of $\varphi$. This fact is used to define the reconstructions of modelled distributions which are given in $(0,t)\times\mathbb{R}^d$ with $t\in(0,\infty]$, not in $\bbR\times\bbR^d$. See Theorem \ref{thm singular reconstruction} and Section \ref{SubsectionNonAnticipative}.

\medskip

\begin{cor} \label{CorLocalDependenceReconstruction}
Pick $\gamma$ positive. If $\bsf\in\mcD^\gamma(T,{\sf g})$ is null on an open set $U\subset\bbR\times\bbR^d$, then $\textbf{\textsf{R}}^{\sf M}\bsf=0$ on $U$.
\end{cor}

\medskip

\begin{Dem}
Since the mapping $x\mapsto{\sf \Pi}^{\sf g}_x\bsf(x)$ is null on $U$, it follows from estimate \eqref{EqReconstructionCondition} that 
$$
\vert \big\langle \textbf{\textsf{R}}^{\sf M}\bsf, p_t(x,\cdot) \big\rangle \vert \lesssim t^{\gamma/4},
$$ 
for all $x\in U$. For a smooth function $\varphi$ with compact support in $U$, one can use the convergence of $\int\varphi(x)p_t(x,z)dx$ to $\varphi$ as $t>0$ goes to $0$ in an appropriate $C^k$ space, to get
$$
\big\langle \textbf{\textsf{R}}^{\sf M}\bsf, \varphi \big\rangle = \lim_{t\to0} \int \varphi(x)\,\big\langle \textbf{\textsf{R}}^{\sf M}\bsf, p_t(x,\cdot) \big\rangle \,dx = 0.
$$
\end{Dem}

\medskip

The following fact is an immediate consequence of uniqueness in the reconstruction theorem; it is used in Section \ref{SectionSolvingPDEs} and implies in particular that the the lift \eqref{EqLiftHolder} of a $\gamma$-H\"older function $f$ in the polynomial regularity structure has indeed $f$ as a reconstruction.

\medskip

\begin{cor}\label{Cor reconstruction of smooth model}
Pick $\gamma$ positive and $\bsf\in\mcD^\gamma(T,{\sf g})$.
If the model $({\sf g}, {\sf \Pi})$ takes values in the space of smooth functions on $\bbR\times\bbR^d$, then the mapping $x\mapsto\big({\sf \Pi}_x^{\sf g}\bsf(x)\big)(x)$ is itself a continuous function and
\begin{equation} \label{EqReconstructionFunctionLike}
\big(\textsf{\textbf{R}}^{\sf M}\bsf\big)(x) = \Big({\sf \Pi}_x^{\sf g}\bsf(x)\Big)(x).
\end{equation}
\end{cor}

\medskip

{\small {\sl \begin{Rem} One may wonder for which class of $T$-valued functions the reconstruction theorem holds. Caravenna and Zambotti proved in \cite{CaravennaZambotti} a new notion of \textbf{\textsf{germs}}. A family of distributions $(f_x)_{x\in\bbR\times\bbR^d}\subset\mcS'(\bbR\times\bbR^d)$ which is measurable in $x$ is called a germ. For any $\gamma\in\bbR$, the germ $(f_x)_x$ is called $\gamma$-\textbf{\textsf{coherent}} if there exists $\beta_0\le\gamma\wedge0$ such that
$$
\big\vert\big\langle f_y-f_x,p_t(x,\cdot)\big\rangle\big\vert\lesssim t^{\frac{\beta_0}4}\big(d(y,x)+t^{\frac14}\big)^{\gamma-\beta_0}.
$$
(Here we present a simplified definition rather than the original definition in \cite{CaravennaZambotti}.) Caravenna and Zambotti \cite{CaravennaZambotti} stated a more general reconstruction theorem at the level of coherent germs (Theorem 5.1 therein), and also stated that the coherence is actually necessary for the existence of the reconstruction (Theorem 6.1 therein). For any models $\sf M=(\Pi,g)$ and modelled distributions $\bsf\in\mcD^\gamma(T,{\sf g})$, the $x$-dependent  distributions ${\sf\Pi}_x^{\sf g}\bsf(x)$ define a coherent germ.
\end{Rem} }}

\bigskip

\hfill \textsf{\textbf{ {\Large \S} Function-like comodules}}

\medskip

Throughout we will work with regularity structures satisfying the following assumption saying that $T$ and $T^+$ contain the polynomial regularity 
structure.

\medskip

\begin{assumA}\label{A1}
The concrete regularity structure $\big((T^+,\Delta^+),(T,\Delta)\big)$ contains the polynomial regularity structure in the following sense. One has $\mcB_\alpha^{(+)} = \big\{X^n_{(+)}\,;\,\vert 
n\vert_{\frak{s}}=\alpha\big\}$, for a symbol $X_{(+)}\in T^{(+)}$ and any integer $\alpha\in\bbN$, and 
$$
\Delta^{(+)} X_{(+)}^n = \sum_{\ell\leq n} \binom{n}{\ell} X_{(+)}^\ell\otimes \Xplus^{n-\ell}.
$$
\end{assumA}

\ssk

The notation $\Xplus^n$ allows to distinguish the elements in $\mcB^+$ and $\mcB$. Note that we always have $X_+$ in the right hand side of the above tensor product while we have $X$ or $X_+$ in the left hand side depending on whether we work on $T$ or $T^+$. Set 
$$
\mcB_X^+ \defeq  \big\{\Xplus^n\,;\,n\in \bbN\times\bbN^d\big\}, \qquad \mcB_{X} \defeq  \big\{X^n\,;\,n\in\bbN\times\bbN^d\big\},
$$ 
and write 
$$
{\bf1}_+ = \Xplus^0, \qquad {\bf1}=X^0.
$$ 
By assumption, for any {\it integer} $\alpha$, the basis $\mcB_\alpha^{(+)}$ consists {\it exactly} of the elements $X^n_{(+)}$ with $\vert n\vert_{\frak{s}}=\alpha$. Therefore these polynomials are the 
only elements of $T^{(+)}$ with integer homogeneities. In particular, $T_0^{(+)}$ is a one-dimensional vector space.

Note the use of $X_+$ in the formula for $\Delta X^n$. The space 
$$
T_X^+ \defeq  \text{span}(\mcB_X^+)
$$ 
with $\Delta^+$ is isomorphic to a polynomial regularity structure, while 
the space 
$$
T_X \defeq  \text{span}(\mcB_X)
$$ 
with $\Delta$ is a right comodule over $T_X^+$. One defines a {\it canonical model} over the polynomial regularity structure 
$$
\mathscr{T}_X \defeq  \big((T_X^+,\Delta^+),(T_{X},\Delta)\big)
$$ 
setting for all $x,y\in\bbR\times\bbR^d$
$$
{\sf g}_x(\Xplus^n)\defeq x^n,\qquad
({\sf \Pi}X^n)(y) \defeq  y^n.
$$
We see that ${\sf g}_{yx}(\Xplus^n)=(y-x)^n$ and $({\sf\Pi}_x^{\sf g}X^n)(y)=(y-x)^n$, so $({\sf g}, {\sf \Pi})$ is indeed a model over $\mathscr{T}_X$.

\medskip

\begin{assumA}\label{A2}
Under Assumption \refA{A1}, we only consider models $({\sf\Pi},{\sf g})$ whose restriction to $\mathscr{T}_X$ is the canonical model.
\end{assumA}

\medskip

We only work from now on with regularity structures satisfying Assumptions \refA{A1} and \refA{A2}. It is useful, to deal with sub-regularity structures of a given regularity structure, to introduce the following notion. A linear subspace $V$ of $T$ is called a \textsf{\textbf{subcomodule}} if
$$
\Delta V\subset 
V\otimes T^+,
$$
that is, defining $V_\alpha\defeq V\cap T_\alpha$, the pair $\big((T^+,\Delta^+),(V,\Delta)\big)$ is a regularity structure. A subcomodule $V$ is said to be \textbf{\textsf{function-like}}, if $V$ satisfies Assumptions \refA{A1} and \refA{A2} and if $V_\beta=0$ whenever $\beta<0$. Given a subcomodule $V$, set
$$
\alpha_0(V) \defeq  \min\Big\{\alpha\in A\,;\, V_\alpha\neq0,\, \alpha\notin\bbN\Big\}
$$
if there is $\alpha\notin\bbN$ such that $V_\alpha\neq0$, and otherwise $\alpha_0(V)\defeq\infty$.

\medskip

\begin{cor} \label{CorEasyReconstruction}
Let $V$ be a function-like comodule. For a positive regularity exponent $\gamma$ and $\bsf\in\mcD^\gamma(V,{\sf g})$, one has $\textsf{\textbf{R}}^{\sf M}\bsf\in\mcC^{\alpha_0(V) \wedge\gamma}(\bbR\times\bbR^d)$ and for all $x\in\bbR\times\bbR^d$
$$
\big(\textsf{\textbf{R}}^{\sf M}\bsf\big)(x) = f_{\bf1}(x).
$$
\end{cor}

\medskip

\begin{Dem} 
Set $\beta=\alpha_0(V)\wedge\gamma$. To see the regularity of $f_{\bf1}$, we write
$$
\bsf(x)=\sum_{|n|_{\mathfrak{s}}<\beta}\frac{f_n(x)}{n!}X^n+\sum_{\tau\in\mcB,\,|\tau|\geq\beta}f_\tau(x)\tau,\qquad
\bsg(x)\defeq\sum_{\tau\in\mcB,\,|\tau|\geq\beta}f_\tau(x)\tau.
$$
By expanding $(\bsf(y)-\widehat{{\sf g}_{yx}}\bsf(x))_{X^n}=O(d(x,y)^{\gamma-|n|_{\mathfrak{s}}})$ for any $|n|_{\mathfrak{s}}<\beta$, we have
\begin{align*}
f_n(y)-\sum_{|\ell|_{\mathfrak{s}}<\gamma-|n|_{\mathfrak{s}}}\frac{f_{n+\ell}(x)}{\ell!}(y-x)^\ell-\ell!\big(\widehat{{\sf g}_{yx}}\,\bsg(x)\big)_{X^n}
=O\big(d(x,y)^{\gamma-|n|_{\mathfrak{s}}}\big).
\end{align*}
Since $\big(\widehat{{\sf g}_{yx}}\,\bsg(x)\big)_{X^n}=O\big(d(x,y)^{\beta-|n|_{\mathfrak{s}}}\big)$, we have
$$
f_n(y)-\sum_{|\ell|_{\mathfrak{s}}<\gamma-|n|_{\mathfrak{s}}}\frac{f_{n+\ell}(x)}{\ell!}(y-x)^\ell
=O\big(d(x,y)^{\beta-|n|_{\mathfrak{s}}}\big).
$$
This implies that $f_0=f_{\bf1}\in\mcC^\beta$ and $f_n=\partial^nf_{\bf1}$.
Therefore
\begin{align*}
\big\langle f_{\bf1}-{\sf\Pi}_x^{\sf g}\bsf(x),p_t(x,\cdot)\big\rangle
&=\bigg\langle f_{\bf1}-\sum_{|n|_{\mathfrak{s}}<\beta}\frac{f_n(x)}{n!}(\cdot-x)^n+{\sf\Pi}_x^{\sf g}\,\bsg(x),p_t(x,\cdot)\bigg\rangle\\
&=\big\langle O\big(d(x,\cdot)^\beta\big)+{\sf\Pi}_x^{\sf g}\,\bsg(x),p_t(x,\cdot)\big\rangle
=O(t^{\frac\beta4}).
\end{align*}
The uniqueness part of the proof of the reconstruction theorem, Theorem \ref{ThmReconstructionRS}, makes it clear that the reconstruction $\textbf{\textsf{R}}^{\sf M}\bsf$ of $\bsf\in\mcD^\gamma(T,\sf g)$, with $\gamma>0$, is characterized by the estimate
$$
\Big|\big\langle \textsf{\textbf{R}}^{\sf M}\bsf - {\sf\Pi}_x^{\sf g}\bsf(x) , p_t(x,\cdot)\big\rangle\Big|  \lesssim t^{\gamma'},
$$
whatever positive exponent $\gamma'$ appears in the upper bound. Hence $\textsf{\textbf{R}}^{\sf M}\bsf=f_{\bf1}$.
\end{Dem}

\bigskip

\subsection{Products and derivatives}
\label{section products and derivatives}

Other regularity structures than the polynomial regularity structure can be used to `model' functions. In good cases, they come equipped with a bilinear operation that plays the role plaid by multiplication in the usual 
setting, and allows to define the image of a modelled distribution by a nonlinear map. This is what this section is about.   

\ssk

Let $V,W$ be subcomodules of $T$ and set 
$$
V_\alpha\defeq V\cap T_\alpha,\qquad W_\alpha\defeq W\cap T_\alpha.
$$

\medskip

\begin{defn*}
A \textbf{\textsf{product}} on $V\times W$ is a continuous bilinear map $\star : V\times W\rightarrow T$, such that $V_\alpha\star W_\beta \subset 
T_{\alpha+\beta}$, for all $\alpha,\beta\in A$. 
The product is said to be \textbf{\textsf{regular}} if
$$
\Delta(\tau\star\sigma)=(\Delta\tau)(\Delta\sigma)
$$
for all $\tau\in V$ and $\sigma\in W$. In the right hand side, the product $(V\otimes T^+)\times (W\otimes T^+)\to T\otimes T^+$ is canonically defined from $\star$ and the product of $T^+$ setting 
$$
(\tau\otimes\mu)(\sigma\otimes\nu) \defeq (\tau\star\sigma)\otimes(\mu\nu).
$$   
\end{defn*}

\medskip

The regularity structures used in the study of singular PDEs have elements that are decorated rooted trees. The product is given as a tree product 
in that setting, and such a product is regular in the above sense. The details will be found in Section \ref{SectionBuildingRS}. A regular product 
$\star$ satisfies   
\begin{align}\label{section2: consequence of regularity of product}
\widehat{g}(\tau\star\sigma) = \widehat{g}(\tau)\star\widehat{g}(\sigma),
\end{align}
for any character $g$ on $T^+$. For regularity structures containing the polynomial regularity structure one asks the following consistency assumption.   

\medskip

\begin{assumA}\label{A3}
Under Assumption \refA{A1}, the product between $T_X$ and $T$ is always defined and satisfies
$$
{\bf1}\star\tau=\tau\star{\bf1}=\tau, \text{ for all } \tau\in T,\qquad X^k\star X^\ell=X^{k+\ell}, \text{ for all } k,\ell\in\bbN\times\bbN^d.
$$
\end{assumA}

\medskip

We remark that Assumption \refA{A3} is not contained in Hairer's general definition \cite{Hai14} (of course always assumed in the specific regularity structure of decorated rooted trees \cite{BHZ, ChandraHairer, BCCH18}). We make this assumption here because it is used in the proof of `Whitney extension theorem' (Theorem \ref{thm Whitney extension} in appendix) by Martin \cite{Martin}. Assumptions \refA{A1}, \refA{A2}, and \refA{A3} are jointly called Assumption \REFA. The proof of the next statement is elementary and left to the 
reader. See the proof of Theorem 4.6 in \cite{Hai14} if needed. For $\alpha\le0<\gamma$, denote by $\mcD_\alpha^\gamma$ the space of modelled distributions of the form
$$
\bsf=\sum_{\alpha\le\vert\tau\vert<\gamma}f_\tau\tau
$$
and write 
$$
\mcQ_{<\gamma} : T\to T_{<\gamma}
$$
for the canonical projection.

\medskip

\begin{prop} \label{PropRegularityProduct}
Let $\alpha_1,\alpha_2\le0<\gamma_1,\gamma_2$, and set $\gamma=(\gamma_1+\alpha_2)\wedge(\gamma_2+\alpha_1)$. Let $\star:V\times W\to T$ be a regular product. Given $\bsf_1\in\mcD_{\alpha_1}^{\gamma_1}(V,\sf g)$ and $\bsf_2\in\mcD_{\alpha_2}^{\gamma_2}(W,\sf g)$, one has
$$
\mcQ_{<\gamma}
(\bsf_1\star \bsf_2)\in\mcD_{\alpha_1+\alpha_2}^{\gamma}(T,\sf g).
$$
The mapping $(\bsf_1,\bsf_2)\mapsto \mcQ_{<\gamma}(\bsf_1\star \bsf_2)$ is continuous.
\end{prop}

\medskip

Let $V$ be a {\it function-like comodule} of $T$ equipped with an associative product
$$
\star:V\times V\to V.
$$
Then $\star$ is naturally extended to the multilinear map from $V^n$ to $V$, for any $n\geq 1$. For any $\bsf\in\mcD^\gamma(V,{\sf g})$ with $\gamma>0$ and a smooth function $F:\bbR\to\bbR$, we define
$$
F^\star(\bsf) \defeq  \mcQ_{<\gamma}\left(\sum_{n=0}^\infty \frac{F^{(n)}(f_{\bf1})}{n!}\,\overline{\bsf}^{\,\star n}\right), \qquad
\overline{\bsf} \defeq  \bsf-f_{\bf1}{\bf1}.
$$
The sum contains only finitely many terms since the sector $V$ is function-like.
Indeed, since $\overline{\bsf}=\sum_{\alpha\le|\tau|<\gamma}f_\tau\tau$ for an $\alpha>0$, we have $\overline{\bsf}^{\,\star n}=\sum_{n\alpha\le |\tau| <n\gamma}g_\tau\tau$ for some coefficients $g_\tau$. The proof of the next proposition is elementary and left to the reader; see Theorem 4.15 in \cite{Hai14} for a proof.

\medskip

\begin{prop} \label{PropFModelled}
Pick a positive regularity exponent $\gamma$. For any $\bsf\in\mcD^\gamma(V,{\sf g})$ and a smooth function $F$, one has $F^\star(\bsf)\in\mcD^\gamma(V,{\sf g})$. Moreover, the mapping $\bsf\mapsto F^\star(\bsf)$ is locally Lipschitz continuous.
\end{prop}

\medskip

Finally we introduce a linear operator playing the role of the derivative.

\medskip

\begin{defn*}
A \textbf{\textsf{derivative}} is a continuous linear map $D:T\to T$, such that $DT_\alpha\subset T_{\alpha-1}$ for all $\alpha\in A$, and
$$
\Delta(D\tau) = (D\otimes\iden)\Delta\tau
$$
for any $\tau\in T$ -- by an abuse of notation, we mean $T_{\alpha-1}=\{0\}$ if $\alpha-1\notin A$.
\end{defn*}

The assumption on $D$ implies
$$
\widehat{g}(D\tau)=D\widehat{g}(\tau)
$$
for any character $g$ on $T^+$. From this property it is straightforward to show the following statement.

\medskip

\begin{prop}\label{prop abstract derivative}
The mapping $\mcD^\gamma(T,{\sf g})\ni \bsf\mapsto D\bsf\in \mcD^{\gamma-1}(T,{\sf g})$ is continuous.
Moreover, if ${\sf\Pi}\circ D=\mathscr{D}\circ {\sf\Pi}$, holds for a first order differential operator $\mathscr{D}$ then
$$
\textbf{\textsf{R}}^{\sf M}(D\bsf) = \mathscr{D}\big(\textbf{\textsf{R}}^{\sf M}\bsf\big)
$$
for any $\bsf\in\mcD^\gamma(T,{\sf g})$ with $\gamma>1$.
\end{prop}

\bigskip

\section{Regularity structures built from integration operators}
\label{SectionIntegration}

We describe in this section a setting where one can lift a given singular PDE into an equation set on a space of modelled distributions. The lift depends on the arbitrary choice of a model on the regularity structure and we will see in Section \ref{SectionSolvingPDEs} that the lifted equation has a unique local in time solution $\bsu=\bsu(\sf M)$ for every model. The solution of the initial singular equation will be defined as the reconstruction of this unique $\bsu$.

The regularity structures used for the study of singular stochastic PDEs have a particular structure that comes from the fixed point formulation of the (system of) PDE(s) under study. We concentrate here on the case where only one second order differential operator is involved, typically $\partial_t-\Delta_x$. (See Section \ref{SectionBuildingRS} and  Appendix {\sf \ref{SectionAppendixComments}} for comments on the general case.) We work then with regularity structures equipped with an operator $\mcI$ that plays the role of the convolution operator $(\partial_t-\Delta_x)^{-1}$, involved in the fixed point formulation of the equation under study. This operator is called an abstract integration map; it is introduced in Section \ref{SubsectionIntegrationOperators}. One associates in Section \ref{SubsectionAdmissible} to an abstract integration operator $\mcI$ a notion of admissible $\sf\Pi$-maps, which roughly means that ${\sf \Pi}(\mcI\tau) = (\partial_t-\Delta_x)^{-1}({\sf \Pi}\tau)$. For some models $\sf M=(g, \Pi)$ with $\sf \Pi$ admissible, the map $\mcI$ is (essentially) intertwined to the operator $(\partial_t-\Delta_x)^{-1}$ via $\sf \Pi$. These particular models play a crucial role in Section \ref{SubsectionLiftingK}. We construct therein a model-dependent operator $\mcK^{\sf M}$ that is (essentially) intertwined to $(\partial_t-\Delta_x)^{-1}$ via the reconstruction map ${\sf R}^{\sf M}$, and which also has a regularizing property analogue to a similar property enjoyed by $(\partial_t-\Delta_x)^{-1}$. The operator $\mcK^{\sf M}$ will be used in Section \ref{SectionSolvingPDEs} to lift the singular PDE into an equation on a space of modelled distributions. The condition that $\sf \Pi$ is admissible has some far reaching consequences and it is not obvious in the first place that one can construct a non-trivial model that is admissible. Section \ref{SubsectionAdmissibleModels} is dedicated to constructing a large class of admissible models. 

The main result of this section is Theorem \ref{thm wellposedness of K}, which gives the continuity property of the operator $\mcK^{\sf M}$. 

\medskip

A remark is in order before we set the scene in Section \ref{SubsectionOperators}. We will restrict our study to regularity structures for which the minimum homogeneity of its elements satisfies
\begin{equation} \label{EqConditionBeta0}
\beta_0 = \min A > -2.
\end{equation}
This condition ensures that elements of the form $\mcI(\tau)$ that appear in the expansion of solutions to the regularity structure lift of the considered class of singular stochastic PDE are of positive homogeneity. (So if $\mcI$ were increasing the homogeneity of all symbols by $a$ we would require $\beta_0>-a$.) While the generalized (KPZ) equation \eqref{EqGKPZ} satisfies for instance this assumption, not all singular stochastic PDEs satisfy it. This is for instance the case of the $\Phi^4_2, \Phi^4_3$ or sine-Gordon equations. This kind of equations can nonetheless be studied within the setting of regularity structures by writing their solutions as the sum of an explicit functional of the noise and a remainder term that solves an equation that can be formulated in a regularity structure satisfying condition \eqref{EqConditionBeta0} -- the so-called da Prato-Debussche trick, after similar operation was used in their work \cite{DPD}. We consider as an example the case of the $\Phi^4_3$ equation
$$
(\partial_t-\Delta_x)u = -u^3 + \zeta,
$$
set on the $3$-dimensional torus, with $\zeta$ a spacetime white noise of H\"older regularity $(-5/2)^-$. We decompose a priori the solution $u$ into $u=X+v$, where
\begin{align*}
(\partial_t-\Delta_x)X &= \zeta,   \\
(\partial_t-\Delta_x)v &= -v^3-3v^2X-3vX^2-X^3.
\end{align*}
The polynomial functions $(X^n)_{1\leq n\leq 3}$ can directly be defined as elements of $\mcC^{(-n/2)^-}$ by probabilistic means. The above equation for $v$ can then be formulated in a regularity structure with three noise symbols for $X,X^2$ and $X^3$, in which $\beta_0=(-3/2)^->-2$. The interested reader will find more details on this matter for a general class of singular stochastic PDEs in Section 5 of Bruned, Chandra, Chevyrev and Hairer's work \cite{BCCH18}.

\vfill \pagebreak

\subsection{Operators on $\bbR\times\bbR^d$}
\label{SubsectionOperators}

\ssk

We will be interested in (systems of) singular stochastic PDEs that involve possibly two types of differential operators. The derivatives $\partial_i$ in the directions of the canonical basis of $\bbR\times\bbR^d$, and the second order differential operator
$$
{\bf L} \defeq  \partial_{x_0}-\Delta_{x'}+1.
$$
Denote by ${\bf L}^{-1}$ the resolution operator associating to a Schwartz function $v\in\mcS(\bbR\times\bbR^d)$ the solution $u\in\mcS(\bbR\times\bbR^d)$ to the equation
$$
{\bf L} u = v.
$$
The strict positiveness of $-\Delta_{x'}+1$ ensures uniqueness of a solution $u\in\mcS(\bbR\times\bbR^d)$ to the preceding equation. (This is the reason why we work with ${\bf L}$ rather than with the heat operator.) The operator ${\bf L}^{-1}$ can be represented using a variant of the elliptic operators $\mcG$ introduced in Section \ref{SubsectionModelsAndCo}. Indeed, we have
\begin{align*}
\big(\partial_{x_0}-\Delta_{x'}+1\big)^{-1} &= -\big(\partial_{x_0}+\Delta_{x'}-1\big)\,\big(-\partial_{x_0}^2+(\Delta_{x'}-1)^2\big)^{-1}   \\
&= - \int_0^\infty (\partial_{x_0}+\Delta_{x'}-1)\, e^{r(\partial_{x_0}^2-(\Delta_{x'}-1)^2)} dr.
\end{align*}
Write
$$
{\bf L}^{-1} = -\int _0^\infty (\partial_{x_0}+\Delta_{x'}-1)\, e^{r(\partial_{x_0}^2-(\Delta_{x'}-1)^2)}\,dr  \eqdef \int_0^\infty K_r\,dr.
$$ 
We thus use the inhomogeneous operator
$$
\widetilde\mcG \defeq  \partial_{x_0}^2 - (\Delta_{x'}-1)^2
$$
instead of the operator $\mcG$ considered in Section \ref{SubsectionModelsAndCo}. The contents in the previous section holds similarly even if we redefine $p_t$ as the kernel of $e^{t\widetilde{\mcG}}$. Since $q_t(x,y)=-(\partial_{x_0}+\Delta_{x'}-1)p_t(x,y)$, by Proposition \ref{PropModelsOnKernels}, the kernel $q_r(x,y)$ of the operator $K_r$ satisfies the $x$-uniform bounds
\begin{align}\label{Bound zeta_x(K_n)}
\left\vert \big\langle {\sf \Pi}^{\sf g}_x\tau , \partial_x^n q_r (x,\cdot)\big\rangle \right\vert
\lesssim
r^{\frac{|\tau|-|n|_{\frak{s}}-2}{4}}
\end{align}
for any $\tau\in \mcB$, $r\in(0,1]$, and $n\in\bbN\times\bbN^d$, with the 
exponent $|n|_{\frak{s}}+2$ coming from the derivative operators $\partial_x^n$ and $(\partial_{x_0}+\Delta_{x'}-1)$ applied to $e^{r\widetilde\mcG}$. 
It is convenient, for technical purposes, we replace $q_r$ with the function
$$
\widetilde{q}_r(x,y)\defeq q_r(x,y)-P_r(\partial_x)q_1(x,y)
$$
for some $r$-dependent polynomial $P_r(\xi)=P_r(\xi_0,\xi_1,\dots,\xi_d)$ whose coefficients are bounded over $r\in[0,1]$, chosen to satisfy the property that

\begin{align}\label{eq: integral of K vanishes}
\int_{\bbR\times\bbR^d}y^n\widetilde{q}_r(x,y)dy=0
\end{align}
for any $n\in\bbN\times\bbN^d$ such that $|n|_{\mathfrak{s}}<N$ for fixed positive number $N$.
For instance, when $N=2$, we choose
$$
\widetilde{q}_r(x,y)\defeq q_r(x,y)-e^{-(r-1)}q_1(x,y).
$$
In general, by noting that Fourier transform of $Q_r(\cdot)\defeq q_r(\cdot,0)$ is given by $\widehat{Q_r}(\xi)=f(\xi)e^{rg(\xi)}$ for some polynomials $f$ and $g$, and that
$$
Q_r(\xi)-P_r(\xi)Q_1(\xi)=Q_1(\xi)\big(e^{(r-1)g(\xi)}-P_r(\xi)\big),
$$
we can choose a polynomial $P_r(\xi)$ such that $\partial_\xi^n\big(e^{(r-1)g(\xi)}-P_r(\xi)\big)\vert_{\xi=0}$ for any $|n|_{\mathfrak{s}}<N$.
By using this modified kernel, we decompose ${\bf L}^{-1}$ under the form
$$
{\bf L}^{-1} = {\bf K}+{\bf K}',
$$ 
with
\begin{align}\label{section 3: decomposition of K into K_r}
{\bf K} \defeq  \int_0^1K_r\,dr-\bigg(\int_0^1P_r(\partial)K_1dr\bigg),\qquad {\bf K}' \defeq  \int_1^\infty K_r\,dr+\bigg(\int_0^1P_r(\partial)K_1dr\bigg).
\end{align}
It is elementary to see that the operators ${\bf K}'$ maps $\mcC^\gamma(\bbR\times\bbR^d)$ into $\mcC^\infty(\bbR\times\bbR^d)$, for any regularity exponent $\gamma\in\bbR$. This decomposition is similar but different to the one $\bar{K}=K+R$ as in \cite[Lemma 5.5]{Hai14}, where $\bar{K}$ is the Green function, $K$ is a singular kernel with a bounded support, and $R$ is a smooth remainder. We concentrate on the operator ${\bf K}$ in the remainder of this section. 
Denote by 
\begin{align}\label{section 3: decomposition of K into q_r}
K(x,y) \defeq  \int_0^1 \widetilde{q}_r(x,y)\,dr
\end{align}
its kernel. The compensation of $K$ ensures that
$$
\int_{\bbR\times\bbR^d}y^nK(x,y)dy=0
$$
for any $|n|_{\mathfrak{s}}<N$ with fixed positive number $N$. This property is used in the proof of Theorem \ref{thm wellposedness of K} in this section, and the proof of Theorem \ref{ThmPropertyBPHZRenorm} later.

\medskip

\begin{prop} \label{PropSchauderEstimates}
\textsf{\textbf{(Schauder estimates for ${\bf L}^{-1}$)}} The operators ${\bf L}^{-1}$ and $\bf K$ are continuous operators from $\mcC^\gamma(\bbR\times\bbR^d)$ into $\mcC^{\gamma+2}(\bbR\times\bbR^d)$, for all non-integer regularity exponents $\gamma\in\bbR$. 
\end{prop}

\medskip

\begin{Dem}
It is sufficient to show the estimate for $\widetilde\bfK=\int_0^1K_r\,dr$. Note that $K_t=-e^{t\widetilde\mcG}(\partial_{x_0}+\Delta_{x'}-1)$, as $\widetilde\mcG$ and $(\partial_{x_0}+\Delta_{x'}-1)$ commute. We use the freedom on the choice of the ($(2,1,\dots,1)$-scaling) elliptic operator used to define the H\"older spaces, while giving equivalent norms, to work with the norm associated with the operator $\widetilde\mcG$ rather than the operator $\mcG$. We emphasize that fact by writing $\widetilde{\mcQ}^{(N)}_s$ for the operators built from $\widetilde\mcG$ in the same way as $\mcQ^{(N)}_s$ is built from $\mcG$. Given a distribution $\Lambda\in\mcC^\gamma(\bbR\times\bbR^d)$, with $\gamma\in\bbR$ non-integer, we read on the identity
$$
\widetilde{\mcQ}^{(N)}_s\big(\widetilde{\bf K}(\Lambda)\big) = \int_0^1 
\widetilde{\mcQ}^{(N)}_s\big(K_r(\Lambda)\big)dr 
=
-\int_0^1 \Big(\frac{s}{s+r}\Big)^N \widetilde{\mcQ}^{(N)}_{s+r}
\Big((\partial_{x_0}+\Delta_{x'}-1)\Lambda\Big)\,
dr,
$$
the estimate
$$
\Big\|Q_s^{(N)}\big(\widetilde{\bf K}(\Lambda)\big)\Big\|_\infty \lesssim 
s^{\frac{\gamma-2}{4}+1} + O(s^N). 
$$
The result follows for all $\gamma+2<4N$. The equivalence of the different H\"older norms corresponding to different choices of $N$ gives the conclusion.
\end{Dem}

\medskip

(We refer the reader to Section 14.3 of the second edition of Friz \& Hairer's lecture notes \cite{FrizHairerBook} for a particularly nice proof of the classical Schauder estimates using different tools.) For a regularity structure $\mathscr{T}$ for which $\beta_0=\min A > -2$, and a model $({\sf\Pi},{\sf 
g})$ on it, Schauder estimates imply in particular that all the distributions $\textbf{K}({\sf \Pi}\tau)(\cdot)$, hence all the distributions $\textbf{K}({\sf \Pi}^{\sf g}_x\tau)$, are actually defined pointwise, for any $x\in\bbR\times\bbR^d$, making sense of $\textbf{K}({\sf \Pi}^{\sf g}_x\tau)(x)$, or even $\partial^n\textbf{K}({\sf \Pi}^{\sf g}_x\tau)(x)$, for $\vert n\vert_\frak{s}<\beta_0+2$. The following lemma allows to take profit from the fact that ${\sf \Pi}^{\sf g}_x\tau$ behaves near $x$ ``as'' an element of $\mcC^{\vert\tau\vert}(\bbR\times\bbR^d)$, to give meaning to $\partial^n\textbf{K}({\sf \Pi}^{\sf g}_x\tau)(x)$, for all multiindices $n$ such that $\vert n\vert_\frak{s}<\vert\tau\vert+2$.

\medskip

\begin{lem} \label{LemmeDefnIDistribution}
Assume $\beta_0>-2$. Given $\tau\in \mcB$ and $n\in\bbN\times\bbN^d$, the 
integral
\begin{equation} \label{EqDefnIDistribution}
\big(\partial^n{\bf K}({\sf \Pi}^{\sf g}_x\tau)\big)(x) \defeq  \big\langle {\sf \Pi}^{\sf g}_x\tau , \partial_x^nK(x,\cdot)\big\rangle 
\defeq   
\int_0^1 \big\langle {\sf \Pi}^{\sf g}_x\tau , \partial_x^n \widetilde{q}_r (x,\cdot)\big\rangle\,dr
\end{equation}
converges for all $x\in\bbR\times\bbR^d$, provided $|n|_\frak{s}<\vert\tau\vert + 2$.
\end{lem}

\medskip

\begin{Dem}
It follows from \eqref{Bound zeta_x(K_n)} that the first term of the right hand side of \eqref{EqDefnIDistribution} is integrable over $r\in(0,1)$ 
if
$$
|\tau|-|n|_{\frak{s}} > -2.
$$
\end{Dem}

\medskip

The $\bbR\times\bbR^d$-indexed distributions $\textbf{\textsf{R}}^{\sf M}\bsf - {\sf \Pi}^{\sf g}_x\bsf(x)$ satisfy a similar bound to \eqref{Bound zeta_x(K_n)} for any modelled distribution $\bsf\in\mcD^\gamma(T,\sf g)$. We can then define properly $\partial^n{\bf K}\big(\textbf{\textsf{R}}^{\sf M}\bsf - {\sf \Pi}^{\sf g}_x\bsf(x)\big)(x)$, for all multiindices $n$ such that $\vert n\vert_\frak{s}<\gamma+2$, as in the preceding lemma.

\bigskip

\subsection{Regularity structures with some abstract integration operators}
\label{SubsectionIntegrationOperators}

Recall Assumption {\REFA} essentially says that we consider regularity structures containing the canonical polynomial structure and models that behave naturally on the latter. In this section, we consider a regularity structure that can represent more specific functions/distributions. To explain the motivation of the next assumption, pick a $\gamma$-H\"older function $f$ with some $\gamma>0$ and let us try to add a new symbol $F$ with homogeneity $\gamma$ which represents $f$ to a polynomial regularity structure $\mathscr{T}_X$. The new model space is $T=\text{span}(\mcB_X\cup\{F\})$.
On the other hand, $T^+$ is required to be rich enough to define the application of the model to $F$ as
\begin{equation}\label{eq:motivativeAssumpB}
{\sf \Pi}F=f,\qquad
({\sf\Pi}_x^{\sf g}F)(\cdot)=f(\cdot)-\sum_{|n|_{\mathfrak{s}}<\gamma}\frac{(\cdot-x)^n}{n!}\partial^nf(x)
=O\big(d(\cdot,x)^\gamma\big).
\end{equation}
To represent coefficients $\partial^nf$, it would be natural to introduce symbols $\{F_n\}_{|n|_{\mathfrak{s}}<\gamma}$ in $T^+$. Then a good definition of $\Delta$ is given by
$$
\Delta F=F\otimes{\bf1}+\sum_{|n|_{\mathfrak{s}}<\gamma}\frac{X^n}{n!}\otimes F_n.
$$
Indeed, by applying ${\sf\Pi}_x^{\sf g}\otimes{\sf g}_x$ on both sides and noting that ${\sf \Pi}=({\sf\Pi}_x^{\sf g}\otimes{\sf g}_x)\Delta$, we have the second equality of \eqref{eq:motivativeAssumpB} under the choice ${\sf g}_x(F_n)=\partial^nf(x)$.
Similarly, we also need to define the application of $\Delta^+$ to each $F_n$ by
$$
\Delta^+ F_n=F_n\otimes{\bf1}+\sum_{|m|_{\mathfrak{s}}<\gamma-|n|_{\mathfrak{s}}}\frac{X^m}{m!}\otimes F_{n+m}
$$
to define ${\sf g}_{yx}(F_n)$ as a remainder term of the Taylor expansion of $\partial^nf$. For this reason, it would be natural to consider the regularity structure that satisfy in addition to Assumption {\REFA} the following set of assumptions. Recall we denote by $\{e_0,e_1,\dots,e_d\}$ the canonical basis of $\bbN\times\bbN^d$.

\medskip

\begin{assumB}\label{B1}
\begin{enumerate}
\setlength{\itemsep}{0.1cm}
   \item The basis $\mcB^+$ of $T^+$ is a commutative monoid with unit ${\bf 1}_+$, freely generated by the symbols
$$
\{\Xplus^{e_i}\}_{0\leq i\leq d} \cup \big\{\mcI_n^+\tau\big\}_{\tau\in\mcB,\, n\in\bbN\times\bbN^d,\,|\tau|+2-|n|_{\frak{s}}>0},
$$
Each element has homogeneity
$$
\big\vert \Xplus^{e_i} \big\vert \defeq  \frak{s}_i,\quad
\big| \mcI^+_n\tau\big| \defeq  |\tau|+2-|n|_\frak{s}.
$$
The operators $\Delta$ and $\Delta^+$ are related by the \emph{intertwining relations}  
   \begin{equation}
   \label{EqIntertwiningDeltaDelta+}
   \copro(\mcI_n^+\tau) = (\mcI_n^+\otimes\iden)\comul\tau + \sum_{\ell\in\bbN\times\bbN^d,\, |\ell|_{\frak{s}}<|\tau|+2-|n|_{\frak{s}}} \frac{\Xplus^\ell}{\ell!}\otimes \mcI_{n+\ell}^+ \tau, 
\end{equation}   
for any $\tau\in\mcB$.   \vspace{0.15cm}

   \item For $n\in\{0,e_1,\dots,e_d\}$ -- i.e., $n\in\bbN\times\bbN^d$ such that $|n|_{\frak{s}}\le1$, there are operators $\mcI_n: T\to T$, with 
   $$
   \big|\mcI_n\tau\big| \defeq  |\tau|+2-|n|_{\frak{s}}, \quad \tau\in
   \mcB.
   $$
One has for any $\tau\in \mcB$
   \begin{equation}
   \label{EqDefnDeltaIkTau}
   \Delta(\mcI_n\tau) = (\mcI_n\otimes\iden)\Delta\tau + \sum_{\ell\in\bbN\times\bbN^d,\, |\ell|_{\frak{s}}<|\tau|+2-|n|_{\frak{s}}} 
   \frac{X^\ell}{\ell!}\otimes \mcI_{n+\ell}^+\tau.
   \end{equation}
\end{enumerate}
For simplicity we write
$$
\mcI^{(+)}\defeq\mcI_0^{(+)}.
$$
\end{assumB}

\medskip

The operator $\mcI_n$ is an abstract version of the convolution operator $\partial^n{\bf K}$. The restriction $|n|_{\frak{s}}\le1$ on $n$ means that we only consider ${\bf K}$ or $\partial_i{\bf K}$; this is sufficient for the study of all (systems of) singular stochastic PDEs whose solutions are functions involving second order differential operators satisfying the above classical Schauder estimate. Two remarks are in order.

\begin{enumerate}
	\item The main point of Assumption \refB{B1} is the introduction of the operator $\mcI_n$. It brings the supplementary operator $\mcI_n^+$ for the reasons mentioned before. The first term of the right hand side of the identity \eqref{EqDefnDeltaIkTau} seems different from the definition of $\Delta F$ mentioned before, but it specifies the action of the recentering operations on $\tau$. Indeed, by applying ${\sf\Pi}_x^{\sf g}\otimes{\sf g}_x$ on both sides and noting that $\widehat{{\sf g}_x}(\tau) = \tau+\sum_{\sigma<\tau}{\sf g}_x(\tau/\sigma)\sigma$, we have
\begin{align*}
{\sf\Pi}_x^{\sf g}(\mcI_n\tau) = {\sf\Pi}\mcI_n\tau-\sum_{\sigma<\tau}{\sf g}_x(\tau/\sigma){\sf\Pi}_x^{\sf g}\mcI_n\sigma-\sum_{\ell}{\sf g}_x(\mcI_{n+\ell}^+\tau)\frac{(\cdot-x)^n}{n!}.
\end{align*}  

	\item We assume that $T^+$ is entirely constructed from the $\mcI^+_n$ operators and the polynomials, and has no other elements. This is because these elements are rich enough to build the regularity structure including polynomials, operators $\mcI_n$, and their (possible) products. The regularity structures that are used for the study of singular stochastic PDEs have the structure as above that is described in details in Section \ref{SectionBuildingRS}.
\end{enumerate}

\smallskip

We notice here some algebraic formulas about $\Delta$ and $\Delta^+$. We let the readers check that identity \eqref{EqIntertwiningDeltaDelta+} ensures the co associativities 
\begin{align*}
(\Delta^{(+)}\otimes\iden)\Delta^{(+)}(\mcI_n^{(+)}\tau) = (\iden\otimes\Delta^+)\Delta^{(+)}(\mcI_n^{(+)}\tau)
\end{align*}
on elements of $T^{(+)}$ of the form $\mcI_n^{(+)}\tau$. Next using identity \eqref{EqIntertwiningDeltaDelta+} we check that the antipode $S_+$ on $T^+$ satisfies the inductive relation 
\begin{equation} \label{EqInductiveS+}
S_+\big(\mcI_n^+\tau\big) = -\sum_  {\ell\in\bbN\times\bbN^d}\frac{(-X)^\ell}{\ell!}\,\mcM_+\big(\mcI^+_{n+\ell}\otimes S_+\big)\Delta\tau,
\end{equation}
where we denote by $\mcM_+$ the multiplication operator in the algebra $T^+$. Together with the relation $S_+(X_+)=-X_+$, such a formula defines indeed a unique algebra morphism. Recall from Appendix {\sf \ref{SectionAppendixAlgebra}} the defining property \eqref{EqDefnAntipode} of the antipode $S_+$ on $T^+$. As $T^+$ is a Hopf algebra by Proposition \ref{PropBialgebraHopf}, it suffices to see that
$$
\mcM_+\big(\textrm{Id}\otimes S_+\big)\Delta^+\big(\mcI^+_n\tau\big) = 0
$$
for all $n\in\bbN\times\bbN^d$ and $\tau\in T$. This relation follows from \eqref{EqInductiveS+} and \eqref{EqIntertwiningDeltaDelta+} writing
\begin{equation*} \begin{split}
\mcM_+\big(\textrm{Id}\otimes S_+\big)\Delta^+\big(\mcI^+_n\tau\big) &= 
\mcM_+\big(\mcI^+_n\otimes S_+\big)\Delta\tau + \sum_{\ell\in\bbN\times\bbN^d}\frac{X^\ell}{\ell!}\,S_+\big(\mcI^+_{n+\ell}\tau\big)  \\
&= \Big\{\mcM_+\big(\mcI^+_n\otimes S_+\big) -  \sum_{\ell,k\in\bbN\times\bbN^d}\frac{X^\ell}{\ell!}\,\frac{(-X)^k}{k!}\,\mcM_+\big(\mcI^+_{n+\ell+k}\otimes S_+\big)\Big\}\Delta\tau   \\
&= \Big\{\mcM_+\big(\mcI^+_n\otimes S_+\big) -  \sum_{\ell\in\bbN\times\bbN^d}\frac{(X-X)^\ell}{\ell!}\,\mcM_+\big(\mcI^+_{n+\ell}\otimes S_+\big)\Big\}\Delta\tau   \\
&= 0.
\end{split} \end{equation*}

\smallskip

Remark that the image by the operator $\mcI$ of a modelled distribution is not a modelled distribution. The next two sections are dedicated to 
constructing a model-dependent map $\mcK^{\sf M}$ that maps continuously all $\mcD^\gamma(T,\sf g)$ into $\mcD^{\gamma+2}(T,\sf g)$ when $\gamma$ is a positive non-integer real number -- the analogue of part of Schauder 
estimates, and is intertwined to the convolution operator $\bf K$
$$
{\bf K}\circ\textbf{\textsf{R}}^{\sf M} = \textbf{\textsf{R}}^{\sf M}\circ \mcK^{\sf M}
$$
via the reconstruction operator $\textbf{\textsf{R}}^{\sf M}$ associated with $\sf M$. We say that $\mcK^{\sf M}$ is a {\it lift of} $\bf K$. The construction of this operator requires the introduction of the notion of admissible model.

\bigskip

\subsection{Admissible models}
\label{SubsectionAdmissible}

In this section we consider only the operator $\mcI=\mcI_0:T\to T$. The following notion plays a key role in the proof of the existence of a lift of the convolution operator $\bf K$. Recall that ${\bf K}(\zeta)$ is well-defined pointwise for any distribution $\zeta\in\mcC^\beta(\bbR\times\bbR^d)$ with $\beta>-2$, by Lemma \ref{LemmeDefnIDistribution}. The following definition is directly from \cite[Definition 6.9]{BHZ}.

\medskip

\begin{defn*}
Let $\mathscr{T}$ be a regularity structure satisfying Assumptions \refA{A3} and \refB{B1}. Assume $\beta_0>-2$. A ${\sf \Pi}$-map on $T$ is said to be \textbf{\textsf{${\bf K}$-admissible}} if it satisfies 
\begin{equation} \label{EqConditionAdmissibility} 
{\sf \Pi}(\mcI\tau) = {\bf K}({\sf \Pi}\tau),
\end{equation}
and
\begin{equation} \label{EqConditionAdmissibility2} 
{\sf \Pi}(X^n\star\tau)(x) = x^n({\sf \Pi}\tau)(x),
\end{equation}
for any $\tau\in\mcB$ and $n\in\bbN\times\bbN^d$. A model $({\sf\Pi},{\sf 
g})$ on $\mathscr{T}$ is said to be ${\bf K}$-admissible if its $\sf \Pi$-map is ${\bf K}$-admissible.
\end{defn*}

\medskip

The notion of admissible $\sf \Pi$-map gives flesh to the idea that the operator $\mcI$ is the regularity structure counterpart of the convolution 
operator ${\bf K}$. The importance of the notion of $\bf K$-admissible model comes from Theorem \ref{thm wellposedness of K} in the next section. It shows that when working with $\bf K$-admissible models $\sf M$, one can upgrade the intertwining relation  \eqref{EqConditionAdmissibility} into an intertwining relation between $\bf K$ and an operator $\mcK^{\sf M}$ 
on modelled distributions, via the reconstruction map $\textbf{\textsf{R}}^{\sf M}$ associated with $\sf M$. 

While for a general model as in Definition \ref{DefnModel} defining a $\sf g$-map satisfying the constraint \eqref{EqEstimateGammayx} is decorrelated from the task of defining a $\sf \Pi$-map satisfying the constraint \eqref{EqestimatePix} it will turn out that imposing the intertwining relation \eqref{EqConditionAdmissibility} on $\sf \Pi$ will constrain strongly $\sf g$. Unlike general models, admissible models on a regularity structure satisfying Assumptions {\REFA} and \refB{B1} will turn out to be partly defined by their $\sf\Pi$ map; this will be proved in Proposition \ref{PropLinkGPiForAdmissibleModels}. The $\sf g$-map of an admissible model on a regularity structure satisfying further a mild Assumption \refB{B2} will turn out to be {\it entirely} defined by its $\sf\Pi$ map.

\ssk

We worked so far with models that are not constrained by anything else than their defining properties \eqref{EqEstimateGammayx} and \eqref{EqestimatePix} and it is not clear that one can further impose additional conditions like \eqref{EqConditionAdmissibility}. We will construct in Section \ref{SubsectionAdmissibleModels} a whole class of admissible models with values in the set of smooth functions. This is all we need for the study of singular stochastic PDEs, as the nonsmooth admissible models involved in this setting are limits of smooth admissible models, and limits of admissible models are admissible. As for now, we keep going and see what can 
be done with admissible models.

\medskip

Recall from Lemma \ref{LemmeDefnIDistribution} the definition of $\partial^n {\bf K}({\sf\Pi}^{\sf g}_x\tau)(x)$, for any $x\in\bbR\times\bbR^d$ 
and $n\in\bbN\times\bbN^d$ such that $\vert n\vert_\frak{s}<\vert\tau\vert+2$, and define the model-dependent polynomial-valued function on $T$
$$
\mcJ^{\sf M}(x)\tau \defeq  \sum_{|n|_\frak{s}<|\tau|+2} \frac{X^n}{n!} \,\partial^n {\bf K}({\sf\Pi}^{\sf g}_x\tau)(x) \in T_X,
$$
for any $\tau\in\mcB$ and $x\in\bbR\times\bbR^d$.

\medskip

\begin{prop} \label{PropUsefulJx}
For a $\bf K$-admissible model $\sf M$ on $\mathscr{T}$ one has, for any $x\in\bbR\times\bbR^d$ and $\tau\in T$,
$$
{\sf \Pi}^{\sf g}_x\big(\mcI\tau + \mcJ^{\sf M}(x)\tau\big) = {\bf K}({\sf \Pi}^{\sf g}_x\tau)
$$
\end{prop}

\medskip

\begin{Dem} 
Using identity \eqref{EqDefnDeltaIkTau} and the admissibility of the model, one has indeed
\begin{equation*} \begin{split}
{\sf \Pi}^{\sf g}_x(\mcI\tau) = \big({\sf \Pi}\otimes {\sf g}_x^{-1}\big)(\Delta\mcI\tau) &= \sum_{\sigma\leq\tau} \big({\sf \Pi}\otimes {\sf g}_x^{-1}\big) (\mcI\sigma\otimes\tau/\sigma) + \sum_{\vert\ell\vert_{\frak{s}}<\vert\tau\vert+2} \frac{(\cdot)^\ell}{\ell!}\,{\sf g}_x^{-1}(\mcI^+_\ell\tau)   \\
&\eqdef \sum_{\sigma\leq\tau} {\bf K}({\sf \Pi}\sigma){\sf g}_x^{-1}(\tau/\sigma) + P_x(\tau, \cdot) = {\bf K}\big({\sf \Pi}_x^{\sf g}\tau\big) + P_x(\tau, \cdot),
\end{split} \end{equation*}
so ${\sf \Pi}^{\sf g}_x(\mcI\tau)$ and ${\bf K}\big({\sf \Pi}_x^{\sf g}\tau\big)$ differ by a polynomial $P_x(\tau, \cdot)$ of degree at most $\vert\tau\vert+2$. We identify this polynomial with $-{\sf \Pi}_x^{\sf g}\big(\mcJ^{\sf M}(x)\tau\big)$ noting that ${\sf \Pi}^{\sf g}_x(\mcI\tau)$ has null derivatives at $x$ up to the order the integer part of $\vert\tau\vert+2$, from condition \eqref{EqestimatePix}.
\end{Dem}

\medskip

\begin{cor} \label{CorIdentityAdmissibility}
For a $\bf K$-admissible model $\sf M$ on $\mathscr{T}$ one has, for any $x,y\in\bbR\times\bbR^d$, 
\begin{equation} \label{EqCommutationRelation}
\widehat{{\sf g}_{yx}}\big(\mcI + \mcJ^{\sf M}(x)\big) = \big(\mcI + \mcJ^{\sf M}(y)\big)\widehat{{\sf g}_{yx}}.
\end{equation}
\end{cor}

\medskip

\begin{Dem}
Given $\tau\in T$ and $x,y\in\bbR\times\bbR^d$, one has both $\big(\widehat{{\sf g}_{yx}}\big(\mcI + \mcJ^{\sf M}(x)\big) - \big(\mcI + \mcJ^{\sf M}(y)\big)\widehat{{\sf g}_{yx}}\big)\tau\in T_X$, and 
\begin{equation*} \begin{split}
{\sf \Pi}^{\sf g}_y\Big(\widehat{{\sf g}_{yx}}\big(\mcI + \mcJ^{\sf M}(x)\big)\tau - \big(\mcI + \mcJ^{\sf M}(y)\big)\widehat{{\sf g}_{yx}}\tau\Big) &= {\sf \Pi}^{\sf g}_x\big(\mcI\tau + \mcJ^{\sf M}(x)\tau\big) - {\sf \Pi}^{\sf g}_y\Big(\big(\mcI + \mcJ^{\sf M}(y)\big)\widehat{{\sf g}_{yx}}\tau\Big)  \\
&= {\bf K}\big({\sf \Pi}_x^{\sf g}\tau\big) - {\bf K}\big({\sf \Pi}_y^{\sf g}\widehat{{\sf g}_{yx}}\tau\big) = 0,
\end{split} \end{equation*}
from Proposition \ref{PropUsefulJx}. As ${\sf \Pi}_y^{\sf g}$ is injective on $T_X$ this implies that 
$$
\widehat{{\sf g}_{yx}}\big(\mcI + \mcJ^{\sf M}(x)\big) - \big(\mcI + \mcJ^{\sf M}(y)\big)\widehat{{\sf g}_{yx}} = 0.
$$
\end{Dem}

\medskip

\begin{prop} \label{PropLinkGPiForAdmissibleModels}
Let $\mathscr{T}$ be a regularity structure satisfying Assumptions {\REFA} and \refB{B1}. The $\sf g$-map of a $\bf K$-admissible model $\sf (\Pi, g)$ on $\mathscr{T}$ satisfies
$$
{\sf g}_x(\mcI^+\tau) = {\bf K}({\sf\Pi}\tau)(x),
$$
for any $\tau\in\mcB$ and all $x\in\bbR\times\bbR^d$.
\end{prop}

\medskip

\begin{Dem}
We show that there is at most one choice of ${\sf g}(\mcI^+\tau)$ such that $(\sf \Pi, g)$ is an admissible model. Applying ${\sf\Pi}\otimes{\sf 
g}_x^{-1}$  to the identity \eqref{EqDefnDeltaIkTau} giving $\Delta(\mcI_n\tau)$, with $n=0$, one gets from the $\bf K$-admissibility of ${\sf \Pi}$
\begin{align} \label{eq Pi I tau into K Pi tau + polynomialsAlpha}
\begin{aligned}
{\sf\Pi}_x^{\sf g}(\mcI\tau) &= {\bf K}\big({\sf\Pi}_x^{{\sf g}}\tau) + 
\sum_{|\ell|_{\frak{s}}<|\tau|+2}\frac{(\cdot)^\ell}{\ell!} \, {\sf g}_x^{-1}(\mcI_\ell^+\tau)   \\
&= {\bf K}\big({\sf\Pi}_x^{{\sf g}}\tau) + \sum_{|\ell|_{\frak{s}}<|\tau|+2}\frac{(\cdot-x)^\ell}{\ell!} \, {\sf f}_x(\mcI_\ell^+\tau),
\end{aligned}
\end{align}
where ${\sf f}$ and ${\sf g}$ are related by the formulas
\begin{align*}
{\sf f}_x(\mcI_\ell^+\tau) \defeq  \sum_{ |m|_{\frak{s}}<|\tau|+2-|\ell|_{\frak{s}} }\frac{x^m}{m!}\,{\sf g}_x^{-1}(\mcI_{\ell+m}^+\tau),\quad
{\sf g}_x^{-1}(\mcI_n^+\tau)=\sum_{ |m|_{\frak{s}}<|\tau|+2-|n|_{\frak{s}} } \frac{(-x)^m}{m!}{\sf f}_x(\mcI_{n+m}^+\tau).
\end{align*}
As in the proof of Proposition \ref{PropUsefulJx}, since the derivatives of ${\sf\Pi}_x^{{\sf g}}(\mcI\tau)$ up to order $|\tau|+2$ vanish at $x$, 
we have
\begin{align*}
{\sf f}_x(\mcI_n^+\tau)=-\partial^n{\bf K}({\sf\Pi}_x^{{\sf g}}\tau)(x),
\end{align*}
hence 
\begin{align} \label{EqCanonicalChoice}
{\sf g}_x^{-1}(\mcI_n^+\tau)=-\sum_{|m|_{\frak{s}}<|\tau|+2-|n|_{\frak{s}}}\frac{(-x)^m}{m!}\partial^{n+m}{\bf K}\big({\sf\Pi}_x^{{\sf g}}\tau\big)(x).
\end{align}   
This implies another inductive formula
\begin{align} \label{eq inductive g(I_k^+tau)}
{\sf g}_x(\mcI^+_k\tau) = \sum_{\sigma\le\tau; |k|_{\frak{s}}<|\sigma|+2}\, {\sf g}_x(\tau/\sigma)\, \partial^k\bfK({\sf\Pi}^{{\sf g}}_x\sigma)(x),
\end{align}
which is proved by applying ${\sf g}_x^{-1}\otimes{\sf g}_x$ to the identity \eqref{EqIntertwiningDeltaDelta+} describing $\Delta^+(\mcI^+_k\tau)$ and using \eqref{EqCanonicalChoice}. Since $\beta_0>-2$, if $k=0$ the condition $|k|_{\frak{s}}<|\sigma|+2$ can be removed.
Hence we have ${\sf g}_x(\mcI^+\tau)=\bfK({\sf\Pi}\tau)(x)$, by the comodule identity \eqref{eq T is comodule}.  
\end{Dem}   

\medskip

Set 
$$
\mcJ^{\sf M+}(x)\tau \defeq  \sum_{|n|_\frak{s}<|\tau|+2} \frac{\Xplus^n}{n!} \,\partial^n {\bf K}({\sf\Pi}^{\sf g}_x\tau)(x) \in T_X^+\subset T^+,
$$
for any $\tau\in\mcB$ and $x\in\bbR\times\bbR^d$. The following statement 
is proved exactly as Proposition \ref{PropUsefulJx} and Corollary \ref{CorIdentityAdmissibility}; it will be used in the proof of Theorem \ref{ThmSmoothAdmissibleModels}.   

\medskip

\begin{prop} \label{PropCommutationRelationPlus}
Given a regularity structure satisfying Assumptions {\REFA} and \refB{B1} and a $\bf K$-admissible model $({\sf\Pi},{\sf g})$ on it, one has, for any $x,y\in\bbR\times\bbR^d$,
$$
{\sf g}_{yx} \Big(\mcI^+\tau + \mcJ^{\sf M+}(x)\tau\Big) = {\bf K}({\sf \Pi}^{\sf g}_x\tau)(y),
$$
and  
\begin{equation}
\label{EqCommutationRelationPlus}
\widehat{{\sf g}_{yx}}^+\Big(\mcI^+ + \mcJ^{\sf M+}(x)\Big) = \Big(\mcI^+ + \mcJ^{\sf M+}(y)\Big)\widehat{{\sf g}_{yx}}.
\end{equation}
Recall here that $\widehat{g}^+\defeq (\iden\otimes g)\Delta^+:T^+\to T^+$ denotes the action of $g\in G^+$ on $T^+$ defined in the same way as $\widehat{g}:T\to T$.
\end{prop}

\bigskip

\subsection{Lifting ${\bf K}$ as a continuous map from $\mcD^\gamma(T,\sf 
g)$ into $\mcD^{\gamma+2}(T,\sf g)$}
\label{SubsectionLiftingK}

For a given ${\bf K}$-admissible model ${\sf M}=({\sf\Pi},{\sf g})$ we define in this section a continuous map $\mcK^{\sf M}$ from $\mcD^\gamma(T,\sf g)$ into $\mcD^{\gamma+2}(T,\sf g)$, for any positive non-integer regularity exponent $\gamma$, intertwined to ${\bf K}$ via the reconstruction operator
\begin{align} \label{eq KR=RK}
{\bf K}\circ\textbf{\textsf{R}}^{\sf M} = \textbf{\textsf{R}}^{\sf M}\circ\mcK^{\sf M}.
\end{align}
To get a grasp on what $\mcK^{\sf M}$ could be one keeps from the reconstruction theorem, Theorem \ref{ThmReconstructionRS}, the image that for $\bsg\in\mcD^{\gamma+2}(T,\sf g)$ and $x\in\bbR\times\bbR^d$, the distribution $\textbf{\textsf{R}}^{\sf M}\bsg - {\sf \Pi}_x^{\sf g}\bsg(x)$ behaves near $x$ like the function $d(\cdot,x)^{\gamma+2}$. For $\bsf\in\mcD^\gamma(T,\sf g)$, since we have
$$
{\bf K}\big(\textbf{\textsf{R}}^{\sf M}\bsf\big) - {\sf\Pi}_x^{\sf g}\Big(\big(\mcI+\mcJ^{\sf M}(x)\big)\bsf(x)\Big) = {\bf K}\Big(\textbf{\textsf{R}}^{\sf M}\bsf - {\sf\Pi}_x^{\sf g}\bsf(x)\Big),
$$
from Proposition \ref{PropUsefulJx}, it then looks natural to add to $\big(\mcI+\mcJ^{\sf M}(x)\big)\bsf(x)$ the polynomial expansion 
$$
\big(\mcN^{\sf M}\bsf\big)(x) \defeq  \sum_{\vert\ell\vert_{\frak{s}} < \gamma+2} \frac{X^\ell}{\ell!}\,\big(\partial^\ell{\bf K}\big)\Big(\textbf{\textsf{R}}^{\sf M}\bsf - {\sf\Pi}_x^{\sf g}\bsf(x)\Big)(x)
$$ 
of ${\bf K}\big(\textbf{\textsf{R}}^{\sf M}\bsf - {\sf\Pi}_x^{\sf g}\bsf(x)\big)$ at point $x$, at order $\gamma+2$, and expect that 
$$
{\bf K}\big(\textbf{\textsf{R}}^{\sf M}\bsf\big) - {\sf\Pi}_x^{\sf g}\Big(\big(\mcI+\mcJ^{\sf M}(x)\big)\bsf(x) + \big(\mcN^{\sf M}\bsf\big)(x)\Big)
$$
behaves like $d(\cdot,x)^{\gamma+2}$ near $x$. (The remark after Lemma \ref{LemmeDefnIDistribution} justifies the good definition of the quantities $(\partial^\ell{\bf K})\big(\textbf{\textsf{R}}^{\sf M}\bsf - {\sf\Pi}_x^{\sf g}\bsf(x)\big)$ in $\big(\mcN^{\sf M}\bsf\big)(x)$, for $\vert\ell\vert_{\frak{s}}<\gamma+2$.) This does not guarantee that the $T$-valued map 
\begin{equation} \label{EqDefnMcK}
\big(\mcK^{\sf M}\bsf\big)(x) \defeq  \big(\mcI+\mcJ^{\sf M}(x)\big)\bsf(x) + \big(\mcN^{\sf M}\bsf\big)(x),\qquad x\in\bbR\times\bbR^d,
\end{equation}
is a modelled distribution, but this turns out to be the case! Note that unlike $\mcI$ or $\mcJ^{\sf M}(x)$, the $T_X$-valued function $\mcN^{\sf M}\bsf$ is a non-local function of $\bsf$ -- i.e. $(\mcN^{\sf M}\bsf)(x)$ is not a function of $\bsf(x)$ only. Note also that one has \emph{formally}
$$
\big(\mcK^{\sf M}\bsf\big)(x) = \mcI\bsf(x) + \sum_{\vert\ell\vert_{\frak{s}} < \gamma+2} \frac{X^\ell}{\ell!}\,\partial^\ell{\bf K}\big(\textbf{\textsf{R}}^{\sf M}\bsf\big)(x).
$$
This identity gives the intuitive meaning of the polynomial part of $\big(\mcK^{\sf M}\bsf\big)(x)$. Decomposition \eqref{EqDefnMcK} is needed to make sense of $\big(\mcK^{\sf M}\bsf\big)(x)$ in a rigorous way. We can prove the following theorem by the same way as \cite[Theorem 5.12]{Hai14} except the use of a decomposition of $K$ by integration rather than the dyadic decomposition as in \cite{Hai14}.

\medskip

\begin{thm} \label{thm wellposedness of K}
Let the regularity structure $\mathscr{T}$ satisfy Assumptions \textbf{\textsf{(A-B1)}} and $\sf (g,\Pi)$ be a $\bf K$-admissible model on it. Let $\gamma$ be a positive non-integer regularity exponent $\gamma$, and choose an integer $N$ such that $\gamma+2<N$ and that the property \eqref{eq: integral of K vanishes} holds for any $|n|_{\mathfrak{s}}<N$. Then the map $\mcK^{\sf M}$ sends continuously $\mcD^\gamma(T,\sf g)$ into $\mcD^{\gamma+2}(T,\sf g)$, and satisfies the intertwining identity \eqref{eq KR=RK}.
\end{thm}

\medskip

Before the proof, we recall from Proposition A.1 of \cite{Hai14} the (anisotropic) integral Taylor formula for the remainder.

\medskip

\begin{lem}\label{lem:anisoTaylor}
There exists a family of Borel probability measures $\{m_k\}_{k\in\bbN\times\bbN^d}$ on $[0,1]^{d+1}$ satisfying the following properties.
For any smooth functions $f$ on $\bbR\times\bbR^d$ and $\gamma>0$, one has the identity
$$
f(y)-\sum_{|k|_\mfs<\gamma}\frac{\partial^kf(x)}{k!}(y-x)^k
=\sum_{|\ell|_\mfs\ge\gamma}\frac{(y-x)^\ell}{\ell!}\int_{[0,1]^{d+1}}\partial^\ell f(x_t)m_\ell(dt),
$$
where $\ell$ runs over a finite set and $x_t\defeq (x_i+t_i(y_i-x_i))_{i=0}^d$.
\end{lem}

\medskip

\begin{Dem}[of Theorem \ref{thm wellposedness of K}] 
{\small 
We use the interwining relation \eqref{EqCommutationRelation} to write

\begin{align*}
(\mcK^{\sf M}\bsf)(y) - \widehat{{\sf g}_{yx}}(\mcK^{\sf M}\bsf)(x)
&= (\mcK^{\sf M}\bsf)(y) - \widehat{{\sf g}_{yx}} \big(\mcI+\mcJ^{\sf M}(x)\big)\bsf(x) - \widehat{{\sf g}_{yx}}(\mcN^{\sf M}\bsf)(x)   \\
&=  (\mcK^{\sf M}\bsf)(y) - \big(\mcI+\mcJ^{\sf M}(y)\big)\widehat{{\sf 
g}_{yx}}\bsf(x) - \widehat{{\sf g}_{yx}}(\mcN^{\sf M}\bsf)(x)   \\
&= \mcI\Big(\bsf(y) - \widehat{{\sf g}_{yx}}\bsf(x)\Big) \hspace{-0.05cm}+\hspace{-0.05cm} \mcJ^{\sf M}(y)\Big(\bsf(y) - \widehat{{\sf g}_{yx}}\bsf(x)\Big) \hspace{-0.05cm}+\hspace{-0.05cm} \Big((\mcN^{\sf M}\bsf)(y) - \widehat{{\sf g}_{yx}}(\mcN^{\sf M}\bsf)(x)\Big).
\end{align*}
For the $\mcI$ term, from the continuity of $\mcI$, one has the estimate
$$
\big\|\mcI\big(\bsf(y) - \widehat{{\sf g}_{yx}}\bsf(x)\big)\big\|_\beta \lesssim \big\| \bsf(y) - \widehat{{\sf g}_{yx}}\bsf(x) \big\|_{\beta-2} \leq \|\bsf\|_{\mcD^\gamma} \, d(y,x)^{\gamma+2-\beta}
$$
for any $\beta\in A$. The $\mcJ^{\sf M}$ and $\mcN^{\sf M}$ terms take values in the polynomial part $T_X$ of $T$. Decompose $K(x,y)$ into the integral of $\widetilde{q}_r$ by \eqref{section 3: decomposition of K into q_r}, and let 
$$
\admi^{\sf M} \eqdef \int_0^1 \admi^{\sf M}_rdr, \quad \textrm{and} \quad \remind^{\sf M} = :\int_0^1\remind^{\sf M}_rdr,
$$ 
stand for the corresponding operators. Since $q_1$ has a smoothing property, we replace $\widetilde{q}_r$ with $q_r$ in the following calculations. Fix $n\in\bbN\times\bbN^d$, and write $(\tau)_{X^n}$ for the component 
of $\tau\in T$ in the direction of $X^n$.
We have for 
$$
(\circleddash)_r \defeq  
\Big( \mcJ^{\sf M}_r(y)\big(\bsf(y) - \widehat{{\sf g}_{yx}}\bsf(x)\big) + \big(\mcN^{\sf M}_r\bsf)(y) - \widehat{{\sf g}_{yx}}(\mcN^{\sf M}_r\bsf)(x)\Big)_{X^n}
$$
the two decompositions
\begin{align*}
(\circleddash)_r &= \sum_{\beta\in A, |n|_{\frak{s}}<\beta+2}
\frac1{n!} \Big\langle {\sf\Pi}^{\sf g}_y \big(\bsf(y) - \widehat{{\sf g}_{yx}}\bsf(x)\big)_\beta, \partial_y^n q_r(y,\cdot) \Big\rangle\\
&\quad+ \bigg\{
\frac1{n!} \Big\langle \textbf{\textsf{R}}^{\sf M}\bsf - {\sf\Pi}_y^{\sf g}\bsf(y) , \partial_y^n q_r(y,\cdot) \Big\rangle
-\sum_{|k|_{\frak{s}}<\gamma+2-|n|_{\frak{s}}} \frac{(y-x)^k}{k!} \Big\langle \textbf{\textsf{R}}^{\sf M}\bsf - {\sf\Pi}_x^{\sf g}\bsf(x) , \partial_x^{n+k}q_r(x,\cdot) \Big\rangle \bigg\}   \\
&\eqdef (\Asterisk)^1_r + (\Asterisk)^2_r,
\end{align*}
and

\begin{align*}
(\circleddash)_r &= \sum_{\beta\in A, |n|_{\frak{s}}<\beta+2}
\frac1{n!} \Big\langle {\sf\Pi}^{\sf g}_y \big(\bsf(y) - \widehat{{\sf g}_{yx}}\bsf(x)\big)_\beta, \partial_y^n q_r(y,\cdot) \Big\rangle  \\
&\quad+ 
\frac1{n!} \Big\langle \textbf{\textsf{R}}^{\sf M}\bsf - {\sf\Pi}_x^{\sf g}\bsf(x) , (\partial^n q_r)_{y,x}^{\gamma+2-|n|_{\frak{s}}} \Big\rangle
+ \frac1{n!} \Big\langle {\sf\Pi}^{\sf g}_x \bsf(x) - {\sf\Pi}^{\sf g}_y\bsf(y) , \partial_y^n q_r(y,\cdot) \Big\rangle\\
&=
\frac1{n!} \Big\langle \textbf{\textsf{R}}^{\sf M}\bsf - {\sf\Pi}_x^{\sf g}\bsf(x) , (\partial^n q_r)_{y,x}^{\gamma+2-|n|_{\frak{s}}} \Big\rangle  
-\sum_{\beta\in A, |n|_{\frak{s}}\ge\beta+2} \frac1{n!} \Big\langle {\sf\Pi}^{\sf g}_y \big(\bsf(y) - \widehat{{\sf g}_{yx}}\bsf(x)\big)_\beta, \partial_y^n q_r(y,\cdot) \Big\rangle
\\
&\eqdef (\bigstar)_r^1+(\bigstar)_r^2,
\end{align*}
where
\begin{align*}
(\partial^n q_r)_{y,x}^{\gamma+2-|n|_{\frak{s}}}(z) \defeq  \partial_y^n q_r(y,z) - \sum_{|k|_{\frak{s}} < \gamma+2-|n|_{\frak{s}}}\frac{(y-x)^k}{k!}\,\partial_x^{n+k}q_r(x,z).
\end{align*}
Choose $r_0\in(0,1]$ such that $r_0^{\frac1{4}}\simeq d(y,x)\wedge1$. We use the $(\Asterisk)$-decomposition to estimate the integral over $0<r<r_0$, and the $(\bigstar)$-decomposition to estimate the integral over $r_0\le r\leq 1$.   \vspace{0.15cm} 

$\bullet$ For $r\in(0,r_0]$, we have from the bound \eqref{Bound zeta_x(K_n)} the estimate 

\begin{align*}
\int_0^{r_0} \big|(\Asterisk)^1_r\big| dr
&\lesssim \sum_{\beta\in A,|n|_{\frak{s}}<\beta+2} d(y,x)^{\gamma-\beta}
\int_0^{r_0}r^{\frac{\beta-|n|_{\frak{s}}-2}{4}} dr   \\
&\lesssim \sum_{\beta\in A,|n|_{\frak{s}}<\beta+2} d(y,x)^{\gamma-\beta} r_0^{\frac{\beta-|n|_{\frak{s}}+2}{4}}
\lesssim d(y,x)^{\gamma+2-|n|_{\frak{s}}}.
\end{align*}
Since $|n|_{\frak{s}}<\gamma+2$, from the bound  \eqref{EqImprovedReconstructionEstimate} in the reconstruction theorem, we get
\begin{align*}
\int_0^{r_0} \big|(\Asterisk)^2_r\big| dr
&\lesssim \int_0^{r_0} r^{\frac{\gamma-|n|_{\frak{s}}-2}{4}}\,dr
+ \sum_{|k|_{\frak{s}}<\gamma+2-|n|_{\frak{s}}} d(y,x)^{|k|_{\frak{s}}}
\int_0^{r_0} r^{\frac{\gamma-|n|_{\frak{s}}-|k|_{\frak{s}}-2}{4}}\,dr\\
&\lesssim r_0^{\frac{\gamma-|n|_{\frak{s}}+2}{4}}
+ \sum_{|k|_{\frak{s}}<\gamma+2-|n|_{\frak{s}}} d(y,x)^{|k|_{\frak{s}}}
r_0^{\frac{\gamma-|n|_{\frak{s}}-|k|_{\frak{s}}+2}{4}}
\lesssim  d(y,x)^{\gamma+2-|n|_{\frak{s}}}.
\end{align*}

$\bullet$ To deal with the integral over $r\in(r_0,1]$, we use the $(\bigstar)$-decomposition. Since this integral does not make sense if $r_0\ge1$, we assume $d(y,x)\le1$. For $(\bigstar)^1_r$, we apply Lemma \ref{lem:anisoTaylor} to write
\begin{align}\label{lem:anisoTaylorAppqr}
(\partial^n q_r)_{y,x}^{\gamma+2-|n|_{\frak{s}}}(z) =
\sum_{\gamma+2-|n|_{\frak{s}} < |\ell|_{\frak{s}} } \frac{(y-x)^\ell}{\ell!} \int_{[0,1]^{d+1}}\partial^{n+\ell} q_r(x_t,z)m_\ell(dt),
\end{align}
where {\it $\ell$ runs over a finite set}. Note that no index $n$ with $\gamma+2-|n|_{\frak{s}} = |\ell|_{\frak{s}}$ exists, because $\gamma\notin\bbZ$. By decomposing 
$$
\textbf{\textsf{R}}^{\sf M}\bsf - {\sf\Pi}_x^{\sf g}\bsf(x)=\textbf{\textsf{R}}^{\sf M}\bsf - {\sf\Pi}_{x_t}^{\sf g}\bsf(x_t)+{\sf\Pi}_{x_t}^{\sf g}\big(\bsf(x_t)-\widehat{{\sf g}_{x_t x}}\bsf(x)\big),
$$ 
and using the bounds \eqref{EqImprovedReconstructionEstimate} and \eqref{Bound zeta_x(K_n)}, we have
\begin{align*}
\int_{r_0}^1\big|(\bigstar)^1_r\big|\,dr
&\lesssim\sum_{\gamma+2-|n|_{\frak{s}} < |\ell|_{\frak{s}} }
d(y,x)^{|\ell|_{\frak{s}}}
\int_{r_0}^1\left\{r^{\frac{\gamma-|n|_{\frak{s}}-|\ell|_{\frak{s}}-2}{4}}
+\sum_{\beta\in A,\,\beta<\gamma} d(y,x)^{\gamma-\beta} r^{\frac{\beta-|n|_{\frak{s}}-|\ell|_{\frak{s}}-2}{4}}\right\} dr   \\
&\lesssim\sum_{\gamma+2-|n|_{\frak{s}} < |\ell|_{\frak{s}} }
d(y,x)^{|\ell|_{\frak{s}}}
\left\{r_0^{\frac{\gamma-|n|_{\frak{s}}-|\ell|_{\frak{s}}+2}{4}}
+\sum_{\beta\in A,\,\beta<\gamma} d(y,x)^{\gamma-\beta} r_0^{\frac{\beta-|n|_{\frak{s}}-|\ell|_{\frak{s}}+2}{4}}\right\}   \\
&\lesssim d(y,x)^{\gamma+2-|n|_{\frak{s}}}.
\end{align*}
We obtain the same bound for the $(\bigstar)^2_r$-term by a similar argument. Note that the terms with integer $\beta$ can be excluded. Indeed, the only elements of $\mcB$ with integer homogeneity are the polynomials, and $\langle{\sf\Pi}_y^{\sf g}X^k,\partial_y^nq_r(y,\cdot)\rangle=0$ for any $|k|_{\mathfrak{s}}<N$ by the property \eqref{eq: integral of K vanishes}. Therefore
\begin{align*}
\int_{r_0}^1\big|(\bigstar)^2_r\big|\,dr
&\lesssim\sum_{\beta\in A,\, |n|_{\frak{s}} >\beta+2 }
d(y,x)^{\gamma-\beta}
\int_{r_0}^1r^{\frac{\beta-|n|_{\frak{s}}-2}4}dr   \\
&\lesssim\sum_{\beta\in A,\, |n|_{\frak{s}} >\beta+2 }
d(y,x)^{\gamma-\beta}r_0^{\frac{\beta-|n|_{\frak{s}}+2}4}   \lesssim d(y,x)^{\gamma+2-|n|_{\frak{s}}}.
\end{align*}

To show the intertwining identity \eqref{eq KR=RK}, it is sufficient to obtain the estimate
\begin{align}\label{section3: KR=RK}
\bigg\vert\bigg\langle{\bf K}\big(\textbf{\textsf{R}}^{\sf M}\bsf-{\sf\Pi}_x^{\sf g}\bsf(x)\big)
-\sum_{|n|_{\frak{s}}<\gamma+2}\frac{(\cdot-x)^n}{n!}\big(\partial^n{\bf K}\big)\big(\textbf{\textsf{R}}^{\sf M}\bsf - {\sf\Pi}_x^{\sf g}\bsf(x)\big)(x),p_t(x,\cdot)\bigg\rangle\bigg\vert
\lesssim t^{\frac{\gamma+2}4}
\end{align}
for any $\bsf\in\mcD^\gamma(F,{\sf g})$.
Then the uniqueness of the reconstruction operator gives ${\bf K}(\textbf{\textsf{R}}^{\sf M}\bsf) = \textbf{\textsf{R}}^{\sf M}(\mcK^{\sf M}\bsf)$. To prove \eqref{section3: KR=RK}, we write
\begin{align*}
{\bf K}\big(\textbf{\textsf{R}}^{\sf M}\bsf\big)-{\sf\Pi}_x^{\sf g}(\mcK^{\sf M}\bsf)(x)
&={\bf K}\big(\textbf{\textsf{R}}^{\sf M}\bsf\big)-{\sf\Pi}_x^{\sf g}\big(\mcI+\mcJ^{\sf M}(x)\big)\bsf(x)-{\sf\Pi}_x^{\sf g}(\mcN^{\sf M}\bsf)(x)\\
&={\bf K}\big(\textbf{\textsf{R}}^{\sf M}\bsf\big)-{\bf K}\big({\sf\Pi}_x^{\sf g}\bsf(x)\big)
-{\sf\Pi}_x^{\sf g}(\mcN^{\sf M}\bsf)(x)\\
&={\bf K}\big(\textbf{\textsf{R}}^{\sf M}\bsf-{\sf\Pi}_x^{\sf g}\bsf(x)\big)
-\sum_{|n|_{\frak{s}}<\gamma+2}\frac{(\cdot-x)^n}{n!}\big(\partial^n{\bf K}\big)\big(\textbf{\textsf{R}}^{\sf M}\bsf - {\sf\Pi}_x^{\sf g}\bsf(x)\big)(x).
\end{align*}
We decompose ${\bf K}$ into the integral of $K_r$ over $r\in[0,1]$ by \eqref{section 3: decomposition of K into q_r} and ignore $P(\partial)K_1$ as it has a smoothing property. Since $\int_{\bbR\times\bbR^d}q_r(y,\cdot)p_t(x,y)dy=q_{r+t}(x,\cdot)$ by definition, we have
\begin{align*}
&\bigg\vert\bigg\langle K_r\big(\textbf{\textsf{R}}^{\sf M}\bsf-{\sf\Pi}_x^{\sf g}\bsf(x)\big)
-\sum_{|n|_{\frak{s}}<\gamma+2}\frac{(\cdot-x)^n}{n!}\big(\partial^nK_r\big)\big(\textbf{\textsf{R}}^{\sf M}\bsf - {\sf\Pi}_x^{\sf g}\bsf(x)\big)(x),p_t(x,\cdot)\bigg\rangle\bigg\vert\\
&\lesssim \big\vert K_{r+t}\big(\textbf{\textsf{R}}^{\sf M}\bsf-{\sf\Pi}_x^{\sf g}\bsf(x)\big)(x)\big\vert
+\sum_{|n|_{\frak{s}}<\gamma+2}t^{\frac{|n|_{\frak{s}}}4}\big\vert \big(\partial^nK_r\big)\big(\textbf{\textsf{R}}^{\sf M}\bsf - {\sf\Pi}_x^{\sf g}\bsf(x)\big)(x)\big\vert\\
&\lesssim (r+t)^{\frac{\gamma-2}4}+\sum_{|n|_{\frak{s}}<\gamma+2}t^{\frac{|n|_{\frak{s}}}4}r^{\frac{\gamma-|n|_{\frak{s}}-2}4}.
\end{align*}
Since $\frac{\gamma-|n|_{\frak{s}}-2}4>-1$, the integral over $r\in[0,t]$ gives the upper bound $t^{\frac{\gamma+2}4}$. For the integral over $r\in[t,1]$, we use the representation $(\bigstar)^1_r$ wth $n=0$. By the estimate of $(\bigstar)^1_r$ obtained before, we have
\begin{align*}
&\big\vert\big\langle (\bigstar)^1_r(\cdot),p_t(x,\cdot)\big\rangle\big\vert\\
&\lesssim\sum_{\gamma+2< |\ell|_{\frak{s}} }
\int_{\bbR\times\bbR^d}d(y,x)^{|\ell|_{\frak{s}}}\left\{r^{\frac{\gamma-|\ell|_{\frak{s}}-2}{4}}
+\sum_{\beta\in A,\,\beta<\gamma} d(y,x)^{\gamma-\beta} r^{\frac{\beta-|\ell|_{\frak{s}}-2}{4}}\right\} p_t(x,y)dy\\
&\lesssim\sum_{\gamma+2< |\ell|_{\frak{s}} }
\left\{t^{\frac{|\ell|_{\frak{s}}}4}r^{\frac{\gamma-|\ell|_{\frak{s}}-2}{4}}
+\sum_{\beta\in A,\,\beta<\gamma} t^{\frac{\gamma-\beta+|\ell|_{\frak{s}}}4} r^{\frac{\beta-|\ell|_{\frak{s}}-2}{4}}\right\}.
\end{align*}
Then the integral over $r\in[t,1]$ gives the upper bound $t^{\frac{\gamma+2}4}$.
}
\end{Dem}

\medskip

Note that the intertwining relation \eqref{eq KR=RK} between $\bf K$ and $\mcK^{\sf M}$ provides indeed an `upgraded' version of the defining identity \eqref{EqConditionAdmissibility} for a $\bf K$-admissible model in 
so far as the former reduces to the latter when applied to the modelled distribution $\bsf(x) = \bsh^{\tau}(x) = \sum_{\sigma<\tau}{\sf g}_x(\tau/\sigma)\sigma$. Indeed, on the one hand we have $\textbf{\textsf{R}}^{\sf M}\bsh^{\tau} = {\sf \Pi}\tau$. On the other hand $\mcK^{\sf M}\bsh^\tau$ has positive regularity and takes its values in a function-like sector when the model takes values in the space of continuous functions, so Corollary \ref{Cor reconstruction of smooth model} applies and identifies the reconstruction of $\mcK^{\sf M}\bsh^\tau$ as ${\sf \Pi}^{\sf g}_x\big(\mcK^{\sf M}\bsh^\tau(x)\big)(x)$, equal to ${\sf \Pi}(\mcI\tau)$, as all the $x$-indexed polynomial terms are null  when evaluated at $x$. 

\bigskip

\subsection{Building admissible smooth models}
\label{SubsectionAdmissibleModels}

We left aside in Section \ref{SubsectionAdmissible} the non-elementary question of existence of non-trivial admissible models to concentrate on their properties. We construct in this section a large class of admissible models for which all $({\sf \Pi}\tau)_{\tau\in T}$ and $\big({\sf g}(\sigma)\big)_{\sigma\in T^+}$ are smooth functions. In applications to singular stochastic PDEs such models can be built from realizations of the noise(s) in the equation.

\ssk

Recall Assumption \textbf{\textsf{(B1)}} in Section \ref{SubsectionIntegrationOperators} describes the action of the `recentering operator' $\Delta$ on elements of $T$ of the form $\mcI_k\tau$. We single out for our needs an assumption on $\Delta$ that provides a crucial induction structure.

\medskip

\begin{assumB}\label{B2}
There exists an increasing sequence $\{\mcB^{(n)}\}_{n=0}^\infty$ of subsets of $\mcB$ such that $\mcB^{(0)}=\{X^k\}_{k\in\bbN\times\bbN^d}$, $\mcB=\bigcup_{n=0}^\infty\mcB^{(n)}$, and
$$
\Delta\tau-\tau\otimes{\bf1}\in T^{(n-1)}\otimes T^{(n-1)+}
$$
for any $n$ and $\tau\in\mcB^{(n)}$, where $T^{(n-1)}$ is the vector space spanned by $\mcB^{(n-1)}$, and $T^{(n-1)+}$ is the subalgebra of $T^+$ generated by the symbols
$$
\{\Xplus^{e_i}\}_{0\leq i\leq d} \cup \big\{\mcI_k^+\sigma\big\}_{\sigma\in \mcB^{(n-1)},\, k\in\bbN\times\bbN^d}.
$$
Here $\mcB^{(-1)}\defeq\emptyset$.
\end{assumB}

\medskip

Assumptions \refB{B1} and \refB{B2} are jointly called Assumption \REFB. Assumption \refB{B2} is satisfied by the regularity structures of decorated trees as defined in Section \ref{SectionBuildingRS}. Indeed, for this case, we can define $\mcB^{(n)}$ as the set of decorated trees $\tau$ such that the sum of the number of edges in $\tau$ and the total number of $\frak{n}$ decorations is $n$. That is, $n$ means the `complicatedness' of the tree.
It should be noted that, by the definition \eqref{EqIntertwiningDeltaDelta+}, we have
$$
\Delta^+ \mcI_k^+\tau\in T^{(n)+}\otimes T^{(n)+}
$$
for any $\tau\in\mcB^{(n)}$, and thus $\Delta^+$ is closed within $T^{(n)+}$. Therefore the pair of subspaces
$$
\mathscr{T}^{(n)} \defeq  \big(T^{(n)+},T^{(n)}\big),
$$ 
equipped with the restrictions of the $\Delta^+$ and $\Delta$ maps, is also a regularity structure.

\ssk

Under Assumptions \textbf{\textsf{(A-B)}}, formula \eqref{eq inductive g(I_k^+tau)} in the proof of Proposition \ref{PropLinkGPiForAdmissibleModels} shows that the $\sf g$-map of an admissible model $\sf (\Pi,g)$ is uniquely determined by its $\sf \Pi$-map. The following theorem essentially comes from \cite[Proposition 3.32]{Hai14}, which is a generalzation of the `Lyons' extension theorem'.

\ssk

\begin{thm} \label{ThmSmoothAdmissibleModels}
Let $\mathscr{T}$ be a regularity structure satisfying Assumption \textbf{\textsf{(A-B)}} and \eqref{EqConditionBeta0}. One can associate to any family $\big([\tau]; \tau\in\mcB, \vert\tau\vert<0\big)$ of smooth functions on $\bbR\times\bbR^d$ a unique $\bf K$-admissible model $({\sf g}, {\sf \Pi})$ on $\mathscr{T}$ such that ${\sf \Pi}\tau = [\tau]$, for all $\tau\in\mcB$ with $\vert\tau\vert<0$.
\end{thm}   

\medskip

\begin{Dem}
We set the scene for an inductive proof of the statement, taking profit of the induction structure given by Assumption \refB{B2}. 
We will define inductively on $n\in\bbN$ the $\bf K$-admissible models ${\sf M}^{(n)}=({\sf g}^{(n)}, {\sf \Pi}^{(n)})$ on $\mathscr{T}^{(n)}$ over such that
$$
{\sf g}^{(n)} : T^{(n)+} \to \mcC^\infty(\bbR\times\bbR^d),\qquad
{\sf\Pi}^{(n)} : T^{(n)} \to \mcC^\infty(\bbR\times\bbR^d),
$$
and
$$
{\sf g}^{(m)}\vert_{T^{(n)+}}={\sf g}^{(n)},\qquad
{\sf\Pi}^{(m)}\vert_{T^{(n)}}={\sf\Pi}^{(n)}
$$
for any $m>n$. We denote by $\textsf{\textbf{R}}^{{\sf M}^{(n)}}$ the reconstruction operator associated with the model ${\sf M}^{(n)}$.
Note that the model ${\sf M}^{(0)}$ on $\mathscr{T}^{(0)}$ is canonically defined by Assumption \refA{A2}.

\ssk

\textcolor{gray}{$\bullet$} 
We now define an extension ${\sf M}^{(n)}$ of ${\sf M}^{(n-1)}$ on $\mathscr{T}^{(n)}$. It is sufficient to define ${\sf \Pi}^{(n)}\tau$ and ${\sf g}^{(n)}(\mcI_k^+\tau)$ for $\tau\in\mcB^{(n)}$.
Recall from the sentences after Theorem \ref{ThmReconstructionRS} that the function 
$$
\bsh^\tau(x) \defeq  \sum_{\sigma<\tau} {\sf g}_x^{(n-1)}(\tau/\sigma)\,\sigma
$$ 
is an element of $\mcD^{|\tau|}(T^{(n-1)},{\sf g}^{(n-1)})$ and its reconstruction is a candidate of ${\sf\Pi}^{(n)}\tau$.
Given that ${\sf g}^{(n-1)}$ and ${\sf \Pi}^{(n-1)}$ take values in smooth functions, any smooth function is a reconstruction of $\bsh^\tau$ for the model ${\sf M}^{(n-1)}$, if $\vert\tau\vert<0$. (Recall the reconstruction operator is defined uniquely only when acting on modelled distributions of positive regularity. We are working here with a modelled distribution of negative regularity when $|\tau|<0$.) 
Define
$$
{\sf\Pi}^{(n)}\tau=
\begin{cases}
[\tau]&(|\tau|<0),\\
\textsf{\textbf{R}}^{{\sf M}^{(n-1)}}(\bsh^\tau)&(|\tau|>0).
\end{cases}
$$
This is a smooth function in both cases. (Recall that ${\bf1}$ is the only element of $T$ of null homogeneity by Assumption \refA{A1}.) The map ${\sf \Pi}^{(n)}$ coincides with ${\sf \Pi}^{(n-1)}$ on $T^{(n-1)}$.   \vspace{0.1cm}


\textcolor{gray}{$\bullet$} Define then an extension ${\sf g}^{(n)}$ of ${\sf g}^{(n-1)}$ to $T^{(n)+}$ by requiring that it is multiplicative, and by setting 
$$
{\sf g}^{(n)}_x(\mcI^+_k\tau) \defeq  \sum_{\sigma\le\tau; |k|_{\frak{s}} < |\sigma|+2}\, {\sf g}^{(n-1)}_x(\tau/\sigma)\, \partial^k\bfK\Big(\big({\sf\Pi}^{(n)}\big)^{{\sf g}^{(n-1)}}_x\sigma\Big)(x),
$$ 
for all $\tau\in\mcB^{(n)}$, in view of \eqref{eq inductive g(I_k^+tau)}. 
Note that $(T^{(n-1)+},T^{(n)})$ is a regularity structure and $({\sf g}^{(n-1)},{\sf\Pi}^{(n)})$ is a model over it.
Closing the induction step amounts to proving that 
\begin{equation} \label{EqBoundConstructionAdmissible} 
\big|{\sf g}^{(n)}_{yx}(\mcI^+_k\tau)\big| \lesssim d(y,x)^{|\tau| + 2 - |k|_\frak{s}},
\end{equation} 
for every $k\in\bbN\times\bbN^d$ with $\vert k\vert_\frak{s}<\vert\tau\vert+2$. Look for that purpose at the $T^{(n-1)+}$-valued function 
$$
\mcK^{{\sf M}^{(n-1)}+} \bsh^\tau\defeq\big(\mcI^++\mcJ^{{\sf M}^{(n-1)}+}(x)\big)\bsh^\tau(x)+(\mcN^{{\sf M}^{(n-1)}+}\bsh^\tau)(x),
$$
where $\mcI^+:\mcB\to\mcB^+$ is linearly extended by imposing $\mcI^+(X^k)=0$, the linear map $\mcJ^{{\sf M}^{(n-1)}+}:T^{(n-1)}\to T^+$ is defined by
$$
\mcJ^{{\sf M}^{(n-1)}+}\sigma\defeq\sum_{|k|_{\frak{s}}<|\tau|+2}\frac{X_+^k}{k!}\partial^k\bfK\Big(\big({\sf\Pi}^{(n-1)}\big)^{{\sf g}^{(n-1)}}_x\sigma\Big)(x),\qquad\sigma\in\mcB^{(n-1)}
$$
in the same way as $\mcJ^{\sf M}$ by replacing ${\sf M}$ with ${\sf M}^{(n-1)}$ and $X^k$ with $X_+^k$, and the function $\mcN^{{\sf M}^{(n-1)}+}\bsh^\tau$ is defined by
$$
(\mcN^{{\sf M}^{(n-1)}}\bsh^\tau)(x)\defeq\sum_{|k|_{\frak{s}}<|\tau|+2}\frac{X_+^k}{k!}\partial^k{\bf K}
\Big(\textbf{\textsf{R}}^{{\sf M}^{(n-1)}}\bsh^\tau(x)-\big({\sf\Pi}^{(n-1)}\big)^{{\sf g}^{(n-1)}}_x\bsh^\tau(x)\Big)(x)
$$
in the same way as $\mcN^{\sf M}$ by replacing ${\sf M}$ with ${\sf M}^{(n-1)}$ and $X^k$ with $X_+^k$, where $\textbf{\textsf{R}}^{{\sf M}^{(n-1)}}\bsh^\tau$ is defined as ${\sf\Pi}^{(n)}\tau$.
Then in the same proof as Theorem \ref{thm wellposedness of K}, we can prove that $\mcK^{{\sf M}^{(n-1)}+} \bsh^\tau$ is an element of $\mcD^{|\tau|+2}(T^{(n-1)+},{\sf g}^{(n-1)})$ and that, denoting by $\textbf{\textsf{R}}^{{\sf g}^{(n-1)}}$ the reconstruction operator associated with the model $({\sf g}^{(n-1)},{\sf g}^{(n-1)})$ on $(T^{(n-1)+},T^{(n-1)+})$, we have $\textbf{\textsf{R}}^{{\sf g}^{(n-1)}}\mcK^{{\sf M}^{(n-1)}+} \bsh^\tau={\bf K}(\textbf{\textsf{R}}^{{\sf M}^{(n-1)}}\bsh^\tau)$. 
In the proof, Proposition \ref{PropCommutationRelationPlus} has the role of Corollary \ref{CorIdentityAdmissibility}.
Since $\textbf{\textsf{R}}^{{\sf M}^{(n-1)}}\bsh^\tau(x)-\big({\sf\Pi}^{(n-1)}\big)^{{\sf g}^{(n-1)}}_x\bsh^\tau(x)=\big({\sf\Pi}^{(n)}\big)^{{\sf g}^{(n-1)}}_x\tau$ in the definition of $(\mcN^{{\sf M}^{(n-1)}}\bsh^\tau)(x)$, we have
$$
\big(\mcK^{{\sf M}^{(n-1)}+} \bsh^\tau\big)(x) = \mcI^+\big(\bsh^\tau(x)\big) + \sum_{ |k|_{\frak{s}}<|\tau|+2} {\sf g}^{(n)}_x(\mcI^+_k\tau)\,\frac{\Xplus^k}{k!}.
$$
The $\Xplus^k$-component of 
$$
\Big(\mcK^{{\sf M}^{(n-1)}+} \bsh^\tau\Big)(y) - \widehat{{\sf g}^{(n-1)}_{yx}}^+\Big(\mcK^{{\sf M}^{(n-1)}+} \bsh^\tau\Big)(x)
$$ 
is then equal to 

$$
{\sf g}^{(n)}_y(\mcI_k^+\tau)  -  \sum_{\eta<\tau} {\sf g}^{(n)}_x(\tau/\eta){\sf g}^{(n)}_{yx}(\mcI_k^+\eta)  -  \sum_m {\sf g}^{(n)}_x(\mcI_{k+m}^+\tau)\frac{(y - x)^m}{m!} = {\sf g}^{(n)}_{yx}(\mcI^+_k\tau),
$$
and of size $d(y,x)^{|\tau|+2-|k|_\frak{s}}$, since $\mcK^{{\sf M}^{(n-1)}+} \bsh^\tau\in\mcD^{|\tau|+2}\big(T^{(n-1)+}, {\sf g}^{(n-1)}\big)$. This shows the bound \eqref{EqBoundConstructionAdmissible}.

\textcolor{gray}{$\bullet$} It remains to show that ${\sf\Pi}$ is $\bfK$-admissible. Given that we assume $\beta_0=\min A>-2$, the elements of $T$ of the form $\mcI\tau$ have positive homogeneity. So the definition of ${\sf\Pi}$ on $\mcI\tau$ comes under the form of the reconstruction of a modelled distribution $\bsh^{\mcI\tau}$. Since $\bsh^{\mcI\tau}$ is function-like with the ${\bf1}$-component ${\sf g}_x(\mcI^+\tau)$, by Corollary \ref{CorEasyReconstruction}, it follows that
\begin{align*}
{\sf\Pi}(\mcI\tau)(x) = \textbf{\textsf{R}}^{\sf M} (\bsh^{\mcI\tau})(x) = {\sf g}_x(\mcI^+\tau) = \bfK({\sf\Pi}\tau)(x).
\end{align*}

\ssk

\textcolor{gray}{$\bullet$} The above construction makes it clear that the map $\sf \Pi$ is entirely determined from its restriction to the elements of negative homogeneity. The uniqueness part of the statement of the theorem follows then from formula \eqref{eq inductive g(I_k^+tau)} giving ${\sf g}_x(\mcI^+_k\tau)$ as it shows that the map $\sf g$ is entirely determined by the $\sf \Pi$ map under the Assumption \textbf{\textsf{(B1-B2)}}.
\end{Dem}

\bigskip

\section{Solving singular PDEs within regularity structures}
\label{SectionSolvingPDEs}

In this section, we formulate singular stochastic PDEs in the sense of modelled distributions. We trade in this section the generality of the above results for the simplicity of an example that contains the main difficulties of the general case. The reader can consult \cite{Hai14} or \cite{BCCH18} for a description of the general case. We consider the generalized 
(KPZ) equation
\begin{equation}   \label{EqgKPZ*}   \begin{split}
\big(\partial_{x_0}-\Delta_{x'}+1\big) u &= f(u)\zeta + 
\sum_{i,j=1}^d g_2^{ij}(u)(\partial_{x_i}u)(\partial_{x_j}u) + \sum_{i=1}^d g_1^i(u)(\partial_{x_i} u) +g_0(u)  \\
&\eqdef f(u)\zeta + g_2(u)(\partial_{x'} u)^2+g_1(u)\partial_{x'}u+g_0(u)\\
&\eqdef f(u)\zeta +g(u,\partial_{x'} u)
\end{split} \end{equation}
with a noise $\zeta\in\mcC^{\beta_0}$. (Remember that the minimum homogeneity in a regularity structure associated with a singular stochastic PDE coincides with the minimum of the regularities of the noises in the equation.) This type of equation appears in a number of problems. If $d=1$ and $\zeta$ is a space-time white noise then \eqref{EqgKPZ*} contains the KPZ equation, which appears in the large scale picture of one-dimensional random interface evolutions. Here $u$ is scalar valued but a vector valued case is used in the description of the random motion of a rubber on a manifold \cite{HairerString}, a random perturbation of the harmonic flow map on loops. If $d=2,3$ and $\zeta$ is a space white noise, then \eqref{EqgKPZ*} contains the generalized PAM
$$
(\partial_t-\Delta_x)u=f(u)\zeta.
$$

\ssk

The differential equation \eqref{EqgKPZ*} with the initial value $u(0,x')=u_0(x')$ has an equivalent integral form
$$
u(x)=e^{x_0(\Delta_{x'}-1)}u_0(x')+{\bf L}^{-1}\big(f(u)\zeta+g(u,\partial_{x'}u)\big)(x).
$$
Under an appropriate setting, the generalized (KPZ) equation \eqref{EqgKPZ*} will be lifted to the following equation on modelled distributions $\bsv\in\mcD^\gamma(T,\sf g)$
\begin{equation}   \label{EqPDERS}
\bsv = \bsh+\mcP^{\sf M}\Big(f^\star\big(\bsv\big) \Xi + g^\star\big(\bsv, D\bsv \big)\Big),
\end{equation}
for some $T_X$-valued modelled distribution $\bsh$ and an operator $\mcP^{\sf M}$ having the role of ${\bf L}^{-1}$. This section is dedicated to giving the meaning to the equation \eqref{EqPDERS} and showing that this equation has a unique solution on a small time interval $(0,t_0)$; this is the content of Theorem \ref{ThmWellPosedness}, which is the main result of this section. At the end of this section we will be in a position to define the model-dependent solution of the singular equation \eqref{EqgKPZ*} as the model-dependent reconstruction $u^{\sf M} = {\sf R}^{\sf M}(\bsu^{\sf M})$ of the unique solution $\bsu^{\sf M}$ of \eqref{EqPDERS}. The function $u^{\sf M}$ will then appear as a continuous function of $\sf M$.  

\ssk

The restriction to each band $[0,t_0]\times [-R,R]$ of spacetime white noise has a norm growing indefinitely as $R$ goes to infinity for each fixed $t_0>0$. To avoid working with unbounded spacial domains and functional spaces involving spacial weights we will assume that all the objects are $\bbZ$-periodic in space -- they would be $\bbZ^d$-periodic in space in a more general setting. The function $\bsh$ in \eqref{EqPDERS} plays the role of the regularity structure lift of the propagator of the initial condition $u_0$. The use of time weights to take care of the free propagation $(e^{ t(\Delta-1) }u_0)_{t>0}$ of the initial condition in a regularity structures setting is made necessary by the classical sharp estimate
\begin{equation} \label{EqSchauderHeatSemigroup}
\big\| \partial_x^k e^{ t(\Delta-1) }u_0\big\|_\infty \lesssim 
t^{-\frac{(\vert k\vert_\frak{s}-\alpha)\vee0}{2}}
\,\|u_0\|_{\mcC^\alpha(\bbR^d)}.
\end{equation}
Theorem \ref{ThmWellPosedness} is proved under spacial periodic boundary conditions and in the space of modelled distributions involving temporal weights exploding in $t=0^+$. We introduce the former in Section \ref{SubsectionPeriodicModel} and the latter in Section \ref{SubsectionTimeWeight}. We examine in Section \ref{SubsectionNonAnticipative} the notion of non-anticipative operator, involved in the analysis of equation \eqref{EqPDERS}. We prove in Section \ref{SubsectionFixedPoint} that \eqref{EqPDERS} is locally in time well-posed.

\bigskip

\subsection{Spatially periodic models}
\label{SubsectionPeriodicModel}

We work on the models and modelled distributions that are spacially $\bbZ^d$-periodic, with $d=1$ here -- we give the definitions for an arbitrary space dimension $d$. All the results and estimates proved above hold true in the periodic case. For any $x=(x_0,x')\in\bbR\times\bbR^d$ and $m\in\bbZ^d$, denote by $x+m\defeq (x_0,x'+m)$.

\medskip

\begin{defn*}
A \textbf{\textsf{model}} $\sf M = (g,\Pi)$ is said to be $\bbZ^d$\textbf{\textsf{-periodic}} if 
for any $m\in\bbZ^d$,
\begin{align*}
{\sf g}_{y+m,x+m} = {\sf g}_{yx},\qquad
\big\langle {\sf \Pi}^{\sf g}_{x+m}\tau \,,\,\varphi(\cdot+m) \big\rangle 
= \big\langle {\sf \Pi}^{\sf g}_x\tau \,,\,\varphi(\cdot) \big\rangle,
\end{align*}
for all $x,y\in\bbR\times\bbR^d, \tau\in T$ and all $\varphi\in\mcS(\bbR\times\bbR^d)$.   
\end{defn*}

\medskip

The canonical model $({\sf g},{\sf \Pi})$ on the polynomial regularity structure $(T_X^+,T_X)$ is $\bbZ^d$-periodic in the above sense. Note that ${\sf g}_x(\Xplus^n)=x^n$ and $({\sf\Pi}X^n)(x)=x^n$ are \emph{not} $\bbZ^d$-periodic functions.
This is the reason why we do not impose periodic conditions on ${\sf g}$ and ${\sf\Pi}$.
It is elementary to see that if $\sf M$ is a $\bbZ^d$-periodic model on $\mathscr{T}$ and $\bsf\in\mcD^\gamma(T,\sf g)$ is $\bbZ^d$-periodic, with 
$\gamma$ positive, then $\textsf{\textbf{R}}^{\sf M}\bsf$ is also $\bbZ^d$-periodic, in the sense that 
$$
\big\langle \textsf{\textbf{R}}^{\sf M}\bsf, \varphi(\cdot+m)\big\rangle = \big\langle \textsf{\textbf{R}}^{\sf M}\bsf, \varphi(\cdot)\big\rangle,
$$
for all $m\in\bbZ^d$ and $\varphi\in\mcS(\bbR\times\bbR^d)$ -- see Proposition 3.38 in \cite{Hai14}. {\it All objects in remainder of this section 
are implicitly assumed to be $\bbZ^d$-periodic.}

\bigskip

\subsection{Modelled distributions with singularity at $x_0=0$}
\label{SubsectionTimeWeight}

We use time weights
$$
\omega(x)\defeq|x_0|^{\frac12}\wedge1,\qquad
\omega(x,y)\defeq\omega(x)\wedge\omega(y)
$$
to treat the boundary condition at $x_0=0$.

\medskip

\begin{defn*}
Fix two exponents $\eta\le\gamma\in\bbR$. One defines the space $\mcD^{\gamma,\eta}(T,{\sf g})$ of \textbf{\textsf{modelled distributions with singularity of weight $\eta$ at $x_0=0$}}, as the space of functions $\bsf$ from $\bbR\times\bbR^d\setminus\{x_0=0\}$ into $T_{<\gamma}$ such that
\begin{align}\label{Def singular MD}
\begin{aligned}
\brarb{\bsf}_{\mcD^{\gamma,\eta}}
&\defeq  \max_{\beta<\gamma}\sup_{x\in(\bbR\setminus\{0\})\times\bbR^d}
\frac{\|\bsf(x)\|_\beta}{\omega(x)^{(\eta-\beta)\wedge0}}<\infty,   \\
\|\bsf\|_{\mcD^{\gamma,\eta}}
&\defeq  \max_{\beta<\gamma}\sup_{x,y\in(\bbR\setminus\{0\})\times\bbR^d,\, d(x,y)\le\omega(x,y)}
\frac{\big\|\bsf(y)-\widehat{{\sf g}_{yx}}\bsf(x)\big\|_\beta}{\omega(x,y)^{\eta-\gamma}d(y,x)^{\gamma-\beta}}<\infty.
\end{aligned}
\end{align}
Set $\trino{\bsf}_{\mcD^{\gamma,\eta}}\defeq \brarb{\bsf}_{\mcD^{\gamma,\eta}}+\|\bsf\|_{\mcD^{\gamma,\eta}}$.
\end{defn*}

\medskip

One also talks of singular modelled distributions. An example of singular modelled distributions is obtained as follows. Given $\eta\in\bbR$ and $v\in\mcC^\eta(\bbT^d)$, the $T_X$-valued function
\begin{align}\label{eq heat propagator}  
(P_\gamma v)(x) \defeq  {\bf1}_{x_0>0}\sum_{|k|_{\frak{s}}<\gamma} \partial^k\big(e^{x_0(\Delta_{x'}-1)}v\big)(x)\,\frac{X^k}{k!}
\end{align}
belongs to $\mcD^{\gamma,\eta}(T,\sf g)$ for any $\gamma\ge\eta$. This is 
a consequence of the Schauder estimates satisfied by the heat semigroup recalled in \eqref{EqSchauderHeatSemigroup} -- see e.g. Lemma 7.5 of \cite{Hai14}. 

\medskip

We recall some embedding theorems. It is easy to see that $|\!|\!| \bsf |\!|\!|_{\mcD^{\gamma,\eta'}}\le|\!|\!| \bsf |\!|\!|_{\mcD^{\gamma,\eta}}$ if $\eta'\le\eta$. If $\eta\le\gamma'\le\gamma$, we also have
$$
|\!|\!| \mcQ_{<\gamma'}\bsf |\!|\!|_{\mcD^{\gamma',\eta}}\lesssim|\!|\!| \bsf |\!|\!|_{\mcD^{\gamma,\eta}}
$$
with an implicit constant depending on the model $\sf M$ -- see e.g. Proposition 3.5 of \cite{HT24}. Instead of $\brarb{\bsf}_{\mcD^{\gamma,\eta}}$ and $|\!|\!| \bsf |\!|\!|_{\mcD^{\gamma,\eta}}$, it will be convenient to consider the seminorms
$$
\brarb{\bsf}_{\mcD^{\gamma,\eta}}'
\defeq 
\max_{\beta<\gamma}\sup_{x\in(\bbR\setminus\{0\})\times\bbR^d}
\frac{\|\bsf(x)\|_\beta}{\omega(x)^{\eta-\beta}}
$$
and $|\!|\!| \bsf |\!|\!|_{\mcD^{\gamma,\eta}}' \defeq \brarb{\bsf}_{\mcD^{\gamma,\eta}}'+ \|\bsf\|_{\mcD^{\gamma,\eta}}$.
In general $|\!|\!| \bsf |\!|\!|_{\mcD^{\gamma,\eta}}\le|\!|\!| \bsf |\!|\!|_{\mcD^{\gamma,\eta}}'$ but the reverse inequality fails. However for any $\bsf\in\mcD^{\gamma,\eta}$ such that
$$
\lim_{x_0\to0}\mcQ_\beta \bsf(x)=0
$$
for any $\beta<\eta$, the reverse inequality $|\!|\!| \bsf |\!|\!|_{\mcD^{\gamma,\eta}}'\lesssim|\!|\!| \bsf |\!|\!|_{\mcD^{\gamma,\eta}}$ holds with an implicit constant depending on the model $\sf M$ -- see e.g. Lemma 6.5 of \cite{Hai14}.

\medskip

The reconstruction theorem, Theorem \ref{ThmReconstructionRS}, is extended to singular modelled distributions as follows. See Appendix {\sf \ref{Proof of singular reconstruction}} for a detailed proof. An extension to the inhomogeneous integral kernels can be found in \cite{HT24}.

\medskip

\begin{thm}\label{thm singular reconstruction}
Let ${\sf M=(g,\Pi)}$ be a model over $\mathscr{T}$ such that $-2<\beta_0<0$. Assume that $-2<\eta\le\gamma$, with $\gamma>0$. Then there exists a continuous linear operator
$$
\textbf{\textsf{R}}^{\sf M}:\mcD^{\gamma,\eta}(T,{\sf g}) \to \mcC^{\eta\wedge\beta_0}(\bbR\times\bbR^d)
$$
such that, for any $\bsf\in\mcD^{\gamma,\eta}(T,{\sf g})$ and $n\in\bbN\times\bbN^d$, the bound
\begin{equation}
\label{EqSingularReconstructionCondition}
\Big|\big\langle \textbf{\textsf{R}}^{\sf M}\bsf - {\sf\Pi}_x^{\sf g}\bsf(x) , \partial_x^np_t(x,\cdot)\big\rangle\Big|  \lesssim  \trino{\bsf}_{\mcD^{\gamma,\eta}} \, \omega(x)^{(\eta\wedge\beta_0)-\gamma} t^{\frac{\gamma-|n|_\mfs}4},
\end{equation}
holds uniformly over $\bsf\in\mcD^{\gamma,\eta}(T, {\sf g}), x\in\bbR\times\bbR^d$ and $0<t\le\omega(x)^4$, where the implicit proportional constant depends polynomially on $\|{\sf g}\|_\gamma+\|{\sf\Pi}^{\sf g}\|_\gamma$, and is independent of $\bsf$. Such an operator is unique if the exponent $\gamma$ is positive.
\end{thm}

\medskip

The operators discussed in previous sections can be extended to the spaces $\mcD^{\gamma,\eta}(T,\sf g)$, as follows. All of the following maps are locally Lipschitz continuous. For the detailed proofs, see \cite[Propositions 6.12, 6.13, 6.15, and 6.16]{Hai14}. Denote by $\mcD_\alpha^{\gamma,\eta}(T,\sf g)$ the space of modelled distributions $\bsf\in\mcD^{\gamma,\eta}(T,\sf g)$ of the form
\begin{align}\label{def D_a^ge}
\bsf=\sum_{\alpha\le|\tau|<\gamma}f_\tau\tau.   \vspace{0.1cm}
\end{align}

\begin{itemize}
\setlength{\itemsep}{0.1cm}
	\item (Proposition \ref{PropRegularityProduct}')
	Let $\alpha_1,\alpha_2\le0<\gamma_1,\gamma_2$, and set $\gamma=(\gamma_1+\alpha_2)\wedge(\gamma_2+\alpha_1)$ and $\eta=(\eta_1+\alpha_2)\wedge(\eta_2+\alpha_1)\wedge(\eta_1+\eta_2)$.
If a regular product $\star:V\times W\to T$ is given, then
\begin{align*}
\mcD_{\alpha_1}^{\gamma_1,\eta_1}(V,{\sf g})\times\mcD_{\alpha_2}^{\gamma_2,\eta_2}(W,{\sf g})
\ni(\bsf_1,\bsf_2)
\mapsto \mcQ_{<\gamma}(\bsf_1\star\bsf_2)\in
\mcD_{\alpha_1+\alpha_2}^{\gamma, \eta}(T,{\sf g}).
\end{align*}
	
	\item (Proposition \ref{PropFModelled}')
Let $\gamma>0$ and $0\le\eta\le\gamma$. If an associative regular product $\star:V\times V\to V$ and a smooth function $F$ is given, then
$$
\mcD^{\gamma,\eta}(V,{\sf g})\ni\bsf \mapsto F^\star(\bsf)\in\mcD^{\gamma,\eta}(V,{\sf g}).   \vspace{0.1cm}
$$

	\item (Proposition \ref{prop abstract derivative}')
Let $\gamma>1$. If a derivative $D:T\to T$ is given, then
$$
\mcD^{\gamma,\eta}(T,{\sf g})\ni\bsf \mapsto D\bsf\in\mcD^{\gamma-1,\eta-1}(T,{\sf g}).   \vspace{0.1cm}
$$

	\item (Theorem \ref{thm wellposedness of K}')
Let $\gamma>0$ and $-2<\eta\wedge\beta_0$. If ${\sf\Pi}$ is ${\bf K}$-admissible, 
$$
\mcD^{\gamma,\eta}(T,{\sf g})\ni\bsf \mapsto \mcK^{\sf M}\bsf\in\mcD^{\gamma+2,\eta\wedge\beta_0+2}(T,{\sf g}).
$$
\end{itemize}
A sketch of the proof of Theorem \ref{thm wellposedness of K}' is given in Appendix {\sf \ref{Proof of singular reconstruction}}. About this statement, note here the gain in the explosion exponent after we applied the operator $\mcK^{\sf M}$. We will use this gain in Section \ref{SubsectionFixedPoint} to gain a small contraction factor in the fixed point formulation of equation \eqref{EqgKPZ*} as an equation on a space of modelled distributions.

\bigskip

\subsection{Non-anticipative operators}
\label{SubsectionNonAnticipative}

A function $f$ on $(\bbR\times\bbT^d)^2$ is said to be {\it non-anticipative} if $f\big((x_0,x'),(y_0,y')\big)=0$, whenever $x_0<y_0$. The kernel {\color{red}$P$} of the resolution operator ${\bf L}^{-1}$ (in the sense that ${\bf L}^{-1}f(x)=\int_{\bbR\times\bbT^d}P(x-y)f(y)dy$) is of the form
$$
P(x,y)={\bf1}_{x_0>y_0}\,p_{x_0-y_0}(x'-y'),
$$
where $p_t$ is the kernel of $e^{t(\Delta-1)}$, thus $P$ is non-anticipative. The aim of this section is to prove the refined multilevel Schauder estimate (Proposition \ref{PropNonAnticipative}) associated with the non-anticipative operator ${\bf L}^{-1}$.

\medskip

We consider the modelled distributions defined on the domain $(0,t)\times\bbT^d$, for a given positive time $t$. Denote by $\mcD_{(0,t)}^{\gamma,\eta}(T,\sf g)$ the set of functions $\bsf : (0,t)\times\bbT^d \to T_{<\gamma}$ such that the bounds \eqref{Def singular MD} hold with the domain of $x,y$ restricted to $(0,t)\times\bbT^d $. Denote by
$$
\trino{\bsf}_{\mcD_{(0,t)}^{\gamma,\eta}} \defeq  \brarb{\bsf}_{\mcD_{(0,t)}^{\gamma,\eta}} + \|\bsf\|_{\mcD_{(0,t)}^{\gamma,\eta}}
$$ 
the associated norms. It is also useful to consider $\brarb{\bsf}_{\mcD_{(0,t)}^{\gamma,\eta}}'$ and $|\!|\!| \bsf |\!|\!|_{\mcD_{(0,t)}^{\gamma,\eta}}' \defeq \brarb{\bsf}_{\mcD_{(0,t)}^{\gamma,\eta}}'+ \|\bsf\|_{\mcD_{(0,t)}^{\gamma,\eta}}$. Since $\omega(x),\omega(y)\le t^{\frac12}$ if $x,y\in(0,t)\times\bbT^d$, we have
$$
|\!|\!| \bsf |\!|\!|_{\mcD_{(0,t)}^{\gamma,\eta-\kappa}}'\lesssim t^{\frac\kappa2}|\!|\!| \bsf |\!|\!|_{\mcD_{(0,t)}^{\gamma,\eta}}'
$$
for any $\kappa>0$ small enough such that $[\eta-\kappa,\eta)\cap A=\emptyset$. The small factor $t^{\kappa/2}$ is used in the fixed point problem in the next section. To apply the reconstruction operator to locally defined modelled distributions, we use the cut-off operator.
The following result is obtained from Proposition \ref{PropRegularityProduct}' and Lemma \ref{App:lem:connectionMD}. See Lemma 5.7 of \cite{HT24} for the detailed proof.

\medskip

\begin{prop}
Let $\sf M=(g,\Pi)$ be a model over $\mathscr{T}$ such that $-2<\beta_0<0$ and let $\gamma>0$ and $\eta\le\gamma\wedge\beta_0$.
Fix a smooth non-increasing function $\chi:(0,\infty)\to[0,1]$ such that $\chi(t)=1$ if $0<t\le\frac12$ and $\chi(t)=0$ if $t\ge1$.
For each $t>0$ and $x\in\bbR\times\bbT^d$, we set $\chi_t(x)={\bf1}_{x_0>0}\chi(x_0/t)$ and define
$$
\boldsymbol{\chi}_t(x)\defeq\sum_{|k|_\mfs<\gamma-\beta_0}\frac{\partial_x^k\chi_t(x)}{k!}X^k.
$$
Then one can define the linear operator $C_t:\mcD^{\gamma,\eta}_{(0,t)}(T,{\sf g})\to\mcD^{\gamma,\eta}(T,{\sf g})$ by
$$
C_t\bsf(x)={\bf1}_{(0,t)\times\bbT^d}(x)\mcQ_{<\gamma}\big(\boldsymbol{\chi}_t(x)\star\bsf(x)\big),
$$
and $C_t$ is uniformly bounded over $t\in(0,1]$ and satisfies $(C_t\bsf)\vert_{(0,\frac{t}2)\times\bbT^d}= \bsf\vert_{(0,\frac{t}2)\times\bbT^d}$.
\end{prop}

\medskip

\begin{prop} \label{PropNonAnticipative}
Pick $\gamma>0$, $-2<\eta\le\beta_0<0$, and $0<\rho\le\gamma+2$. For any $\bf K$-admissible model $\sf M=(g,\Pi)$ and $t\in(0,1]$, there exists a continuous linear map $\mcP_t^{\sf M}:\mcD_{(0,t)}^{\gamma,\eta}(T,{\sf g})\to\mcD^{\rho,\eta+2}(T,{\sf g})$ such that the following properties hold for any $\bsf\in\mcD_{(0,t)}^{\gamma,\eta}(T,{\sf g})$.
\begin{itemize}
\setlength{\itemsep}{0.1cm}
\item[(1)]
For any $x\in(0,\frac{t}2)\times\bbT^d$, one has $\mcP_t^{\sf M}\bsf(x)-\mcI\bsf(x)\in T_X$.
\item[(2)]
For any $\kappa>0$, one has
\begin{align}\label{ineq small factor improved}
\trino{\mcP_t^{\sf M}\bsf}_{\mcD^{\rho,\eta+2-\kappa}_{(0,t)}} 
\lesssim t^{\kappa/2}\,\trino{\bsf}_{\mcD^{\rho,\eta}_{(0,t)}},
\end{align}
where the implicit proportional constant is independent of $\bsf$ and $t$. 
\item[(3)]
If ${\sf\Pi}_x^{\sf g}\tau$ happens to be continuous for any $\tau\in T$, then the function $\textbf{\textsf{R}}^{\sf M}\bsf\defeq\textbf{\textsf{R}}^{\sf M}(C_t\bsf)$ satisfies
\begin{align}\label{eq KR=RK half plane}
\textbf{\textsf{R}}^{\sf M}\big(\mcP_t^{\sf M}\bsf\big)(x) = \int_{[0,x_0]\times\bbT^d }P(x,y)\,\textbf{\textsf{R}}^{\sf M}\bsf(y)\,dy.
\end{align}
\end{itemize}
\end{prop}

\medskip

\begin{Dem}
We provide only a sketch here. See Theorem 5.9 of \cite{HT24} for the detailed proof.
Recall the decomposition ${\bf L}^{-1}={\bf K}+{\bf K}'$.
Noting that $\textbf{\textsf{R}}^{\sf M}\bsf\defeq\textbf{\textsf{R}}^{\sf M}(C_t\bsf)\in\mcC^\eta$, we denote by $(\mcK')^{\sf M}$ the lift of the ${\bf K}'$ operator in the polynomial part of the regularity structure
$$
\big((\mcK')^{\sf M}\bsf\big) \defeq  \sum_{\vert\ell\vert_{\frak{s}} 
< \rho} \partial^\ell{\bf K}'\big(\textbf{\textsf{R}}^{\sf M}\bsf\big)\,
\frac{X^\ell}{\ell!} \in T_X.
$$
Since $\bfK'$ maps $\mcC^\eta$ into $\mcC^\infty$, we can show that $(\mcK')^{\sf M}\bsf\in\mcD^{\rho,\eta+2}$.
Then we can also show that
$$
\mcP_t^{\sf M}\bsf=\mcQ_{<\rho}\mcK^{\sf M}(C_tf)+(\mcK')^{\sf M}\bsf\in\mcD^{\rho,\eta+2}
$$
by Theorem \ref{thm wellposedness of K}'.
The property (1) is obvious from the definition. To show the property (3), note that
\begin{align*}
\textbf{\textsf{R}}^{\sf M}\mcP_t^{\sf M}\bsf
={\bf K}\big(\textbf{\textsf{R}}^{\sf M}\bsf\big)+{\bf K}'\big(\textbf{\textsf{R}}^{\sf M}\bsf\big)
={\bf L}^{-1}\big(\textbf{\textsf{R}}^{\sf M}\bsf\big).
\end{align*}
Since $\textbf{\textsf{R}}^{\sf M}\bsf$ vanishes on $(-\infty,0)\times\bbT^d$ by Corollary \ref{CorLocalDependenceReconstruction}, we have \eqref{eq KR=RK half plane}.
To show the property (2), recall the sufficient condition for the equivalence between two norms $\trino{\cdot}_{\mcD^{\gamma,\eta}}$ and $\trino{\cdot}_{\mcD^{\gamma,\eta}}'$.
Set $\bsg\defeq\mcP_t^{\sf M}\bsf$.
By definition, $\bsg$ takes values in the function-like sector $V$ with $\alpha_0(V)=\beta_0+2$. 
If $\eta\le-1$, it is sufficient to show that
$$
\lim_{t\downarrow0}\bsg_{\bf1}(x)=0
$$
for the equivalence between $\trino{\bsg}_{\mcD^{\rho,\eta+2}}$ and $\trino{\bsg}_{\mcD^{\rho,\eta+2}}'$.
Recall that $\bsg_{\bf1}=\textbf{\textsf{R}}^{\sf M}\bsg={\bf L}^{-1}\big(\textbf{\textsf{R}}^{\sf M}\bsf\big)$.
Since $\textbf{\textsf{R}}^{\sf M}\bsf\in\mcC^\eta$, we have $\bsg_{\bf1}\in\mcC^{\eta+2}$ by Schauder estimate, so it is H\"older continuous. Since $\bsg_{\bf1}=0$ on $(-\infty,0)\times\bbT^d$ from the non-anticipativity of $K_{\bf L}$, it also vanishes at $x_0=0$.
If $-1\le\eta\le0$, we also have $\lim_{t\downarrow0}\bsg_{X^k}(x)=0$ for any $|k|_\mfs=1$ by a similar argument. By the equivalence between $\trino{\bsg}_{\mcD^{\rho,\eta+2}}$ and $\trino{\bsg}_{\mcD^{\rho,\eta+2}}'$, we have \eqref{ineq small factor improved} as follows.
\begin{align*}
|\!|\!|\bsg |\!|\!|_{\mcD_{(0,t)}^{\rho,\eta+2-\kappa}(0,T)}
&\lesssim
|\!|\!|\bsg |\!|\!|_{\mcD_{(0,t)}^{\rho,\eta+2-\kappa}(0,T)}'
\lesssim
t^{\frac\kappa2}|\!|\!|\bsg |\!|\!|_{\mcD_{(0,t)}^{\rho,\eta+2}(0,T)}'\\
&\lesssim
t^{\frac\kappa2}|\!|\!|\bsg |\!|\!|_{\mcD_{(0,t)}^{\rho,\eta+2}(0,T)}
\lesssim
t^{\frac\kappa2}|\!|\!|\bsf |\!|\!|_{\mcD_{(0,t)}^{\gamma,\eta}(0,T)}.
\end{align*}
\end{Dem}

\bigskip

\subsection{Fixed point solution}
\label{SubsectionFixedPoint}

Finally we make sense of the equation \eqref{EqPDERS} and show its local well-posedness.

\medskip

\begin{defn*}
A \textbf{\textsf{regularity structure}} $\mathscr{T}$ is said to be \textbf{\textsf{associated with equation \eqref{EqgKPZ*}}} if it satisfies Assumptions \textbf{\textsf{(A-B)}} and contains subcomodules $S, DS, F, N$ of $T$ satisfying Assumption \refA{A1} and the following constraints. 
\begin{itemize}
\setlength{\itemsep}{0.1cm}
	\item The symbol $\Xi$ and the set $\partial S$ are contained in $N$.
	
	\item The sector $S$ is function-like and regular products
$$
S\times\cdots\times S\to F,\qquad DS\times DS\to N,\qquad F\times N\to T
$$
are given and satisfy Assumption \refA{A3}. We denote them all by the same symbol $\star$.
	
	\item Abstract integration operators
$$
\mcI:T\to S,\qquad \mcI_{e_i}:T\to DS,\qquad (1\leq i\leq d)
$$
are given and satisfy Assumption \REFB.
	
	\item Derivative operators
$$
D_i:S\to DS, \qquad (1\leq i\leq d)
$$
are given and satisfy 
$$
{\sf \Pi}\circ D_i = \partial_{x_i}\circ{\sf \Pi}, \quad\textrm{and}\quad D_iX^k = k_iX^{k-e_i}{\bf1}_{k\ge e_i},\quad\text{and}\quad D_i\mcI\tau=\mcI_{e_i}\tau.
$$
\end{itemize}
\end{defn*}

\medskip

The element $\Xi$ represents the noise $\zeta$. The spaces $S$ and $DS$ are used to represent the solution $u$ and its derivative $\partial_{x'}u$, respectively. (The letter $S$ is chosen for `solution'.) The space $F$ are used to represent $f(u)$ and $g_n(u)$, with $n=0,1,2$. (The 
letter $F$ is chosen for `function'.) The space $N$ is used to represent the `singular' elements $\zeta$, $\partial_{x'}u$, and $(\partial_{x'}u)^2$. (The letter $N$ is chosen for `noise'.) The only role played by the intermediate spaces $DS,F, N$ is to clarify on which spaces the product $\star$ is defined; they play no other role. We will see in Section \ref{SectionBuildingRS} how to construct explicitly a regularity structure associated with the generalized (KPZ) equation. The $\star$ product is used to define nonlinear images of singular modelled distributions as in Section \ref{section products and derivatives}. In this setting, the regularity structure lift of the generalized (KPZ) equation is formulated under the form 
\begin{equation}   \label{EqPDERSPrime} \begin{split}
\bsv &= \bsh+\mcP_t^{\sf M}\Big(f^\star\big(\bsv\big) \Xi + g^\star\big(\bsv, D\bsv \big)\Big)  \\
       &\eqdef \Phi_t^{\bsh,{\sf M}}(\bsv),
\end{split} \end{equation}   
for appropriate choices of $t\in(0,1]$ and $\rho>0$, where $f^\star(\bsv)$ denotes the composition operator (Proposition \ref{PropFModelled}), and $f^\star(\bsv)\Xi\defeq f^\star(\bsv)\star\Xi$, and
\begin{align*}
g^\star(\bsv,D\bsv)
&\defeq\mcQ_{<\gamma'}\big\{g_2^\star(\bsv)\star(D\bsv)^{\star2}+g_1^\star(\bsv)\star (D\bsv)+g_0^\star(\bsv)\big\}\\
&\defeq\mcQ_{<\gamma'}\bigg\{\sum_{i,j=1}^d(g_2^{ij})^\star(\bsv)\star(D_i\bsv)\star(D_j\bsv)+\sum_{i=1}^d(g_1^i)^\star(\bsv)\star (D_i\bsv)+g_0^\star(\bsv)\bigg\}
\end{align*}
for an appropriate $\gamma'\in\bbR$ -- we will choose $\gamma'=\gamma+\beta_0$ in the proof of Theorem \ref{ThmWellPosedness}.

\medskip

Pick a $\bf K$-admissible model $\sf M=(g,\Pi)$ over $\mathscr{T}$ and $\bsh\in\mcD^{\gamma,\eta}(T_X,\sf g)$. Assume that $\Phi^{\bsh,{\sf M}}$ sends $\mcD^{\gamma,\eta}(S,\sf g)$ into itself, which turns out to be the 
case as proved below under the conditions of Theorem \ref{ThmWellPosedness}.

\medskip

\begin{defn*}
A \textbf{\textsf{solution to equation \eqref{EqPDERS} on the time interval $(0,t_0)$}} is a fixed point of the map $\Phi_{t_0}^{\bsh,{\sf M}} : \mcD^{\gamma,\eta}_{(0,t_0)}(S,{\sf g})\rightarrow\mcD^{\gamma,\eta}_{(0,t_0)}(S,{\sf g})$.
\end{defn*}

\medskip

\begin{thm} \label{ThmWellPosedness}
Assume that $f$ and $g$ are smooth functions. Let $\mathscr{T}$ be a regularity structure associated with equation \eqref{EqgKPZ*} and satisfying Assumptions \textbf{\textsf{(A-B)}}, with $\beta_0\in(-2,-1)$. Pick $\eta\in(0,\beta_0+2]$ and $\gamma>-\beta_0$. Then for any $\bf K$-admissible model $\sf M=(g,\Pi)$ and any $\bsh\in\mcD^{\gamma,\eta}(T_X,{\sf g})$, there exists a positive time $t_0=t_0(\bsh,{\sf M})$ such that equation \eqref{EqPDERS} has a unique solution $\bsu$ on the time interval $(0,t_0)$. The time $t_0$ can be chosen to be a lower semicontinuous function of $\bsh$ and $\sf M$.
\end{thm}

\medskip

\begin{Dem}
Recall that $\mcD_\alpha^{\gamma,\eta}(T,{\sf g})$ denotes the set of modelled distributions of the form \eqref{def D_a^ge}.
Starting from $\bsv\in\mcD^{\gamma,\eta}(S,{\sf g})$, we show that
$$
f^\star(\bsv)\Xi + g^\star(\bsv, D\bsv )\in\mcD^{\gamma+\beta_0,2\eta-2}(T,{\sf g}).
$$
From the `singular/exploding' version of Proposition \ref{PropFModelled} given at the end of Section \ref{SubsectionTimeWeight}, one has $f^\star(\bsv),(g_2^{ij})^\star(\bsv),(g_1^i)^\star(\bsv),g_0^\star(\bsv)\in\mcD_0^{\gamma,\eta}(F,{\sf g})$. Since $\Xi\in\mcD_{\beta_0}^{\infty,\infty}(N,{\sf g})$, one has 
$$
f^\star(\bsv)\Xi\in \mcD^{\gamma+\beta_0,\eta+\beta_0}(T,{\sf g})\subset \mcD^{\gamma+\beta_0,2\eta-2}(T,{\sf g})
$$ 
from the singular version of Proposition \ref{PropRegularityProduct} given at the end of Section \ref{SubsectionTimeWeight}. Noting that the smallest homogeneity in the subcomodule $\partial S$ 
is $\beta_0+1<0$, which is the homogeneity of $\mcI_{e_i} \Xi$, one has $ D_i\bsv \in\mcD_{\beta_0+1}^{\gamma-1,\eta-1}(\partial S,{\sf g})$ 
and $( D_i\bsv )\star(D_j\bsv)\in\mcD_{2\beta_0+2}^{\gamma+\beta_0,2\eta-2}(N,{\sf g})$. Thus $g^\star(\bsv, D\bsv )\in\mcD^{\gamma+\beta_0,2\eta-2}(T,{\sf g})$. From Proposition \ref{PropNonAnticipative}, one has
\begin{align*}
\trino{\Phi_t^{\bsh,\sf M}(v)}_{\mcD_{(0,t)}^{\gamma,\eta}}
&\lesssim \trino{\bsh}_{\mcD_{(0,t)}^{\gamma,\eta}}+t^{\eta/2}\trino{f^\star(\bsv)\Xi + g^\star(\bsv, D\bsv )}_{\mcD_{(0,t)}^{\gamma+\beta_0,2\eta-2}}   \\
&\lesssim \trino{\bsh}_{\mcD_{(0,t)}^{\gamma,\eta}}+t^{\eta/2}C\Big(\trino{\bsv}_{\mcD_{(0,t)}^{\gamma,\eta}}\Big)
\end{align*}
for some locally bounded function $C$. Then one can associate with each positive radius $\lambda\gtrsim\trino{\bsh}_{\mcD^{\gamma,\eta}}$ a time horizon $t(\lambda)$ such that $\Phi_{t(\lambda)}^{\bsh,\sf M}$ sends the ball of $\mcD_{(0,t(\lambda))}^{\gamma,\eta}(S,{\sf g})$ of radius $\lambda$ into itself. From the local Lipschitz continuity result, the map $\Phi_{t(\lambda)}^{\bsh,\sf M}$ is also a contraction on the ball of $\mcD_{(0,t(\lambda))}^{\gamma,\eta}(S,{\sf g})$ of radius $\lambda$. As such, it has a unique fixed point on the ball of radius $\lambda$. An elementary argument gives the uniqueness of a fixed point within $\mcD_{(0,t(\lambda))}^{\gamma,\eta}(S,{\sf g})$, as in the proof of Theorem 4.7 in \cite{HairerKPZ}.
\end{Dem}

\medskip

The proof makes it clear that one can ask $f$ and $g$ to have finite regularity rather than being smooth. We do not try to optimize the regularity assumptions on $f$ and $g$ here. Thinking of $\bsh$ as the regularity structure lift of the free propagation of an initial condition on $\bbT^d$, assuming in $\bsh\in\mcD^{\gamma,\eta}(T_X,\sf g)$ allows us to work with an initial condition of H\"older regularity $\eta$ -- recall the constraint $\eta\in(0,\beta_0+2]$. Note that the map is uniformly contracting on a small enough time interval for $f$ and $g$ ranging in a bounded set. In order to compare fixed points of $\Phi^{\bsh,{\sf M}}$ associated with different admissible models over $\mathscr{T}$ -- hence different maps on different spaces, we use the metric $d_\gamma$ in Definition \ref{DefnModelledDistribution} with a slight modification to $\mcD_{(0,t)}^{\gamma,\eta}$ norms.
One can then prove the following statement in terms of this metric by making explicit in the reconstruction theorem and the lifting theorem that the operators $\textbf{\textsf{R}}^{\sf M}$ and $\mcP_t^{\sf M}$ depend in a locally Lipschitz way on $\sf M$ with respect to the pseudo-distance ${\sf d}_\gamma$ on the space of models over $\mathscr{T}$ introduced in \eqref{EqPseudoDistanceModels}. We do not give the details here and refer the reader to the corresponding results in \cite{Hai14}, Theorem 3.10 and Theorem 5.12 therein.

\medskip

\begin{prop}
Given any time $t_0'<t_0(\bsh,\sf M)$, the restriction to $(0,t_0']\times\bbT^d$ of $\bsu$ defines locally a continuous function of $\bsh\in\mcD^{\gamma,\eta}(T_X,{\sf g})$ and the $\bf K$-admissible model $\sf M$.
\end{prop}

\medskip

Together with Theorem \ref{ThmRenormPDEs} in Section \ref{SectionrenormalizationMultiPreLie} below and Chandra \& Hairer's convergence result \cite{ChandraHairer} this continuity result allows to give meaning of the solution to a singular stochastic PDE as a limit in probability of solutions of renormalized equations driven by a mollified noise. This result holds more generally for all the equations that can be treated using regularity structures. Emphasize that this continuity result is fundamental. In a random setting where the noise is random and the models of interest are constructed as measurable functionals of the noise the continuity allows to transport automatically support theorems or large deviation results about random models into corresponding results about the solutions of the regularity structure lifts of the equations under study. See Hairer \& Sch\"onbauer's work \cite{HairerSchonbauer} on support theorems, Hairer \& Weber's work \cite{HairerWeber} on large deviation results, or Hairer \& Mattingly's work \cite{HairerMattingly} on the strong Feller property for solutions of singular stochastic PDEs, for a sample.   

\medskip

The last statement of this section makes the link between solving equation \eqref{EqgKPZ*} with a smooth noise $\zeta$ and the corresponding problem in the regularity structure equipped with the canonical model ${\sf M}^\zeta$ associated with the smooth noise. The latter is constructed in Section \ref{SectionMultiPreLie} and the only thing we presently need to know about it is that its reconstruction map $\textsf{\textbf{R}}^{{\sf M}^\zeta}$ is multiplicative with respect to the $\star$-product of modelled distributions, and sends the noise symbol $\Xi$ on the smooth function $\zeta$. For positive exponents $\gamma\in(-\beta_0,2)$ and $\eta\in(0,\beta_0+2]$, pick $v\in \mcC^\eta(\bbT^d)$ and denote by $P_\gamma v$ the lift in the polynomial structure of the heat propagator acting 
on $v$, defined by \eqref{eq heat propagator}.

\medskip

\begin{prop}   \label{PropRSCharacterizationSmoothPDE}
Let $\bsu\in\mcD_{(0,t_0)}^{\gamma,\eta}(T,\sf g)$ stand for the solution 
in a sufficiently small time interval $(0,t_0)$ of the fixed point problem
\begin{equation}   \label{EqRSgKPZ}
\bsu = P_\gamma v+\mcP_{t_0}^{{\sf M}^\zeta}\Big(f^\star(\bsu)\,\Xi + g^\star\big(\bsu,D\bsu\big)\Big).
\end{equation}
Then on the domain $(0,\frac{t_0}2)\times\bbT^d$, the function $u\defeq \textbf{\textsf{R}}^{{\sf M}^\zeta}(C_{t_0}\bsu)$ coincides with the solution to the well-posed equation \eqref{EqgKPZ*} with initial condition $v$.
\end{prop}

\medskip

\begin{Dem}
As in \eqref{eq KR=RK half plane}, the function $u$ satisfies the equation
\begin{align*}
u(x) = Pv(x)+\int_{(0,x_0)\times\bbT^d}P(x,y) \textbf{\textsf{R}}^{{\sf 
M}^\zeta} C_{t_0}\Big(f^\star(\bsu)\Xi + g^\star\big(\bsu,D\bsu\big)\Big) (y)\,dy,
\end{align*}
with $Pv$ the free propagation of the initial condition. We take advantage of the fact that ${\sf M}^\zeta$ is  a smooth model to write
$$
\textbf{\textsf{R}}^{{\sf M}^\zeta}(\bsw) (x) = {\sf \Pi}_x^\zeta\big(\bsw(x)\big)(x)
$$
for any modelled distribution $\bsw\in\mcD^{\alpha,\eta}(T,{\sf g}^\zeta)$ with $\alpha>0$ -- see identity \eqref{EqReconstructionFunctionLike}. Moreover for any $\bsf\in\mcD_{(0,t_0)}^{\gamma.\eta}(T,{\sf g}^\zeta)$, we have $C_{t_0}\bsf(x)=\bsf(x)$ if $x\in(0,\frac{t_0}2)\times\bbT^d$. Thus, in this domain, we can use the multiplicative character of the map $\textbf{\textsf{R}}^{{\sf M}^\zeta}$ and write 
$$
\textbf{\textsf{R}}^{{\sf M}^\zeta}\big(f(\bsu)\big) = f\big(\textbf{\textsf{R}}^{{\sf M}^\zeta}\bsu\big),\quad
\textbf{\textsf{R}}^{{\sf M}^\zeta}( D\bsu  ) = \partial_x\textbf{\textsf{R}}^{{\sf M}^\zeta}\bsu.
$$
as a consequence of Corollary \ref{CorEasyReconstruction} and Proposition \ref{prop abstract derivative}, and
\begin{align*}
\textbf{\textsf{R}}^{{\sf M}^\zeta}
\Big(f(\bsu)\Xi + g\big(\bsu, D\bsu  \big)\Big)
= f\Big(\textbf{\textsf{R}}^{{\sf M}^\zeta}\bsu\Big)\zeta
+g\Big(\textbf{\textsf{R}}^{{\sf M}^\zeta}\bsu,\partial_{x'}\textbf{\textsf{R}}^{{\sf M}^\zeta}\bsu\Big).
\end{align*}   
This finishes the proof.
\end{Dem}

\medskip

Arrived at that stage we have a model-dependent notion of solution $u^{\sf M}$ to 
the generalized (KPZ) equation, under the form
$$
u^{\sf M} = \textsf{\textbf{R}}^{\sf M}(\bsu^{\sf M}),\quad \bsu^{\sf M} = \Phi^{\bsh,\sf M}(\bsu),
$$
indexed by the set of $\bf K$-admissible models $\sf M$ on the regularity structure $\mathscr{T}$ associated with the equation. Theorem \ref{ThmSmoothAdmissibleModels} gives us a whole family of smooth $\bf K$-admissible models which we can use. However the $\bf K$-admissible models of interest are not smooth as we wish they satisfy the identity ${\sf \Pi}\Xi=\zeta$ for a non-smooth noise $\zeta$. The combinatorial structure of the elements of $\mathscr{T}$ detailled in Section \ref{SectionBuildingRS} allows to associate to any regularized version $\zeta_\epsilon$ of $\zeta$ a $\bf K$-admissible model ${\sf \Pi}^\epsilon$ such that ${\sf \Pi}^\epsilon\Xi=\zeta_\epsilon$, and ${\sf \Pi}^\epsilon$ is multiplicative for the $\star$-product on $T$. We will talk of ${\sf \Pi}^\epsilon$ as the naive interpretation map. However, these models diverge as $\epsilon>0$ goes to $0$. The tools needed to construct some $\epsilon$-dependent smooth models that have a limit as $\epsilon$ goes to $0$ are developed in the next section at the same level of generality as Section \ref{SectionBasicsRS} and Section \ref{SectionIntegration}. The so-called renormalization operation involved in the construction of these converging $\bf K$-admissible models will be given a dynamical meaning in Section \ref{SectionMultiAndrenormalizedEquations}. It is only after Theorem \ref{ThmRenormPDEs} in Section \ref{SectionMultiAndrenormalizedEquations} that we will be able to give an answer to the question ``{\it What dynamics does $u^{\sf M}$ follow?}'' when ${\sf M}$ is the renormalized naive model associated with the $\sf \Pi$-map ${^{k_\epsilon}}{\sf \Pi}^\epsilon$ introduced in the next section.

\bigskip

\section{Renormalization structures}
\label{SectionConcreteRenormStructure}

We introduce in this section the fundamental notion of renormalization structure, and a notion of compatibility between some regularity and renormalization structures. We emphasized in the previous paragraph that it is generically not possible to define a canonical $\bf K$-admissible model as a limit of the canonical $\bf K$-admissible models associated with regularized noises $\zeta_\epsilon$ if the noise(s) $\zeta$ is (are) not sufficiently regular. On a technical level, the non-convergence of the models ${\sf M}^\epsilon$ is related to the fact that the canonical model is defined by some intricate convolution of kernels that explode on the diagonal. Limit models need to be constructed by probabilistic means as limits in probability of models built from regularized noises, \emph{using a moving window}, as in Meta-Theorem $1$ in Section \ref{SectionIntro}. The implementation of this moving window picture involves the renormalization structures that we introduce in this section. Note that we do not need to know the details of the renormalization operation; the only properties that we need are encoded in the definition of a   renormalization structure and the compatibility condition with a regularity structure given below in Definition \ref{DefnRenormalizationStructure}. An example of renormalization structure will be given in Section \ref{SectionBuildingRS}, where the renormalization operation will be intimately related to the Taylor expansion procedure. 

\medskip

Renormalization structures are defined in Section \ref{SectionConcreteRenorS}. If we call the concrete regularity structures from Section \ref{SectionConcreteRS} right regularity structures, then renormalization structures 
$$
\mathscr{U}=\big((U,\delta),(U^-,\delta^-)\big)
$$ 
look like left regularity structures, with the difference that elements of the space $U^-$ have non-positive homogeneities. A fundamental notion of compatibility between renormalization and regularity structures is introduced in Section \ref{SubsectionCompatibleStructures}; it accounts for the fact that the renormalization operation induces a renormalization operation on $T^+$ and `commutes' with the recentering operators $\Delta$ and $\Delta^+$. This property allows to associate with each model $\sf M$ over $\mathscr{T}$ and each character $k$ on $U^-$ a new model $^k{\sf M}$ on $\mathscr{T}$. This is the main result of Section \ref{SubsectionCompatibleStructures}, Theorem \ref{ThmrenormalizationActionModels}. A large class of characters $k$ produces some $\bf K$-admissible models $^k{\sf M}$ if $\sf M$ is $\bf K$-admissible.

\bigskip

\subsection{Definition}
\label{SectionConcreteRenorS}

A renormalization structure is made up of two ingredients. First, it is a vector space $U$ with a basis whose elements are built by induction from elementary elements and multilinear operators giving new elements. The use of the symbol $\tau$ for a generic basis vector emphasizes this recursive, tree-like, definition. Each basis vector $\tau$ is a placeholder for a function $[\tau]$ from $(0,1]$ into a Banach space, typically $\bbR,\bbC$, a H\"older space or an algebra, whose structure as an element of the target space is encoded in the structure of $\tau$. In the cases of interest, the functions $[\tau]$ have no limit in $0^+$ and the basic problem is to remove in a `consistent' way the diverging pieces of these $[\tau]$ so as to end up with a collection of functions parametrized by $\epsilon>0$ having a limit where $\varepsilon$ goes to $0$. The functions $[\tau]$ are then said to have been renormalized. What `consistent' means is part of what follows.

\ssk

Roughly speaking, the basic operation for renormalizing a placeholder $\tau$ consists in removing from $\tau$ its different diverging pieces, in all possible sensible ways. This is the second ingredient of a renormalization structure. Tuples of pieces of elements of $U$ are not necessarily elements of $U$; we store them in a side space $U^-$. Endowing $U^-$ with an algebra structure allows to store the removed pieces of $\tau$ as an element of $U^-$ under the form of a product. We require nonetheless that any $\tau$ amputated from diverging pieces is an element of $U$; this is a restriction on which pieces of any $\tau\in U$ can be removed. We thus have a splitting map 
$$
\delta: U\rightarrow U^-\otimes U,
$$
with $\delta\tau$ the sum of all the elements from $U^-\otimes U$ corresponding to removing from $\tau$ all possible diverging allowed pieces, possibly several at a time. The removed pieces may themselves have diverging 
subpieces, and it makes sense to assume we have another splitting map 
$$
\delta^-: U^-\rightarrow U^-\otimes U^-,
$$
that extracts them on the left hand side of the tensor product $U^-\otimes U^-$. That the remaining piece is still in $U^-$ rather than in another 
space is a consistency requirement. 

\medskip

\begin{defn} \label{DefnRenormalizationStructure}
A \textsf{\textbf{renormalization structure}} is a pair of graded vector spaces  
$$
U \eqdef \bigoplus_{\beta\in B} U_\beta, \qquad U^- \eqdef \bigoplus_{\alpha\in B^-} U_\alpha^-  
$$
such that the following holds.   \vspace{0.1cm}

\begin{itemize}
\setlength{\itemsep}{0.1cm}
   \item The vector spaces $U_\alpha^-$ and $U_\beta$ are finite dimensional.

   \item The space $U^-$ is a connected graded bialgebra with unit ${\bf1}_-$, counit ${\bf1}_-'$, coproduct 
$$
\delta^-:U^-\to U^-\otimes U^-,
$$ 
and grading $B^-\subset(-\infty,0]$, with $0\in B^-$.

   \item The index set $B$ for $U$ is a locally finite subset of $\bbR$ bounded below. The space $U$ is a left comodule over $U^-$, that is $U$ is 
equipped with a splitting map $\delta : U \rightarrow U^-\otimes U$, which satisfies 
\begin{align} \label{U is comodule over U^-}
(\textrm{\emph{Id}}\otimes\delta)\delta = (\delta^-\otimes\textrm{\emph{Id}})\delta,\quad\textrm{and}\quad ({\bf1}_-'\otimes\iden)\delta=\iden.
\end{align}
Moreover, for any $\beta\in B$, one has
\begin{align}\label{EqSpecialRenorm}
\delta U_\beta \subset {\bigoplus_{\alpha\le0} } U_{\alpha}^-\otimes U_{\beta-\alpha}.
\end{align}
\end{itemize}
We denote by 
$$
\mathscr{U} \defeq  \Big((U,\delta), (U^-,\delta^-)\Big)
$$
a renormalization structure.
\end{defn} 

\medskip

Similarly to the regularity structure, let $\mcU_\alpha^-$  and $\mcU_\beta$ be bases of $U_\alpha^-$ and $U_\beta$, respectively, and set
$$
\mcU^- \defeq  \bigcup_{\alpha\in B^-} \mcU_\alpha^-, \qquad \mcU \defeq  \bigcup_{\beta\in B} \mcU_\beta.
$$ 
Note that, unlike in the definition of a concrete regularity structure satisfying Assumption \refA{A1}, we do not require that $\mcU_0$ is one dimensional in the definition of a renormalization structure. Since all $\alpha\in B^-$ are non-positive, one has $\beta-\alpha\geq \beta$ in \eqref{EqSpecialRenorm}. Proposition \ref{PropBialgebraHopf} in Appendix {\sf \ref{SectionAppendixAlgebra}} can be applied to the negative grading $B^-$ of $U^-$, and says that $U^-$ is a Hopf algebra; we denote by $S_-$ its antipode. Choosing a basis of $U^-$ provides an associated decomposition of $\delta\tau$ of the form
\begin{align}\label{section 5: trianglelefteq}
\delta\tau \eqdef \sum_{\varphi\trianglelefteq\tau} \varphi\otimes \tau/^-\varphi,
\end{align}
where the $\varphi$ are distinct elements of the chosen basis. The notation $\varphi\trianglelefteq\tau$ means that $\varphi$ is a basis element that appears as one of the left hand side members of the finite sum giving $\delta\tau$. We call $\delta$ a \textit{renormalization splitting} and fix throughout a basis of $U^-$. The results we prove in the sequel do not depend on that arbitrary choice. Similarly to what we saw in Section \ref{SectionConcreteRS} for the Hopf algebra $(T^+,\Delta^+)$, the $\delta^-$ splitting of the Hopf algebra $(U^-,\delta^-)$ induces a convolution group law on the set $G^-$ of characters on $U^-$
$$
(k_1*k_2)\tau
\defeq  (k_1\otimes k_2)\delta^-\tau, \qquad (\tau\in U^-).
$$
The inverse of a character $k$ for the convolution product is explicit and given by $k\circ S_-$. Given a character $k$ on $U^-$, we define a linear map $\widetilde{k}: U\rightarrow U$, setting   
$$
\widetilde{k} \defeq  (k\otimes\textrm{Id})\delta.
$$
The group $G^-$ acts on $U$ from right. Indeed, as a direct consequence of the comodule property in \eqref{U is comodule over U^-}, one has
\begin{align*}
\widetilde{k_1*k_2} &= ((k_1*k_2)\otimes\iden)\delta
=(k_1\otimes k_2\otimes\iden)(\delta^-\otimes\iden)\delta\\
&=(k_1\otimes k_2\otimes\iden)(\iden\otimes\delta)\delta=
(k_1\otimes\widetilde{k_2})\delta\\
&=\widetilde{k_2} \circ \widetilde{k_1}
\end{align*}
for any $k_1,k_2\in G^-$.

\bigskip

\subsection{Compatible renormalization and regularity structures}
\label{SubsectionCompatibleStructures}

We introduce a `compatibility' property between regularity and renormalization structures. We use the notations from Appendix {\sf \ref{SectionAppendixAlgebra}}. In particular, given an algebra $A$ and two spaces $E,F$, we define a linear map $\mcM^{(13)}$ from the algebraic tensor product $A\otimes E\otimes A\otimes F$ to the algebraic tensor product $A\otimes E\otimes F$ setting
$$
\mcM^{(13)}\Big(a_1\otimes e\otimes a_2\otimes f\Big) \defeq  (a_1a_2)\otimes e\otimes f.
$$
Recall we write $\mathscr{T}=\big((T^+,\Delta^+),(T,\Delta)\big)$ for a 
regularity structure and $S_+$ for the antipode map on $T^+$.

\medskip

\begin{defn}   \label{*DefnCompatibility}
A regularity structure $\mathscr{T}$ is said to be \textbf{\textsf{compatible}} with a   renormalization structure $\mathscr{U}$ if the following three compatibility conditions hold true.
\begin{enumerate}
\setlength{\itemsep}{0.5cm}

   \item[\textsf{\textbf{(a)}}] 
The spaces $T$ and $U$ coincide as linear spaces and the bases $\mcB$ and 
$\mcU$ coincide. (Each element $\tau\in\mcB$ is in particular homogeneous 
in both $T$ and $U$, but it may belong to $\mcB_{\beta_1}$ and $\mcU_{\beta_2}$ with $\beta_1\neq\beta_2$.) Moreover,
\begin{equation}   \label{delta T to U-T}
\delta T_\beta \subset U^-\otimes T_\beta,\quad \textrm{ for all }\beta\in A.
\end{equation}   
   
   \item[\textsf{\textbf{(b)}}] There exists an algebra morphism 
$$
\delta^+:T^+\to U^-\otimes T^+
$$ 
such that
\begin{equation} \label{T+ is a comodule over U^-}
\big(\textrm{\emph{Id}}\otimes\delta^+\big)\delta^+ = \big(\delta^-\otimes\textrm{\emph{Id}}\big)\delta^+,\quad\textrm{and}\quad
\big({\bf1}_-'\otimes\iden\big)\delta^+ = \iden
\end{equation}
and
\begin{equation}   \label{delta+ T+ to U-T+}
\delta^+ T_\alpha^+ \subset U^-\otimes T_\alpha^+, \quad \textrm{ for all 
}\alpha\in A^+.
\end{equation}   
   
   \item[\textsf{\textbf{(c)}}]  The compatibility conditions
\begin{align}\label{*EqCompatibilityCondition}
\big(\iden \otimes \Delta^{(+)}\big)\delta^{(+)} = \mcM^{(13)}\big(\delta^{(+)}\otimes\delta^+\big)\Delta^{(+)}
\end{align}
and
\begin{align}\label{*EqCompatibilityCondition2}
\big(\iden\otimes{\bf1}_+'\big)\delta^+={\bf1}_+'(\cdot){\bf1}_-
\end{align}
hold.
\end{enumerate}
\end{defn} 

\medskip

Emphasize the fact that the homogeneity notion in $T$ captures the notion of regularity of the associated analytic objects encoded by elements of $T$ while the homogeneity notion in $U$ captures the diverging behavior of the corresponding regularized objects, as the regularization parameter goes to $0$. It makes sense that the two notions of homogeneities are unrelated. Definition \ref{*DefnCompatibility} also captures the fact that the renormalization procedure encoded in $\mathscr{U}$ induces a renormalization operation on $T^+$ and commutes with the recentering operators $\Delta$ and $\Delta^+$. We will see 
in Proposition \ref{prop simplification of compatibility} below that the six conditions from Definition \ref{*DefnCompatibility} hold iff condition \eqref{delta T to U-T} and condition \eqref{*EqCompatibilityCondition}, in its form without the $+$ labels, hold, under a reasonable assumption on $\delta^+$ that holds true for the regularity and renormalization structures associated with (systems of) 
singular stochastic PDEs.

Compare conditions \eqref{delta T to U-T} and \eqref{EqSpecialRenorm}. Emphasize here as in item \textsf{\textbf{(a)}} that the notion of homogeneity is relative to the grading used to define it. An element of $T=U$ may thus have different homogeneities, depending on whether it is considered as an element of $T$ or $U$. By condition \textsf{\textbf{(a)}}, the space $T$ is a left $U^-$-comodule. The map $\delta^+$ in \textsf{\textbf{(b)}} accounts for the effect in $T^+$ of the renormalization process. By \eqref{T+ is a comodule over U^-}, the space $T^+$ is also a left $U^-$-comodule. Hence for given a character $k$ on $U^-$, we can define linear maps $\widetilde k:T\to T$ and $\widetilde k^+:T^+\to T^+$, by
$$
\widetilde k = \big(k\otimes\iden\big)\delta,\quad\textrm{and}\quad 
\widetilde k^+ = \big(k\otimes\iden\big)\delta^+.
$$
Properties \eqref{delta T to U-T} and \eqref{delta+ T+ to U-T+} ensure that homogeneities of elements of $T$ and $T^+$ are stable under these actions. 
Condition \textbf{\textsf{(c)}}, read with the $+$ labels, somehow says that the renormalization operation encoded in $\widetilde{k}$ commutes with the Taylor expansion operation on the coefficients of any modelled distribution, encoded in $\Delta^+$. Condition \textbf{\textsf{(c)}}, read without the $+$ labels, says something similar for modelled distributions. Note that the Hopf algebra $T^+$ is a left $U^-$-comodule bialgebra. By Proposition \ref{*compatibility antipode}, we have the following compatibility condition on the antipode
\begin{align}\label{*EqCompatibilityAntipode}
\delta^+ \circ S_+ = \big(\iden\otimes S_+\big)\circ \delta^+.
\end{align}

\ssk

Recall that given a model $\sf  M = (g,\Pi)$ on $\mathscr{T}$, the anchored interpretation operator ${\sf \Pi}_x^{\sf g}$ associated with $\sf M$ is given for any $x\in\bbR^d$, by 
$$
{\sf \Pi}_x^{\sf g} = ({\sf \Pi}\otimes {\sf g}_x^{-1})\Delta.
$$
The next statement and its proof are part of Theorem 6.15 in Bruned, Hairer and Zambotti's work \cite{BHZ} on the algebraic renormalization of regularity structures. It tells us that the $\widetilde{k}$ and $\widetilde{k}^+$ maps have jointly a natural and simple action on the space of models on $\mathscr{T}$.

\medskip

\begin{thm}   \label{ThmrenormalizationActionModels}
Let a renormalization structure $\mathscr{U}=(U,U^-)$ be compatible with a regularity structure $\mathscr{T}=(T^+,T)$. Given any character $k$ 
on $U^-$, and any model ${\sf M = (g,\Pi)}$ on $\mathscr{T}$, define $^k{\sf M}=({}^k{\sf g}, {}^k{\sf\Pi})$, on $\mathscr{T}$ setting 
\begin{equation}   \label{EqrenormalizedModel}
^k{\sf M} \defeq  \Big({\sf g}\circ \widetilde{k}^+, {\sf \Pi}\circ \widetilde{k}\Big).
\end{equation}
One has
\begin{equation}   \label{EqgRenorm}
\big({\sf g}_y\circ\widetilde{k}^+\big)*\big({\sf g}_x\circ\widetilde{k}^+\big)^{-1} = {\sf g}_{yx}\circ\widetilde{k}^+,
\end{equation}
and 
\begin{equation}   \label{EqPiXRenorm}
\Big(\big({\sf \Pi}\circ\widetilde{k}\big)\otimes\big({\sf g}_x\circ\widetilde{k}^+\big)^{-1}\Big)\Delta = {\sf \Pi}_x^{\sf g}\circ\widetilde{k},
\end{equation}
for any $x,y\in\bbR^d$. Moreover, the size conditions \eqref{EqEstimateGammayx} and \eqref{EqestimatePix} hold for $^k{\sf M}=({}^k{\sf g},{}^k{\sf\Pi})$, so $^k{\sf M}$ is a model.
\end{thm} 

\medskip

\begin{Dem} 
The proof is short and simple because the notion of compatibility between some regularity and renormalization structures is tailored for that purpose. One has 
\begin{equation*}
\begin{split}
\big({\sf g}_y\circ\widetilde{k}^+\big)*\big({\sf g}_x\circ\widetilde{k}^+\big)^{-1} 
&= \big((k\otimes {\sf g}_y)\delta^+\big)\otimes\Big((k\otimes {\sf g}_x) \delta^+\circ S_+\Big)\Delta^+   \\
&\hspace{-0.15cm}\overset{\eqref{*EqCompatibilityAntipode}}{=} 
\big((k\otimes {\sf g}_y)\delta^+\big)\otimes\big((k\otimes {\sf g}_x^{-1})\delta^+\big)\Delta^+   \\
&= \big({\sf g}_y\otimes {\sf g}_x^{-1}\big)\circ\big(\widetilde{k}^+\otimes\widetilde{k}^+\big)\Delta^+   \\
&\hspace{-0.15cm}\overset{\eqref{*EqCompatibilityCondition}}{=} 
\big({\sf g}_y\otimes {\sf g}_x^{-1}\big)\circ(k\otimes\Delta^+)\delta^+  
 \\
&= {\sf g}_{yx}\circ\widetilde{k}^+.
\end{split}
\end{equation*}
and 
\begin{equation*}
\begin{split}
\Big(\big({\sf \Pi}\circ\widetilde{k}\big)\otimes\big({\sf g}_x\circ\widetilde{k}^+\big)^{-1}\Big)\Delta \ \ 
&\hspace{-0.15cm}\overset{\eqref{*EqCompatibilityAntipode}}{=} 
\Big((k\otimes{\sf \Pi})\delta \otimes (k\otimes {\sf g}_x^{-1})\delta^+\Big)\Delta   \\
&= \big({\sf \Pi}\otimes {\sf g}_x^{-1}\big)\circ\big(\widetilde{k}\otimes\widetilde{k}^+\big)\Delta   \\
&\hspace{-0.15cm}\overset{\eqref{*EqCompatibilityCondition}}{=} 
\big({\sf \Pi}\otimes {\sf g}_x^{-1}\big)\circ(k\otimes\Delta)\delta   \\
&= {\sf\Pi}^{\sf g}_x\circ\widetilde{k}.
\end{split}
\end{equation*}
The size conditions \eqref{EqEstimateGammayx} and \eqref{EqestimatePix} on $^k{\sf M}$ follow now from formulas \eqref{EqgRenorm} and \eqref{EqPiXRenorm}, and from the fact that the maps $\widetilde{k}$ and $\widetilde{k}^+$ preserve the spaces $T_\beta$ and $T_\alpha^+$, respectively, as a consequence of the stability conditions \eqref{delta T to U-T} and \eqref{delta+ T+ to U-T+}.
\end{Dem}

\medskip

Together with Corollary \ref{Cor reconstruction of smooth model} this statement implies in particular that if the model $\sf M$ takes values in the space of {\it continuous functions} then the reconstruction operator $^{k}\textbf{\textsf{R}}$ associated with the renormalized model is related 
to the reconstruction operator $\textbf{\textsf{R}}$ associated with the unrenormalized model by the relation 
$$
^{k}\textbf{\textsf{R}} = \textbf{\textsf{R}}\circ\widetilde{k}.
$$
This point will be used crucially in Section \ref{SectionrenormalizationMultiPreLie}, where we will give a dynamical picture of the renormalization of models.   

\medskip

We consider in the remainder of this section the case of interest for the 
study of (systems of) singular stochastic PDE(s) where the regularity structure $\mathscr{T}$ is built from integration operators and satisfies Assumption \REFB. Unfortunately, even if a model $\sf M$ is ${\bf K}$-admissible, ${}^k{\sf M}$ is not always ${\bf K}$-admissible for any $k\in G^-$. We put forward an assumption under which one builds ${\bf K}$-admissible models using elements $k$ of a non-trivial subgroup $G_\textrm{ad}^-$ of $G^-$. Assume $\mcB=\mcU$ and let $\mcF$ stand for the family of operators
$$
\mcF \defeq  \{\mcI_p\}_{|p|_{\frak{s}}\le1}\cup\{X^n\star\}_{n\in\bbN^{d+1}\setminus\{0\}}, 
$$
acting on the basis $\mcB=\mcU$, where $X^n\star$ denotes the linear operator on $T$ defined by $\tau\mapsto X^n\star\tau$. Recall that such a multiplication is always given by Assumption \refA{A3}. 

\medskip

\begin{assumC}\label{C1}
The regularity structure $\mathscr{T}$ is built from integration operators and satisfies Assumption \refB{B1} and the renormalization structure $\mathscr{U}$ is compatible with $\mathscr{T}$. Moreover, the following holds.
\begin{itemize}
\setlength{\itemsep}{0.1cm}
	\item The algebra $U^-$ is generated by the basis elements $\mcU_{<0}\defeq \bigcup_{\alpha<0}\mcU_\alpha$ and the unit $\textbf{\textsf{1}}_-$.   
	\item 
	Let $\frak{J}^-$ be the ideal of $U^-$ generated by the set $\big(\mcF(\mcU)\big)\cap\mcU_{<0}$. The linear map $\delta:U\to U^-\otimes U$ satisfies, for any operator $F\in\mcF$ and $\tau\in\mcU$,
	\begin{align}\label{eq relation delta and K}
	\delta(F\tau)-(\text{\rm Id}\otimes F)\delta\tau\in \frak{J}^-\otimes U.
	\end{align}
	
	\item 
	We define a projection operator $P_-:U\to U^-$ setting $P_-\tau \defeq  \tau{\bf1}_{\tau\in\mcU_{<0}}$, for any $\tau\in\mcU$. The linear map $\delta^-:U^-\to U^-\otimes U^-$ is defined by $\delta^- = (\text{\rm Id}\otimes P_-)\delta$ on $\mcU_{<0}$ and its multiplicative extension.
\end{itemize}
\end{assumC}      

\medskip

Define the subset $G_{\text{\rm ad}}^-$ of $G^-$ by
$$
G_{\text{\rm ad}}^- \defeq  \Big\{k\in G^-\, ;\, k\big(F(\tau)\big)=0\ \text{for any}\ F(\tau)\in\big(\mcF(\mcU)\big)\cap\mcU_{<0}\Big\}.
$$
\medskip

\begin{prop}
The set $G_{\text{\rm ad}}^-$ is a subgroup of $G^-$, and for any $k\in G_{\text{\rm ad}}^-$ and any $\bf K$-admissible model $\sf M$, one has ${}^k{\sf M}$ is also $\bf K$-admissible. The group $G^-_\textrm{\emph{ad}}$ 
is called the \textbf{\textsf{renormalization group}}.
\end{prop}

\medskip

The definition of the group $G_{\text{\rm ad}}^-$ gives the meaning to assumption \eqref{eq relation delta and K}. Up to irrelevant terms for $k\in G_{\text{\rm ad}}^-$, the renormalization operations in $U^-$ or $U$ of 
an `integral' is the integral of its renormalized integrand, and multiplication by a polynomial has no effect on the renormalization process.   

\medskip

\begin{Dem}
Note that $k(\frak{J}^-)=0$ for any $k\in G_{\text{\rm ad}}^-$. Let $\tau$ be an element of $\mcU_\alpha$ such that $F\tau\in\mcU_{<0}$ for some $F\in\mcF$.

\ssk

\textit{\textsf{(a)}}
Given $k,h\in G_{\text{\rm ad}}^-$ since identity \eqref{eq relation delta and K} and the third property of Assumption \refC{C1} ensures that $(k*h)(F\tau)=(k\otimes h)\delta^-(F\tau)=0$, for all $F\in\mcF$, we have $k*h\in G_{\text{\rm ad}}^-$. Next we show that $k^{-1}=k\circ S_-\in G_{\text{\rm ad}}^-$. Denote by $\mcM^-$ the multiplication operator in $U^-$ and pick $\sigma\in U_\beta$ with $\beta<0$. Since $\delta\sigma\in {\bf1}_-\otimes\sigma+ \sum_{\alpha<0}U_{\alpha}^-\otimes U_{\beta-\alpha}$, by applying the operator $\mcM^-(\text{\rm Id}\otimes S_-P_-)$ to \eqref{eq relation delta and K} we have from Assumption \textsf{\textbf{(B1)}} and the fact that $\frak{J}^-$ is an ideal
$$
S_-(F\sigma)\in \sum_{\alpha<0}\mcM^-\Big(U_{\alpha}^-\otimes S_-\big(P_-FU_{\beta-\alpha})\big)\Big)+\frak{J}^-,
$$
which implies $k^{-1}(FU_\beta)=k(S_-FU_\beta)=0$, by an induction on $\beta$.

\ssk

\textit{\textsf{(b)}}
Let $F=\mcI_p$. By \eqref{eq relation delta and K},
$$
{}^k{\sf \Pi}(\mcI_p\tau)
=(k\otimes{\sf \Pi})\delta \mcI_p\tau
=(k\otimes{\sf \Pi}\mcI_p)\delta\tau
=\partial^p{\bf K}(k\otimes{\sf \Pi})\delta\tau
=\partial^p{\bf K}({}^k{\sf \Pi}\tau).
$$
Since we have a similar identity for $F=X^n\star$, we obtain that ${}^k{\sf M}$ is admissible.
\end{Dem}   

\medskip

We end this section by showing that the definition of compatible renormalization and regularity structures takes then a simple form under the following additional mild assumption. It essentially says that multiplications by a polynomial and integrations are not the sources of renormalization problems.

\ssk

\begin{assumC}\label{C2}
The algebra morphism $\delta^+:T^+\to U^-\otimes T^+$, is determined by the identities
\begin{align}\label{eq relation delta and delta^+}
\delta^+\Xplus^\ell={\bf1}_-\otimes \Xplus^\ell,\qquad
\delta^+(\mcI_p^+\tau)=(\text{\rm Id}\otimes\mcI_p^+)\delta\tau.
\end{align}
\end{assumC}

\ssk

\begin{prop}\label{prop simplification of compatibility}
Under Assumption \refC{C2}, assume $\mathscr{T}$ satisfies the property \eqref{delta T to U-T} and the version of identity \eqref{*EqCompatibilityCondition} without the $+$ labels. Then the other conditions in Definition \ref{*DefnCompatibility} follow automatically.
\end{prop}

\medskip

\begin{Dem}
The comodule property \eqref{T+ is a comodule over U^-} follows from \eqref{U is comodule over U^-} and the definition \eqref{eq relation delta and delta^+}.
Indeed,
\begin{align*}
\big(\iden\otimes\delta^+\big)\delta^+(\mcI_n^+\tau)
&=\big(\iden\otimes\delta^+\mcI_n^+\big)\delta\tau
=\big(\iden\otimes\iden\otimes\mcI_n^+\big)(\iden\otimes\delta)\delta\tau   \\
&=\big(\iden\otimes\iden\otimes\mcI_n^+\big)(\delta^-\otimes\iden)\delta\tau
=\big(\delta^-\otimes\iden\big)\big(\iden\otimes\mcI_n^+\big)\delta\tau 
  \\
&=\big(\delta^-\otimes\iden\big)\delta^+(\mcI_n^+\tau).
\end{align*}
The counit part of \eqref{T+ is a comodule over U^-} and \eqref{*EqCompatibilityCondition2} are left to readers.
The condition \eqref{delta+ T+ to U-T+} follows from \eqref{delta T to U-T} and the definition \eqref{eq relation delta and delta^+}.
The $(+)$-labelled version of \eqref{*EqCompatibilityCondition} is checked for $\mcI_n^+\tau\in\mcB^+$ as follows.
\begin{align*}
\mcM^{(13)}\big(\delta^+\otimes\delta^+\big)\Delta^+(\mcI_n^+\tau)
&= \mcM^{(13)}\left(\big(\delta^+\mcI_n^+\otimes\delta^+\big)\Delta\tau
+\sum_{\ell\in\bbN\times\bbN^d}\delta^+\frac{\Xplus^\ell}{\ell!}\otimes\delta^+(\mcI_{n+\ell}^+\tau)\right)   \\
&= \mcM^{(13)}\left(\big((\iden\otimes\mcI_n^+)\delta\otimes\delta^+\big)\Delta\tau
+\sum_{\ell\in\bbN\times\bbN^d}{\bf1}_-\otimes\frac{\Xplus^\ell}{\ell!}\otimes(\iden\otimes\mcI_{n+\ell}^+)\delta\tau\right)   \\
&= \big(\iden\otimes\mcI_n^+\otimes\iden\big)\mcM^{(13)}(\delta\otimes\delta^+)\Delta\tau
+\sum_{\ell\in\bbN\times\bbN^d, \varphi\trianglelefteq\tau}\varphi\otimes\frac{\Xplus^\ell}{\ell!}\otimes\mcI_{n+\ell}^+(\tau/^-\varphi)   \\
&= \big(\iden\otimes\mcI_n^+\otimes\iden\big)(\iden\otimes\Delta)\delta\tau
+\sum_{\ell\in\bbN\times\bbN^d, \varphi\trianglelefteq\tau}\varphi\otimes\frac{\Xplus^\ell}{\ell!}\otimes\mcI_{n+\ell}^+(\tau/^-\varphi),
\end{align*}
and
\begin{align*}
\big(\iden\otimes\Delta^+\big)\delta^+(\mcI_n^+\tau)
&=(\iden\otimes\Delta^+\mcI_n^+)\delta\tau   \\
&=\sum_{\varphi\trianglelefteq\tau}\varphi\otimes
\left((\mcI_n^+\otimes\iden)\Delta(\tau/^-\varphi)
+\sum_{\ell\in\bbN\times\bbN^d}\frac{\Xplus^\ell}{\ell!}\otimes\mcI_{n+\ell}^+(\tau/^-\varphi)\right),
\end{align*}
hence we have
$$
\mcM^{(13)}\big(\delta^+\otimes\delta^+\big)\Delta^+(\mcI_n^+\tau) =
\big(\iden\otimes\Delta^+\big)\,\delta^+(\mcI_n^+\tau).
$$
It remains to prove \eqref{T+ is a comodule over U^-} and \eqref{*EqCompatibilityCondition} for elements of the form $X_+^n\tau$. This is elementary using the multiplicative property of $S^+$ and $\Delta^+$.
\end{Dem}

\medskip

\section{Multi-pre-Lie structure and renormalized equations}
\label{SectionMultiAndrenormalizedEquations}

Let us summarize the successive steps that we have followed after the formalism of regularity structures was set up in Section \ref{SectionBasicsRS}. We described in Section \ref{SectionIntegration} a particular class of regularity structures, and the class of admissible models on them, that have the property that one can lift the integral operator $(\partial_t-\Delta_x)^{-1}$ into an operator on some spaces of modelled distributions that has a Schauder-type continuity property given in Theorem \ref{thm wellposedness of K}. This result played a crucial role in the local in time well-posedness result proved in Theorem \ref{ThmWellPosedness} of Section \ref{SectionSolvingPDEs}. It gives us a modelled distribution $\bsu^{\sf M}$ that solves a well-defined regularity structure formulation of a fixed point formulation of an ill-defined singular (stochastic) PDE. We define the (model-dependent) solution of this ill-defined equation as the (model-dependent) reconstruction $u^{\sf M}= {\sf R}^{\sf M}(\bsu^{\sf M})$ of $\bsu^{\sf M}$. To make this definition consistent with the initial objective we would like to use some models $\sf M$ for which ${\sf \Pi}(\Xi)=\zeta$. The construction of some admissible model that has this property is made very non-trivial by the fact that $\zeta$ has low regularity. As a matter of fact, this cannot be done in a deterministic reasonable way but one can construct some random admissible models that are limits in a probabilistic sense of some smooth models built from the canonical admissible model ${\sf M}^\epsilon$ associated with a regularized noise $\zeta_\epsilon$. The construction recipe $({\sf M},k)\mapsto {}^k{\sf M}$ for these renormalized models was given in Theorem \ref{ThmrenormalizationActionModels} in Section \ref{SectionConcreteRenormStructure}. Denote by ${\sf M}^\epsilon=({\sf g} ^\epsilon, {\sf \Pi}^\epsilon)$ the canonical admissible model associated with a regularized noise $\zeta_\epsilon$. We will see in Section \ref{SectionBHZCharacter} that there is a particular choice of character $k_\epsilon$ for which the ${}^{k_\epsilon}{\sf \Pi}^\epsilon(\tau)(x)$ are centered for all the $\tau$ of negative homogeneity and all state space point $x$. Chandra \& Hairer first proved the probabilistic convergence of the renormalized admissible models ${}^{k_\epsilon}{\sf M}^\epsilon$ to some limit admissible random model $\overline{\sf M}$. We will not prove this result in the tourist guide and refer the reader to the review \cite{BHReview} for some information on this matter. Rather we will see in the present section that one can give a somewhat explicit description of the dynamics of $u^{\overline{\sf M}}$ based on the following facts.
\begin{enumerate}
	\item The continuity of the map ${\sf M}\mapsto u^{\sf M}$ ensures that the dynamics of $u^{\overline{\sf M}}$ is the limit of the dynamics of the $\overline{u}_\epsilon \defeq u^{{}^{k_\epsilon}{\sf M}^\epsilon}$.
	
	\item We will see in the present section that $\overline{u}_\epsilon$ is actually the solution of an explicit stochastic PDE driven by the regularized noise $\zeta_\epsilon$, called the {\it renormalized equation}. This is the main result of this section, stated in Theorem \ref{ThmRenormPDEs}.
\end{enumerate}

We will concentrate in this section on the study of the generalized (KPZ) equation
\begin{equation}   \label{EqgKPZ}   \begin{split}
\big(\partial_{x_0}-\Delta_{x'}+1\big) u &= f(u)\zeta + g_2(u)(\partial_{x'} u)^2 + g_1(u)(\partial_{x'} u) +g_0(u)  \\
&= f(u)\zeta +g(u,\partial_{x'} u),
\end{split} \end{equation}
with a given initial condition. It already involves the main difficulties 
of the most general situation, with the advantage of leaving aside a number of purely technical and notational matters compared to the most general situation. 

\ssk

\noindent {\it 1. Picard iteration and decorated trees.} We saw in Section \ref{SectionSolvingPDEs} that a there is a unique modelled distribution 
$$
\bsu = \sum_{\tau\in\mcB} u_\tau \tau \in \mcD_{(0,t_0)}^{\gamma,\eta}(T,\sf g),
$$ 
with $\gamma\in(-\beta_0,2)$ and $\eta\in(0,\beta_0+2]$ solving the lift \eqref{EqPDERSPrime} in the regularity structure $\mathscr{T}$ of equation \eqref{EqgKPZ}. It satisfies on the domain $(0,\frac{t_0}2)\times\bbT^d$ the fixed point problem
\begin{equation}   \label{EqIdentification1} \begin{split}
\bsu &\simeq \mcI\Big(f^\star(\bsu)\Xi + g_2^\star(\bsu) \star ( D\bsu  )^{\star2} + g_1^\star(\bsu) \star D\bsu   + g_0^\star(\bsu) \Big)\\
   &\simeq \frac{f^{(k)}(u)}{k!}\,u_{\tau_1}\cdots u_{\tau_k}\,\mcI\big(\tau_1\cdots\tau_k \Xi\big)   
   + \frac{g_2^{(k)}(u)}{k!}\,u_{\tau_1}\cdots u_{\tau_k}u_{\sigma_1}u_{\sigma_2} \mcI\big(\tau_1\cdots\tau_kD_i\sigma_1D_j\sigma_2\big)\\
   &\quad+ \frac{g_1^{(k)}(u)}{k!}\,u_{\tau_1}\cdots u_{\tau_k}u_{\sigma_1}\,\mcI\big(\tau_1\cdots\tau_kD_i\sigma_1\big)
   + \frac{g_0^{(k)}(u)}{k!}\,u_{\tau_1}\cdots u_{\tau_k}\,\mcI\big(\tau_1\cdots\tau_k\big),
\end{split}   \end{equation}
up to some model-dependent non-trivial polynomial components, with $\tau_k, \sigma_\ell\in\mcB$, and implicit sums over $\mcB$ and $i,j\in\{1,\dots,d\}$. We see on this identity that $T$ needs at least to be stable by the operations 
$$
\big(\tau_1,\dots,\tau_k,\sigma_1,\sigma_2\big) \mapsto
\mcI\big(\tau_1\cdots\tau_k\big),\
\mcI\big(\tau_1\cdots\tau_k\Xi\big),\
\mcI\big(\tau_1\cdots\tau_kD_i\sigma_1\big),\
\mcI\big(\tau_1\cdots\tau_kD_i\sigma_1D_j\sigma_2\big);
$$
this naturally endows the elements of $T$ with a tree/inductive structure. This fact is common to all the equations that can be treated by the methods of regularity structures. This leads us in Section \ref{SectionMultiPreLie} to set the framework of {\it rooted decorated trees} as a convenient encoding of the elements of $T$. 

\ssk

\noindent {\it 2. Decorated trees and pre-Lie algebras.} The importance of this algebraic setting comes from the fact that the vector space $V$ spanned by the set of all rooted trees with vertex and edge decorations in some given sets happens to be a universal object in a class of algebraic structures called multi-pre-Lie algebras. Morphisms of such multi-pre-Lie algebras defined on $V$ are thus determined by their restrictions to a set of generators. We show in Section \ref{SectionCoherence} that the modelled distribution solution of the regularity structure lift of equation \eqref{EqgKPZ} involves precisely such a morphism, with values in the space of vector fields; see Proposition \ref{PropPreLieMorphism}. 

\ssk

\noindent {\it 3. Pre-Lie algebras and renormalization.} The regularity structure associated with equation \eqref{EqgKPZ} is built from $V$, with $T$ and $T^+$ subsets of $V$. Building $\mathscr{T}$ within $V$, any renormalization structure $\mathscr{U}$ compatible with $\mathscr{T}$ and satisfying Assumption \textbf{\textsf{(C)}} will also be built within $V$, with $U$ and $U^-$ some subsets of $V$. Theorem \ref{ThmRenormPDEs} below shows that $\overline{u}_\epsilon$ is the solution of an explicit equation driven by $\zeta_\epsilon$. This result was first proved in the seminal work \cite{BCCH18} of Bruned, Chandra, Chevyrev and Hairer. The proof builds on the fact that the dual renormalization map $\widetilde{k}^*$ that one can associate to any $k\in G^-_\textrm{ad}$ happens to be a multi-pre-Lie morphism under a {\it compatibility condition} on the multi-pre-Lie structure and the renormalization operator $\delta$, found here under the form of Assumption \textbf{\textsf{(D3)}}. 

\medskip

Assumptions \textbf{\textsf{(D1-D3)}} to be found in this section are 
all met in the case of a general subcritical system of singular stochastic PDEs, and we verify them by hand in Section \ref{SectionBuildingRS} where we construct the regularity and renormalization structures associated with the generalized (KPZ) equation. We emphasize them here as `assumptions' to stress the mechanics at work in the most general case.

\bigskip

\subsection{Free $\sf E$-multi-pre-Lie algebra generated by $\sf N$}
\label{SectionMultiPreLie}

We introduce in the first paragraph the space of edge and node decorated trees. Decoration spaces are associated to any given system of singular stochastic PDEs, and the associated space of decorated trees provides the background scene from which one can define the regularity and renormalization structures associated with the system. The multi-pre-Lie structure of the space of decorated trees is introduced in another paragraph and its dual operator described explicitly.

\bigskip

\hfill \textcolor{black}{\textsf{\textbf{ {\Large \S} Decorated trees}}}

\medskip

\begin{defn*}
Let $\frak{T}_{\rm n}$ (called a \emph{node type set}) and $\frak{T}_{\rm 
e}$ (called an \emph{edge type set}) be abstract sets.
\begin{itemize}
\setlength{\itemsep}{0.1cm}

   \item A \textbf{\textsf{rooted tree}} $\tau$ is a finite connected non-planar graph without loops, with a node set $N_\tau$ and an edge set $E_\tau$, and with a distinguished node $\rho_\tau$, called the \textbf{\textsf{root}}. The root defines a natural order on each edge, from the root to the leaves. In particular, each edge $e\in E_\tau$ is written as the form $e=(u,v)$, where $u,v\in N_\tau$ are endpoints of $e$ and $u$ is closer to the root. $u$ is called a parent of $v$, and $v$ is called a child of $u$.
   
We identify two trees $\tau$ and $\sigma$ if they are graph isomorphic, so we always write a graph by putting ancestors lower and descendants upper. The root is put at the bottom. Here is an example.
$$
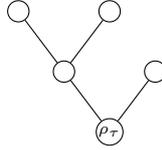
\begin{tikzpicture}
\coordinate (A1) at (0,0);
\coordinate (A2) at (-0.6,0.8);
\coordinate (A3) at (0.6,0.8);
\coordinate (A4) at (0,1.6);
\coordinate (A5) at (-1.2,1.6);
\draw (A1)--(A5);
\draw (A1)--(A3);
\draw (A2)--(A4);
\foreach \n in {2,3,4,5} \filldraw[white] (A\n) circle (4pt);
\foreach \n in {2,3,4,5} \draw (A\n) circle (4pt);
\filldraw[white] (A1) circle (5pt);
\draw (A1) circle (5pt);
\node at (A1) {\tiny $\rho_\tau$};
\end{tikzpicture}   
$$

   \item A \textbf{\textsf{typed rooted tree}} is a rooted tree with type 
maps $\frak{t}_{\rm n}:N_\tau\to\frak{T}_{\rm n}$ and $\frak{t}_{\rm e}:E_\tau\to\frak{T}_{\rm e}$. 
   Moreover, a \textbf{\textsf{rooted decorated tree}} is a typed rooted tree $\tau$ with two maps
   $$
   \frak{n}:N_\tau\to\bbN^{d+1},\qquad \frak{e}:E_\tau\to\bbN^{d+1}.
   $$
   We denote a generic typed rooted tree by Greek letters like $\tau$, and a generic rooted decorated trees with two decorations $\frak{n},\frak{e}$ by $\tau_{\frak{e}}^{\frak{n}}$ or a bold letter $\bstau$.
      
\end{itemize}
\end{defn*}   

\medskip

We will consider later rooted trees $\tau$ equipped with three decorations $\frak{n},\frak{o},\frak{e}$ -- see Section \ref{SectionDecoratedTrees} 
for the precise definitions. In this section, we hide the $\frak{o}$-decoration in the node type map, so we consider the type sets
$$
\frak{T}_{\rm n}=\{\bullet,\circ\}\cup\{{\color{red}\bullet}^{\cdot,\alpha}\}_{\alpha\in\bbR}.
$$
The node type $\bullet$ represents the monomial ${\bf1}=X^0$, and $\circ$ represents the noise $\Xi$. The third node type ${\color{red}\bullet}^{\cdot,\alpha}$ is a node with the $\frak{o}$-decoration $\alpha$. The set $\frak{T}_{\rm e}$ labels the set of differential operators involved in 
the system of equations under study. There is a single operator $\partial_{x_0}-\Delta_{x'}+1$ in the example of the generalized (KPZ) equation \eqref{EqgKPZ}, so the set $\frak{T}_{\rm e}$ consists of only one element, associated with the integration operator $\mcI$ in that case. If we consider a system of singular stochastic PDEs involving different operators, different operators $\mcI$'s would be associated with each of them and the set $\frak{T}_{\rm e}$ would collect them all. 

An element $X^n\in T$ is denoted by $\bullet^n$, that is a graph with only one node with the type $\bullet$ and the $\frak{n}$-decoration $n\in\bbN^{d+1}$. An edge with $\frak{e}$-decoration $p\in\bbN^{d+1}$ represents the operator $\mcI_p$, with the notations of Section \ref{SubsectionIntegrationOperators}, for one of the operators $\mcI$ involved in the equation.

All operations appearing in the equation \eqref{EqIdentification1} are graphically defined as follows. In the following pictures, types and decorations are omitted unless necessary, and the root of $\bstau$ in the first bullet and of $\bstau_1,\dots,\bstau_m$ in the third bullet is denoted by a square.

\begin{itemize}
\item The integration $\bstau\mapsto\mcI_p(\bstau)$ is given by the map connecting the root of $\bstau$ with a new node, which becomes a root of the tree $\mcI_p(\bstau)$, and giving the $\frak{e}$-decoration $p\in\bbN^{d+1}$ to the connecting edge.
$$
\mcI_p\big(
\begin{tikzpicture}[baseline=5]
\coordinate (A1) at (0,0);
\coordinate (A2) at (0,0.3);
\filldraw[black!30] (A2) circle [x radius = 0.4, y radius = 0.3];
\node at (A2) {$\bstau$};
\filldraw (-0.07,-0.07) rectangle (0.07,0.07);
\end{tikzpicture}
\big)
=
\begin{tikzpicture}[baseline=-5]
\coordinate (A1) at (0,0);
\coordinate (A2) at (0,0.3);
\filldraw[black!30] (A2) circle [x radius = 0.4, y radius = 0.3];
\node at (A2) {$\bstau$};
\filldraw (-0.07,-0.07) rectangle (0.07,0.07);
\coordinate (A4) at (0,-0.5);
\draw (A1)--(A4);
\fill (A4) circle (2pt);
\node at ($0.5*(A1)+0.5*(A4)+(0.1,0)$) {\tiny $p$};
\end{tikzpicture}
$$
\item The product $\frak{b}^n\star\bstau$ for $\frak{b}\in\frak{T}_{\rm n}$, $n\in\bbN^{d+1}$, and $\bstau$ with $\frak{t}_{\rm n}(\rho_\tau)=\bullet$ and $\frak{n}(\rho_\tau)=0$ is given changing the node type of $\rho_\tau$ to $\frak{b}$ and $\frak{n}$-decoration to $n$. For example, if $\frak{b}=\circ$ and $n=0$,
$$
\circ\star
\begin{tikzpicture}[baseline=5]
\coordinate (A1) at (0,0);
\coordinate (A2) at (0,0.3);
\filldraw[black!30] (A2) circle [x radius = 0.4, y radius = 0.3];
\node at (A2) {$\bstau$};
\fill (A1) circle (2pt);
\end{tikzpicture}
=
\begin{tikzpicture}[baseline=5]
\coordinate (A1) at (0,0);
\coordinate (A2) at (0,0.3);
\filldraw[black!30] (A2) circle [x radius = 0.4, y radius = 0.3];
\node at (A2) {$\bstau$};
\fill[white] (A1) circle (2pt);
\draw (A1) circle (2pt);
\end{tikzpicture}
$$
\item The product of trees $\mcI_{p_j}(\bstau_j)$ ($j=1,\dots,m$) is given by the \emph{tree product}, that is joining their roots.
$$
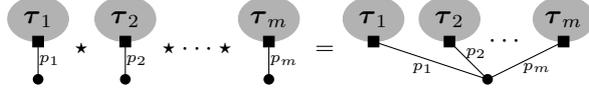
\begin{tikzpicture}[baseline=-5]
\coordinate (A1) at (0,0);
\coordinate (A2) at (0,0.3);
\filldraw[black!30] (A2) circle [x radius = 0.4, y radius = 0.3];
\node at (A2) {$\bstau_1$};
\filldraw (-0.07,-0.07) rectangle (0.07,0.07);
\coordinate (A4) at (0,-0.5);
\draw (A1)--(A4);
\fill (A4) circle (2pt);
\node at ($0.5*(A1)+0.5*(A4)+(0.15,0)$) {\tiny $p_1$};
\end{tikzpicture}
\star
\begin{tikzpicture}[baseline=-5]
\coordinate (A1) at (0,0);
\coordinate (A2) at (0,0.3);
\filldraw[black!30] (A2) circle [x radius = 0.4, y radius = 0.3];
\node at (A2) {$\bstau_2$};
\filldraw (-0.07,-0.07) rectangle (0.07,0.07);
\coordinate (A4) at (0,-0.5);
\draw (A1)--(A4);
\fill (A4) circle (2pt);
\node at ($0.5*(A1)+0.5*(A4)+(0.15,0)$) {\tiny $p_2$};
\end{tikzpicture}
\star
\cdots
\star
\begin{tikzpicture}[baseline=-5]
\coordinate (A1) at (0,0);
\coordinate (A2) at (0,0.3);
\filldraw[black!30] (A2) circle [x radius = 0.4, y radius = 0.3];
\node at (A2) {$\bstau_m$};
\filldraw (-0.07,-0.07) rectangle (0.07,0.07);
\coordinate (A4) at (0,-0.5);
\draw (A1)--(A4);
\fill (A4) circle (2pt);
\node at ($0.5*(A1)+0.5*(A4)+(0.2,0)$) {\tiny $p_m$};
\end{tikzpicture}
=
\begin{tikzpicture}[baseline=-5]
\coordinate (O) at (1.5,-0.5);
\fill (O) circle (2pt);
\coordinate (A1) at (0,0);
\coordinate (A2) at (0,0.3);
\filldraw[black!30] (A2) circle [x radius = 0.4, y radius = 0.3];
\node at (A2) {$\bstau_1$};
\filldraw (-0.07,-0.07) rectangle (0.07,0.07);
\draw (A1)--(O);
\node at ($0.5*(A1)+0.5*(O)+(-0.1,-0.1)$) {\tiny $p_1$};
\coordinate (B1) at (1,0);
\coordinate (B2) at (1,0.3);
\filldraw[black!30] (B2) circle [x radius = 0.4, y radius = 0.3];
\node at (B2) {$\bstau_2$};
\filldraw (0.93,-0.07) rectangle (1.07,0.07);
\draw (B1)--(O);
\node at ($0.5*(B1)+0.5*(O)+(0.1,0.1)$) {\tiny $p_2$};
\coordinate (D) at (1.75,0);
\node at (D) {$\cdots$};
\coordinate (C1) at (2.5,0);
\coordinate (C2) at (2.5,0.3);
\filldraw[black!30] (C2) circle [x radius = 0.4, y radius = 0.3];
\node at (C2) {$\bstau_m$};
\filldraw (2.43,-0.07) rectangle (2.57,0.07);
\draw (C1)--(O);
\node at ($0.5*(C1)+0.5*(O)+(0.15,-0.1)$) {\tiny $p_m$};
\end{tikzpicture}
$$
\end{itemize}   

Thus we see that the rooted trees obtained by the above operations are sufficient to describe the fixed point problem \eqref{EqIdentification1}.
The symbol $\textcolor{red}{\bullet}^{\cdot,\alpha}$ does not come from the fixed point problem \eqref{EqIdentification1}, but its use is made clear in Section \ref{SectionrenormalizationMultiPreLie}. As we concentrate in this section on the generalized (KPZ) equation the edge type set $\frak{T}_{\rm e}$ will consist of a single element, suggestively denoted by $\mcI$. There is no difficulty in working with a finite edge type set.

\medskip

\begin{defn*}
Let $\mcV$ be the set of \emph{all} rooted decorated trees with type sets 
$\frak{T}_{\rm n}=\{\bullet,\circ\}\cup\{{\color{red}\bullet}^{\cdot,\alpha}\}_{\alpha\in\bbR}$ and $\frak{T}_{\rm e}=\{\mcI\}$, and let $V$ be the vector space spanned by $\mcV$.
Moreover, denote by $\big(\bstau^*:V\to\bbR\big)_{\tau\in\mcV}$ the dual basis of $\mcV$ and let $V^*$ be the vector space spanned by $\{\bstau^*\}_{\bstau\in\mcV}$.
\end{defn*}

\medskip

Throughout this section we view each element of $\mcV$ as the rooted tree $\tau$ with the composite decorations $(\frak{t}_{\rm n},\frak{n}):N_\tau\to{\sf N}$ and $(\frak{t}_{\rm e},\frak{e}):E_\tau\to{\sf E}$, where
$$
{\sf E}\defeq \frak{T}_{\rm e}\times\bbN^{d+1}\simeq\bbN^{d+1},\quad
{\sf N}\defeq \frak{T}_{\rm n}\times\bbN^{d+1}.
$$
The set ${\sf N}$ is considered as a subset of $\mcV$ consisting of simple trees
$$
{\sf N}=\{\frak{b}^n\}_{\frak{b}\in\frak{T}_{\rm n},n\in\bbN^{d+1}} = \big\{\circ^\ell,\bullet^m,\textcolor{red}{\bullet}^{n,\alpha}\big\}_{\ell,m,n\in\bbN^{d+1}, \alpha\in\bbR}.
$$
Write ${\sf N}^0\defeq \{\frak{b}^0\}_{\frak{b}\in\frak{T}_{\rm n}}\simeq\frak{T}_{\rm n}$. We introduce a few notations. Note that each $\bstau\in\mcV$ has a decomposition of the form
\begin{align}\label{section 6: generic element of U}
\bstau = \frak{b}^n\star\Bigstar_{i=1}^a\mcI_{p_i}(\bstau_i) = \frak{b}^n\star\mcI_{p_1}(\bstau_1)\star\cdots\star\mcI_{p_a}(\bstau_a)
\end{align}
with $\frak{b}^n\in{\sf N}$, $p_1,\dots,p_a\in\bbN^{d+1}$, and $\bstau_1,\dots,\bstau_a\in\mcV$. Taking care of the number of automorphisms of $\bstau$ that leave it fixed, for $\bstau$ of the form 
\begin{align}\label{section 6: generic element of U overlap} 
\bstau=\frak{b}^n\star\Bigstar_{j=1}^b \big(\mcI_{q_j}(\bssigma_j)\big)^{\star m_j},
\end{align}
with $(q_i,\bssigma_i)\neq(q_j,\bssigma_j)$ for any $i\neq j$, define inductively
$$
S(\bstau) \defeq  n!\prod_{j=1}^b S(\bssigma_j)^{m_j}\,m_j!.
$$
Then we define the paring $\wangle{\cdot,\cdot}$ between $V$ and $V^*$ by
\begin{align}\label{Eq normalized pairing}
\wangle{\bstau,\bssigma^*} \defeq  S(\bssigma)\,\bssigma^*(\bstau)
\end{align}
for $\bstau,\bssigma\in\mcV$.    We see $V^*$ as a part of the algebraic dual of $V$. (As $V$ is infinite-dimensional, $V^*$ is not equal to the full algebraic dual of $V$.) The `copy' space $V^*$ will play an important role in the second half part of this section.

\bigskip

\hfill \textcolor{black}{\textsf{\textbf{ {\Large \S} Canonical model}}}

\medskip

Given a {\it smooth} noise $\zeta\in\mcC^\infty(\bbR\times\bbR^d)$, we define the canonical operator ${\sf \Pi}^\zeta$ on the whole of $V$ requiring that it is multiplicative with respect to the $\star$ product and setting, for all $x\in\bbR\times\bbR^d$,
$$
{\sf \Pi}^\zeta(\circ^n)(x) = x^n\zeta(x), 
\qquad {\sf \Pi}^\zeta(\bullet^n)(x) =
{\sf \Pi}^\zeta(\textcolor{red}{\bullet}^{n,\alpha})(x)  = x^n,
$$
and 
$$
{\sf \Pi}^\zeta(\mcI_p\bstau) = \partial^p{\bf K}({\sf \Pi}^\zeta\bstau),
$$
for all $n,\alpha, \bstau, p$. The regularity structures we will work with have spaces $T$ and $T^+$ that are subsets of $V$. Since all functions ${\sf\Pi}^\zeta\bstau$ are smooth, the restriction of ${\sf\Pi}^\zeta$ to 
$T_{<0}$ defines the \textbf{\textsf{canonical model}}
$$
{\sf M}^\zeta = \sf \big(g^\zeta,{\Pi}^\zeta\big)
$$ 
on the regularity structure $\mathscr{T}$, from Theorem \ref{ThmSmoothAdmissibleModels}. Things are explicit here as the multiplicativity and the $\bf K$-admissibility properties fix the definition of ${\sf \Pi}^\zeta$ on all decorated trees in $V$. Emphasize the fact that since the map ${\sf \Pi}^\zeta$ is multiplicative its associated reconstruction map is also 
multiplicative.

\bigskip

\hfill \textcolor{black}{\textsf{\textbf{ {\Large \S} Multi-pre-Lie algebras}}}

\medskip

We first recall the definition of a multi-pre-Lie algebra, referring the reader to Foissy's article \cite{FoissyTypedTrees} for basics on multi-pre-Lie algebras. All we need to know on the subject is the following definition and the result of Proposition \ref{PropUniversalPreLie2} below. 
 
\medskip

\begin{defn*}
Let $\sf{E}$ be a set.
A vector space $W$, equiped with a family $(\rightslice_{\sf {e}})_{\sf{e}\in\sf{E}}$ of bilinear maps from $W\times W$ into $W$, is called an \textbf{\textsf{$\sf{E}$-multi-pre-Lie algebra}} if one has
$$
(a\rightslice_{\sf{e}}b)\rightslice_{\sf{e}'}c - a\rightslice_{\sf{e}}(b\rightslice_{\sf{e}'}c)  = (b\rightslice_{\sf{e}'}a)\rightslice_{\sf{e}}c - b\rightslice_{\sf{e}'}(a\rightslice_{\sf{e}}c),
$$
for all $a,b,c\in W$, and $\sf{e}, \sf{e}'\in\sf{E}$.
\end{defn*}

\medskip

The two arguments of a pre-Lie product $a\rightslice_{\sf{e}}b$ do not play a symmetric role, and we think here of $a$ as acting on $b$ via the operator $\rightslice_{\sf{e}}$; we read $a\rightslice_{\sf{e}}b$ from left 
to right. Here is an example of $\sf E$-multi-pre-Lie algebra. Take $\sf E$ finite, identified with $\{1,\dots,\vert \sf E\vert\}$, and consider the space of smooth functions on $\bbR^{\vert \sf E\vert}$. Then the family of differentiation operators
$$
G \triangleright_{\sf e} H \defeq  G\partial_{x_{\sf e}}H
$$
defines an $\sf E$-multi-pre-Lie algebra. If $\sf E$ consists of a single 
element $\rightslice$, this operator is called a {\it pre-Lie product}, and a vector space equipped with a pre-Lie product is called a {\it pre-Lie algebra}. Any pre-Lie algebra is Lie-admissible, in the sense that the map $(a,b)\mapsto a\rightslice b - b\rightslice a$ defines a Lie bracket. 
The relevance of the multi-pre-Lie structure in the study of singular stochastic PDEs comes from Proposition \ref{PropPreLieMorphism} in the next section, as it identifies the components $u_\tau$ of solutions $\bsu=\sum u_\tau\tau$ regularity structures lifts of a singular stochastic PDEs as $\sf E$-multi-pre-Lie algebra morphisms.

\medskip

We define now the multi-pre-Lie structure in the space $V^*$. The reason for working on $V^*$ rather than on $V$ will appear clearly in Section \ref{SectionCoherence} and Section \ref{SectionrenormalizationMultiPreLie}. The spaces $V$ and $V^*$ being infinite dimensional, the symbol $\otimes$ denotes below the algebraic tensor product of these spaces with themselves, without any completion. 

\medskip

\begin{defn*}
Given ${\sf e}\in{\sf E}\simeq\bbN^{d+1}$, a node $v$ of a decorated tree $\bssigma\in \mcV$ and $\bstau\in \mcV$, denote by 
$$
\bstau\overset{\sf e}{\rightarrow}_{(v)}\bssigma,
$$
the element of $\mcV$ obtained by grafting $\bstau$ on the node $v$ of $\bssigma$, along an edge of $\frak{e}$-decoration ${\sf e}$.
Define also
\begin{align*}
\bstau\graftingatv\sigma_{\frak{e}}^{\frak{n}}
&\defeq  \sum_{m\in\bbN^{d+1}; m\le \frak{n}(v)\wedge p_{\sf e}}\binom{\frak{n}(v)}{m}\,
\bstau\xrightarrow{{\sf e}-m}_{(v)}\sigma_{\frak{e}}^{\frak{n}-m{\bf1}_v}
\in V,\\
\bstau\grafting\sigma_{\frak{e}}^{\frak{n}}
&\defeq \sum_{v\in N_\sigma}\bstau\graftingatv\sigma_{\frak{e}}^{\frak{n}}\in V,
\end{align*}
where ${\bf1}_v$ is the indicator function of $v$. Recall that the binomial coefficient of multiindices is defined at the end of Section \ref{SectionIntro}. Finally, define a linear map 
$$
\grafting : V^*\otimes V^*\rightarrow V^*
$$ 
by
$$
\bstau^*\grafting\bssigma^* \defeq  (\bstau\grafting\bssigma)^*,\qquad
\bstau,\bssigma\in\mcV,
$$
where the map $(\cdot)^*:V\to V^*$ is the linear extension of the map $\mcV\ni\bstau\mapsto\bstau^*\in V^*$.
\end{defn*}

\medskip

Here is an example
$$
\begin{tikzpicture}
\coordinate (A1) at (0,0);
\fill (A1) circle (2pt);
\end{tikzpicture}
\
\graftingatv
\begin{tikzpicture}[baseline=10]
\coordinate (A1) at (0,0);
\coordinate (A2) at (-0.33,0.57);
\coordinate (A3) at (0.33,0.57);
\draw (A1)--(A3);
\draw (A2)--(A1);
\foreach \n in {1,3} \fill (A\n) circle (2pt);
\fill[green] (A2) circle (2pt);
\node at ($(A2)+(-0.2,0)$) {\tiny $n$};
\end{tikzpicture}
=
\sum_{m\le {\sf e}\wedge n}
\binom{n}{m}
\begin{tikzpicture}[baseline=15]
\coordinate (A1) at (0,0);
\coordinate (A2) at (-0.33,0.57);
\coordinate (A3) at (0.33,0.57);
\coordinate (A4) at (-0.33,1.24);
\draw (A1)--(A3);
\draw (A2)--(A1);
\draw (A2)--(A4);
\foreach \n in {1,3,4} \fill (A\n) circle (2pt);
\fill[green] (A2) circle (2pt);
\node at ($0.5*(A2)+0.5*(A4)+(0.4,0)$) {\tiny ${\sf e}-m$};
\node at ($(A2)+(-0.5,0)$) {\tiny $n-m$};
\end{tikzpicture},
$$
where $v$ is colored in green. The next statement is fundamental and can be proved as Corollary 9 in Foissy's work \cite{FoissyTypedTrees} -- it was first proved in Proposition 4.21 of Bruned, Chandra, Chevyrev and Hairer's work \cite{BCCH18}. A proof can be found in Appendix {\sf \ref{Appendix Proof of mpL}}.

\medskip

\begin{prop} \label{PropUniversalPreLie}
The space $V^*$ with the operators $\{\grafting\}_{{\sf e}\in{\sf E}}$ is the free ${\sf E}$-multi-pre-Lie algebra generated by ${\sf N}$, in the sense that the universal property \textit{\textsf{(b)}} in Appendix {\sf \ref{Appendix Proof of mpL}} holds.
\end{prop}

\medskip

Any morphism from $V^*$ into an $\sf E$-multi-pre-Lie algebra is thus determined by its restriction to the generators ${\sf N}$ of $V^*$. This is the universal property of the free $\sf{E}$-multi-pre-Lie algebra with generators $\sf N$. In particular if two $\sf{E}$-multi-pre-Lie morphisms from $V^*$ into the same $\sf{E}$-multi-pre-Lie algebra coincide on the generators of $V^*$ then they are equal.

\medskip

The space $T=U$ of the regularity and renormalization structures associated with the generalized (KPZ) equation is a subspace of $V$ with each space $T_\beta, T^+_\alpha, U_{\beta'}, U^-_{\alpha'}$ spanned by finitely many rooted decorated trees. Denote by $\pi_U:V\to U$ the canonical projection. The next assumption is a piece of properties to be satisfied by the basis $\mcB$ of $T$ and $U$. In Section \ref{SectionDecoratedTrees}, $\mcB$ is defined as the set of all trees \textbf{\textsf{strongly conforming}} to the rule. The first one means that the $\frak{n}$-decoration is independent of the rule and the second one means that the rule is local. The last one describes that the projection  map $\pi_U$ behaves consistently with respect to all the grafting products $\grafting$.

\medskip

\begin{assumD}\label{D1}
The homogeneous basis $\mcB$ of $T$ and $U$ is a subset of $\mcV$ with the following properties. (Recall that the notions of homogeneity in $T$ and $U$ are different.)
\begin{itemize}
\setlength{\itemsep}{0.1cm}
\item If $\bstau=\tau_{\frak{e}}^{\frak{n}}\in\mcB$, then $\tau_{\frak{e}}^{\frak{m}}\in\mcB$ for any $\frak{m}:N_\tau\to\bbN^{d+1}$.
\item If $\bstau=\frak{b}^n\star\Bigstar_{i=1}^a \mcI_{p_i}(\bstau_i)\in\mcB$, then $\frak{b},\bstau_1,\dots,\bstau_a\in\mcB$.
\item For any $\bstau,\bssigma\in\mcV$ and ${\sf e}\in{\sf E}$,
$$
\pi_U\big(\bstau\grafting(\pi_U\bssigma)\big) = \pi_U\big((\pi_U\bstau)\grafting\bssigma\big) = \pi_U\big(\bstau\grafting\bssigma\big).
$$
\end{itemize}
\end{assumD}

\medskip

Set 
$$
U^* \defeq  \text{\rm span}\big\{\bstau^*\,;\, \bstau\in\mcB\big\}
$$ 
and denote by $\pi_{U^*}:V^*\to U^*$ be the canonical projection. Then we 
define the map 
$$
\flatgrafting\, : U^*\otimes U^*\to U^*
$$
setting
$$
\bstau^*\flatgrafting\bssigma^* \defeq  \pi_{U^*}\big(\bstau^*\grafting\bssigma^*\big).
$$
The following statement is proved in Appendix {\sf \ref{Appendix Proof of mpL}}. 

\medskip

\begin{prop} \label{PropUniversalPreLie2}
Under Assumption \refD{D1}, the space $U^*$ with the operators $\{\flatgrafting\}_{{\sf e}\in{\sf E}_{\le1}}$ is the $\sf E$-multi-pre-Lie algebra generated by ${\sf N}\cap \mcB$.
\end{prop}   

\medskip

Finally we define an operator playing the role of `antiderivative'.

\medskip

\begin{defn*}
For each $i\in\{0,1,\dots,d\}$, define the linear map $\uparrow_i:V^*\to V^*$ by
$$
\uparrow_i(\tau_{\frak{e}}^{\frak{n}})^*
=\sum_{v\in N_\tau}(\tau_{\frak{e}}^{\frak{n}+e_i\mathbf{1}_v})^*.
$$
\end{defn*}

\medskip

The map $\uparrow_i$ sends $U^*$ into itself under Assumption \refD{D1}. We denote by
$$
\downarrow_i:V\to V
$$
the dual map of $\uparrow_i:V^*\to V^*$ under the pairing \eqref{Eq normalized pairing}, that is, 
$$
\wangle{\,\downarrow_i\bstau,\bssigma^*}=\wangle{\bstau,\uparrow_i\bssigma^*},
$$
for any $\bstau,\bssigma\in\mcV$.
Moreover, we extend the pairing \eqref{Eq normalized pairing} into a pairing between $V\otimes V$ and $V^*\otimes V^*$ setting
$$
\wangle{\bstau_1\otimes\bstau_2,\bssigma_1^*\otimes\bssigma_2^*} \defeq  \wangle{\bstau_1,\bssigma_1^*}\,\wangle{\bstau_2,\bssigma_2^*}.
$$
Under such pairings, denote by 
$$
\updownharpoons_{\sf e} : V\rightarrow V\otimes V,
$$ 
the dual map of $\grafting\,:V^*\otimes V^*\to V^*$, that is,
\begin{align}\label{section 6:dual}
\wangle{\updownharpoons_{\sf e}\bseta, \bstau^*\otimes\bssigma^*} \defeq  \wangle{\bseta, \bstau^*\grafting\bssigma^*},
\end{align}
for any $\bstau,\bssigma,\bseta\in\mcV$ and ${\sf e}\in \sf E$. The following explicit formulas for $\downarrow_i$ and $\updownharpoons_{\sf e}$ are 
helpful to get a graphical image. It is used only in the proof of Theorem 
\ref{thm: all requirements of regul and renor str} giving an explicit construction of the regularity and renormalization structures associated with the generalized (KPZ) equation.

\medskip

\begin{lem}\label{lem:dual between graft and cut}
For any $i\in\{0,1,\dots,d\}$ and any $\bstau=\tau_{\frak{e}}^{\frak{n}}\in\mcV$, one has
\begin{align*}
\downarrow_i\big(\tau_{\frak{e}}^{\frak{n}}\big) = \sum_{v\in N_\tau,\, 
e_i\le \frak{n}(v)}\frak{n}(v)\,\tau_{\frak{e}}^{\frak{n}-e_i\mathbf{1}_v}.
\end{align*}
Moreover for any $\bstau=\tau_{\frak{e}}^{\frak{n}}\in\mcV$ and any $\sf e\in E$ one has
\begin{align}\label{section6:defofcut}
\updownharpoons_{\sf e}(\tau_{\frak{e}}^{\frak{n}})
&=\sum_{e=(v,w)\in E_\tau; \frak{e}(e)\le {\sf e}}
\frac1{({\sf e}-\frak{e}(e))!}
(C_e\tau)_{\frak{e}}^{\frak{n}}
\otimes (P_e\tau)_{\frak{e}}^{\frak{n}+({\sf e}-\frak{e}(e)){\bf1}_v},
\end{align}
where $C_e\tau$ and $P_e\tau$ are the two connected components of the graph $\tau\setminus\{e\}$, with $P_e\tau$ containing the root of $\tau$.
(Again, recall that the factorial of a multiindex is defined at the end of Section \ref{SectionIntro}.)
\end{lem}

\ssk

\begin{Dem}
We show that the equation \eqref{section 6:dual} holds for the map $\updownharpoons_{\sf e}$ defined by the second formula \eqref{section6:defofcut}. The first formula is proved by a similar argument. Note that, for any elements $\bstau=\frak{b}^n\star\Bigstar_{i=1}^a\mcI_{p_i}(\bstau_i)\in \mcV$ and $\bssigma=\frak{u}^{m}\star\Bigstar_{j=1}^b\mcI_{q_j}(\bssigma_j)\in \mcV$, one has
\begin{align}\label{section 6:general inner product}
\wangle{\bstau,\bssigma^*}=
{\bf1}_{\frak{b}=\frak{u},n=m,a=b}\, n!\sum_{s\in S_a}\prod_{i=1}^a{\bf1}_{p_i=q_{s(i)}}\wangle{\bstau_i,\bssigma_{s(i)}^*},
\end{align}
where $S_a$ is the symmetric group of the set $\{1,2,\dots,a\}$. For $\bseta$ of the form $\bseta={\frak{b}}^n\star\Bigstar_{i=1}^a\mcI_{p_i}(\bseta_i)$, 
we divide the right hand side of \eqref{section6:defofcut} according to the edge $e$ is connected to the root or not and have
\begin{align*}
\updownharpoons_{\sf e}\bseta
&=\sum_i\frac1{({\sf e}-p_i)!}\bseta_i
\otimes \frak{b}^{n+{\sf e}-p_i}
\star\Bigstar_{j:j\neq i}\mcI_{p_j}(\bseta_j)
+\sum_i\sum_{(\bseta_i)} \bseta_i^1\otimes \frak{b}^n\star\mcI_{p_i}(\bseta_i^2)
\star\Bigstar_{j:j\neq i}\mcI_{p_j}(\bseta_j)\\
&\eqdef \updownharpoons_{\sf e}^1\bseta+\updownharpoons_{\sf e}^2\bseta,
\end{align*}
where we write $\updownharpoons_{\sf e}\bseta_i=\sum_{(\bseta_i)}\bseta_i^1\otimes\bseta_i^2$ following Sweedler's notation in the first equality. Similarly for any $\bstau\in\mcV$ and $\bssigma=\frak{u}^{m}\star\Bigstar_{j=1}^b\mcI_{q_j}(\bssigma_j)\in \mcV$ one has
\begin{align*}
\bstau\grafting\bssigma
&= \sum_\ell\binom{m}{\ell}
\frak{u}^{m-\ell}\star\mcI_{{\sf e}-\ell}(\bstau) \star\Bigstar_{j=1}^b \mcI_{q_j}(\bssigma_j)
+ \sum_{j=1}^b\frak{u}^m\star\mcI_{q_j}\big(\bstau\grafting\bssigma_j\big)
\star\Bigstar_{k;k\neq j}\mcI_{q_k}(\bssigma_k)  \\
&\eqdef \bstau\rootgrafting\bssigma + \bstau\nonrootgrafting\bssigma.
\end{align*}
Hence it is sufficient to show that
\begin{align}
\label{section 6:dual1}\wangle{\bseta, \bstau^*\rootgrafting\bssigma^*} 
&= \wangle{\updownharpoons_{\sf e}^1\bseta, \bstau^*\otimes\bssigma^*}, 
  \\
\label{section 6:dual2}\wangle{\bseta, \bstau^*\nonrootgrafting\bssigma^*} 
&= \wangle{\updownharpoons_{\sf e}^2\bseta, \bstau^*\otimes\bssigma^*}.
\end{align}
It is not difficult to show \eqref{section 6:dual1} directly from \eqref{section 6:general inner product}. For \eqref{section 6:dual2}, it is sufficient to consider $\bssigma = \frak{b}^n\star\Bigstar_{i=1}^a\mcI_{q_i}(\bssigma_i)$, and for such $\bssigma$ one has
\begin{align*}
\wangle{\bseta, \bstau^*\nonrootgrafting\bssigma^*}
= n! \sum_{i=1}^a\sum_{s\in S_a}{\bf1}_{q_{s(i)}=p_i}
\wangle{\bseta_i, \bstau^*\grafting\bssigma_{s(i)}^*}
\prod_{j;j\neq i}{\bf1}_{q_{s(j)}=p_j}\wangle{\bseta_j, \bssigma_{s(j)}^*}
\end{align*}
and
\begin{align*}
\wangle{\updownharpoons_{\sf e}^2\bseta, \bstau^*\otimes\bssigma^*}
&= n! \sum_{i=1}^a\sum_{(\bseta_i)}\sum_{s\in S_a}\wangle{\bseta_i^1, \bstau^*}
{\bf1}_{q_{s(i)}=p_i} \wangle{\bseta_i^2, \bssigma_{s(i)}^*}
\prod_{j;j\neq i}{\bf1}_{q_{s(j)}=p_j}\wangle{\bseta_j, \bssigma_{s(j)}^*}.
\end{align*}
Since 
$$
\sum_{(\bseta_i)}\wangle{\bseta_i^1, \bstau^*} \,\wangle{\bseta_i^2, \bssigma_{s(i)}^*} 
= \wangle{\updownharpoons_{\sf e}\bseta, \bstau^*\otimes\bssigma_{s(i)}^*},
$$ 
identity \eqref{section 6:dual2} follows if \eqref{section 6:dual} holds for $\bssigma=\bssigma_i$, which leads to an induction on the number of edges contained in $\bssigma$. The case $\bssigma=\frak{b}^n\in{\sf N}$ is an easy exercise.
\end{Dem}

\bigskip

\subsection{Modelled distributions solutions of singular PDEs}
\label{SectionCoherence}

The approximate description \eqref{EqIdentification1} of the fixed point problem \eqref{EqRSgKPZ} leads to an explicit formula for the coefficients of the solution $\bsu$. Noting that $\gamma\in(-\beta_0,2)$ can be arbitrarily chosen, the solution $\bsu$ of \eqref{EqIdentification1} is of the form
\begin{align}\label{Section 6: expansion of bsu}
\bsu=\sum_{|k|_\mfs<\gamma}\frac{u_k}{k!}X^k+\sum_{\bstau\in\mcB,\, |\bstau|<\gamma-2}u_{\mcI(\bstau)}\,\mcI(\bstau).
\end{align}
Inserting such an expansion into \eqref{EqIdentification1}, we see that all coefficients $u_{\mcI(\bstau)}$ are cylindrical functions of
$$
{\sf u}\defeq(u_k)_{k\in\bbN^{1+d}}\in\bbR^{\bbN^{1+d}}.
$$
Here we say that a function of $\sf u$ is cylindrical if it depends only on a finite number of entries among $(u_k)_{k\in\bbN^{1+d}}$. For any smooth cylindrical function $F$, we denote by $\partial_kf\defeq\frac{\partial}{\partial u_k}F$ the derivative with respect to $u_k$. Moreover, we define the derivative operators $(D_i)_{i=0}^d$ by setting
$$
D_iF\defeq\sum_{k\in\bbN^{1+d}}u_{k+e_i}\partial_kF,
$$
and $D^n\defeq \prod_{i=0}^dD_i^{n_i}$ for $n=(n_i)_{i=0}^d\in\bbN^{d+1}$.

\medskip

\begin{defn*}
Set ${\sf u}_0=u_0$ and ${\sf u}_1=(u_{e_i})_{i=1}^d$. Define the linear map $F$ from $V^*$ to the space of cylindrical functions of $\sf u$ as follows. For the primitive trees in ${\sf N}^0\defeq \{\frak{b}^0\}_{\frak{b}\in\frak{T}_{\rm n}}\simeq\frak{T}_{\rm n}$, set

\begin{equation}  \label{EqDefnFG} \begin{split}
F(\circ^*)({\sf u}) &\defeq  f({\sf u}_0),   \\
F(\bullet^*)({\sf u}) &\defeq  g({\sf u}_0,{\sf u}_1)\defeq g_2({\sf u}_0)({\sf u}_1)^2+g_1({\sf u}_0){\sf u}_1+g_0({\sf u}_0)   \\
&\defeq\sum_{i,j=1}^dg_2^{ij}(u_0)u_{e_i}u_{e_j}+\sum_{i=1}^dg_1^i(u_0)u_{e_i}+g_0(u_0),\\
F\big((\textcolor{red}{\bullet}^{0,\alpha})^*\big)({\sf u}) &\defeq  0.
\end{split} \end{equation}
For a generic tree 
$$
\bstau = \frak{b}^n\star\Bigstar_{i=1}^a\mcI_{p_i}(\bstau_i),
$$
define inductively
\begin{align}\label{defn: Ffg}
F(\bstau^*)({\sf u}) \defeq  \left\{\prod_{i=1}^a F(\bstau_i^*)({\sf u})\right\} \left\{D^n\prod_{i=1}^a \partial_{p_i}\right\}F\big(\frak{b}^*\big)({\sf u}).
\end{align}
\end{defn*}

\medskip

Here are some examples. Recall that $\bullet$ represents ${\bf1}=X^0$, and $\circ$ represents the noise $\Xi$.
\begin{align*}
F\big((\Xi\star\mcI(\Xi))^*\big)({\sf u})&=f(u_0)f'(u_0),\\
F\big((\mcI_{e_i}(\Xi)\star\mcI_{e_j}(\Xi))^*\big)({\sf u})&=
\begin{cases}
g_2^{ij}(u_0)f(u_0)^2&(i\neq j),\\
2g_2^{ii}(u_0)f(u_0)^2&(i=j).
\end{cases}
\end{align*}

We recall some useful formulas for the derivative operators.

\medskip

\begin{lem}
For any smooth cylindrical function $F$ of $\sf u$ and any $n\in\bbN^{1+d}$, one has
\begin{align}\label{lem:FaadiBruno}
\frac{D^nF}{n!}=\sum_{\substack{m:\bbN^{1+d}\times(\bbN^{1+d}\setminus\{0\})\to\bbN \\ \sum_q(\sum_km(k,q))q=n}}
\bigg\{\prod_{\substack{k\in\bbN^{1+d} \\ q\in\bbN^{1+d}\setminus\{0\}}}\frac1{m(k,q)!}\Big(\frac{u_{k+q}}{q!}\Big)^{m(k,q)}\bigg\}
\bigg(\prod_{\substack{k\in\bbN^{1+d} \\ q\in\bbN^{1+d}\setminus\{0\}}}\partial_k^{m(k,q)}\bigg)F
\end{align}
(Fa\`a di Bruno formula from Lemma A.1 of \cite{BCCH18}).
Here $0!=1$, $u_k^0=1$, and $\partial_k^0=\iden$ by convention, so the sum and the multiplications in the right hand side are over only finitely many parameters.
Moreover, for any $k,n\in\bbN^{1+d}$, one has
\begin{equation}\label{lem:commutationpartialD}
\partial_kD^nF=\sum_{\ell\in\bbN^{1+d}}\binom{n}{\ell}D^{n-\ell}\partial_{k-\ell}F,
\end{equation}
where $\binom{n}{\ell}=0$ if $\ell>n$ and $\partial_{k-\ell}=0$ if $\ell>k$ by convention, so the sums are over $\ell\le n\wedge k$.
\end{lem}

\medskip

\begin{Dem}
The proof of \eqref{lem:FaadiBruno} is an induction on $n$. In the case $n=e_i$, the right hand side of \eqref{lem:FaadiBruno} coincides with $D_iF$ by definition.
Next assume that \eqref{lem:FaadiBruno} holds for $n$ and consider $n+e_i$.
For simplicity, we denote by $M_n$ the set of all maps $m:\bbN^{1+d}\times(\bbN^{1+d}\setminus\{0\})\to\bbN $ such that $\sum_q(\sum_km(k,q))q=n$, and write $A(m)=\prod_{k,q}\frac1{m(k,q)!}\big(\frac{u_{k+q}}{q!}\big)^{m(k,q)}$ and $B(m)=\big(\prod_{k,q}\partial_k^{m(k,q)}\big)F$.
By Leibniz rule, we can divide $D_i\frac{D^nF}{n!}$ into two terms according to that $D_i$ is applied to $A(m)$ or $B(m)$.
By definition of $D_i=\sum_{\ell\in\bbN^{1+d}}u_{\ell+e_i}\partial_\ell$, the latter part is reorganized as
\begin{equation}\label{proof:lem:FaadiBruno1}
\sum_{m\in M_n}\sum_{\ell\in\bbN^{1+d}}\big(m(\ell,e_i)+1\big)A(\widetilde{m}_\ell)B(\widetilde{m}_\ell),
\end{equation}
where $\widetilde{m}_\ell\in M_{n+e_i}$ is defined by $\widetilde{m}_\ell(k,q)=m(k,q)+{\bf1}_{(k,q)=(\ell,e_i)}$. On the other hand, since $D_iu_k^m=mu_k^{m-1}u_{k+e_i}$, the former part is also reorganized as
\begin{equation}\label{proof:lem:FaadiBruno2}
\sum_{m\in M_n}\sum_{\ell\in\bbN^{1+d},\,p\in\bbN^{1+d}\setminus\{0\};m(\ell,p)\ge1}(p_i+1)\big(m(\ell,p+e_i)+1\big)A(\widetilde{m}_{\ell,p})B(\widetilde{m}_{\ell,p}),
\end{equation}
where $\widetilde{m}_{\ell,p}\in M_{n+e_i}$ is defined by $\widetilde{m}_{\ell,p}(k,q)=m(k,q)-{\bf1}_{(k,q)=(\ell,p)}+{\bf1}_{(k,q)=(\ell,p+e_i)}$. 
The term \eqref{proof:lem:FaadiBruno1} can be absorbed into the sum \eqref{proof:lem:FaadiBruno2} where the condition on $p$ is replaced by `$p\in\bbN^{1+d}$', by setting $\widetilde{m}_{k,0}\defeq\widetilde{m}_k$ and $m(\ell,0)=1$ for any $\ell$. Conversely, for any $\mu\in M_{n+e_i}$, if $\mu(\ell,p)\ge1$ for some $\ell\in\bbN^{1+d}$ and $p\in\bbN^{1+d}\setminus\{0\}$, then there exists a unique $m\in M_n$ such that $\widetilde{m}_{\ell,p-e_i}=\mu$. Therefore the sum \eqref{proof:lem:FaadiBruno2} is reorganized as
\begin{align*}
D_i\frac{D^nF}{n!}&=\sum_{m\in M_n}\sum_{\ell,p\in\bbN^{1+d};m(\ell,p)\ge1}(p_i+1)\big(m(\ell,p+e_i)+1\big)A(\widetilde{m}_{\ell,p})B(\widetilde{m}_{\ell,p})\\
&=\sum_{\mu\in M_{n+e_i}}\bigg(\sum_{m\in M_n,\,\ell,p\in\bbN^{1+d};\widetilde{m}_{\ell,p}=\mu}(p_i+1)\big(m(\ell,p+e_i)+1\big)\bigg)A(\mu)B(\mu).
\end{align*}
It turns out that the quantity inside the large parentheses is equal to
$$
\sum_{\ell\in\bbN^{1+d},q\in\bbN^{1+d}\setminus\{0\}}q_i\mu(\ell,q)=n_i+1
$$
because of the condition that $\mu\in M_{n+e_i}$. Thus we have $\frac{D_i}{n_i+1}\frac{D^nF}{n!}=\sum_{\mu\in M_{n+e_i}}A(\mu)B(\mu)$. This yields that \eqref{lem:FaadiBruno} holds for any $n$.

The proof of \eqref{lem:commutationpartialD} is also an induction on $n$. The case $n=e_i$ follows from Leibniz rule:
\begin{align*}
\partial_kD_iF&=\sum_{\ell\in\bbN^{1+d}}\partial_k(u_{\ell+e_i}\partial_\ell F)
=\sum_{\ell\in\bbN^{1+d}}(\partial_ku_{\ell+e_i})\partial_\ell F+\sum_{\ell\in\bbN^{1+d}}u_{\ell+e_i}\partial_k\partial_\ell F\\
&={\bf1}_{k\ge e_i}\partial_{k-e_i}F+D_i\partial_kF.
\end{align*}
Assuming that \eqref{lem:commutationpartialD} holds for $n$, we have for $n+e_i$
\begin{align*}
\partial_kD^{n+e_i}F&=\partial_kD^n(D_iF)=\sum_{\ell}{\textstyle\binom{n}{\ell}}D^{n-\ell}\partial_{k-\ell}D_iF
=\sum_{\ell}{\textstyle\binom{n}{\ell}}D^{n-\ell}({\bf1}_{k-\ell\ge e_i}\partial_{k-\ell-e_i}F+D_i\partial_{k-\ell}F)
\\
&=\sum_{m\in\bbN^{1+d}}{\textstyle\binom{n}{m-e_i}}D^{n+e_i-m}\partial_{k-m} F+\sum_{m\in\bbN^{1+d}}{\textstyle\binom{n}{m}}D^{n+e_i-m}\partial_{k-m} F.
\end{align*}
Since $\binom{n}{m-e_i}+\binom{n}{m}=\binom{n+e_i}{m}$, it turns out that \eqref{lem:commutationpartialD} also holds for $n+e_i$.
\end{Dem}

\medskip

Using the Fa\`a di Bruno formula \eqref{lem:FaadiBruno}, we can give in the following lemma a representation of the nonlinear terms of \eqref{EqRSgKPZ}.
Given a modelled distribution $\bsu\in\mcD^{\gamma,\eta}(T,\sf g)$ of the form \eqref{Section 6: expansion of bsu}, set
$$
\mcF(\bsu) \defeq  \mathcal{Q}_{<\gamma+\beta_0}\left(\sum_{\frak{b}\in{\sf N}^0} \big(F(\frak{b}^*)\big)^\star(\bsu, D\bsu  )\star\frak{b}\right) = f^\star(\bsu)\star \Xi + g^\star(\bsu, D\bsu  ).
$$

\begin{lem}\label{lem ExpSol bsu}
Consider a generic tree of the form $\bstau = \frak{b}^n\star\Bigstar_{j=1}^b \big(\mcI_{p_j}(\bssigma_j)\big)^{\star m_j}$ such that $|\bstau|<\gamma+\beta_0$ and $(p_i,\bssigma_i)\neq(p_j,\bssigma_j)$ for any $i\neq j$.
The $\bstau$-component of $\mcF(\bsu(x))$ is given by
\begin{equation} \label{EqCoherence*}
\left\{\prod_{j=1}^b \frac{u_{\mcI(\bssigma_j)}(x)^{m_j}}{m_j!}\right\} \left\{\frac{D^n}{n!}\prod_{j=1}^b \partial_{p_j}^{m_j}\right\}F\big(\frak{b}^*\big)\big((u_k(x))_{|k|_\mfs\le1}\big).
\end{equation}
Consequently, if $\bsu$ solves the fixed point problem \eqref{EqRSgKPZ} then
for any $\bstau\in\mcB$ with $|\bstau|<\gamma-2$ and any $x\in(0,\frac{t_0}2)\times\bbT^d$, one has
\begin{equation} \label{EqCoherence}
u_{\mcI(\bstau)}(x)=\frac1{S(\bstau)}F(\bstau^*)\big((u_k(x))_{k\in\bbN^{1+d}}\big),
\end{equation}
where the right hand side depends only on $u_k(x)$ with $|k|_\mfs<\gamma$.
\end{lem}

\medskip

\begin{Dem}
We consider here
$$
\bstau = \circ^n\star\Bigstar_{j=1}^b \big(\mcI(\bssigma_j)\big)^{\star m_j};
$$
the other cases are proved by similar arguments. The element $\bstau$ appears in the term $f^\star(\bsu)\star \Xi$. Inserting the expansion \eqref{Section 6: expansion of bsu} into $f^\star(\bsu)\star \Xi$, its $\bstau$-component is given by
$$
\sum_{\bstau_1,\dots,\bstau_a\,;\, \Bigstar_{i=1}^a\bstau_i=\bstau}
\frac{f^{(a)}(u_0)}{a!} \, u_{\bstau_1}\dots u_{\bstau_a},
$$
where $\bstau_1,\dots,\bstau_a$ are elements of $\mcB$ of the forms $X^k$ with $k\neq0$ or $\mcI(\bssigma_j)$, and $u_{X^k}=\frac{u_k}{k!}$ by definition.
Rearranging $\bstau_1,\dots,\bstau_a$ so that duplicate elements are grouped together, the above quantity is reorganized as
\begin{align*}
\sum_{\substack{q_1,\dots,q_r\in\bbN^{1+d}\setminus\{0\},\,i\neq j\Rightarrow q_i\neq q_j, \\ 
n_1,\dots,n_r\in\bbN,\, n_1q_1+\cdots+n_rq_r=n}}
\frac{f^{(n_1+\cdots+n_r+m_1+\cdots+m_b)}(u_0)}{n_1!\cdots n_r!m_1!\cdots m_b!}
\Big(\frac{u_{q_1}}{q_1!}\Big)^{n_1}\cdots\Big(\frac{u_{q_r}}{q_r!}\Big)^{n_r}
u_{\mcI(\bssigma_1)}^{m_1}\cdots u_{\mcI(\bssigma_b)}^{m_b}
\end{align*}
Applying the formula \eqref{lem:FaadiBruno} to the sum over $q_1,\dots,q_r$ and $n_1,\dots,n_r$, the above quantity is equal to
\begin{align*}
\frac{1}{n!\,m_1!\cdots m_b!}\,
u_{\mcI(\bssigma_1)}^{m_1}\cdots u_{\mcI(\bssigma_b)}^{m_b}
D^n\partial_0^{m_1+\cdots+m_b}f(u_0).
\end{align*}
This is a particular case of \eqref{EqCoherence*} with $\frak{b}=\circ$ and $p_j=0$.

If $\bsu$ solves the fixed point problem \eqref{EqRSgKPZ}, then
$$
\bsu(x)=\sum_{|k|_\mfs<\gamma}\frac{u_k(x)}{k!}X^k+\mcQ_{<\gamma}\mcI\big(\mcF(\bsu(x))\big)
$$
for any $x\in(0,\frac{t_0}2)\times\bbT^d$. Therefore the quantity \eqref{EqCoherence*} should be equal to $u_{\mcI(\bstau)}(x)$ for any $\bstau\in\mcB$ with $|\bstau|<\gamma-2$. Assuming $u_{\mcI(\bssigma_j)}=F(\bssigma_j^*)/S(\bssigma_j)$ inductively, we have
$$
u_{\mcI(\bstau)}=\frac1{S(\bstau)}\,F(\bstau^*).
$$
This concludes the proof.
\end{Dem}

\ssk

Modelled distributions satisfying identity \eqref{EqCoherence} are called `\textit{coherent}' in \cite{BCCH18}.

\medskip

Recall ${\sf E}\simeq \bbN^{d+1}$, and define the family of differential operators
$$
G\triangleright_{\sf e} H \defeq  G\,\partial_{u_{\sf e}}H \qquad ({\sf e}\in{\sf E}),
$$
acting on smooth functions of $(u_0,u_1)$, with $u_1=(u_{X_i})_{i=1}^d$. The family $\{\triangleright_{\sf e}\}_{{\sf e}\in{\sf E}}$ defines an $\sf E$-multi-pre-Lie algebra structure. 

\medskip

\begin{prop}   \label{PropPreLieMorphism}
The map $F$ is an ${\sf E}$-multi-pre-Lie algebra morphism: For any ${\sf e}\in{\sf E}$ and any decorated trees $\bstau, \bssigma$ in $\mcB$, one has
\begin{equation}   \label{EqCharacterizationGrafting*}
F\big(\bstau^*\grafting\bssigma^*\big) = F(\bstau^*)\triangleright_{\sf 
e} F(\bssigma^*).
\end{equation}
\end{prop}

\medskip

\begin{Dem}
Assume $\bssigma$ is of the form
$$
\bssigma = \frak{b}^n\star\Bigstar_{i=1}^a \mcI_{p_i}(\bssigma_i).
$$
Then by definition,
\begin{align*}
\bstau\grafting\bssigma
&= \sum_\ell\binom{n}{\ell}
\frak{b}^{n-\ell}\star\mcI_{{\sf e}-\ell}(\bstau) \star\Bigstar_{i=1}^a 
\mcI_{p_i}(\bssigma_i)
+ \sum_{i=1}^a
\frak{b}^n\star\mcI_{p_i}\big(\bstau\grafting\bssigma_i\big)
\star\Bigstar_{j;j\neq i}\mcI_{p_j}(\bssigma_j).  
\end{align*}
Hence
\begin{align*}
F\big(\bstau^*\grafting\bssigma^*\big) 
&= \sum_\ell\binom{n}{\ell}
F(\bstau^*)\left\{\prod_{i=1}^a F(\bssigma_i^*)\right\}
D^{n-\ell}\partial_{{\sf e}-\ell}\left\{\prod_{i=1}^a \partial_{{p_i}}\right\}F(\frak{b}^*) \\
&\quad+\sum_{i=1}^a F(\bstau^*\grafting\bssigma_i^*)
\left\{\prod_{j;j\neq i}F(\bssigma_i^*)\right\}
D^n\left\{\prod_{i=1}^a \partial_{{p_i}}\right\}F(\frak{b}^*).
\end{align*}
On the other hand, by Leibniz rule,
\begin{align*}
F(\bstau^*)\triangleright_{\sf e} F(\bssigma^*)
&= F(\bstau^*)\left\{\prod_{i=1}^a F(\bssigma_i^*)\right\}
\partial_{{\sf e}}D^n\left\{\prod_{i=1}^a \partial_{{p_i}}\right\}F(\frak{b}^*)   \\
&\quad+F(\bstau^*)\sum_{i=1}^a \partial_{{\sf e}}F(\bssigma_i^*)
\left\{\prod_{j;j\neq i}F(\bssigma_i^*)\right\}
D^n\left\{\prod_{i=1}^a \partial_{{p_i}}\right\}F(\frak{b}^*).
\end{align*}
The first terms in the expansions of $F\big(\bstau^*\grafting\bssigma^*\big)$ and $F(\bstau^*)\triangleright_{\sf e} F(\bssigma^*)$ coincide, because of the identity \eqref{lem:commutationpartialD}. The second terms turn out to coincide if \eqref{EqCharacterizationGrafting*} holds for $\bstau^*$ and $\bssigma_i^*$, which leads an induction on the number of edges contained in $\bssigma^*$.
\end{Dem}

\medskip

Assumption \refD{D1} is a necessary condition for the basis $\mcB$. The next assumption means that $\mcB$ is sufficiently large to describe all terms in the right hand side of \eqref{EqIdentification1}.

\medskip

\begin{assumD}\label{D2}
One has $F(\bstau^*)=0$ for any $\bstau\in\mcV\setminus\mcB$.
\end{assumD}

\medskip

In particular, Assumption \refD{D2} holds if $\mcB$ contains all trees strongly conforming to the rule as in Section \ref{SectionDecoratedTrees}. Indeed, if $\bstau$ is not strongly conforming and does not have any node with $\textcolor{red}{\bullet}^{\cdot,\alpha}$ decoration, then $\bstau$ have an edge $\mcI_p$ with $|p|_{\frak{s}}\ge2$ or have a node with at least three leaving edges $\mcI_p$ with $|p|_{\frak{s}}=1$.
Since $F(\bullet^*)$ is at most quadratic with respect to $u_1$, we have $F(\bstau^*)=0$.
We define
$$
\Upsilon\defeq F\vert_{U^*}.
$$
By Assumption \refD{D2}, we can conclude that $\Upsilon$ is an ${\sf E}$-multi-pre-Lie algebra morphism on the ${\sf E}$-multi-pre-Lie algebra $(U^*,\{\grafting_\flat\}_{{\sf e}\in {\sf E}})$.   

\medskip

\begin{prop}
Under Assumption \refD{D2}, the map $\Upsilon$ is an ${\sf E}$-multi-pre-Lie algebra morphism: For any ${\sf e}\in{\sf E}$ and any decorated trees $\bstau, \bssigma$ in $\mcB$, one has
\begin{equation}   \label{EqCharacterizationGrafting}
\Upsilon\big(\bstau^*\flatgrafting\bssigma^*\big) = \Upsilon(\bstau^*)\triangleright_{\sf e} \Upsilon(\bssigma^*).
\end{equation}
\end{prop}   

\medskip

\begin{Dem}
Since $\Upsilon\circ\pi_{U^*}=F\circ\pi_{U^*}=F$,
$$
\Upsilon\big(\bstau^*\flatgrafting\bssigma^*\big) 
=F\big(\bstau^*\grafting\bssigma^*\big)
=F(\bstau^*)\triangleright_{\sf e} F(\bssigma^*)
=\Upsilon(\bstau^*)\triangleright_{\sf e} \Upsilon(\bssigma^*).
$$
\end{Dem}

\medskip

The next proposition is proved by an induction similar to the induction used in the proof of Proposition \ref{PropPreLieMorphism}, noting that $D_i$ satisfies Leibniz rule.

\ssk

\begin{prop}   \label{PropPreLieMorphism2}
Under Assumption \refD{D2}, for any $i\in\{0,1,\dots,d\}$ and $\bstau\in\mcV$, one has
\begin{equation}   \label{EqCharacterizationGrafting2*}
F(\,\uparrow_i\bstau^*) = D_iF(\bstau^*).
\end{equation}
and for $\bstau\in\mcB$,
\begin{equation}   \label{EqCharacterizationGrafting2}
\Upsilon(\,\uparrow_i\bstau^*) = D_i\Upsilon(\bstau^*).
\end{equation}
\end{prop}

\medskip

\begin{Dem}
Assume $\bstau$ is of the form
$$
\bstau = \frak{b}^n\star\Bigstar_{k=1}^a \mcI_{p_k}(\bstau_k).
$$
Then by definition,
\begin{align*}
\uparrow_i\tau
&= \frak{b}^{n+e_i}\star\Bigstar_{k=1}^a \mcI_{p_k}(\bstau_k)
+ \sum_{k=1}^a
\frak{b}^n\star\mcI_{p_k}(\, \uparrow_i\bstau_k)
\star\Bigstar_{j;j\neq k}\mcI_{p_j}(\bstau_j).  Z
\end{align*}
Hence
\begin{align*}
F(\,\uparrow_i\bstau^*)
&= \left\{\prod_{k=1}^a F(\bstau_k^*)\right\}D^{n+e_i}\left\{\prod_{k=1}^a \partial_{{p_k}}\right\}F(\frak{b}^*)
+ \sum_{k=1}^aF(\,\uparrow_i\bstau_k^*)\left\{\prod_{j;j\neq k} F(\bstau_j^*)\right\}D^n\left\{\prod_{k=1}^a \partial_{{p_k}}\right\}F(\frak{b}^*).  
\end{align*}
If \eqref{EqCharacterizationGrafting2*} holds for $\bstau_i$, then the above quantity is equal to $D_iF(\bstau^*)$ by Leibniz rule. Thus the proof is reduced to an induction on the number of edges contained in $\bstau$.   
One also obtains \eqref{EqCharacterizationGrafting2} by Assumption \refD{D2}.
\end{Dem}

\bigskip

\subsection{Renormalization structure over a multi-pre-Lie algebra}
\label{SectionrenormalizationMultiPreLie}

We now come to the main result of \cite{BCCH18} giving a dynamical meaning to the renormalization operations on models associated with elements $k\in G^-_\textrm{ad}$ of the renormalization group and more generally to elements $k\in G^-$. We keep working on the example of the generalized (KPZ) equation.

\medskip

In Theorem \ref{thm: all requirements of regul and renor str} in Section \ref{SectionBuildingRS}, we show that one can choose $U$ stable under all 
the splitting maps $\updownharpoons_{\sf e}$, that is
\begin{equation} \label{EqSplittingU}
\updownharpoons_{\sf e}(U)\subset U\otimes U
\end{equation}
for any ${\sf e}\in {\sf E}$. The restricted map 
$$
{\updownharpoons_{\sf e}}\vert_{U} : U\to U\otimes U
$$ 
is then the dual of the map $\flatgrafting$ for any ${\sf e}\in{\sf E}$. The following assumption is thus to be understood as a constraint on which renormalization schemes $\delta$ can be used.

\medskip

\begin{assumD}\label{D3}
\begin{enumerate}
\setlength{\itemsep}{0.5cm}
	\item For any ${\sf e}\in {\sf E}$, the space $U$ is stable under $\updownharpoons_{\sf 
e}$,
and one has
\begin{equation}   \label{EqCompatibilityConditionPreLiRenorm}
\big(\textrm{\emph{Id}}\,\otimes ({\updownharpoons_{\sf e}}\vert_U)\big)\delta 
= \mcM^{(13)}\big(\delta\otimes \delta\big){\updownharpoons_{\sf e}}\vert_U. 
\end{equation}
and
\begin{equation}   \label{EqCompatibilityConditionPreLiRenorm2}
\delta\, \circ \downarrow_i\, =\big(\textrm{\emph{Id}}\, \otimes\downarrow_i\big)\delta.
\end{equation}
	\item When $\bstau=\tau_{\frak{e}}^{\frak{n}}$ is an element of $\mcB$ without $\textcolor{red}{\bullet}^{\cdot,\alpha}$ decorations, then for any subforest $\varphi$ of $\tau$, $\frak{n}_\varphi:N_\varphi\to\bbN^{d+1}$ with $\frak{n}_\varphi\le\frak{n}$, and $\frak{e}_{\partial\varphi}':\partial\varphi\to\bbN^{d+1}$, the tree
\begin{equation}\label{D3:contraction}
(\tau/^{\text{\rm red}}\varphi)_{\frak{e}+\frak{e}_{\partial\varphi}'}^{[\frak{n}-\frak{n}_\varphi]_\varphi,\frak{o}(\frak{n}_\varphi+\pi\frak{e}_{\partial\varphi}',\frak{e})}
\end{equation}
is also contained in $\mcB$ -- see Section \ref{SectionDecoratedTreesCoproducts} for the notation.
\end{enumerate}
\end{assumD}

\medskip

Assumptions \refD{D1}, \refD{D2}, and \refD{D3} are jointly called Assumption \REFD.
Identity \eqref{EqCompatibilityConditionPreLiRenorm} is the ${\sf E}$-multi-pre-Lie version of the compatibility condition \eqref{*EqCompatibilityCondition} between the splitting map $\Delta$ of a regularity structure and 
a renormalization splitting $\delta$. 
Recall that any character $k$ of $U^-$ defines a linear map $\widetilde{k}=(k\otimes\iden)\delta:U\to U$. Denote by $\widetilde{k}^*: U^*\rightarrow U^*$ the dual map of $\widetilde{k}$ under the pairing \eqref{Eq normalized pairing}. Anticipating over Section \ref{SectionBuildingRS}, say here that $\textcolor{red}{\bullet}^{0,\alpha}$ is used to denote the result of extracting from a decorated tree $\bstau$ the entire tree, but keeping track of the homogeneity $\alpha=\vert\bstau\vert$ of the tree that was removed. Using the duality relation defining $\widetilde{k}^*$ and the definition of $\widetilde{k}$ we see that 
$$
\widetilde{k}^*(\circ^*)=\circ^*,\quad \widetilde{k}^*(\bullet^*)=\bullet^*,
$$
and
$$
\widetilde{k}^*\big((\textcolor{red}{\bullet}^{0,\alpha})^*\big) = 0
$$
for $\alpha>0$, and 
$$
\widetilde{k}^*\big((\textcolor{red}{\bullet}^{0,\alpha})^*\big) = \sum_{\bstau\in\mcB,\vert\bstau\vert=\alpha} \frac{k(\bstau)}{S(\bstau)}\,\bstau^*,
$$
for $\alpha<0$. The following result is part of Proposition 4.18 in Bruned, Chandra, Chevyrev and Hairer's work \cite{BCCH18}. It is the reason why we insisted on making a difference between $U$ and $U^*$, to emphasize the dual action of $\widetilde{k}$.

\medskip

\begin{prop}   \label{PropRenormMultiPreLie}
Under the compatibility Assumption \refD{D3}, given any character $k$ on $U^-$, the map $\widetilde{k}^*$ is an ${\sf E}$-multi-pre-Lie morphism: For any edge type ${\sf e}\in {\sf E}$, and any $\bstau, \bssigma\in \mcB$, one has
$$
\widetilde{k}^*(\bstau^*)\flatgrafting\widetilde{k}^*(\bssigma^*) = \widetilde{k}^*(\bstau^*\flatgrafting\bssigma^*),
$$
and
\begin{align}\label{PropRenormMultiPreLieEq2}
\widetilde{k}^*\circ\uparrow_i\, =\, \uparrow_i\circ\,\widetilde{k}^* \qquad (1\leq i\leq d).
\end{align}
\end{prop}

\medskip

\begin{Dem} 
We prove the dual identities writing
\begin{equation*} \begin{split}
\updownharpoons_{\sf e}\circ\,\widetilde{k} 
&= (k\,\otimes\updownharpoons_{\sf e})\delta 
\overset{\eqref{EqCompatibilityConditionPreLiRenorm}}{=} \Big(k\otimes \textrm{Id}\otimes \textrm{Id}\Big) \mcM^{(13)}\big(\delta\otimes \delta\big)\updownharpoons_{\sf e}
= \Big(\big(k\otimes \textrm{Id}\big)\delta\otimes \big(k\otimes \textrm{Id}\big)\delta\Big)\updownharpoons_{\sf e}
= \big(\,\widetilde{k}\otimes \widetilde{k}\,\big)\updownharpoons_{\sf e},
\end{split} \end{equation*}
and
$$
\downarrow_i\circ\,\widetilde{k}
=(k\, \otimes\downarrow_i)\delta
\overset{\eqref{EqCompatibilityConditionPreLiRenorm2}}{=}
(k\otimes\iden)\delta\,\circ\downarrow_i
=\widetilde{k}\,\circ\downarrow_i.
$$
\end{Dem}

\medskip

Pick a character $k$ on $U^-$. For primitive trees $\frak{b}\in{\sf N}^0$ 
define
\begin{align} \label{EqDefnFk}
F^{(k)}(\frak{b}^*) \defeq  F\Big(\widetilde{k}^*(\frak{b}^*)\Big),
\end{align}
so we have
\begin{align*}
F^{(k)}(\circ^*) = f({\sf u}_0),\quad
F^{(k)}(\bullet^*) = g({\sf u}_0,{\sf u}_1),\quad
F^{(k)}\big((\textcolor{red}{\bullet}^{0,\alpha})^*\big) 
= {\bf 1}_{\alpha<0}\sum_{\bstau\in\mcB\cap U_\alpha}\frac{k(\bstau)}{S(\bstau)}\,F(\bstau^*).
\end{align*} 
By Lemma \ref{lem ExpSol bsu}, $F^{(k)}\big((\textcolor{red}{\bullet}^{0,\alpha})^*\big)$ is also a function of only ${\sf u}_0$ and ${\sf u}_1$. For a tree 
$$
\bstau = \frak{b}^n\star\Bigstar_{i=1}^a\mcI_{p_i}(\bstau_i),
$$ 
define inductively the functions of $(u_0,u_1)$
\begin{align}\label{defn: Fkfg}
F^{(k)}(\bstau^*)({\sf u}) \defeq  \left\{\prod_{i=1}^a F^{(k)}(\bstau_i^*)(u_0,u_1)\right\} \left\{D^n\prod_{i=1}^a \partial_{{p_i}}\right\}F^{(k)}(\frak{b}^*)({\sf u}),
\end{align}
similarly to \eqref{defn: Ffg}. 
Note that $F^{(k)}$ can be defined for all elements $\bstau$ in $\mcV$. The map $F^{(k)}$ is an ${\sf E}$-multi-pre-Lie morphism with respect to $\{\grafting\}_{{\sf e}\in {\sf E}}$ by the same proof as the proof of Proposition \ref{PropPreLieMorphism}. We prove the following proposition to ensure that $F^{(k)}(\bstau^*)=0$ for any $\bstau\in\mcV\setminus\mcB$. The notations are all defined in Section \ref{SectionDecoratedTreesCoproducts}; the reader can skip it now and come back to it later.

\medskip

\begin{prop}\label{Prop origin of red node}
Let $k$ be a character of $U^-$.
For any $\bstau\in\mcV$, the function $F^{(k)}(\bstau^*)$ is represented as a linear combination of the functions $F(\bssigma^*)$,
where $\bssigma=\sigma_{\frak{e}}^{\frak{n}}$ runs over all elements of 
$\mcV$
without $\textcolor{red}{\bullet}^{\cdot,\alpha}$ decorations and such that $\bstau$ is obtained by contracting $\bssigma$ by its subforest $\varphi$ with maps $\frak{n}_\varphi$ and $\frak{e}_{\partial\varphi}'$ as the form \eqref{D3:contraction}.
\end{prop}

\medskip

\begin{Dem}
The proof is an induction on the number of edges contained in $\bstau$. It is sufficient to consider $\bstau$ of the form
$$
\bstau = \textcolor{red}{\bullet}^{n,\alpha}\star\Bigstar_{i=1}^a \mcI_{p_i}(\bstau_i).
$$
Assume that the result holds for each $\bstau_i$. Then by definition, $F^{(k)}(\bstau^*)$ is a linear combination of the functions of the form
$$
\left\{\prod_{i=1}^aF(\bssigma_i^*)\right\}
\left\{D^n\prod_{i=1}^a \partial_{{p_i}}\right\}F(\bseta^*),
$$
where $\bssigma_i$ is an element of $\mcV$ without $\textcolor{red}{\bullet}^{\cdot,\alpha}$ decorations such that $\bstau_i$ is obtained as a contraction of $\bssigma_i$, and $\bseta\in\mcB\cap U_\alpha$. By the formula $D^n\partial_kF=\sum_{\ell\in\bbN^{1+d}}(-1)^\ell\binom{n}{\ell}\partial_{k-\ell}D^{n-\ell}F$ obtained similarly to \eqref{lem:commutationpartialD} and by Proposition \ref{PropPreLieMorphism2}, the above function is a linear combination of the functions
$$
\left\{\prod_{i=1}^aF(\bssigma_i^*)\right\}
\left\{\prod_{i=1}^a \partial_{q_i}\right\}F(\boldsymbol{\mu}^*).
$$
where $q_i\le p_i$ and $\boldsymbol{\mu}\in\mcB$ is a tree appearing in the expansion of $\uparrow^m\bseta\defeq(\prod_{i=0}^d\uparrow_i^{m_i})\bseta$ for some $m\in\bbN^{1+d}$. By an argument similar to the proof of Proposition \ref{PropPreLieMorphism}, we can show that the above function is equal to
$$
F\Big(\sum_{v_1,\dots,v_a\in N_\mu}
\big(\bssigma_1\overset{q_1}{\curvearrowright}_{(v_1)}
(\bssigma_2\overset{q_2}{\curvearrowright}_{(v_2)}
\cdots
(\bssigma_a\overset{q_a}{\curvearrowright}_{(v_a)}\boldsymbol{\mu})\cdots)
\big)^*\Big)
$$
The trees inside $F$ produce $\bstau$ when we contract $\bssigma_1,\dots,\bssigma_a,\boldsymbol{\mu}$ as in the explicit formula of ${\sf D}^-$ in Section \ref{SectionDecoratedTreesCoproducts}.
\end{Dem}

\medskip

As a result of Proposition \ref{Prop origin of red node}, under Assumptions \refD{D2} and \refD{D3}, we have that $F^{(k)}(\bstau^*)=0$ for any $\bstau\in\mcV\setminus\mcB$,
since $\bstau$ is an element of $\mcB$ if and only if $\bstau$ is a contraction of a tree in $\mcB$ without $\textcolor{red}{\bullet}^{\cdot,\alpha}$ decoration. This ensures that the map
$$
\Upsilon^{(k)}\defeq F^{(k)}\vert_{U^*},
$$
is an $\sf E$-multi-pre-Lie morphism on $U^*$ with respect to $\{\flatgrafting\}_{\sf e\in E}$.
Moreover, denoting by $S$ the subspace of $T$ spanned by $\{\bullet^n\}_{n\in\bbN^{d+1}}\cup\mcI(\mcB)$,
we have that for any $\frak{b}\in{\sf N}^0$ and any function $\bsu:\bbR^{d+1}\to S$, the function
$$
\big(F^{(k)}(\frak{b}^*)\big)^\star(\bsu, D\bsu  )\star\frak{b}:\bbR^{d+1}\to V
$$
is actually $T$-valued.

\medskip

\begin{cor}   \label{CorrenormalizedVectorFields}
Under Assumption \REFD, one has $\Upsilon\circ\widetilde{k}^* = \Upsilon^{(k)}$, for all $k\in G^-$.
\end{cor}

\medskip

\begin{Dem} 
It follows from Propositions \ref{PropPreLieMorphism} and \ref{PropRenormMultiPreLie} that the map $\Upsilon\circ\widetilde{k}^*$ is an ($\flatgrafting$ vs $\triangleright_{\sf e}$) $\sf E$-multi-pre-Lie morphism. 
Because of Proposition \ref{Prop origin of red node}, the map $\Upsilon^{(k)}$ is also an ($\flatgrafting$ vs $\triangleright_{\sf e}$) $\sf E$-multi-pre-Lie morphism. Hence it is sufficient to show that they are equal on the generators of $U^*$, that is,
$$
F\big(\widetilde{k}^*(\frak{b}^n)^*\big)=F^{(k)}\big((\frak{b}^n)^*\big)
$$
for any $(\frak{b},n)\in\frak{T}_{\rm n}\times\bbN^{d+1}$. The case $n=0$ is given by definition \eqref{EqDefnFk}. For $n=(n_i)_{i=0}^d\in\bbN^{d+1}$, by writing $\uparrow^n\defeq \prod_{i=0}^d\uparrow_i^{n_i}$, we 
have
\begin{align*}
F\big(\widetilde{k}^*(\frak{b}^n)^*\big)
&=F\big(\widetilde{k}^*(\uparrow^n\frak{b}^*)\big)
\overset{\eqref{PropRenormMultiPreLieEq2}}{=}
F\big(\uparrow^n(\widetilde{k}^*\frak{b}^*)\big)
\overset{\eqref{EqCharacterizationGrafting2}}{=}
D^nF\big(\widetilde{k}^*\frak{b}^*\big)
\overset{\eqref{EqDefnFk}}{=}
D^nF^{(k)}(\frak{b}^*)   \\
&=F^{(k)}\big((\frak{b}^n)^*\big),
\end{align*}
using \eqref{EqCharacterizationGrafting2} and \eqref{PropRenormMultiPreLieEq2} in the last equality.
\end{Dem}

\medskip

Similarly to the definition of $\mcF(\bsu)$, given a modelled distribution $\bsu\in\mcD^{\gamma,\eta}(T,\sf g)$ with $\gamma\in(-\beta_0,2)$ such that $(0,\gamma+\beta_0)\cap A=\emptyset$, set
\begin{align}\label{section 6: renormalized force term in MD sense}
\mcF^{(k)}(\bsu) \defeq  \mathcal{Q}_{\le 0}\left(\sum_{\frak{b}\in{\sf N}^0}\big(F^{(k)}(\frak{b}^*)\big)^\star(\bsu, D\bsu  )\star\frak{b}\right).
\end{align}
Note the appearance in \eqref{section 6: renormalized force term in MD sense} of a number of symbols ${\color{red} \bullet}^{0,\alpha}$, with $\alpha<0$, that have no counterpart in $\mcF(\bsu)$.

\medskip

\begin{lem}\label{lem kF=Fk}
If $\bsu$ is a solution of equation \eqref{EqRSgKPZ} then
\begin{equation} \label{EqIdentityTildeKUpsilon}
\widetilde{k}\big(\mcF(\bsu)\big) = \mcF^{(k)}\big(\widetilde{k}(\bsu)\big)
\end{equation}
on the domain $(0,\frac{t_0}2)\times\bbT^d$.
\end{lem}

\medskip

\begin{Dem}
Let $\bsu$ be evaluated at the fixed $x$ in the domain $(0,\frac{t_0}2)\times\bbT^d$.
Lemma \ref{lem ExpSol bsu} implies that 
$$
\mcF(\bsu)=\sum_{|\bstau|<\gamma+\beta_0}\frac{\Upsilon(\bstau^*)\big((u_k)_{k\in\bbN^{1+d}}\big)}{S(\bstau)}\bstau,
$$
or equivalently, $\wangle{\mcF(\bsu),\bstau^*} = 
\Upsilon(\bstau^*)$, for any $\bstau\in \mcB$ with $|\bstau|<\gamma+\beta_0$. Noting that $\widetilde{k}$ and $\widetilde{k}^*$ preserve the grading of $T$, as a consequence of the compatibility condition \eqref{delta T to U-T}, we have
\begin{align*}
\wangle{\widetilde{k}\big(\mcF(\bsu)\big),\bstau^*}
= \wangle{\mcF(\bsu),\widetilde{k}^*(\bstau^*)}
= \Upsilon\big(\widetilde{k}^*(\bstau^*)\big)
= \Upsilon^{(k)}(\bstau^*)
\end{align*}
for any $\bstau\in \mcB$ with $|\bstau|<\gamma+\beta_0$, or equivalently,
$$
\widetilde{k}\big(\mcF(\bsu)\big)=\sum_{|\bstau|<\gamma+\beta_0}\frac{\Upsilon^{(k)}(\bstau^*)\big((u_k)_{k\in\bbN^{1+d}}\big)}{S(\bstau)}\bstau.
$$
Next we consider the $\bstau$-component of $\mcF^{(k)}\big(\widetilde{k}(\bsu)\big)$ by an argument similar to Lemma \ref{lem ExpSol bsu}. Note that Assumption \textbf{\textsf{(C)}} yields $\widetilde{k}(X^n) = 
X^n$, and $\widetilde{k}\big(\mcI(\bstau)\big) = \mcI\big(\widetilde{k}(\bstau)\big)$. Hence
\begin{align*}
\widetilde{k}(\bsu)
= \widetilde{k}\left\{\sum_{|k|_\mfs<\gamma}\frac{u_k}{k!}X^k+\mcQ_{<\gamma}\mcI\big(\mcF(\bsu)\big)\right\}
= \sum_{|k|_\mfs<\gamma}\frac{u_k}{k!}X^k+ \sum_{\bstau}\frac{\Upsilon^{(k)}(\bstau^*)}{S(\bstau)}\,\mcI(\bstau).
\end{align*}
Thus by an argument similar to the former half part of Lemma \ref{lem ExpSol bsu}, we see that the $\bstau$-component of $\mcF^{(k)}\big(\widetilde{k}(\bsu)\big)$ is equal to \eqref{EqCoherence*}, where $u_{\mcI(\bssigma_j)}$ is replaced by $\frac{\Upsilon^{(k)}(\bssigma_j^*)}{S(\bssigma_j)}$, and $F(\frak{b}^*)$ is replaced by $\Upsilon^{(k)}(\frak{b}^*)$.
By definition of $\Upsilon^{(k)}(\bstau^*)$, we see that the $\bstau$-component of $\mcF^{(k)}\big(\widetilde{k}(\bsu)\big)$ is equal to $\frac{\Upsilon^{(k)}(\bstau^*)}{S(\bstau)}$.
Therefore
$$
\wangle{\widetilde{k}\big(\mcF(\bsu)\big),\bstau^*}=\wangle{\mcF^{(k)}\big(\widetilde{k}(\bsu)\big),\bstau^*}
$$ 
for any $\bstau\in \mcB$ with $|\bstau|<\gamma+\beta_0$.
\end{Dem}

\medskip

The next statement provides a dynamical picture of the renormalization operation on models. As its proof will make it clear, it is a consequence of identity \eqref{EqIdentityTildeKUpsilon} and Theorem \ref{ThmrenormalizationActionModels}, giving in particular the reconstruction operator of a 
renormalized smooth model in terms of the unrenormalized smooth model, together with the multiplicativity property of the canonical model associated with a smooth noise.  

\medskip

\begin{thm}   \label{ThmRenormPDEs}
Assume \REFD. Let $\zeta$ be a smooth noise with canonical model $\sf M^\zeta = (\Pi^\zeta, g^\zeta)$. Given a character $k\in G^-_\textrm{\emph{ad}}$, denote 
by $^k{\sf M}^\zeta = \big({\sf g}^\zeta\circ\widetilde{k}^+,{\sf\Pi}^\zeta\circ\widetilde{k}\big)$ its associated renormalized $\bfK$-admissible model -- see Theorem \ref{ThmrenormalizationActionModels}. Pick $\eta\in(0,\beta_0+2]$ and $\gamma\in(-\beta_0,2)$. Given an initial condition $v\in\mcC^\eta(\bbT^d)$, let $\bsu^{(k)}\in\mcD_{(0,t_0)}^{\gamma,\eta}\big(T, {\sf g}^\zeta\circ\widetilde{k}^+\big)$ stand for the solution on $(0,t_0)$ to the equation
$$
\bsu^{(k)} = P_\gamma v+\mcP_{t_0}^{{}^k{\sf M}^\zeta}
\Big(f^\star\big({\bsu^{(k)}}\big)\Xi + g^\star\big({\bsu^{(k)}},D{\bsu^{(k)}}\big)\Big).
$$
Then 
$$
u^{(k)} \defeq  \textbf{\textsf{R}}^{^k{\sf M}^\zeta}(C_{t_0}{\bsu^{(k)}})
$$
is the solution on $(0,\frac{t_0}2)$ to the well-posed equation
$$
\big(\partial_{x_0} - \Delta_{x'}+1\big) u^{(k)} = f\big({u^{(k)}}\big)\zeta + 
g\big({u^{(k)}},\partial_{x'}{u^{(k)}}\big) + \sum_{\tau\in \mcB,\, |\tau|<0} \frac{k(\bstau)}{S(\bstau)}\,F(\bstau^*)\big({u^{(k)}},\partial_{x'}{u^{(k)}}\big)
$$
started from $v$.
\end{thm}

\medskip

\begin{Dem}
The proof is similar to the proof of Proposition \ref{PropRSCharacterizationSmoothPDE}. The function $u^{(k)}$ satisfies the equation
\begin{align*}
u^{(k)}(x) = Pv(x)+\int_{(0,x_0)\times\bbT^d} P(x,y) \textbf{\textsf{R}}^{^k{\sf M}^\zeta} C_{t_0}\Big(f^\star(\bsu^{(k)})\Xi + g^\star\big(\bsu^{(k)},D\bsu^{(k)}\big)\Big) (y)dy.
\end{align*}
Since $^k{\sf M}^\zeta$ is a smooth model one has
$$
\textbf{\textsf{R}}^{^k{\sf M}^\zeta}(\bsw) (x) = {\sf \Pi}_x^\zeta\Big(\widetilde{k}\big(\bsw(x)\big)\Big)(x).
$$
for any modelled distribution $\bsw\in\mcD^\alpha\big(T, {\sf g}^\zeta\circ\widetilde{k}^+\big)$ with $\alpha>0$. Applying Lemma \ref{lem kF=Fk} 
to $\bsu^{(k)}$ one has
\begin{align*}
&\textbf{\textsf{R}}^{^k{\sf M}^\zeta}
C_{t_0}\Big[f^\star(\bsu^{(k)})\Xi + g^\star\big(\bsu^{(k)},D\bsu^{(k)}\big)\Big](x) \\
&= {\sf \Pi}_x^\zeta\left[\widetilde{k}\Big(\mcF(\bsu^{(k)}(x)\Big)\right](x) = {\sf\Pi}_x^\zeta\left[\mcF^{(k)}\Big(\widetilde{k}\big(\bsu^{(k)}(x)\big)\Big)\right](x).
\end{align*}
for any $x\in(0,\frac{t_0}2)\times\bbT^d$.
We see from the definition of the $F^{(k)}$ that the term $\mcF^{(k)}(\bsw)$ is a sum of functions of the form
$$
H^\star(\bsw)R^\star(D\bsw),
$$
for smooth functions $H:\bbR\to\bbR$ and polynomials $R$ that are at most 
quadratic. We now use the fact that since the map ${\sf \Pi}^\zeta$ is {\it multiplicative} so is its associated reconstruction operator. 
The latter has value ${\sf\Pi}^\zeta_x(\cdot)(x)$ at point $x$, so we have
$$
{\sf\Pi}_x^\zeta\Big[H^\star\big(\bsw(x)\big)R^\star\big(D\bsw(x)\big)\Big](x) = 
H\Big[\big({\sf\Pi}_x^\zeta\bsw(x)\big)(x)\Big]
R\Big[\big({\sf\Pi}_x^\zeta D\bsw(x)\big)(x)\Big],
$$
with $\bsw(x) = \widetilde{k}\big(\bsu^{(k)}(x)\big)$.
Since $\big({\sf\Pi}_x^\zeta\bsw(x)\big)(x)=\textbf{\textsf{R}}^{^k{\sf M}^\zeta}(C_{t_0}\bsu^{(k)})(x)=u^{(k)}(x)$
and $\big({\sf\Pi}_x^\zeta D\bsw(x)\big)(x)=\textbf{\textsf{R}}^{^k{\sf M}^\zeta}(C_{t_0}D\bsu^{(k)})(x)=\partial_{x'}u^{(k)}(x)$ on the domain $(0,\frac{t_0}2)\times\bbT^d$, we have in the end
$$
{\sf\Pi}_x^\zeta\Big[\mcF^{(k)}\Big(\widetilde{k}\big(\bsu^{(k)}(x)\big)\Big)\Big](x) = \mcF^{(k)}\Big(u^{(k)},\partial_{x'} u^{(k)}\Big)(x).
$$
\end{Dem}

\medskip

{\small {\sl \begin{Rem} 
The preceding proof underlines the fundamental role played by the multiplicative property of the centered naive interpretation operators ${\sf \Pi}^\zeta_x$. The canonical smooth model ${\sf \Pi}^\zeta$ is not the only multiplicative model that one can associate with a smooth noise $\zeta$. The class of models associated with `preparation maps' introduced by Bruned in \cite{BrunedRecursive} provides a general setting where to obtain the renormalized equation for a class of renormalization procedures including the procedure implemented here \cite{BailleulBruned}.
\end{Rem}}}

\bigskip

\section{The BHZ character}
\label{SectionBHZCharacter}

Among all the characters $k$ on $U^-$ that can be used to build a renormalization map $\widetilde{k}$, Bruned, Hairer and Zambotti proved in \cite{BHZ} that there is a unique character $k$ whose associated random model is centered and translation invariant, in a probabilistic sense, when the smooth noise $\zeta$ in the preceding section is random, centered and translation invariant. We describe it in this section and name it `BHZ character', after the initials of Bruned, Hairer and Zambotti. We also call the associated renormalized model the BHZ model.

\ssk

Arrived at that stage, the only piece of the story that will be missing to have a complete proof of the meta-theorems from Section \ref{SectionIntro} will be a proof of the fact that one can indeed construct some regularity structures satisfying the different assumptions that we put forward in the course of obtaining the above results, and a proof of convergence of the BHZ smooth renormalized models. We will tackle the first point in Section \ref{SectionBuildingRS}. The second point is the object of Chandra \& Hairer's work \cite{ChandraHairer}; we do not treat it here. We refer the reader to the work \cite{LOTT, HS, BH23} for some alternative proofs of the convergence of BHZ renormalized models in situations where the law of the noise satisfies a spectral gap inequality.

\medskip

{\it We assume throughout this section that we work with regularity and renormalization structures satisfying Assumptions \textbf{\textsf{(A-C)}}.} To have a picture in mind, think of the structures associated with the generalized (KPZ) equation \eqref{EqgKPZ}. Elements of $T=U$ are thus given by node and edge decorated trees. Denote by $\bbR[U]$ the commutative algebra generated by $U$. Recall from Assumption \refC{C1} that $U^-$ is an algebra generated by $\mcU_{<0}$ and a unit $\textbf{\textsf{1}}_-$, and if one extends first the splitting map $\delta:U\to U^-\otimes U$ into an algebra morphism $\hat\delta : \bbR[U]\to U^-\otimes\bbR[U]$, then the splitting map $\delta^-$ satisfies
$$
\delta^- = (\textrm{Id}\otimes P_-){\widehat\delta}\vert_{U^-},
$$
for an algebra morphism projection map $P_- : \mathbb{R}[U]\rightarrow U^-$. (See Section \ref{SectionDecoratedTreesRegRenor} for the details; this projection map sends in particular $\bullet$ and all the ${\color{red} \bullet}^{0,\alpha}$ to $\textbf{\textsf{1}}_-$). Denote by {\color{red}$\widehat{\bf1}_{-}$} the unit of $\mathbb{R}[U]$, seen as the empty graph. (We use a distinct notation for {\color{red}$\widehat{\bf1}_{-}$} and ${\bf1}_-$ to emphasize that they do not live in the same space.) For basis elements $\bstau=\tau_{\frak{e}}^{\frak{n}}$ and $\bssigma=\sigma_{\frak{f}}^\frak{m}$ of $\bbR[U]$, we write
$$
\bssigma\prec\bstau
$$
if $\sigma$ is a strict {\it subgraph} of $\tau$, or $\sigma=\tau$ and $\frak{m}\leq \frak{n}$ with $\frak{m}\neq\frak{n}$. (Do not get misled by the notation, $\sigma$ may be a product of disjoint subtrees of $\tau$.) Recall from Section \ref{SubsectionCompatibleStructures} the notation $\mcF=\{\mcI_p\}_{|p|_\mfs\le1}\cup\{X^n\star\}_{n\in\bbN^{1+d}\setminus\{0\}}$ for the family of integral operators and multiplication by a non-null monomial. The following assumption describes the properties from the splitting map $\delta$ that are relevant here. The renormalization structure built in Section \ref{SectionBuildingRS} for the generalized (KPZ) equation satisfies it.


\medskip

\noindent \textbf{\textsf{Assumption (E) -- }}{\it 
\begin{enumerate}
\setlength{\itemsep}{0.1cm}

   \item[\emph{\textsf{(a)}}] For any $\bstau\in\mcU_{<0}$, one has the splitting formula
   $$
   \delta\bstau \in P_-(\bstau)\otimes \textcolor{red}{\bullet}^{0,\vert\tau\vert}+U_{\bstau}^-\otimes U,
   $$
   where 
   $$
   U_{\bstau}^- \defeq  \text{\rm span}\Big\{P_-(\bssigma)\in U^-\, ;\, \bssigma\prec\bstau\Big\}.
   $$
   
   \item[\emph{\textsf{(b)}}] For any $\bstau\in\mcU_{<0}$ and $F\in\mcF$ 
such that $F(\bstau)\in\mcU_{<0}$, one has
   $$
   \delta\big(F(\bstau)\big)\in (\text{\rm Id}\otimes F)\delta\bstau + P_-\big(F(\bstau)\big)\otimes\textcolor{red}{\bullet}^{0,\vert F(\bstau)\vert} + \frak{J}_{\bstau}^-\otimes U,
   $$
   where $\frak{J}_{\bstau}^-$ is the ideal of $U^-$ generated by $\Big\{P_-\big(F(\bssigma)\big)\, ;\, F\in\mcF,\, \bssigma\in \mcU_{<0},\, \bssigma\prec\bstau\Big\}$.
   
\end{enumerate} }   

\medskip

Property \textsf{(b)} is a refinement of the property \eqref{eq relation delta and K} in Assumption \textbf{\textsf{(C1)}}. 
Note that, under the definitions of gradings of $T$ and $U$ in Section \ref{SectionBuildingRS} later, $\textcolor{red}{\bullet}^{0,\beta}\in\mcB_\beta$ but $\textcolor{red}{\bullet}^{0,\beta}\in\mcU_0$, for any $\beta\in\bbR$. (Basis elements of $U$ with $0$-homogeneity are not necessarily unique, unlike in Assumption \refA{A1} on concrete regularity structures.) Hence the above properties are consistent with the definitions of compatible renormalization and regularity structures. Property \textsf{(b)} also ensures the existence of the following map, defined by induction on the order relation $\prec$. Each element $\bstau$ of $U^-$ has by definition a unique representative $\bstau_U$ in $\mathbb{R}[U]$. Denote by $\mcM$ the multiplication operator on $\mathbb{R}[U]$ and extend it naturally on $U^-\otimes\mathbb{R}[U]$ setting $\mcM(\bstau\otimes\bssigma)=\mcM(\bstau_U\otimes\bssigma)$; it takes values in $\mathbb{R}[U]$.

\medskip

\begin{defn*}
Under Assumption \textbf{\textsf{(E)}}, the \textbf{\textsf{negative twisted antipode}} is an algebra morphism 
$$
S'_-: U^-\hspace{-0.1cm}\rightarrow\bbR[U]
$$
given recursively by $S_-'{\bf1}_- = \widehat{\bf1}_-$ and, for every basis element $\bstau\in\mcU_{<0}$ by
\begin{equation}   \label{EqCharacterizationTwistedNegativeAntipode}
S'_-\big(P_-(\bstau)\big) = -\mcM_-\big(S'_-\otimes \textrm{\emph{Id}}\big)\Big(\delta\bstau - P_-(\bstau)\otimes\textcolor{red}{\bullet}^{0,\vert\bstau\vert}\Big).
\end{equation}
\end{defn*}

\medskip

The $P_-(\bstau)$ generating $U^-$ as an algebra for $\bstau$ ranging in $\mcU_{<0}$, identity \eqref{EqCharacterizationTwistedNegativeAntipode} characterizes indeed uniquely an algebra morphism. The intuitive meaning of this recursive definition should be clear. One extracts from $\bstau$ all possible subdiverging quantities $\varphi_1$, but also extracts from $\varphi_1$ all its subdiverging quantities, and so on. This formula is close to 
the Dyson-Salam renormalization formula for the antipode in Hopf algebras 
\cite{FigueroaGB}; like the latter, it can be rewritten as a sum over forests of diverging sub-forests, as in Zimmermann forest formula. This will 
not be useful here, and the only thing that matters here is property \eqref{EqCharacterizationTwistedNegativeAntipode}. The forest representation is however useful for the analysis of the convergence of renormalized models \cite{ChandraHairer}.

\ssk

Do not be mislead by the name  of $S'_-$: This is {\it not} the antipode of a Hopf algebra structure. Bruned, Hairer and Zambotti named it like that because its defining relation \eqref{EqCharacterizationTwistedNegativeAntipode} looks like the defining relation \eqref{EqDefnAntipode} for the 
antipode in a Hopf algebra.

\ssk

Recall from Section \ref{SectionMultiPreLie} the definition of the naive interpretation operator ${\sf \Pi}^\zeta$ corresponding to a smooth noise $\zeta$ in $\bbR\times\bbR^d$. We consider a random smooth noise $\zeta$, invariant by translation and centered. 
Define the character $h^\zeta$ on $\mathbb{R}[U]$ by setting
$h^\zeta(\widehat{\bf1}_{-})\defeq 1$ and
\begin{equation} \label{EqRenormCharacter}
h^\zeta(\bstau) \defeq  \bbE\big[{\sf \Pi}^\zeta\bstau\big](0)
\end{equation}
for $\bstau\in\mcU$, and define a character on $U^-$ setting 
$$
k^\zeta \defeq  h^\zeta\circ S_-'.
$$
(Keep in mind that $S_-'$ gives back elements of $\bbR[U]$ and that $h^\zeta$ is multiplicative, so $k^\zeta$ is multiplicative on $U^-$.) The associated `{\sf BPHZ} renormalized' interpretation operator $^{k^\zeta}{\sf \Pi}^\zeta$ is defined on $U$ by 
$$
 ^{k^\zeta}{\sf \Pi}^\zeta\bstau = \big(k^\zeta\otimes {\sf \Pi}^\zeta\big)\delta\bstau = \Big((h^\zeta\circ S'_-)\otimes {\sf \Pi}^\zeta\Big)\delta\bstau.
$$
The acronym BPHZ stands for Bogoliubov, Parasiuk, Hepp and Zimmermann, who made deep contributions to the renormalization problem in quantum field theory. We call \textbf{\textsf{BHZ character}}, after Bruned, Hairer and Zambotti, the character $k^\zeta$ on $U^-$. The reason for introducing the negative twisted antipode operator lies entirely in the following simple computations used in the proof of the next statement claiming that the BPHZ renormalization associated with the BHZ character recenters probabilistically the $\sf \Pi$ map at all points in spacetime. Its proof is taken from the proof of Theorem 6.17 in Bruned, Hairer and Zambotti's work \cite{BHZ} on the 
algebraic renormalization of regularity structures. 

\medskip

\begin{thm}   \label{ThmPropertyBPHZRenorm}
Let $\zeta$ stand for a smooth noise that is centered and translation invariant in law, and such that $\partial^n\zeta(0)$ has finite moments of any order for any $n\in\bbN^{d+1}$. We work with compatible regularity and renormalization structures under Assumptions \textbf{\textsf{(A-C)}} and Assumption \textbf{\textsf{(E)}}. The character $k^\zeta$ belongs to $G_{\text{\rm ad}}^-$ and one has
\begin{equation}   \label{EqStationaryModel}
\bbE\Big[\big( {}^{k^\zeta}{\sf \Pi}^\zeta\bstau\big)(x)\Big] = 0
\end{equation}
for any $\bstau\in \mcU_{<0}$ and $x\in\bbR\times\bbR^d$.
\end{thm}

\ssk

\begin{Dem}
First we show that $k^\zeta\in G_{\text{\rm ad}}^-$. Let $F\in\mcF$ and $\bstau\in\mcU_{<0}$ be such that $F\bstau\in\mcU_{<0}$. If $F=X^n\star$ then $h^\zeta(F\bstau)=\bbE\big[{\sf \Pi}^\zeta X^n\star\bstau\big](0)=x^n\bbE\big[{\sf \Pi}^\zeta\bstau(x)\big]|_{x=0}=0$. Next consider $F=\mcI_p$. Recall from Section \ref{SubsectionOperators} that we defined the operator $\bf K$ so that 
$$
\int_{\bbR\times\bbR^d} y^pK(x,y)\,dy = 0
$$
for all $p\in\bbN\times\bbN^d$ such that $\vert p\vert_\frak{s}<N$ for fixed $N\in\bbN$. Now pick $N=2$. Note that $\bbE[{\sf\Pi}^\zeta\bstau](x)=\bbE[{\sf\Pi}^\zeta\tau_{\frak{e}}^{\frak{n}}](x)$ is a polynomial of $x$ with degree $\sum_{v\in N_\tau}|\frak{n}(v)|_{\frak{s}}$, since $K(x,y)$ depends on $x-y$ only. We use it here to have
$$
h^\zeta\big(\mcI_p\bstau\big) = \int_{\bbR\times\bbR^d}\partial_x^pK(x,y)\bbE\big[{\sf \Pi}^\zeta\bstau(y)\big]=0,
$$
since $\sum_{v\in N_\tau}|\frak{n}(v)|_{\frak{s}}\le1$ (otherwise $|\bstau|\ge\beta_0+2>0$). Hence $h^\zeta(F\bstau)=0$ for all $F\in\mcF$. Since Assumption \textbf{\textsf{(E)}} guarantees that we have
\begin{align*}
S_-'(F\bstau) \in -\mcM_-\big(S_-'\otimes F\big)\delta\bstau+\mcM_-\big(S_-'\frak{J}_{\bstau}\otimes U\big),
\end{align*}
we can conclude that $h^\zeta\big(S_-'(F\bstau)\big)=0$, by an induction on the size of the graph $\tau$. Hence $k^\zeta\in G_{\text{\rm ad}}^-$.

\ssk

-- The negative twisteed antipode $S'_-$ is defined so as to have identity \eqref{EqStationaryModel} for $x=0$. Indeed, since ${\sf\Pi}^{\zeta}(\textcolor{red}{\bullet}^{0,\beta})\equiv1$ for all $\beta$, one has from 
the defining relation \eqref{EqCharacterizationTwistedNegativeAntipode} for the twisted antipode, for any $\bstau\in\mcU_{<0}$,
\begin{align*} 
\bbE\Big[\big({}^{k^\zeta}{\sf \Pi}^\zeta\bstau\big)(0)\Big]
&=\sum_{\bsvarphi\trianglelefteq\bstau} h^\zeta\big(S_-'(\bsvarphi)\big)\, \bbE\big[\big({\sf\Pi}^\zeta(\bstau/^-\bsvarphi)\big)(0)\big]\\
&= \sum_{\bsvarphi\trianglelefteq\bstau} h^\zeta\big(S'_-(\bsvarphi)\big)\,h^\zeta(\bstau/^-\bsvarphi)\\
&= h^\zeta\Big( \mcM\big(S'_-\otimes \textrm{Id}\big)\delta\bstau\Big) = h^\zeta\big(S_-'\tau\big)\Big(h^\zeta(\textcolor{red}{\bullet}^{0,\vert\bstau\vert})-1\Big) = 0.
\end{align*}
Recall the homogeneity and grading notions on $U$ and $T$ are different. It is the homogeneity of $\bstau$, seen as an element of $T$, that appears in $\textcolor{red}{\bullet}^{0,\vert\bstau\vert}$. It is elementary to 
go from $\bbE\big[\big({}^{k^\zeta}{\sf \Pi}^\zeta\tau\big)(0)\big]=0$, 
to $\bbE\big[\big({}^{k^\zeta}{\sf \Pi}^\zeta\tau\big)(x)\big]=0$, for all $x\in\bbR\times\bbR^d$, using the probabilistic translation invariance property of ${\sf \Pi}^\zeta$.
\end{Dem}

\medskip

{\small {\sl \begin{Rem} 
Note that the cointeraction identity between $\delta$ and $\delta^-$ implies that we have 
\begin{equation}  \label{EqTworenormalizations}
^{k_1\star k_2}{\sf \Pi}^\zeta = {^{k_1}({^{k_2}{\sf \Pi}^\zeta})},
\end{equation}
for any two characters $k_1,k_2$ on $U^-$. There is no other character $k$ on $U^-$ than $h^\zeta\circ S'_-$ such that the renormalized naive interpretation operator $^k{\sf \Pi}^\zeta \defeq  
(k\otimes {\sf \Pi}^\zeta)\delta$, has property \eqref{EqStationaryModel} 
of Theorem \ref{ThmPropertyBPHZRenorm}. The uniqueness claim amounts to proving that for any non-null character $k\neq 1$ there exists an 
element $\bstau\in U$ such that $\bbE\big[\big(^{k\star k^\zeta}{\sf \Pi}^\zeta\bstau\big)(0)\big]\neq 0$. See the second part of the proof of Theorem 6.18 in \cite{BHZ}. 
\end{Rem}}}

\medskip

Assume now that $\zeta=\xi_\epsilon$ is the regularized version of a random irregular noise $\xi$, centered and translation invariant, and write 
${\sf \Pi}^\epsilon$ for ${\sf \Pi}^{\xi_\epsilon}$. The BHZ character $h$ from \eqref{EqRenormCharacter} becomes $\epsilon$-dependent as well. Set 
\begin{equation} \label{EqDefnKEpsilon}
k_\epsilon \defeq  h^{\xi_\epsilon}\circ S'_-.
\end{equation}
Identity \eqref{EqTworenormalizations} tells us that if the maps $^{k_\epsilon}{\sf \Pi}^\epsilon$ converge to a limit when $\epsilon$ goes to zero, then for any character $k$ on $U^-$, the renormalized interpretation map $^{k\star k_\epsilon}{\sf \Pi}^\epsilon$ is also converging. There is thus a whole class of converging renormalization schemes indexed by the group of characters of $U^-$, if there is a single converging renormalization scheme. If we insist on building $\bf K$-admissible models, this provides a family of convergent models indexed by the renormalization group $G^-_\textrm{ad}$.   

\medskip

Recall the arguments in Section \ref{SectionMultiAndrenormalizedEquations}. We say that the family of smooth cylindrical functions $F=\{F(\frak{b}^*)\}_{\frak{b}\in{\sf N}^0}$ is a nonlinearity. Denote by $u=\textsc{Sol}(\xi ; F)$ for the solution to the PDE
$$
(\partial_{x_0}-\Delta_{x'}+1)u=\sum_{\frak{b}\in{\sf N}^0}F(\frak{b}^*)(u,\partial_{x'}u){\sf\Pi}^\xi\frak{b}
$$
driven by  the smooth noise $\xi_\epsilon$, associated with a given initial condition. 
The arguments in Section \ref{SectionrenormalizationMultiPreLie} implies that, the nonlinearity $F=\{F(\frak{b}^*)\}_{\frak{b}\in{\sf N}^0}$ can be extended to smooth cylindrical functions $F(\bstau^*)$ for any $\bstau\in\mcV$ by \eqref{defn: Fkfg}, and the group $G^-_\textrm{ad}$ acts on the set $\frak{F}$ of nonlinearities by $F^{(k)}\defeq F\circ\widetilde{k}^*$ as in Corollary \ref{CorrenormalizedVectorFields}. The renormalization group acquires a dynamical meaning from Theorem \ref{ThmRenormPDEs} if one notices that 
\begin{equation*}
\Big\{ \textsc{Sol}\Big(\xi_\epsilon ; F^{(k)}\Big)\Big\}_{k\in G^-_\textrm{ad}}
=
\Big\{ \textsc{Sol}\Big(\xi_\epsilon ; F^{(k_\epsilon*k)}\Big)\Big\}_{k\in G^-} 
= \Big\{ \textsc{Sol}\Big(\xi_\epsilon ; \big(F^{(k_\epsilon)}\big){}^{(k)}\Big)\Big\}_{k\in G^-}
 \end{equation*}   
for any fixed positive $\epsilon$. This remark tells us that the \emph{family} of solutions of the singular stochastic PDE \eqref{EqGenericPDE} is parametrized by the subset $\big(F^{(k)}\big)_{k\in G^-}$ of the space $\frak{F}$ of nonlinearities. This remains true at the limit when $\epsilon>0$ goes to $0$. We will see in the Section \ref{SectionManifold} that this subset is actually a finite dimensional immersed manifold.

\bigskip

\section{The manifold of solutions}   
\label{SectionManifold}

We take for granted in this section the convergence result of Chandra \& Hairer from \cite{ChandraHairer}, and work with the limit random admissible model ${\sf M}=({\sf g},{\sf\Pi})$, obtained as a limit in probability of the renormalized naive models $^{k_\epsilon}{\sf M}^\epsilon$ when $\epsilon>0$ goes to $0$. Recall from equality \eqref{section 6: renormalized force term in MD sense} the expression of $\mcF^{(k)}(\bsu)$, for $k\in G^-_\textrm{ad}$. Pick $\eta\in(0,\beta_0+2], \gamma\in(-\beta_0,2)$ and an initial condition $u_0\in\mcC^\eta(\bbT^d)$. Write $\bsu(k)\in\mcD^{\gamma,\eta}_{(0,t_0)}(T,\sf g)$ for the solution to the equation
\begin{equation*} \begin{split}
\bsu(k) &= P_\gamma u_0+\mcP_{t_0}^{{\sf M}}\Big( \mcF^{(k)}(\bsu(k)) \Big)   \eqdef \Psi_k\big(\bsu(k)\big),
\end{split} \end{equation*}
and set
$$
u(k) \defeq  \textbf{\textsf{R}}^{\sf M}C_{t_0}\big(\bsu(k)\big).
$$ 
By continuity of the solution map the family of functions $\big\{u(k)\big\}_{k\in G^-_\textrm{ad}}$ coincides with the limit of the family 
$$
\Big\{\textsc{Sol}\big(\xi_\epsilon;F^{(k*k_\epsilon)}\big)\Big\}_{k\in G^-_\textrm{ad}}.
$$   
Note that $\Psi_k$ depends linearly, hence smoothly, on $k$. We saw in Theorem \ref{ThmWellPosedness} in Section \ref{SectionSolvingPDEs} that given a bounded set of nonlinearities in $C^\infty$, there exists a positive time 
horizon $t_0$ such that the `integral' map $\Psi_k$ is a contraction from $\mcD^{\gamma,\eta}_{(0,t_0)}(T,\sf g)$, uniformly with respect to the nonlinearities in the given bounded set. So the continuous linear map $\big(\textrm{Id}-\partial_{\bsu}\Psi_k\big)$, from the Banach space $\mcD^{\gamma,\eta}_{(0,t_0)}(T,\sf g)$ into itself has a continuous inverse, given under the form of the classical Neumann series. The map $\big(\textrm{Id}-\partial_{\bsu}\Psi_k\big)$ is thus a continuous isomorphism of $\mcD^{\gamma,\eta}_{(0,t_0)}(T,\sf g)$ by the open mapping theorem. It is then 
a direct consequence of the implicit function theorem that the unique fixed point $\bsu(k)$ of the equation
$$
\bsu(k) = \Psi_k\big(\bsu(k)\big)
$$
is a smooth function of $k\in G^-_{\textrm{ad}}$.

\medskip

\begin{prop} \label{PropSolutionManifold}
The family $\big\{u(k)\big\}_{k\in G^-_{\textrm{\emph{ad}}}}$ forms a finite dimensional immersed submanifold of $\mcC^\eta(\bbR\times\bbR^d)$, where $\eta$ is the parameter chosen in Theorem \ref{ThmRenormPDEs}.
\end{prop}   

\medskip

\begin{Dem}
It suffices from the implicit function theorem to see that $D_k\bsu(k)$, the derivative of $\bsu(k)$ with respect to $k$, has constant rank; this follows from the linearity of the reconstruction map if we can see that $D_k\bsu(k)$ is injective. (The reconstruction map is not injective without further assumptions.) The linear map $D_k\bsu(k)$ sends the tangent space to $G^-_{\textrm{ad}}$ at the identity into $\mcD^{\gamma,\eta}_{(0,t_0)}(T,\sf g)$. But picking $h$ in that tangent space and setting $\bsv \defeq  (D_k\bsu(k))(h)\in\mcD^{\gamma,\eta}_{(0,t_0)}(T,\sf g)$, the modelled distribution $\bsv$ cannot be null unless $h=0$, since $\bsv$ is the solution to the affine equation
$$
\bsv = \mcP_{t_0}^{\sf M}
\bigg(
D_k(\mcF^{(k)})(\bsu(k))\bsv + \sum_{\bstau\in\mcB,\,|\bstau|<0} \frac{h(\bstau)}{S(\bstau)}F(\bstau^*)(\bsu(k), D\bsu(k)  )
\bigg).
$$
\end{Dem}   

\medskip

\begin{Rem}
{\small The use of the implicit function theorem actually shows that the solution $\bsu$ of the equation
\begin{equation} \label{EqGKPZ-RS} \begin{split}
\bsu &= P_\gamma u_0+\mcP_{t_0}^{\sf M}\Big(f^\star(\bsu)\Xi + g^\star\big(\bsu,  D\bsu  \big) \Big) ,
\end{split} \end{equation}
a smooth function of $f,g\in C^n$, for $n$ large enough. This gives a direct access to Taylor expansions in small noise, where $f$ is replaced by $af$, for a small positive parameter $a$, or if $f$ is the value at $a=0$ of a smooth family $f(a,\cdot)\in C^\infty$, as the solution $\bsu$ happens then to be a smooth function of the expansion parameter $a$. Elementary classical calculus is used to see that the derivatives of $\bsu$ with respect to the parameter $a$ are solutions of affine equations obtained by formal differentiation of equation \eqref{EqGKPZ-RS} with respect to the parameter. This kind of questions has a long history, under the name `stochastic Taylor expansion' in a stochastic calculus setting -- after seminal works by Azencott \cite{Azencott} and Ben Arous \cite{BenArous}, where it was used together with the stationary phase method on Wiener space to get heat kernel estimates for elliptic and sub-elliptic diffusions. Inahama \& Kawabi extended the approach to a rough paths setting in \cite{InahamaKawabi}, and Friz, Gassiat and Pigato made a first use of this type of ideas in a regularity structures setting in \cite{FGP}. The result of Proposition \ref{PropSolutionManifold} holds for all subcritical singular stochastic PDEs, with the above straightforward proof.   }
\end{Rem}

\bigskip

\section{Building regularity and renormalization structures}
\label{SectionBuildingRS}

In the end, for the above results to hold, we require from the regularity structure $\mathscr{T}$ and the renormalization structure $\mathscr{U}$ that they satisfy the different assumptions introduced along the way for different purposes. We summarize them here, with a quick description of what they are useful for.   

\ssk

{\small \begin{center}
\resizebox{\textwidth}{!}{ 
\begin{tabular}{l|c|l}\hline
{\sf Assumption} & {\sf {\sf Section}} & {\sf What it is useful for}   \\ 

\hline \hline
\textbf{\textsf{(A1-2)}} & {\sf \ref{SubsectionModelsAndCo}} & Inclusion of the polynomial structure in our regularity structures.   \\
\textbf{\textsf{(A3)}} & {\sf \ref{section products and derivatives}} & Product between $T_X$ and $T$.   \\
\textbf{\textsf{(B1)}} &  {\sf \ref{SubsectionIntegrationOperators}} & Actions of $\Delta^{(+)}$ on $\mcI_n^{(+)}$.    \\
\textbf{\textsf{(B2)}} &  {\sf \ref{SubsectionAdmissibleModels}} & Induction structure on $\Delta$ for building admissible models.   \\
\textbf{\textsf{(C)}}  &  {\sf \ref{SubsectionCompatibleStructures}} & Compatibility between the maps $\delta^{(+)}$ and $\mcI_n^{(+)}$.   \\
\textbf{\textsf{(D1-2)}}  &  
{\sf \ref{SectionMultiPreLie}, \ref{SectionCoherence}} & 
Largeness of the basis $\mcB$ of $T$ and $U$.  
\\
\textbf{\textsf{(D3)}}  &  {\sf \ref{SectionrenormalizationMultiPreLie}} & 
Compatibility between multi-pre-Lie and renormalization structures.   \\
\textbf{\textsf{(E)}}  &  {\sf \ref{SectionBHZCharacter}} & Structure assumption on $U^-$, and induction structure on $\delta$.   \\
\hline
\end{tabular}   }
\end{center}   } 

\bigskip

Following Bruned, Hairer and Zambotti \cite{BHZ}, we describe in this section a setting tailor made for the study of the generalized (KPZ) equation \eqref{EqgKPZ} where all these conditions hold true. We introduce a homogeneity map on 
the decorated trees from Section \ref{SectionMultiPreLie}.

\medskip

\begin{defn*}
Let $\frak{T}_{\text{\rm n}}$ and $\frak{T}_{\text{\rm e}}$ be abstract finite sets, equipped with homogeneity maps $|\cdot|:\frak{T}_{\text{\rm n}},\frak{T}_{\text{\rm e}}\to\bbR$. 
\begin{itemize}
\setlength{\itemsep}{0.1cm}
   \item On the sets ${\sf N}\defeq \frak{T}_{\text{\rm n}}\times\bbN\times\bbN^d$ and ${\sf E}\defeq \frak{T}_{\text{\rm e}}\times\bbN\times\bbN^d$, the homogeneity maps are extended by
$$
\begin{cases}
|k|_n \defeq  |n|+|k|_{\frak{s}}, &(n,k)\in{\sf N},   \\
|\ell|_e \defeq  |e|-|\ell|_{\frak{s}}, &(e,\ell)\in{\sf E}.
\end{cases}
$$
   
   \item The \textbf{\textsf{naive homogeneity}} of a decorated tree $\tau_{\frak{e}}^{\frak{n}}$ is defined by
$$
\vert\tau_{\frak{e}}^{\frak{n}}\vert' \defeq  \sum_{e\in E_\tau} \vert\frak{e}(e)\vert_{\frak{t}_{\text{\rm e}}(e)} + \sum_{n\in N_\tau} \vert\frak{n}(n)\vert_{\frak{t}_{\text{\rm n}}(n)}.
$$
\end{itemize}
\end{defn*}   

\medskip

We start from the sets
$$
\frak{T}_{\text{n}}=\{\bullet,\circ\},\quad\frak{T}_{\text{e}}=\{\mcI\},
$$
for an abstract symbol $\mcI$ -- we use on purpose the same symbol as the 
abstract integration map from Section \ref{SubsectionIntegrationOperators}. The node type set $\frak{T}_{\text{n}}$ is enlarged later. The two elements $\bullet$ and $\circ$ of $\frak{T}_{\text{n}}$ represent the monomial ${\bf1}=X^0$ and the noise $\Xi$, respectively. The set $\frak{T}_{\rm e}$ consists of only one integration operator $\mcI$. Each element has homogeneity
$$
\vert\bullet\vert=0,\qquad \vert\circ\vert=\beta_0,\qquad \vert\mcI\vert=2,
$$
where $\beta_0\in(-2,0)$ is the regularity of the noise $\zeta$ in the equation. (Would the equation under study involve several noises with different regularities we would introduce several $\circ$ symbols with the corresponding homogeneities.) Given that the polynomial structure is needed to encode at a regularity structure level the term $P_\gamma u_0$ describing the propagation of the initial condition, and the piece of $\mcP_t^{\sf M}$ taking values in the polynomial regularity structure, the use of trees with a node decoration encoding multiplication by polynomials appears as natural. On the other hand, the use of edge decorations for equations that do not involve derivatives of the solution in their formulation, like the generalized (PAM) equation
$$
(\partial_t-\Delta_x) u = f(u)\zeta,
$$
may look strange. The necessity to use edge decorations to encode derivatives of quantities of the form $\mcI(\cdot)$, even in such a case, comes from the renormalization process implemented in this setting, as the latter involves Taylor expansions.

As said in Section \ref{SectionMultiAndrenormalizedEquations}, the final form of a generic element of our regularity structures will be the datum of a decorated tree together with a coloring and an additional decoration 
$\frak{o}:N_\tau\to\bbZ[\beta_0]$, which plays an important role in the compatibility condition between regularity and renormalization structures from Definition \ref{*DefnCompatibility}. In a nutshell, this additional decoration will keep track of the naive homogeneity of the `diverging' trees that will be extracted by the renormalization map $\delta$. This is what will allow to have maps $\delta^{(+)}$ satisfying the fundamental conditions 
$$
\delta^{(+)} T^{(+)}_\beta \subset U^-\otimes T^{(+)}_\beta \qquad (\beta\in A^{(+)})
$$
involved in the definition of compatible regularity and renormalization structures. So one should not be surprised that we will use the naive homogeneity to define the gradings in $U$ and $U^-$ and a different notion of 
homogeneity in $T$ and $T^+$, taking into account the $\frak{o}$-decorations. The discussion will be general enough for the reader to see what needs to be added to deal with the general case.

\bigskip

\subsection{Rules and extended decoration}
\label{SectionDecoratedTrees}

Working with the set of all decorated trees as a candidate for a regularity structure is not reasonable and we first identify a few notions that help clarifying the matter. Recall the abstract self-explaining formulation 
\begin{equation} \label{EqRSgKPZModel}
\bsv = \mcI\Big(f^\star\big(\bsv\big)\Xi + g^\star\big(\bsv,\partial \bsv\big)\Big) + (T_X)
\end{equation}
of the generalized (KPZ) equation. In the present tree setting the $\star$ product is given by the `{\it joining}' operator $\mathscr{J}$ on trees. If one wants to make sense of Picard  iteration within the concrete regularity structure, one needs to make sense of a number of recursive relations -- recall the subcomodules introduced in Section \ref{SubsectionFixedPoint} and see the pictures in Section \ref{SectionMultiPreLie}. General constraints of this type come under the name of \textbf{\textsf{rule}}, that is the definition for each node type $\frak{b}\in\frak{T}_{\rm n}$, of constraints on which kind of tuples of edges $\{e^i=(e^i_-,e^i_+)\}_i$ can have $\frak{t}_{\text{n}}(e^i_-)=\frak{b}$, for all $i$, in a tree allowed by the rule. The choice of a rule is determined by the equation under consideration. Consider the right hand side of equation \eqref{EqRSgKPZModel}. Making sense of the nonlinear term $f^\star(\bsv)\Xi+g_0^\star(\bsv)$ requires that one can find $\mathscr{J}\big(\mcI(\cdot),\dots,\mcI(\cdot)\big)X^n$ or $\mathscr{J}\big(\mcI(\cdot),\dots,\mcI(\cdot)\big)X^n\Xi$, within the trees allowed by the rule, that is the corresponding nodes are of the form
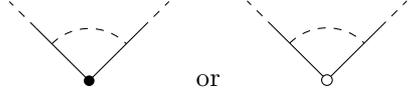
\begin{align*}
\begin{tikzpicture}
\coordinate (A1) at (0,0);
\coordinate (A2) at (-0.7,0.7);
\coordinate (A3) at (0.7,0.7);
\foreach \n in {2,3} \coordinate (B\n) at ($0.7*(A\n)$);
\foreach \n in {2,3} \coordinate (C\n) at ($1.5*(A\n)$);
\foreach \n in {2,3} \draw (A1)--(A\n);
\foreach \n in {2,3} \draw[dashed] (A\n)--(C\n);
\draw[dashed] (B2) to [out=45, in=135] (B3);
\foreach \n in {1} \fill (A\n) circle (2pt);
\end{tikzpicture}
\text{$\quad$or$\quad$}
\begin{tikzpicture}
\coordinate (A1) at (0,0);
\coordinate (A2) at (-0.7,0.7);
\coordinate (A3) at (0.7,0.7);
\foreach \n in {2,3} \coordinate (B\n) at ($0.7*(A\n)$);
\foreach \n in {2,3} \coordinate (C\n) at ($1.5*(A\n)$);
\foreach \n in {2,3} \draw (A1)--(A\n);
\foreach \n in {2,3} \draw[dashed] (A\n)--(C\n);
\draw[dashed] (B2) to [out=45, in=135] (B3);
\filldraw[white] (A1) circle (2pt);
\draw (A1) circle (2pt);
\end{tikzpicture}.
\end{align*}
Making sense of the other terms $g_2^\star(\bsv)\star( D\bsv )^{\star2}+g_1^\star(\bsv)\star( D\bsv )$ requires that one can find \\
$\mathscr{J}\big(\mcI(\cdot),\dots,\mcI(\cdot),\mcI_{e_i}(\cdot)\big)X^n$, or $\mathscr{J}\big(\mcI(\cdot),\dots,\mcI(\cdot),\mcI_{e_i}(\cdot),\mcI_{e_j}(\cdot)\big)X^n$ for some $i,j=1,\dots,d$ within the trees allowed by the rule, so each node of the corresponding elements of $C$ has the 
form
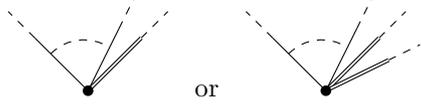
\begin{align*}
\begin{tikzpicture}
\coordinate (A1) at (0,0);
\coordinate (A2) at (-0.7,0.7);
\coordinate (A3) at (0.7,0.7);
\coordinate (A4) at (0.407,0.825);
\coordinate (B1) at ($0.7*(A2)$);
\coordinate (B2) at ($0.7*(A4)$);
\foreach \n in {2,3,4} \coordinate (C\n) at ($1.5*(A\n)$);
\foreach \n in {2,4} \draw (A1)--(A\n);
\foreach \n in {2,3,4} \draw[dashed] (A\n)--(C\n);
\draw[double distance=0.7pt] (A1)--(A3);
\draw[dashed] (B1) to [out=45, in=150] (B2);
\foreach \n in {1} \fill (A\n) circle (2pt);
\end{tikzpicture}
\text{$\quad$or$\quad$}
\begin{tikzpicture}
\coordinate (A1) at (0,0);
\coordinate (A2) at (-0.7,0.7);
\coordinate (A3) at (0.7,0.7);
\coordinate (A4) at (0.407,0.825);  
\coordinate (A5) at (0.825,0.407);
\coordinate (B1) at ($0.7*(A2)$);
\coordinate (B2) at ($0.7*(A4)$);
\foreach \n in {2,3,4,5} \coordinate (C\n) at ($1.5*(A\n)$);
\foreach \n in {2,4} \draw (A1)--(A\n);
\foreach \n in {2,3,4,5} \draw[dashed] (A\n)--(C\n);
\foreach \n in {3,5} \draw[double distance=0.7pt] (A1)--(A\n);
\draw[dashed] (B1) to [out=45, in=150] (B2);
\foreach \n in {1} \fill (A\n) circle (2pt);
\end{tikzpicture}.
\end{align*}
The operators $\mcI_{e_i}$ are represented by the double line in the above picture. Given a rule, a decorated \textsf{\textbf{conforming tree}} is a tree such that all nodes of the tree, except perhaps the root, satisfy the rule. Denote by 
$$
C
$$
the set of conforming trees. If all node of the tree satisfy the rule, the tree is called \textbf{\textsf{strongly conforming}}. We denote by
$$
SC
$$
the set of strongly conforming trees. A rule is said to be \textbf{\textsf{normal}} if any subtree of a strongly conforming tree is also strongly conforming.

\ssk

To construct regularity and renormalization structures, the rooted decorated trees obtained from the above iterations are not sufficient. Another important operation is the \textbf{\textsf{contraction}} of rooted trees, 
involved in the definition of the splitting maps $\Delta$ and $\delta$. Given a typed rooted tree $\tau$ and a family $\varphi$ of disjoint typed subtrees of $\tau$, we use the notation
$$
\tau/^{\text{red}}\varphi
$$
to denote the typed rooted tree obtained by identifying each subtree $\tau_i$ with a single node $\textcolor{red}{\bullet}$ with red color in the quotient tree. Here is an example, with $\varphi$ in {\color{green}{green}},
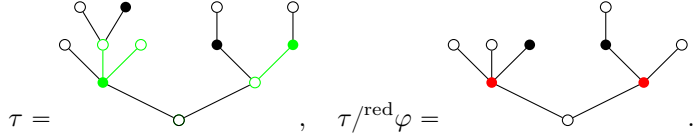
$$
\tau=
\begin{tikzpicture}
\coordinate (A1) at (0,0);
\coordinate (A2) at (-1,0.5);
\coordinate (A3) at (1,0.5);
\coordinate (A4) at (-1.5,1);
\coordinate (A5) at (-1,1);
\coordinate (A6) at (-0.5,1);
\coordinate (A7) at (0.5,1);
\coordinate (A8) at (1.5,1);
\coordinate (A9) at (-1.3,1.5);
\coordinate (A10) at (-0.7,1.5);
\coordinate (A11) at (0.5,1.5);
\coordinate (A12) at (1.5,1.5);
\draw (A2)--(A1)--(A3);
\draw[green] (A5)--(A2)--(A6);
\draw (A4)--(A2);
\draw[green] (A8)--(A3);
\draw (A7)--(A3);
\foreach \n in {9,10} \draw (A\n)--(A5);
\draw (A11)--(A7);
\draw (A12)--(A8);
\foreach \n in {7,10} \fill (A\n) circle (2pt);
\foreach \n in {2,8} \fill[green] (A\n) circle (2pt);
\foreach \n in {1,3,4,5,6,9,11,12} \filldraw[white] (A\n) circle (2pt);
\foreach \n in {1,3,5,6} \draw[green] (A\n) circle (2pt);
\foreach \n in {1,4,9,11,12} \draw (A\n) circle (2pt);
\end{tikzpicture},
\quad
\tau/^{\text{red}}\varphi=
\begin{tikzpicture}
\coordinate (A1) at (0,0);
\coordinate (A2) at (-1,0.5);
\coordinate (A3) at (1,0.5);
\coordinate (A4) at (-1.5,1);
\coordinate (A5) at (-1,1);
\coordinate (A6) at (-0.5,1);
\coordinate (A7) at (0.5,1);
\coordinate (A8) at (1.5,1);
\coordinate (A9) at (0.5,1.5);
\draw (A2)--(A1)--(A3);
\foreach \n in {4,5,6} \draw (A\n)--(A2);
\foreach \n in {7,8} \draw (A\n)--(A3);
\draw (A9)--(A7);
\foreach \n in {6,7} \fill (A\n) circle (2pt);
\foreach \n in {2,3} \fill[red] (A\n) circle (2pt);
\foreach \n in {1,4,5,8,9} \filldraw[white] (A\n) circle (2pt);
\foreach \n in {1,4,5,8,9} \draw (A\n) circle (2pt);
\end{tikzpicture}.
$$
We allow such an operation for the set $SC$ of strongly conforming trees. Precisely, if each connected component of $\varphi$ belongs to $SC$, then we assume that $\tau/^{\text{red}}\varphi\in SC$. Hence each element of $SC$ is a rooted decorated tree with a node type set
$$
\frak{T}_{\text{n}}^{SC}=\{\bullet,\textcolor{red}{\bullet},\circ\}.
$$
The analytic role of $\textcolor{red}{\bullet}$ is the same as that of $\bullet$. In particular, the homogeneity of $\textcolor{red}{\bullet}$ is $0$. This is an example of the coloring of the tree. Only decorated trees 
without red color appear in the analysis of the well-posedness problem \eqref{EqRSgKPZModel}, but colors are used in the definition of the splitting maps in the renormalization structure.

\ssk

Recall from Assumption \textbf{\textsf{(B1)}} and Section \ref{SubsectionAdmissibleModels} that the algebra $T^+$ is spanned by elements of the form
\begin{align}\label{Section 9: generic element of T^+}
X^n\prod_{i=1}^N\mcI_{k_i}^+(\tau_i),
\end{align}
where $n\in\bbN\times\bbN^d,\,k_1,\dots,k_N\in\bbN\times\bbN^d$, and $\tau_1,\dots,\tau_N\in SC$. It is convenient to consider an element like \eqref{Section 9: generic element of T^+} as a tree by interpreting $\mcI_k^+$ as the planting operator like $\mcI_k$ and the product $\prod$ as the tree product 
$\mathscr{J}$. To distinguish such trees from elements of $SC$, we give a 
blue color to their roots, encoding in this way the $+$ sign in $\mcI^+$.
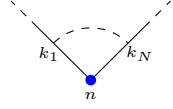
\begin{align}\label{picture generic element of C+}
\begin{tikzpicture}
\coordinate (A1) at (0,0);
\coordinate (A2) at (-0.7,0.7);
\coordinate (A3) at (0.7,0.7);
\foreach \n in {2,3} \coordinate (B\n) at ($0.7*(A\n)$);
\foreach \n in {2,3} \coordinate (C\n) at ($1.5*(A\n)$);
\foreach \n in {2,3} \draw (A1)--(A\n);
\foreach \n in {2,3} \draw[dashed] (A\n)--(C\n);
\draw[dashed] (B2) to [out=45, in=135] (B3);
\fill[blue] (A1) circle (2pt);
\node at ($0.5*(A1)+0.5*(A2)+(-0.2,0)$) {\tiny $k_1$};
\node at ($0.5*(A1)+0.5*(A3)+(0.3,0)$) {\tiny $k_N$};
\node at ($(A1)+(0,-0.2)$) {\tiny $n$};
\end{tikzpicture}
\end{align}
The set $C$ consists of such trees, where we see that the rule is broken at the root. This is because $C$ is only conforming, not strongly conforming. Each element of $C$ is thus a rooted decorated tree with a node type set
$$
\frak{T}_{\text{n}}^C=\{\bullet,\textcolor{red}{\bullet},\textcolor{blue}{\bullet},\circ\}.
$$
A node of a conforming tree has the type $\textcolor{blue}{\bullet}$ if and only if it is a root. The homogeneity of $\textcolor{blue}{\bullet}$ is $0$. The trees with a blue root will only be involved in the description of the space $T^+$.

\ssk

A rule is said to be \textbf{\textsf{subcritical}} if for any $\gamma\in\bbR$, only finitely many elements of $SC$ have naive homogeneity less than $\gamma$.  A \textbf{\textsf{complete}} rule will guarantee that a rooted decorated tree obtained from the contraction of a strongly conforming tree by extracting `diverging' pieces, and changing the decorations accordingly, will still be strongly conforming. Proposition 5.21 in \cite{BHZ} ensures that {\it any normal subcritical rule can be extended into a normal subcritical complete rule}. We take this result for granted and do not reprove it here. The above rule on the set of decorated trees is normal, subcritical and complete.

\ssk

To construct compatible regularity and a renormalization structures we introduce an additional decoration. Denote by $N_\tau^{\text{\rm red}}$ the 
subset of $N_\tau$ consisting of the nodes with type $\textcolor{red}{\bullet}$.

\ssk

\begin{defn*}
A tree with \textbf{\textsf{extended decoration}} is a rooted decorated tree $\tau_{\frak{e}}^{\frak{n}}$ with a map
$$
\frak{o}:N^{{\rm red}}_\tau\to\bbZ[\beta_0].
$$
We write $\bstau=\tau_{\frak{e}}^{\frak{n},\frak{o}}$ for a generic tree with extended decoration.
The \textbf{\textsf{extended homogeneity}} of such a tree is defined by
$$
\big|\tau_{\frak{e}}^{\frak{n},\frak{o}}\big| \defeq  \sum_{e\in E_\tau} \big|\frak{e}(e)\big|_{\frak{t}_{\rm e}(e)} + \sum_{n\in N_\tau} \big|\frak{n}(n)\big|_{\frak{t}_{\rm n}(n)}+ \sum_{n\in N_\tau^{\rm red}}\frak{o}(n).
$$
We extend the naive homogeneity to the set of decorated trees with an extended decoration setting
$$
\big|\tau_{\frak{e}}^{\frak{n},\frak{o}}\big|' \defeq  \big|\tau_{\frak{e}}^{\frak{n}}\big|'.
$$   
\end{defn*}

\ssk

Note here that only trees without $\frak{o}$-decoration, that is $\frak{o}=0$, appear in the analysis of the fixed point problem \eqref{EqRSgKPZModel}. Indeed, the trees without $\frak{o}$-decoration are stable under the coproducts $(\Delta^+,\Delta)$ defined below.
The $\frak{o}$-decoration is only involved in the analysis of the renormalization procedure and the associated convergence problem -- see the second equality of \eqref{EqStabilityPrimed}. Without it, the condition \eqref{delta T to U-T} for the compatibility of $\mathscr{T}$ and $\mathscr{U}$ does not hold. We define the set
$$
\mathcal{SC}
$$
of strongly conforming trees with $\frak{o}$-decoration as the minimal set which contains $SC$ and 
such that the vector space spanned by $\mathcal{SC}$
is stable under all the coproducts defined below. (One could also consider $\mathcal{SC}$ as a set of rooted decorated trees with node type set
$$
\frak{T}_{\rm n}^{\mcS\mcC}=\{\bullet,\circ\}\cup\{\textcolor{red}{\bullet}^{\cdot,\alpha}\}_{\alpha\in\bbZ[\beta_0]}.
$$
We used such an identification in Section \ref{SectionMultiAndrenormalizedEquations}. In the present section we treat $\frak{o}$ as a decoration, rather than as part of a node type.) Similarly, we define
$$
\mcC
$$
as the set of decorated trees with extended decorations of the form \eqref{Section 9: generic element of T^+}, where $\tau_1,\dots,\tau_N\in\mcS\mcC$. We use the bold symbol $\bstau$ to denote a generic element of $\mathcal{SC}$ or $\mcC$. The above rule on the set of extended decorated trees is normal, subcritical and complete. (The subcriticality of the rule on 
this set of trees with extended decorations comes from the fact that, for 
any fixed $\gamma\in\bbR$, the decoration $\alpha$ of trees with extended 
homogeneity less than $\gamma$, will only range in the set of homogeneities of subtrees of strongly conforming trees $\tau_{\frak{e}}^{\frak{n}}$ with homogeneity less than $\gamma$.)

\bigskip

\subsection{Coproducts}
\label{SectionDecoratedTreesCoproducts}

We define coproducts in the spaces of rooted decorated trees. This requires first that we define what we mean by `subtrees' and `subforests'. Recall that the type sets
$$
\frak{T}_{\rm n}^{SC}=\{\bullet,\textcolor{red}{\bullet},\circ\},\qquad
\frak{T}_{\rm n}^{C}=\big\{\bullet,\textcolor{red}{\bullet},\textcolor{blue}{\bullet},\circ\big\},\qquad
\frak{T}_{\rm e}=\{\mcI\}
$$
are fixed. Given a typed rooted tree $\tau$, a nonempty connected subgraph of $\tau$ is called a \textbf{\textsf{subtree}} if it inherits from $\tau$ its type map. Any possibly empty family of disjoint subtrees of $\tau$ is called a \textbf{\textsf{subforest}}. Given a rooted tree $\tau$ and 
a subforest $\varphi=\big\{\tau_1,\dots,\tau_m\big\}$, we use the notation
$$
\tau/\varphi
$$
to denote the rooted tree obtained by identifying each subtree $\tau_i$ with a single node with node type $\bullet$ in the quotient tree. Precisely, writing $y\sim_\varphi z$ if $y$ and $z$ are in the same connected component of $\varphi$, we define $\tau/\varphi$ as the tree consisting of the node set $N_\tau/\sim_\varphi$ and the edge set $E_\tau\setminus E_\varphi$.
Moreover, we write
$$
\tau/^{\text{red}}\varphi,\quad\text{or}\quad \tau/^{\text{blue}}\varphi
$$
if we give a corresponding color to the nodes of $\varphi$ in the quotient tree.

\begin{itemize}
\setlength{\itemsep}{0.1cm}
\item For any function $f:N_\tau\to\bbN\times\bbN^d$, define the function 
$[f]_\varphi$ on $N_{\tau/\varphi}$ by
$$
[f]_\varphi\big([x]\big) \defeq  \sum_{y\sim_\varphi x} f(y),
$$
where $[x]$ denotes the equivalence class of $x\in N_\tau$.

\item Denote by $\partial\varphi$ the leaves of $\varphi$, that is, the set of edges $(x,y)\in E_\tau$ such that $x\in N_\varphi$ and $y\in N_\tau\setminus N_\varphi$.
For any function $g:\partial\varphi\to\bbN\times\bbN^d$, define the function $\pi g$ on $N_\tau$ by setting
$$
(\pi g)(x)\defeq \sum_{e=(x,y)\in\partial\varphi}g(e).
$$

\item For any decorations $\frak{n}_\varphi$ and $\frak{e}_\varphi$ on $\varphi$, define the function $\frak{o}(\varphi,\frak{n}_\varphi,\frak{e}_\varphi):N_{\tau/\varphi}\to\bbZ[\beta_0]$ by
$$
\frak{o}\big(\varphi,\frak{n}_\varphi,\frak{e}_\varphi\big)\big([\tau_j]\big)
=\left\vert(\tau_j)_{\frak{e}_\varphi\vert_{\tau_j}}^{\frak{n}_\varphi\vert_{\tau_j}}\right\vert'
$$
for each $1\leq j\leq m$, and $\frak{o}(\varphi,\frak{n}_\varphi,\frak{e}_\varphi)=0$ outside $[\varphi]$.
\end{itemize}

\medskip

Define
$$
T \defeq  \text{\rm span}(\mathcal{SC}),   \qquad
{\sf T}^+ \defeq  \text{\rm span}(\mathcal{C}),   \qquad
{\sf U}^- \defeq  \bbR[\mathcal{SC}].
$$
Note that ${\sf T}^+$ is an algebra with the tree product and unit ${\bf1}_+ \defeq  \textcolor{blue}{\bullet}^0$, and ${\sf U}^-$ is an algebra with 
the forest product and unit ${\bf1}_-\defeq \emptyset$. The space $T^+$ will 
be built from the side space ${\sf T}^+$ and the space $U^-$ from the side space ${\sf U}^-$. Similarly, the different splitting maps defining a regularity structure and a renormalization structure are built from splitting maps taking values in, or defined on, the spaces $T,{\sf T}^+, {\sf U}^-$.

\medskip

\begin{defn*} 
We introduce three splitting operators.
\begin{itemize}
\setlength{\itemsep}{0.1cm}

	\item[{\sf 1.}] The linear map 
	$$
	{\sf D}: T\to T\otimes {\sf T}^+
	$$ 
	is defined for $\tau_{\frak{e}}^{\frak{n},\frak{o}}\in\mcS\mcC$ by
\begin{equation} \label{EqFormulaDPlus}
\begin{split}
{\sf D}\tau_{\frak{e}}^{\frak{n},\frak{o}} 
\defeq & \sum_{\mu} \sum_{\frak{n}_\mu, \frak{e}'_{\partial\mu}} 
\frac{1}{\frak{e}'_{\partial\mu}!}{{\frak{n}}\choose{\frak{n}_\mu}}
\mu^{\frak{n}_\mu+\pi\frak{e}'_{\partial\mu},\frak{o}\vert_\mu}_{\frak{e}}
\otimes 
\big(\tau/^{\text{\rm blue}}\mu\big)^{[\frak{n}-\frak{n}_\mu]_\mu,\frak{o}\vert_{\tau\setminus\mu}}_{\frak{e}+\frak{e}'_{\partial\mu}},
\end{split}  \end{equation}
where the first sum is over all subtrees $\mu$ of $\tau$ which contains the root of $\tau$,
and the second sum is over functions $\frak{n}:N_\mu\to\bbN\times\bbN^d$, 
with $\frak{n}_\mu\leq \frak{n}$ and functions $\frak{e}'_{\partial\mu}: \partial\mu\rightarrow\bbN\times\bbN^d$.
The algebra morphism 
$$
{\sf D}^{\,+} : {\sf T}^+ \to {\sf T}^+\otimes {\sf T}^+,
$$ 
is defined by the same formula \eqref{EqFormulaDPlus} for $\tau_{\frak{e}}^{\frak{n},\frak{o}}\in{\mcC}$.   \vspace{0.15cm}

	\item[{\sf 2.}] The algebra morphism 
	$$
	{\sf D}^- : {\sf U}^-\to{\sf U}^-\otimes {\sf U}^-
	$$ 
	is defined by ${\sf D}^-({\bf1}_-) \defeq  {\bf1}_-\otimes{\bf1}_-$, and for $\tau_{\frak{e}}^{\frak{n},\frak{o}}\in\mcS\mcC$
\begin{equation*} \begin{split}
{\sf D}^-\big(\tau_{\frak{e}}^{\frak{n},\frak{o}} \big) 
&\defeq  \sum_{\varphi}\sum_{\frak{n}_\varphi, \frak{e}'_{\partial\varphi}} \frac{1}{\frak{e}'_{\partial\varphi}!}{{\frak{n}}\choose{\frak{n}_\varphi}}
\varphi^{\frak{n}_\varphi+\pi\frak{e}'_{\partial\varphi},\frak{o}\vert_\varphi}_{\frak{e}} 
\otimes \big(\tau/^{\text{\rm red}}\varphi\big)^{[\frak{n}-\frak{n_\varphi}]_\varphi,[\frak{o}]_\varphi+\frak{o}(\varphi,\frak{n}_\varphi+\pi\frak{e}_{\partial\varphi}',\frak{e})}_{\frak{e}+\frak{e}'_{\partial\varphi}},
\end{split}  \end{equation*}
where the first sum is over all subforests $\varphi$ of $\tau$ which contains all red nodes of $\tau$, and the sum over $\frak{n}_\varphi$ and $\frak{e}_{\partial\varphi}'$ is taken as in item {\sf 1} of the present definition.   \vspace{0.15cm}

	\item[{\sf 3.}] The algebra morphism 
	$$
	\overline{\sf D}^{\,-} : {\sf T}^+ \to {\sf U}^-\otimes{\sf T}^+
	$$ 
	is defined by the same formula as ${\sf D}^-$, with the first sum restricted to subforests $\varphi$ which are disjoint from the root of $\tau$.
\end{itemize}
\end{defn*}

\medskip

{\sl  \begin{Rem}
Since the right hand side of \eqref{EqFormulaDPlus} may become an infinite series, we have to consider the `bigraded spaces' of rooted decorated trees defined in Section 2.3 of \cite{BHZ}. For any collection of vector spaces $\{V_n\}_{n\in\bbN^2}$, we denote by $V=\bigboxplus_{n\in\bbN^2}V_n$ the space of all sequences $(v_n)_{n=(n_1,n_2)\in\bbN^2}$ with $v_n\in V_n$ such that there exists $k\in\bbN$ and $v_n=0$ unless $n_2\le k$. The tensor product of two bigraded spaces $V=\bigboxplus_{n\in\bbN^2}V_n$ and $W=\bigboxplus_{n\in\bbN^2}W_n$ is defined by
$$
V\mathbin{\hat{\otimes}}W\defeq\bigboxplus_{n\in\bbN^2}\bigg(\bigoplus_{k+\ell=n}V_k\otimes W_\ell\bigg).
$$
For example, the bigraded space ${\sf V}=\bigboxplus_{n\in\bbN^2}{\sf V}_n$ of rooted decorated trees is given by setting ${\sf V}_{(n_1,n_2)}$ as the vector space spanned by all decorated trees $\tau_{\frak{e}}^{\frak{n},\frak{o}}$ such that $\sum_{e\in E_\tau}|\frak{e}(e)|_\mfs=n_1$ and $|N_\tau\setminus\frak{t}^{-1}\{\textcolor{red}{\bullet},\textcolor{blue}{\bullet}\}|=n_2$. The spaces $T$, ${\sf T}^+$, and ${\sf U}^-$ above are defined as sub-bigraded spaces of ${\sf V}$, and their tensor products are also defined as bigraded spaces in the above sense. As in Lemma 2.14 of \cite{BHZ}, triangular maps between bigraded spaces are well-defined. For any bigraded spaces $V$ and $W$, the family of linear maps $A_{mn}:V_n\to W_m$ is called {\em triangular} if $A_{mn}=0$ unless $m_1\ge n_1$ and $m_2\le n_2$. The the linear map from $V$ to $W$
$$
A\big((v_n)_{n\in\bbN^2}\big)\defeq\bigg(\sum_{n\in\bbN^2}A_{mn}v_n\bigg)_{m\in\bbN^2}
$$
is well-defined. The maps ${\sf D}$, ${\sf D}^+$, ${\sf D}^-$, and $\overline{\sf D}^{\,-}$ above are well-defined as triangular maps. In the following, when dealing with infinite series, we use these facts implicitly. In Section \ref{SectionDecoratedTreesRegRenor}, we introduce some truncation maps which reduce infinite series to finite sums.
\end{Rem} }

\medskip

As suggested by the target spaces of the preceding maps, the splitting map $\Delta$ will be constructed from $\sf D$ and the map $\Delta^+$ from ${\sf D}^+$, the maps $\delta$ and $\delta^-$ from ${\sf D}^-$, and the map $\delta^+$ from $\overline{\sf D}^-$. Only trees with blue roots appear 
in the right hand side of the tensor products defining $\sf D$. This is consistent with the fact that the trees with blue roots will represent later elements of $T^+$. The restriction on the choice of $\varphi$ to subforests which are disjoint from the root in the definition of $\overline{\sf D}^-$ ensures that it takes values in ${\sf U}^-\otimes{\sf T}^+$ and that the multiplicative property
$$
\overline{\sf D}^{\,-}(\bstau\bssigma) = \big(\overline{\sf D}^{\,-}\bstau\big)\big(\overline{\sf D}^{\,-}\bssigma\big)
$$ 
holds. This reflects the fact that the product of two functions
$$
{\sf g}(\bstau)\,{\sf g}(\bssigma),\quad \bstau,\bssigma\in \mathcal{C}
$$
does not cause any new renormalization.

\medskip

{\sl  \begin{Rem}
Keep in mind that the elements of $U^-$ are meant to be evaluated by characters of $U^-$, and turned to numbers, while elements of $U$ are meant to be turned to distributions. This is done jointly in a renormalized naive model $(k\otimes{\sf\Pi}^\zeta)\delta$. Recall that the problem of renormalization comes from the fact that the kernel of the operator $\bf K$ explodes on the diagonal. The building block of the renormalization operations $\delta$ and $\delta^-$ is best understood in the light of the following archetype problem. Let $g : ([0,1]^d)^n\rightarrow\bbR$ be a function that is smooth outside the deep diagonal $\textsf{\emph{diag}} \defeq  \big\{{\bf z}=(z_1,\dots,z_n)\in([0,1]^d)^n; z_1=\cdots=z_n\big\}$, near which it behaves as $\vert {\bf z} - (z_1,\dots,z_1)\vert^{-a}$, for an exponent $a>d$. The function $g$ is not integrable in any neighbourhood 
of the deep diagonal, so it only makes sense as a distribution on $([0,1]^d)^n\backslash\textsf{\emph{diag}}$
\begin{equation*}
\int_{([0,1]^d)^n} g({\bf z})f({\bf z})d{\bf z},
\end{equation*}
for $f$ smooth, with support with empty intersection with the deep diagonal. Can we define a distribution $\Lambda$ on $([0,1]^d)^n$ that extends this distribution? This can be done defining $\Lambda$ on $([0,1]^d)^n$

\begin{equation*} 
(\Lambda, \psi) = \int_{([0,1]^d)^n} g({\bf z})\left(\psi({\bf z}) - \psi({\bf z}_1) -\cdots-\frac{({\bf z}-{\bf z}_1)^{[a-d]}}{[a-d]!}\,\psi^{([a-d])}({\bf z}_1)\right)\,d{\bf z},
\end{equation*}
for any smooth function $\psi$ on $([0,1]^d)^n$. This formula defines indeed a distribution, which coincides with the distribution associated with 
$g$ outside the deep diagonal, since the $\psi^{[\ell]}({\bf z}_1)$ are null for functions with compact support with null intersection with $\textsf{\emph{diag}}$. Taylor expansion appears as the building block of this extension procedure. In this parallel, $\tau$ has two pieces, $g$ and $f$, so the role of $\varphi$ in ${\sf D}^-$ would be played by either of them, and the role of the projector $p_-$ in $\delta^-$, defined below, would select only the diverging term. The term $\varphi\otimes(\tau/^-\varphi)$ in $\delta^-\tau$ would precisely correspond to a term $g({\bf z})\,\frac{({\bf z}-{\bf z}_1)^n}{n!}\,f^{(n)}({\bf z}_1)$ in the integral defining $\Lambda$. A formula like the above defining relation for ${\sf D}^-$ appears if one deals with a multiple integral where several subintegrals define functions of their external variables of the same kind as $g$, 
and one uses a similar kind of extension procedure as above.   
\end{Rem} }

\ssk

The following lemma is proved in Appendix {\sf \ref{Appendix Proof of coassociativity}}.

\medskip

\begin{lem}\label{coassosictive cointeraction Dpm}
One has the coassociativity formulas
$$
\begin{aligned}
({\sf D}\otimes\text{\rm Id}){\sf D} &= (\text{\rm Id}\otimes {\sf D}^+){\sf D},   
&({\sf D}^+\otimes\text{\rm Id}){\sf D}^+ &= (\text{\rm Id}\otimes {\sf D}^+){\sf D}^+,   \\
({\sf D}^-\otimes\text{\rm Id}){\sf D}^- &= \big(\text{\rm Id}\otimes{\sf D}^-\big){\sf D}^-,   \quad
&({\sf D}^-\otimes\text{\rm Id})\overline{\sf D}^{\,-} &= \big(\text{\rm 
Id}\otimes\overline{\sf D}^{\,-}\big)\overline{\sf D}^{\,-}.
\end{aligned}
$$
Moreover, one has the cointeraction formulas
\begin{align*}
\mcM^{(13)}\big({\sf D}^-\otimes\overline{\sf D}^{\,-}\big){\sf D} &= (\text{\rm Id}\otimes{\sf D}){\sf D}^-,   \\
\mcM^{(13)}\big(\overline{\sf D}^{\,-}\otimes\overline{\sf D}^{\,-}\big){\sf D}^+ &= (\text{\rm Id}\otimes{\sf D}^+)\overline{\sf D}^{\,-}.
\end{align*}
\end{lem}

\bigskip

\subsection{Regularity and renormalization structures}
\label{SectionDecoratedTreesRegRenor}

We define the Hopf algebra parts of regularity and renormalization structures, from the side spaces ${\sf T}^+$ and ${\sf U}^-$. We use the shorthand notation $\mcI^+_n(\tau)$ to denote the tree $\mcI_n(\tau)$ with a blue root, with $\mcI_n(\tau)$ standing for $\mcI(\tau)$ with decoration $n$ on the edge outgoing from the root. We define subsets $\mathcal{C}^+\subset \mathcal{C}$ and $\mathcal{SC}^-\subset \mathcal{SC}$, by
\begin{align*}
\mathcal{C}^+ &\defeq  \left\{X^n\mathscr{J}\big(\mcI_{m_1}^+(\bstau_1),\dots,\mcI_{m_b}^+(\bstau_b)\big)\in\mathcal{C} \,;\, \text{$\big\vert\mcI_{m_j}^+(\bstau_j)\big\vert > 0$, for any $j=1,\dots,b$}\right\},\\
\mathcal{SC}^-& \defeq  
\Big\{\bstau\in \mathcal{SC}\,;\, \vert\bstau\vert' < 0\Big\}
\end{align*}
and set
$$
T=U = \;\textrm{span}(\mcS\mcC),   \qquad
T^+ \defeq \;\text{\rm span}(\mathcal{C}^+),   \qquad
U^- \defeq \;\bbR[\mathcal{SC}^-].
$$
Note the use of the two notions of homogeneity in these definitions, the extended homogeneity $\vert\cdot\vert$ for $T$ and $T^+$, and the naive homogeneity $\vert\cdot\vert'$ for $U$ and $U^-$. Denote by 
$$
p_+ : {\sf T}^+\to T^+
$$
the canonical projection, and define an algebra morphism 
$$
p_- : {\sf U}^-\to U^-
$$ 
setting
\begin{align*}
p_-(\bstau) \defeq  
\begin{cases}
{\bf1}_-,   &\textrm{for}\;\bstau={\bf1}_-, \textcolor{red}{\bullet}^{0,\alpha},\\
\bstau,   &\textrm{for}\;\bstau\in\mathcal{SC}^-,\\
0,   &\textrm{for}\;\bstau\in\mcS\mathcal{C}\setminus\Big\{\mathcal{SC}^-\cup\{\textcolor{red}{\bullet}^{0,\alpha}\}_{\alpha\in\bbZ[\beta_0]}\Big\}.
\end{cases}
\end{align*}

\medskip

\begin{defn*}
Define the linear maps
$$
\begin{aligned}
\Delta&\defeq  \big(\text{\rm Id} \otimes p_+\big){\sf D} : T\to T\otimes T^+,
&\Delta^+&\defeq  \big(p_+ \otimes p_+\big){\sf D}^+\vert_{T^+} : T^+ \to T^+ \otimes T^+,\\
\delta &\defeq  \big(p_-\otimes\text{\rm Id}\big){\sf D}^-\vert_U : U \to U^-\otimes U,
&\delta^-& \defeq  \big(p_- \otimes p_-\big){\sf D}^-\vert_{U^-} : U^-\to U^-\otimes U^-,\\
\delta^+ &\defeq  (p_-\otimes\text{\rm Id})\overline{\sf D}^{\,-}\vert_{T^+} : T^+\to U^-\otimes T^+.
\end{aligned}
$$
\end{defn*}

\medskip

It follows from the multiplicativity of $p_\pm$ that $\Delta^+$ and $\delta^\pm$ {\it are algebra morphisms}. The assumption $p_-(\textcolor{red}{\bullet}^{0,\mathfrak{o}})={\bf1}_-$, is needed to ensure the formulas
\begin{align*}
\delta\bstau&={\bf1}_-\otimes\bstau+\sum_{\vert\bsvarphi\vert<0}\bsvarphi\otimes(\bstau/\bsvarphi),\\
\delta^-\bssigma&={\bf1}_-\otimes\bssigma+\bssigma\otimes{\bf1}_- + \sum_{\vert\bssigma\vert<\vert\mcFsi\vert<0}\mcFsi\otimes(\bssigma/\mcFsi)
\end{align*}
for $\bstau\in\mathcal{SC}$ and $\bssigma\in\mathcal{SC}^-$. 

\medskip

\begin{thm}\label{thm: all requirements of regul and renor str}
Set
\begin{align*}
\mathscr{T} &\defeq  \big((T^+,\Delta^+), (T,\Delta)\big),   \\
\mathscr{U} &\defeq  \big((U^-,\delta^-), (U,\delta)\big).
\end{align*}

\begin{enumerate}
\setlength{\itemsep}{0.1cm}

\item $\mathscr{T}$ is a regularity structure satisfying Assumptions \textbf{\textsf{(A)}} and \textbf{\textsf{(B)}}, with the grading $|\cdot|$.

\item $\mathscr{U}$ is a renormalization structure satisfying Assumption \textbf{\textsf{(E)}}, with the grading $|\cdot|'$.

\item $\mathscr{T}$ and $\mathscr{U}$ are compatible and satisfy Assumption \textbf{\textsf{(C)}}.

\item Assumption \textbf{\textsf{(D)}}, the compatibility between the splittings ${\updownharpoons_{\sf e}}\vert_U$ and $\delta$ holds true.

\end{enumerate}

\end{thm}   

\medskip

\begin{Dem} 
Write as shorthand
\begin{align*}
{\sf D}\bstau\;\textrm{ or } \;{\sf D}^+\bstau &= \sum_i \bssigma_i\otimes\bseta_i,   \\
{\sf D}^-\bstau\;\textrm{ or } \;\overline{\sf D}^{\,-}\bstau &= \sum_j 
\bsvarphi_j\otimes\mcFsi_j,
\end{align*}
and note that the following stability formulas of the naive and extended homogeneities. One has
\begin{align}\label{stability of extended homogeneity}
&|\bstau| = |\bssigma_i| + |\bseta_i|,   \\
&|\bstau|' = |\bsvarphi_j|' + |\mcFsi_j|',\quad |\bstau| = |\mcFsi_j| 
  \label{EqStabilityPrimed}
\end{align}
for each $i$ and $j$. Here we define $|{\bf1}_-|' \defeq  0$.

\ssk

\textit{\textsf{(a)}} 
By the first identity of \eqref{stability of extended homogeneity},
$$
(p_+\otimes p_+) {\sf D}^+p_+=(p_+\otimes p_+){\sf D}^+
$$
holds on ${\sf T}^+$. Then one has the comodule property of $\Delta$ as follows.
\begin{align*}
(\Delta\otimes\text{\rm Id})\Delta
&= \big(\text{\rm Id}\otimes p_+\otimes p_+\big)({\sf D}\otimes\text{\rm Id}){\sf D}   \\
&= \big(\text{\rm Id}\otimes p_+\otimes p_+\big)(\text{\rm Id}\otimes{\sf D}^+){\sf D}   \\
&= \big(\text{\rm Id}\otimes p_+\otimes p_+\big)(\text{\rm Id}\otimes{\sf D}^+)(\text{\rm Id}\otimes p_+){\sf D}   \\
&= (\text{\rm Id}\otimes\Delta^+)\Delta.
\end{align*}
The coassociativity of $\Delta^+$ is obtained similarly. One gets for free the existence of an antipode on $T^+$ from the fact that $T^+$ is a connected graded bialgebra -- see Proposition \ref{PropBialgebraHopf} in Appendix {\sf \ref{SectionAppendixAlgebra}}.

\ssk

\textit{\textsf{(b)}}
The comodule properties of $\delta$ and $\delta^-$ are obtained by the similar way to {\sf (a)}, since
$$
(p_-\otimes p_-){\sf D}^-p_-=(p_-\otimes p_-){\sf D}^-
$$
holds on ${\sf U}^-$, by identity \eqref{EqStabilityPrimed}.
By definition, ${\bf1}_-$ is the only element in $U^-$ of $0$ homogeneity, so $U^-$ is a connected graded bialgebra.

\ssk

\textit{\textsf{(c)}}
We prove the cointeraction property
$$
\mcM^{(13)}(\delta\otimes\delta^+)\Delta=(\iden\otimes\Delta)\delta;
$$
the proofs of other properties are left to readers. See also Proposition \ref{prop simplification of compatibility}. The second identity of \eqref{EqStabilityPrimed} yields
$$
\delta^+\circ p_+ = (\iden\otimes p_+)\delta^+
$$
on ${\sf T}^+$. Thus we have
\begin{align*}
\mcM^{(13)}(\delta\otimes\delta^+)\Delta
&=\mcM^{(13)}\big(\delta\otimes(\delta^+\circ p_+)\big){\sf D}   \\
&=\mcM^{(13)}\Big(p_-\otimes\iden\otimes p_-\otimes p_+\Big)\big({\sf D}^-\otimes\overline{\sf D}^{\,-}\big){\sf D}   \\
&=\big(p_-\otimes\iden\otimes p_+\big)\mcM^{(13)}\big({\sf D}^-\otimes\overline{\sf D}^{\,-}\big){\sf D}   \\
&=\big(p_-\otimes\iden\otimes p_+\big)(\iden\otimes{\sf D}){\sf D}^-\\
&=\big(\iden\otimes\Delta\big)\delta.
\end{align*}

\ssk

\textit{\textsf{(d)}}
Recall the explicit formula for the map $\updownharpoons_{\sf e}$, from Lemma \ref{lem:dual between graft and cut}. It is obvious that $U$ is stable under $\updownharpoons_{\sf e}$. Define
$$
\updownharpoons(\tau_{\frak{e}}^{\frak{n},\frak{o}})
\defeq  \sum_{\sigma\in A(\tau)}\sum_{\frak{n}_\sigma,\frak{e}_{\partial \sigma}'}
\frac1{\frak{e}_{\partial\sigma}'!}\binom{\frak{n}}{\frak{n}_\sigma}
(\tau/\sigma)_{\frak{e}+\frak{e}_{\partial\sigma}'}^{\frak{n}-\frak{n}_\sigma,\frak{o}\vert_{\tau\setminus\sigma}}
\otimes \sigma_{\frak{e}}^{\frak{n}_\sigma+\frak{e}_{\partial\sigma}',\frak{o}},
$$
where $A(\tau)\defeq \{P_e\tau\}_{e\in E_\tau}$ -- where $P_e\tau$ is a connected component of the graph $\tau/\{e\}$ containing the root of $\tau$, see \eqref{section6:defofcut}. Comparing this with the definition of $\Delta^+$,  it is not difficult to show the equality
$$
\mcM^{(13)}\big(\delta\otimes\delta\big)\updownharpoons
= (\text{\rm Id} \,\otimes \updownharpoons)\delta
$$
proceeding as in the proof of point {\sf (c)}. Note that the contracted tree $\tau/\sigma$ is always planted. Let $p_{\sf e}$ be the canonical projection on the set of planted trees $\eta$ with
$$
\frak{n}(\rho_\eta)=0,\quad \frak{e}(e_\eta)={\sf e},
$$
where $e_\eta$ is the only one edge leaving the root $\rho_\eta$, and let 
$c$ be the map sending the tree of the form $\mcI_n(\tau)$ to $\tau$. Then
$$
{\updownharpoons_{\sf e}}\vert_U=
(c\circ p_{\sf e}\otimes\text{\rm Id})\updownharpoons.
$$
on $U$. Since it is elementary to show
\begin{equation*}  \begin{split}
(\textrm{Id}\otimes c)\delta &= \delta\circ c,   \\
(\textrm{Id}\otimes p_{\sf e})\delta &=\delta\circ p_{\sf e},
\end{split}  \end{equation*}
the compatibility condition follows by writing
\begin{align*}
\mcM^{(13)}(\delta\otimes\delta){\updownharpoons_{\sf e}}\vert_U
&= \mcM^{(13)}\Big(\big(\text{\rm Id}\otimes c\circ p_{\sf e}\big)\delta\otimes\delta\Big)\updownharpoons\\
&= \big(\text{\rm Id}\otimes c\circ p_{\sf e}\otimes\text{\rm Id}\big) \mcM^{(13)}(\delta\otimes\delta)\updownharpoons\\
&= \big(\text{\rm Id}\otimes c\circ p_{\sf e}\otimes\text{\rm Id}\big) \big(\text{\rm Id}\,\otimes\updownharpoons\big)\delta = (\text{\rm Id}\,\otimes ({\updownharpoons_{\sf e}}\vert_U))\delta.
\end{align*}
\end{Dem}

\bigskip

\subsection{Some examples}

\newcommand{\rtd}[1]
{\begin{tikzpicture}
\coordinate (A1) at (0,0);
#1
\end{tikzpicture}
}

\newcommand{\rtdb}[2]
{\begin{tikzpicture}[baseline=#1pt]
\coordinate (A1) at (0,0);
#2
\end{tikzpicture}
}

\newcommand{\pol}[3]
{\coordinate (A#1) at ($(A#2)+0.4*({cos(#3)},{sin(#3)})$);}

\newcommand{\drc}[1]
{\foreach \n in {#1} \filldraw[white] (A\n) circle (1.7pt);
\foreach \n in {#1} \draw (A\n) circle (1.7pt);}
\newcommand{\drcb}[1]
{\foreach \n in {#1} \filldraw[white] (A\n) circle (1.7pt);
\foreach \n in {#1} \draw (A\n) circle (1.7pt);
\foreach \n in {#1} \node at ($(A\n)+(0.17,0)$) {\tiny $1$};}
\newcommand{\drb}[1]
{\foreach \n in {#1} \fill (A\n) circle (1.7pt);}
\newcommand{\drbb}[1]
{\foreach \n in {#1} \fill (A\n) circle (1.7pt);
\foreach \n in {#1} \node at ($(A\n)+(0.17,0)$) {\tiny $1$};}

\newcommand{\drl}[2]
{\foreach \n in {#2} \draw (A#1)--(A\n);}
\newcommand{\drll}[2]
{\foreach \n in {#2} \draw[double distance=0.7pt] (A#1)--(A\n);}

\newcommand{\drr}[2]
{\foreach \n in {#1} \fill[red] (A\n) circle (1.7pt);
\foreach \n in {#1} \node[right] at ($(A\n)+(0,0)$) {\tiny $#2$};}

-- For simplicity we consider the equation
$$
(\partial_t-\Delta_x+1)u=f(u)\zeta+g(u)(\partial_xu)^2
$$
for $d=1$ with the noise $\zeta\in\mcC^{-1-\kappa}$ for sufficiently small $\kappa>0$. Theorem \ref{ThmRenormPDEs} yields that, for any $k\in G^-_\textrm{ad}$ one has the renormalized equation
\begin{align*}
(\partial_t-\Delta_x+1)u^{(k)}
&=f(u^{(k)})\zeta+g(u^{(k)})(\partial_xu^{(k)})^2\\
&\quad+k(\circ)f(u^{(k)})+k(\circ\text{\tiny$1$})f'(u^{(k)})\partial_xu^{(k)}
+2k\left(\rtdb{2}{\pol{2}{1}{90}\drll{1}{2}\drb{1}\drc{2}}\right)
f(u^{(k)})g(u^{(k)})\partial_xu^{(k)}\\
&\quad+k\left(\rtdb{2}{\pol{2}{1}{90}\drl{1}{2}\drc{1,2}}\right)
f(u^{(k)})f'(u^{(k)})
+k\left(\rtdb{2}{\pol{2}{1}{120}\pol{3}{1}{60}
\drll{1}{2,3}\drb{1}\drc{2,3}}  \right)
f^2(u^{(k)})g(u^{(k)}).
\end{align*}
The double line $|\!|$ represents the edge with $\frak{e}$-decoration $(0,1)\in\bbN\times\bbN$. The dot $\circ\text{\tiny$1$}$ represents the node with $\frak{n}$-decoration $(0,1)\in\bbN\times\bbN$. More terms are needed when $\zeta\in\mcC^{-3/2-\kappa}$.

\medskip

-- The table below is the list of strongly conforming trees associated with the generalized (KPZ) equation \eqref{EqgKPZ}, without red nodes. Fix $d=1$ for simplicity. Fix also $\beta_0=-3/2-\kappa$ for sufficiently small $\kappa>0$. The dot $\bullet\text{\tiny1}$ represents the node with $\frak{n}$-decoration $(0,1)\in\bbN\times\bbN$.

\medskip

\begin{center}
\resizebox{\textwidth}{!}{ 
\renewcommand{\arraystretch}{1}
\begin{tabular}{l|c}\hline
{\sf Homogeneity} & {\sf Rooted decorated trees}  \\ 
\hline \hline

$\beta_0=-3/2-\kappa$ & \rtd{\drc{1}}  \\\hline
\vspace{-0.3cm} & \\

$2\beta_0+2=-1-2\kappa$ & 
\rtd{\pol{2}{1}{90}\drl{1}{2}\drc{1,2}}
\
\rtd{\pol{2}{1}{120}\pol{3}{1}{60}
\drll{1}{2,3}\drb{1}\drc{2,3}}  
\\\hline
\vspace{-0.3cm} & \\

$3\beta_0+4=-1/2-3\kappa$ &
\rtd{\pol{2}{1}{90}\pol{3}{2}{90}
\drl{1}{3}\drc{1,2,3}}
\
\rtd{\pol{2}{1}{120}\pol{3}{1}{60}
\drl{1}{2,3}\drc{1,2,3}}
\
\rtd{\pol{2}{1}{90}\pol{3}{2}{120}\pol{4}{2}{60}
\drl{1}{2}\drll{2}{3,4}\drb{2}\drc{1,3,4}}
\
\rtd{\pol{2}{1}{120}\pol{3}{1}{60}\pol{4}{2}{90}
\drl{2}{4}\drll{1}{2,3}\drb{1}\drc{2,3,4}}
\
\rtd{\pol{2}{1}{135}\pol{3}{1}{90}\pol{4}{1}{45}
\drl{1}{2}\drll{1}{3,4}\drb{1}\drc{2,3,4}}
\
\rtd{\pol{2}{1}{120}\pol{3}{1}{60}\pol{4}{2}{120}\pol{5}{2}{60}
\drll{1}{4}\drll{1}{3}\drll{2}{5}\drb{1,2}\drc{3,4,5}}
\\\hline
\vspace{-0.3cm} & \\

$\beta_0+1=-1/2-\kappa$ &
\rtd{\drcb{1}}\ 
\rtd{\pol{2}{1}{90}\drll{1}{2}\drb{1}\drc{2}}
\\\hline
\vspace{-0.3cm} & \\

$4\beta_0+6=-4\kappa$ &
\rtd{\pol{2}{1}{90}\pol{3}{2}{90}\pol{4}{3}{90}
\drl{1}{4}\drc{1,2,3,4}}
\
\rtd{\pol{2}{1}{90}\pol{3}{2}{120}\pol{4}{2}{60}
\drl{1}{2}\drl{2}{3,4}\drc{1,2,3,4}}
\
\rtd{\pol{2}{1}{120}\pol{3}{1}{60}\pol{4}{2}{90}
\drl{2}{4}\drl{1}{2,3}\drc{1,2,3,4}}
\
\rtd{\pol{2}{1}{135}\pol{3}{1}{90}\pol{4}{1}{45}
\drl{1}{2,3,4}\drc{1,2,3,4}}
\
\rtd{\pol{2}{1}{90}\pol{3}{2}{90}\pol{4}{3}{120}\pol{5}{3}{60}
\drl{1}{3}\drll{3}{4,5}\drb{3}\drc{1,2,4,5}}
\
\rtd{\pol{2}{1}{90}\pol{3}{2}{120}\pol{4}{2}{60}\pol{5}{3}{90}
\drl{1}{2}\drll{2}{3,4}\drl{3}{5}\drc{1,3,4,5}\drb{2}}
\
\rtd{\pol{2}{1}{90}\pol{3}{2}{135}\pol{4}{2}{90}\pol{5}{2}{45}
\drl{1}{2}\drl{2}{3}\drll{2}{4,5}\drc{1,3,4,5}\drb{2}}
\
\rtd{\pol{2}{1}{120}\pol{3}{1}{60}\pol{4}{2}{120}\pol{5}{2}{60}
\drl{1}{2,3}\drll{2}{4,5}\drc{1,3,4,5}\drb{2}}
\
\rtd{\pol{2}{1}{120}\pol{3}{1}{60}\pol{4}{2}{90}\pol{5}{4}{90}
\drll{1}{2,3}\drl{2}{5}\drb{1}\drc{2,3,4,5}}
\
\rtd{\pol{2}{1}{120}\pol{3}{1}{60}\pol{4}{2}{120}\pol{5}{2}{60}
\drll{1}{2,3}\drl{2}{4,5}\drc{2,3,4,5}\drb{1}}
\
\rtd{\pol{2}{1}{120}\pol{3}{1}{60}\pol{4}{2}{90}\pol{5}{3}{90}
\drll{1}{2,3}\drl{2}{4}\drl{3}{5}\drb{1}\drc{2,3,4,5}}
\
\rtd{\pol{2}{1}{135}\pol{3}{1}{90}\pol{4}{1}{45}\pol{5}{2}{90}
\drl{1}{2}\drll{1}{3,4}\drl{2}{5}\drb{1}\drc{2,3,4,5}}
\
\rtd{\pol{2}{1}{135}\pol{3}{1}{90}\pol{4}{1}{45}\pol{5}{2}{90}
\drll{1}{2,3}\drl{1}{4}\drl{2}{5}\drb{1}\drc{2,3,4,5}}
\
\rtd{\pol{2}{1}{150}\pol{3}{1}{110}\pol{4}{1}{70}\pol{5}{1}{30}
\drl{1}{2,3}\drll{1}{4,5}\drb{1}\drc{2,3,4,5}}
\\

&
\rtd{\pol{2}{1}{90}\pol{3}{2}{120}\pol{4}{2}{60}\pol{5}{3}{120}\pol{6}{3}{60}
\drl{1}{2}\drll{2}{3,4}\drll{3}{5,6}\drc{1,4,5,6}\drb{2,3}}
\
\rtd{\pol{2}{1}{120}\pol{3}{1}{60}\pol{4}{2}{90}\pol{5}{4}{120}\pol{6}{4}{60}
\drll{1}{2,3}\drl{2}{4}\drll{4}{5,6}\drc{2,3,5,6}\drb{1,4}}
\
\rtd{\pol{2}{1}{120}\pol{3}{1}{60}\pol{4}{2}{120}\pol{5}{2}{60}\pol{6}{4}{90}
\drll{1}{2,3}\drll{2}{4,5}\drl{4}{6}\drb{1,2}\drc{3,4,5,6}}
\
\rtd{\pol{2}{1}{135}\pol{3}{1}{45}\pol{4}{2}{135}\pol{5}{2}{90}\pol{6}{2}{45}
\drll{1}{2,3}\drl{2}{4}\drll{2}{5,6}\drb{1,2}\drc{3,4,5,6}}
\
\rtd{\pol{2}{1}{135}\pol{3}{1}{45}\pol{4}{2}{90}\pol{5}{3}{120}\pol{6}{3}{60}
\drll{1}{2,3}\drl{2}{4}\drll{3}{5,6}\drb{1,3}\drc{2,4,5,6}}
\
\rtd{\pol{2}{1}{135}\pol{3}{1}{90}\pol{4}{1}{45}\pol{5}{2}{120}\pol{6}{2}{60}
\drl{1}{2}\drll{1}{3,4}\drll{2}{5,6}\drb{1,2}\drc{3,4,5,6}}
\
\rtd{\pol{2}{1}{135}\pol{3}{1}{90}\pol{4}{1}{45}\pol{5}{2}{120}\pol{6}{2}{60}
\drll{1}{2,3}\drl{1}{4}\drll{2}{5,6}\drb{1,2}\drc{3,4,5,6}}
\
\rtd{\pol{2}{1}{120}\pol{3}{1}{60}\pol{4}{2}{120}\pol{5}{2}{60}\pol{6}{4}{120}\pol{7}{4}{60}
\drll{1}{3,6}\drll{2}{5}\drll{4}{7}\drb{1,2,4}\drc{3,5,6,7}}
\
\rtd{\pol{2}{1}{135}\pol{3}{1}{45}\pol{4}{2}{110}\pol{5}{2}{70}\pol{6}{3}{110}\pol{7}{3}{70}
\drll{1}{2,3}\drll{2}{4,5}\drll{3}{6,7}\drb{1,2,3}\drc{4,5,6,7}}
\\\hline
\vspace{-0.3cm} & \\ 

$2\beta_0+3=-2\kappa$ &
\rtd{\pol{2}{1}{90}\drl{1}{2}\drc{1}\drcb{2}}
\
\rtdb{-2}{\pol{2}{1}{90}\drl{1}{2}\drcb{1}\drc{2}}
\
\rtd{\pol{2}{1}{90}\pol{3}{2}{90}\drl{1}{2}\drll{2}{3}\drb{2}\drc{1,3}}
\
\rtd{\pol{2}{1}{90}\pol{3}{2}{90}\drll{1}{2}\drl{2}{3}\drb{1}\drc{2,3}}
\
\rtd{\pol{2}{1}{120}\pol{3}{1}{60}\drl{1}{2}\drll{1}{3}\drc{2,3}\drb{1}}
\
\rtdb{-2}{\pol{2}{1}{120}\pol{3}{1}{60}\drll{1}{2,3}\drc{2,3}\drbb{1}}
\
\rtd{\pol{2}{1}{120}\pol{3}{1}{60}\drll{1}{2,3}\drb{1}\drc{2}\drcb{3}}
\
\rtd{\pol{2}{1}{90}\pol{3}{2}{120}\pol{4}{2}{60}
\drll{1}{2}\drll{2}{3,4}\drb{1,2}\drc{3,4}}
\
\rtd{\pol{2}{1}{120}\pol{3}{1}{60}\pol{4}{2}{90}
\drll{1}{2,3}\drll{2}{4}\drb{1,2}\drc{3,4}}
\\

\hline
\end{tabular}   }
\end{center} 

\bigskip

-- Here are some examples of the actions of splitting map $\sf D^-$. The dot $\textcolor{red}{\bullet}\text{\tiny$(\alpha)$}$ represents the node with $\frak{o}$-decoration $\alpha\in\bbZ[\beta_0]$.
\begin{align*}
{\sf D}^-\circ&={\bf1}_-\otimes\circ+\circ\otimes\textcolor{red}{\bullet}\text{\tiny$(\beta_0)$},\\
{\sf D}^-\, \rtdb{1}{\pol{2}{1}{90}\drll{1}{2}\drb{1}\drc{2}}
&={\bf1}_-\otimes\,\rtdb{1}{\pol{2}{1}{90}\drll{1}{2}\drb{1}\drc{2}} 
+\circ\otimes\, \rtdb{1}{\pol{2}{1}{90}\drll{1}{2}\drb{1}\drr{2}{(\beta_0)}}\!
+\bullet\otimes\, \rtdb{1}{\pol{2}{1}{90}\drll{1}{2}\drc{2}\drr{1}{(0)}}\!
+\circ\bullet\otimes\ \rtdb{1}{\pol{2}{1}{90}\drll{1}{2}\drr{2}{(\beta_0)}\drr{1}{(0)}}\!
+\,\rtdb{1}{\pol{2}{1}{90}\drll{1}{2}\drb{1}\drc{2}} \otimes \textcolor{red}{\bullet}\text{\tiny$(\beta_0+1)$}.
\end{align*}
For larger trees, it is inconvenient to write down all possible terms. Note that some of them vanishes by the application of $p_-$ or {\color{red}$k\in G_{\rm ad}^-$}. Omitting them by $(\cdots)$, one has for example

\begin{align*}
{\sf D}^-
\rtdb{1}{\pol{2}{1}{120}\pol{3}{1}{60}\pol{4}{2}{120}\pol{5}{2}{60}\pol{6}{4}{120}\pol{7}{4}{60}
\drll{1}{3,6}\drll{2}{5}\drll{4}{7}\drb{1,2,4}\drc{3,5,6,7}}
&= {\bf1}_-\otimes
\rtdb{1}{\pol{2}{1}{120}\pol{3}{1}{60}\pol{4}{2}{120}\pol{5}{2}{60}\pol{6}{4}{120}\pol{7}{4}{60}
\drll{1}{3,6}\drll{2}{5}\drll{4}{7}\drb{1,2,4}\drc{3,5,6,7}}
+
\rtdb{1}{\pol{2}{1}{120}\pol{3}{1}{60}\pol{4}{2}{90}
\drll{1}{2,3}\drll{2}{4}\drb{1,2}\drc{3,4}}
\otimes
\rtdb{1}{\pol{2}{1}{90}\pol{3}{2}{120}\pol{4}{2}{60}
\drll{1}{2}\drll{2}{3,4}\drb{2}\drc{3,4}\drr{1}{(2\beta_0+3)}}
+2\
\rtdb{1}{\pol{2}{1}{120}\pol{3}{1}{60}\pol{4}{2}{90}
\drll{1}{2,3}\drll{2}{4}\drb{1,2}\drc{3,4}}
\otimes
\rtdb{1}{\pol{2}{1}{120}\pol{3}{1}{60}\pol{4}{2}{90}
\drll{1}{2,3}\drll{2}{4}\drb{1}\drc{3,4}
\fill[red] (A2) circle (1.7pt);
\node[left] at ($(A2)+(0,0)$) {\tiny $(2\beta_0+3)$};}
+
\rtdb{1}{\pol{2}{1}{120}\pol{3}{1}{60}\pol{4}{2}{120}\pol{5}{2}{60}\pol{6}{4}{120}\pol{7}{4}{60}
\drll{1}{3,6}\drll{2}{5}\drll{4}{7}\drb{1,2,4}\drc{3,5,6,7}}
\otimes \textcolor{red}{\bullet}\text{\tiny$(4\beta_0+6)$}
+(\cdots).
\end{align*}


\vfill\pagebreak

\appendix
\section{Summary of notations}
\label{SectionAppendixSummary}

The following is a summary of the notations that we used in several sections.

\begin{center}
\resizebox{\textwidth}{!}{ 
\renewcommand{\arraystretch}{1.2}
\begin{tabular}{l|c|l}\hline
{\sf Notations} & {\sf Section} &  \hfill{\sf Meaning}   \\ 
\hline \hline
$\mathscr{T} = \big((T^+,\Delta^+),(T,\Delta)\big)$ & {\sf \ref{SectionConcreteRS}}  & (Concrete) regularity structure.   \\
$\mcB^+$, $\mcB$ & {\sf\ref{SectionConcreteRS}}  & Bases of $T^+$ and $T$.   \\
$S_+$ & {\sf \ref{SectionConcreteRS}}  & Antipode of $T^+$.   \\
$G^+$, $\widehat{g}$ & {\sf \ref{SectionConcreteRS}} & Character group of 
$T^+$, and an action of $g\in G^+$ on $T$.   \\
$d(\cdot,\cdot)$, $|\cdot|_{\frak{s}}$ & {\sf \ref{SubsectionModelsAndCo}} & Scaled metric and scaled degree. \\
$\big(T_X^+,T_X\big)$, $\big(\mcB_X^+,\mcB_X\big)$ & {\sf \ref{SubsectionModelsAndCo}} & Polynomial regularity structure and their bases. \\

$\mcD_{(0,t)}^{\gamma,\eta}$  & {\sf \ref{SubsectionNonAnticipative}} & Singular modelled distributions on the time interval $(0,t)$. \\

$\mathscr{U} = \big((U,\delta),(U^-,\delta^-)\big)$ & {\sf \ref{SectionConcreteRenorS}}  & Renormalization structure.   \\
$G^-$, $\widetilde{k}$ & {\sf \ref{SectionConcreteRenorS}}  & Character group of $U^-$, and an action of $k\in G^-$ on $U$.   \\
$^k{\sf M}$  & {\sf \ref{SubsectionCompatibleStructures}}  &  Renormalized model.   \\

$\tau$, $N_\tau$, $E_\tau$, $\rho_\tau$  &  {\sf \ref{SectionMultiPreLie}}   & Rooted tree, node set, edge set, and root.  \\
$V$, $V^*$  &  {\sf \ref{SectionMultiPreLie}}   & Vector space spanned by rooted decorated trees, and its copy space.   \\
$\grafting$, $\grafting_\flat$   &   {\sf \ref{SectionMultiPreLie}}   &   
Grafting operator $V^*\otimes V^*\to V^*$, and its projection on $U^*$.   
\\
${\sf M}^\zeta = \big({\sf g}^\zeta, {\sf \Pi}^\zeta\big)$  &  {\sf \ref{SectionMultiPreLie}}   & Canonical model associated with a smooth noise $\zeta$.  \\
$\updownharpoons_{\sf e}$, ${\updownharpoons_{\sf e}}\vert_U$ &  {\sf \ref{SectionrenormalizationMultiPreLie}}  & Dual map of $\curvearrowright_{\sf e}$, and its restriction to $U$.   \\

$S'_- : U^- \rightarrow \bbR[U]$ & {\sf \ref{SectionBHZCharacter}}  &  Twisted negative antipode.   \\

$\mathcal{SC}$  &  {\sf \ref{SectionDecoratedTrees}}   & Set of all strongly conforming decorated trees.   \\
$\mathcal{C}$  &   {\sf \ref{SectionDecoratedTrees}}  & Set of all conforming decorated trees.   \\
\hline
\end{tabular}   }
\end{center}

\bigskip

\section{Basics from algebra}
\label{SectionAppendixAlgebra}

We recall some basics of bialgebras, Hopf algebras, and comodules without 
proofs. See \cite{Sweedler, Manchon, Foissy16} for details. Note that, for any two algebras $A$ and $B$ with units ${\bf1}_A$ and ${\bf1}_B$ respectively, the tensor space $A\otimes B$ is also an algebra with the product
$$
(a_1\otimes b_1)\cdot(a_2\otimes b_2) \defeq  (a_1a_2)\otimes(b_1b_2),
\qquad(a_1,a_2\in A, \ b_1,b_2\in B)
$$
and with unit ${\bf1}_A\otimes{\bf1}_B$.

\medskip

\begin{defn*}

A \textbf{\textsf{bialgebra}} $(B,\mcM,{\bf1},\triangle,\theta)$ is a $5$-tuple of the following components.

\begin{itemize}
\setlength{\itemsep}{1pt}
\item An algebra $B$ with product $\mcM:B\otimes B\to B$, and unit ${\bf1}$.
\item An algebra morphism $\triangle: B\to B\otimes B$ satisfying the \textbf{\textsf{coassociativity}}
$$
(\triangle\otimes\iden)\triangle=(\iden\otimes\triangle)\triangle.
$$
\item An algebra morphism $\theta:B\to\bbR$, satisfying
$$
(\theta\otimes\iden)\triangle = (\iden\otimes\theta)\triangle = \iden,
$$
where we identify $a\otimes\tau=\tau\otimes a=a\tau$, for any $a\in\bbR$ and $\tau\in B$.
\end{itemize}
The map $\triangle$ is called a \textbf{\textsf{coproduct}}, and the map $\theta$ is called a \textbf{\textsf{counit}}. An algebra morphism $S:B\to B$, such that
\begin{equation} \label{EqDefnAntipode}
\mcM(\iden\otimes S)\triangle = \mcM(S\otimes\iden)\triangle = \theta(\cdot){\bf1}
\end{equation}
is called an \textbf{\textsf{antipode}}. A bialgebra equipped with an antipode $S$ is called a \textbf{\textsf{Hopf algebra}}.

\end{defn*}

\medskip

The counit $\theta(\cdot)$ is traditionally denoted $\varepsilon(\cdot)$. 
We use a different letter as $\epsilon$ already stands for a regularization parameter in this work. The following result gives a sufficient condition for a bialgebra to be a Hopf algebra.
A bialgebra $B$ is called \textbf{\textsf{graded}} if it is a direct sum $\bigoplus_{\lambda\in \Lambda} B_\lambda$ of vector spaces such that
\begin{itemize}
\setlength{\itemsep}{1pt}
\item $\Lambda$ be a locally finite subset of $[0,\infty)$ such that $0\in \Lambda$ and $\Lambda+\Lambda\subset \Lambda$.
\item ${\bf1}\in B_0$ and $B_\lambda \cdot B_\mu\subset B_{\lambda+\mu}$, 
for any $\lambda,\mu\in\Lambda$.
\item $\triangle B_\lambda\subset\bigoplus_{\mu,\nu\in \Lambda,\, \mu+\nu=\lambda} B_\mu\otimes B_\nu$.
\end{itemize}
We call $\Lambda$ a \textbf{\textsf{grading}} in this paper. A graded bialgebra with $B_0 = \langle{\bf1}\rangle$ is said to be \textbf{\textsf{connected}}.

\medskip

\begin{prop}
\cite[Exercises pages 228 and 238]{Sweedler},
\cite[Proposition II.1.1 and Corollary II.3.2]{Manchon}   
\label{PropBialgebraHopf}
Any connected graded bialgebra is a Hopf algebra. Moreover, one has the following properties.
\begin{itemize}
\setlength{\itemsep}{1pt}
\item $\theta({\bf1})=1$ and $\theta(\tau)=0$ for any $\tau\in\bigoplus_{\lambda>0}B_\lambda$.
\item $\triangle{\bf1}={\bf1}\otimes{\bf1}$, and for any $\tau\in B_\lambda$ with $\lambda>0$,
$$
\triangle\tau\in \bigg\{\tau\otimes{\bf1}+{\bf1}\otimes\tau + \underset{\mu+\nu=\lambda,\,0<\mu<\lambda}{\sum_{\mu,\nu\in \Lambda} }B_\mu\otimes 
B_\nu \bigg\}.
$$
\end{itemize}
\end{prop}

\noindent Based on the first assertion, we denote by ${\bf1}'$ the counit 
$\theta$ of a connected graded bialgebra. The preceding formula for $\triangle\tau$ gives an inductive formula for the antipode. For $\tau\neq {\bf 1}$ and $\triangle\tau = \tau\otimes{\bf 1} + {\bf 1}\otimes \tau + \sum \tau_1\otimes\tau_2$, one has
$$
S(\tau) = -\tau - \sum S(\tau_1)\tau_2.
$$ 

\medskip

On the dual space $B'$ of the bialgebra $B$, the \textbf{\textsf{convolution product}} is defined by
$$
(f*g)\tau \defeq  (f\otimes g)\triangle\tau,
$$
for all $f,g\in B',\ \tau\in B$, where we identify $a\otimes b=ab$ for any $a,b\in\bbR$. The coassociativity of $\triangle$ implies the associativity of the convolution
$$
(f*g)*h=f*(g*h),
$$
for all $f,g\in B'$, and the counit $\theta$ is indeed a unit of the convolution product
$$
f*\theta = \theta*f=f,
$$
for all $f\in B'$. Hence the triplet $(B',*,\theta)$ is a unital ring. Moreover, the subset $G\subset B'$ of algebra morphisms $g:B\to\bbR$ is stable under the convolution product. The existence of an antipode $S$ implies that $G$ is a group. Indeed, the inverse of $g\in G$ is given by $g^{-1}=g\circ S$. Each element of $G$ is called a \textbf{\textsf{character}}, and when $B$ is a Hopf algebra, the set $G$ is called the \textbf{\textsf{character group}}.

\medskip

We recall comodules and comodule bialgebras.
Given an algebra $A$ and two spaces $E,F$, we define on the algebraic tensor product $A\otimes E\otimes A\otimes F$ the $A\otimes E\otimes F$-valued map
$$
\mcM^{(13)}\big(a_1\otimes e\otimes a_2\otimes f\big) \defeq  (a_1a_2)\otimes e\otimes f.
$$

\begin{defn*}
Let $(B,\mcM,{\bf1},\triangle,\theta)$ be a bialgebra.
\begin{itemize}
\item A linear space $M$ equipped with a linear map $\delta:M\to B\otimes 
M$, with the properties
$$
(\text{\rm Id}_B\otimes\delta)\delta=(\triangle\otimes\text{\rm Id}_M)\delta,\quad\textrm{and}\quad
(\theta\otimes\iden_M)\delta=\iden_M,
$$
is called a left $B$-\textbf{\textsf{comodule}}.
Similarly, a linear space $N$ is called a right $B$-comodule if a linear map $\rho:N\to N\otimes B$, exists and satisfies
$$
(\rho\otimes\iden_B)\rho = (\iden_N\otimes\triangle)\rho,\quad\textrm{and}\quad (\iden_N\otimes\theta)\rho = \iden_N.
$$
\item A bialgebra $M$ is called a left $B$-\textbf{\textsf{comodule bialgebra}} if $M$ is a left $B$-comodule by an algebra morphism $\delta:M\to B\otimes M$, such that
$$
\mcM^{(13)}(\delta\otimes\delta)\triangle_M=(\text{\rm Id}\otimes\triangle_M)\delta,\quad\textrm{and}\quad
(\text{\rm Id}\otimes\theta_M)\delta = \theta_M(\cdot){\bf1},
$$
where $\triangle_M$ is a coproduct of $M$, and $\theta_M$ is a counit of $M$.
\end{itemize}
\end{defn*}

\medskip

\begin{prop} \cite[Proposition 2]{Foissy16} \label{*compatibility antipode}
Let $M$ be a $B$-comodule bialgebra. If $M$ has an antipode $S_M$, then 
$$
\delta\circ S_M=(\iden_B\otimes S_M)\delta.
$$
\end{prop}

\bigskip

\section{Technical proofs}
\label{SectionAppendixProofs}

This section is dedicated to proving Theorem \ref{thm singular reconstruction} and Lemma \ref{coassosictive cointeraction Dpm}.

\bigskip

\subsection{Proof of Theorem \ref{thm singular reconstruction}}
\label{Proof of singular reconstruction}

For any $x\in\bbR\times\bbR^d$, and $\lambda\in(0,1]$, denote by $\varphi\mapsto\varphi_x^\lambda$ the transformation of functions on $\bbR\times\bbR^d$ defined by
$$
\varphi_x^\lambda(y) \defeq  \lambda^{-d-2}\varphi\Big(\lambda^{-2}(y_0-x_0),\lambda^{-1}(y'-x')\Big).
$$
The following bound appears in the Hairer's original paper \cite{Hai14}. Recall $\beta_0=\min A$ and that $p_t$ denotes the heat kernel of the operator $\widetilde{\mcG}=\partial_{x_0}^2-(\Delta_{x'}-1)^2$ on $\bbR\times\bbR^d$.

\medskip

\begin{lem}
Let $\sf M=(\Pi,g)$ be a model over the regularity structure $\mathscr{T}$ and $\bsf\in\mcD^\gamma(T,{\sf g})$ with $\gamma\in\bbR$.
Assume $\beta_0>-4$. Then for any Schwartz function $\varphi\in \mcS(\bbR\times\bbR^d)$, $x\in\bbR\times\bbR^d$, and $\lambda\in(0,1]$, one has the bound
\begin{align*}
\left|\big\langle \textbf{\textsf{R}}^{\sf M}\bsf-{\sf\Pi}_x^{\sf g}\bsf(x),\varphi_x^\lambda \big\rangle\right|
\le C_\varphi\|{\sf\Pi}^{\sf g}\|\left\|\bsf\right\|_{\mcD^\gamma}\lambda^\gamma,
\end{align*}
where the constant $C_\varphi$ depends on the size $\sup_{|k|_{\frak{s}},|\ell|_{\frak{s}}\le N}\|x^k\partial^\ell\varphi\|_{L^\infty(\bbR\times\bbR^d)}$ for $N>0$ large enough.
\end{lem}

\medskip

\begin{Dem}
Write $\Lambda_x \defeq  \textbf{\textsf{R}}^{\sf M}\bsf-{\sf\Pi}_x^{\sf g}\bsf(x)$, to shorten notations. Using $p_0=\int_0^{\lambda^4}\widetilde{\mcG} p_tdt+p_{\lambda^4}$ and the symmetry of $\widetilde{\mcG}$,
\begin{align*}
\langle \Lambda_x,\varphi_x^\lambda \rangle
&= \int_{\bbR\times\bbR^d}\int_0^{\lambda^4}\big\langle \Lambda_x, \widetilde{\mcG}_y 
p_t(y,\cdot)\big\rangle\, \varphi_x^\lambda(y)\,dt dy
+ \int_{\bbR\times\bbR^d} \big\langle \Lambda_x,p_{\lambda^4}(y,\cdot) \big\rangle \,\varphi_x^\lambda(y)\,dy   \\
&= \int_{\bbR\times\bbR^d}\int_0^{\lambda^4}\big\langle \Lambda_x, 
p_t(y,\cdot)\big\rangle\, \widetilde{\mcG}_y \varphi_x^\lambda(y)\,dt dy
+ \int_{\bbR\times\bbR^d} \big\langle \Lambda_x,p_{\lambda^4}(y,\cdot) \big\rangle \,\varphi_x^\lambda(y)\,dy
\eqdef(A)+(B).
\end{align*}
Using the properties of models as in Proposition \ref{PropModelsOnKernels}, one has
\begin{align*}
\left\vert\langle \Lambda_x, p_t(y,\cdot)\rangle\right\vert
\lesssim t^{\frac\gamma4}+\sum_{\beta\in[\beta_0,\gamma)}t^{\frac\beta4}d(y,x)^{\gamma-\beta}.
\end{align*}
This implies $|(B)|\lesssim\lambda^\gamma$. Moreover, since $\beta_0>-4$, one has
\begin{align*}
|(A)|\lesssim
\sum_{\beta\in[\beta_0,\gamma]}\int_0^{\lambda^4}t^{\frac\beta4}dt
\int_{\bbR\times\bbR^d}d(y,x)^{\gamma-\beta}\, |\widetilde{\mcG}_y \varphi_x^\lambda(y)|\,dy
\lesssim
\sum_{\beta\in[\beta_0,\gamma]}\lambda^{\beta+4}\lambda^{\gamma-\beta-4}\lesssim\lambda^\gamma.
\end{align*}
\end{Dem}

\medskip

We prove a fundamental fact on the connection of modelled distributions. For any interval $I\subset\bbR$, denote by $\mcD_I^\gamma(T,{\sf g})$ the set of functions $\bsf:I\times\bbR^d\to T_{<\gamma}$ which satisfies the bounds of $\brarb{\bsf}_{\mcD_I^\gamma}$ and $\|\bsf\|_{\mcD_I^\gamma}$ as in Definition \ref{DefnModelledDistribution} with $\bbR\times\bbR^d$ replaced by $I\times\bbR^d$.

\medskip

\begin{lem}\label{App:lem:connectionMD}
Let $\sf M=(g,\Pi)$ be a model over $\mathscr{T}$ and let $a<b<c$.
If $\bsf:[a,c]\times\bbR^d\to T_{<\gamma}$ satisfies the bounds of $\trino{\bsf}_{\mcD_{[a,b]}^\gamma}$ and $\trino{\bsf}_{\mcD_{[b,c]}^\gamma}$, then $\bsf\in\mcD_{[a,c]}^\gamma(T,{\sf g})$ and one has
$$
\trino{\bsf}_{\mcD_{[a,c]}^\gamma}
\lesssim (1+\|\widehat{\sf g}\|_\gamma)\big(\trino{\bsf}_{\mcD_{[a,b]}^\gamma}+\trino{\bsf}_{\mcD_{[b,c]}^\gamma}\big).
$$
\end{lem}

\medskip

\begin{Dem}
It is sufficient to show the bound of $\|\bsf(y)-\widehat{{\sf g}_{yx}}\bsf(x)\|_\beta$ for $x\in[a,b]\times\bbR^d$ and $y\in[b,c]\times\bbR^d$. Setting $z=(b,y')$, we have
\begin{align*}
\|\bsf(y)-\widehat{{\sf g}_{yx}}\bsf(x)\|_\beta
&\le\|\bsf(y)-\widehat{{\sf g}_{yz}}\bsf(z)\|_\beta+\big\|\widehat{{\sf g}_{yz}}\big(\bsf(z)-\widehat{{\sf g}_{zx}}\bsf(x)\big)\big\|_\beta\\
&\le\|\bsf\|_{\mcD_{[b,c]}^\gamma}d(y,z)^{\gamma-\beta}
+\|\widehat{\sf g}\|_\gamma\|\bsf\|_{\mcD_{[a,b]}^\gamma}\sum_{\alpha\in[\beta,\gamma)}d(y,z)^{\alpha-\beta}d(z,x)^{\gamma-\alpha}\\
&\lesssim(1+\|\widehat{\sf g}\|_\gamma)\big(\trino{\bsf}_{\mcD_{[a,b]}^\gamma}+\trino{\bsf}_{\mcD_{[b,c]}^\gamma}\big)d(y,x)^{\gamma-\beta}.
\end{align*}
\end{Dem}

\medskip

We recall now from J. Martin's work \cite[Theorem 5.3.16]{Martin} the existence of a `Whitney extension' map on locally defined modelled distributions. 

\medskip

\begin{thm}
\label{thm Whitney extension}
Let $\sf M=(g,\Pi)$ be a model over $\mathscr{T}$ with a regular product $\star$ satisfying Assumption \REFA. Then there exists a continuous linear operator $E:\mcD_{(-\infty,t]}^\gamma(T,{\sf g})\to\mcD^\gamma(T,\sf g)$ such that $(E\bsf)\vert_{(-\infty,t]\times\bbR^d}=\bsf$, and the bound
$$
\trino{E\bsf}_{\mcD^\gamma}\le C\trino{\bsf}_{\mcD_{(-\infty,t]}^\gamma}
$$
holds for a positive constant $C=C(\|\widehat{\sf g}\|_\gamma)$ depending polynomially on $\|\widehat{\sf g}\|_\gamma$, and independent of $t>0$ and $\bsf\in\mcD_{(-\infty,t]}^\gamma(T,\sf g)$. A similar result holds for the modelled distributions defined on $[t,\infty)\times\bbR^d$.
\end{thm}

\medskip

\begin{Dem}
For simplicity, we consider $t=0$. It is sufficient to construct the continuous linear extension operator $\widetilde{E}:\mcD_{(-\infty,0]}^\gamma(T,{\sf g})\to \mcD_{(-\infty,1]}^\gamma(T,{\sf g})$. Indeed, once we pick $\chi\in C^\infty(\bbR\times\bbR^d)$ which depends only on $x_0$ such that $\chi(x_0,x')=1$ when $x_0\le0$ and $\chi(x_0,x')=0$ when $x_0\ge1$, and set
$$
\boldsymbol{\chi}(x)=\sum_{|k|_\mfs<\gamma-\beta_0}\frac{\partial^k\chi(x)}{k!}X^k,
$$
then by Proposition \ref{PropRegularityProduct}, we can define the modelled distribution on $\bbR\times\bbR^d$ by
$$
E\bsf\defeq\mcQ_{<\gamma}\big(\boldsymbol{\chi}\star(\widetilde{E}\bsf)\big).
$$
Then the operator $E:\mcD_{(-\infty,0]}^\gamma(T,{\sf g})\to \mcD^\gamma(T,{\sf g})$ satisfies the desired properties.

We construct the extension operator $\widetilde{E}:\mcD_{(-\infty,0]}^\gamma(T,{\sf g})\to \mcD_{(-\infty,1]}^\gamma(T,{\sf g})$. Let $h(x_0,x')=h_{x_0}(x')$ be the kernel of the operator $e^{x_0\Delta_{x'}}$ with $x_0>0$, and define the function
$$
\bsh(x)=\sum_{|k|_{\frak{s}}<\gamma-\beta_0}\frac{\partial^kh(x)}{k!}X^k
$$
on $(0,\infty)\times\bbR^d$. Then we have the properties
\begin{align}\label{7:proof:prop:flatextension3}
\int_{\bbR^d}\bsh(x)dx'={\bf1},\qquad
\int_{\bbR^d}|\partial^kh(x_0,x')||x'|^\alpha\lesssim x_0^{\frac{\alpha-|k|_\mfs}2}
\end{align}
for any $\alpha\ge0$.
For any $\bsf\in\mcD_{(-\infty,0]}^\gamma(T,{\sf g})$, we set $(\widetilde{E}\bsf)(x)\defeq \bsf(x)$ if $(-\infty,0]\times\bbR^d$ and
$$
(\widetilde{E}\bsf)(x) \defeq  \mcQ_{<\gamma}\int_{\bbR^d}\bsh\big(x-(0,y')\big)\star\widehat{{\sf g}_{x(0,y')}}f(0,y')\,dy'
$$
if $x\in(0,\infty)\times\bbR^d$. 
We first prove
\begin{align}
\label{7:prop:flatextension3}
\|(\widetilde{E}\bsf)(x)-\widehat{{\sf g}_{x(0,x')}}\bsf(0,x')\|_\beta
\lesssim\|\widehat{\sf g}\|\|\bsf\|_{\mcD^\gamma}\, x_0^{\frac{\gamma-\beta}2}
\end{align}
for any $\beta<\gamma$ and $x=(x_0,x')\in[0,1]\times\bbR^d$. Since
$$
(\widetilde{E}\bsf)(x)-\widehat{{\sf g}_{x(0,x')}}\bsf(0,x')
=\int_{\bbR^d}\bsh(x_1,x'-y')\star\widehat{{\sf g}_{x(0,x')}}
\big\{\widehat{{\sf g}_{(0,x')(0,y')}}\bsf(0,y')-\bsf(0,x')\big\}dy'
$$
by the first property of \eqref{7:proof:prop:flatextension3}, we have \eqref{7:prop:flatextension3} as follows.
\begin{align*}
&\|(\widetilde{E}\bsf)(x)-\widehat{{\sf g}_{x(0,x')}}\bsf(0,x')\|_\beta\\
&\le\sum_{|k|_\mfs<\gamma-\beta_0}\int_{\bbR^d}|\partial^kh(x_0,x'-y')|
\big\|\widehat{{\sf g}_{x(0,x')}}
\big\{\widehat{{\sf g}_{(0,x')(0,y')}}\bsf(0,y')-\bsf(0,x')\big\}\big\|_{\beta-|k|_\mfs}dy'\\
&\le\|\widehat{\sf g}\|\|\bsf\|_{\mcD^\gamma}\sum_{|k|_\mfs<\gamma-\beta_0}
\sum_{\alpha\in[\beta-|k|_\mfs,\gamma)}
\int_{\bbR^d}|\partial^kh(x_0,x'-y')|\,
x_0^{\frac{\alpha-\beta+|k|_\mfs}2}|x'-y'|^{\gamma-\alpha}dy'\\
&\lesssim\|\widehat{\sf g}\|\|\bsf\|_{\mcD^\gamma}\sum_{|k|_\mfs<\gamma-\beta_0}
\sum_{\alpha\in[\beta-|k|_\mfs,\gamma)}x_0^{\frac{\alpha-\beta+|k|_\mfs}2}x_0^{\frac{\gamma-\alpha-|k|_\mfs}2}
\lesssim\|\widehat{\sf g}\|\|\bsf\|_{\mcD^\gamma}\,x_0^{\frac{\gamma-\beta}2}.
\end{align*}
The bound of $\big\|(\widetilde{E}\bsf)(x)\big\|_\beta$ on $x\in[0,1]\times\bbR^d$ follows from \eqref{7:prop:flatextension3} and the bound of $\bsf(0,x')$.
For the bound of $\|(\widetilde{E}\bsf)(y)-\widehat{{\sf g}_{yx}}(\widetilde{E}\bsf)(x)\|_\beta$, it is sufficient to consider the case that $x,y\in[0,1]\times\bbR^d$ by Lemma \ref{App:lem:connectionMD}. If $x_0^{\frac12}\wedge y_0^{\frac12}\le d(x,y)$, we can bound above $\|(\widetilde{E}\bsf)(y)-\widehat{{\sf g}_{yx}}(\widetilde{E}\bsf)(x)\|_\beta$ by 

\vfill \pagebreak

\begin{align*}
&\|(\widetilde{E}\bsf)(y)-\widehat{{\sf g}_{y(0,y')}}\bsf(0,y')\|_\beta + \big\|\widehat{{\sf g}_{y(0,y')}}\big\{\bsf(0,y')-\widehat{{\sf g}_{(0,y')(0,x')}}\bsf(0,x')\big\}\big\|_\beta   \\
&\quad+\big\|\widehat{{\sf g}_{yx}}\big\{\widehat{{\sf g}_{x(0,x')}}\bsf(0,x')-(\widetilde{E}\bsf)(x)\big\}\big\|_\beta,
\end{align*}
and using \eqref{7:prop:flatextension3}, we obtain the desired bound. If $d(x,y)<x_0^{\frac12}\wedge y_0^{\frac12}$, by using the formula

\begin{align*}
\bsh(y)-\widehat{{\sf g}_{yx}}\bsh(x)
&=\sum_{|k|_\mfs<\gamma-\beta_0}\frac{X^k}{k!}\bigg(\partial^kh(y)-\sum_{|\ell|_\mfs<\gamma-\beta_0-|k|_\mfs}\frac{(y-x)^\ell}{\ell!}\partial^{k+\ell}h(x)\bigg)\\
&=\sum_{|k|_\mfs<\gamma-\beta_0}\frac{X^k}{k!}
\sum_{|n|_\mfs\ge\gamma-\beta_0-|k|_\mfs}\frac{(y-x)^n}{n!}\int_{[0,1]^{d+1}}\partial^{k+n}h(x_t)m_n(dt)
\end{align*}
obtained from Lemma \ref{lem:anisoTaylor}, where $n$ runs over a finite set, we decompose
\begin{align*}
&(\widetilde{E}\bsf)(y)-\widehat{{\sf g}_{yx}}(\widetilde{E}\bsf)(x)\\
&=\mcQ_{<\gamma}\int_{\bbR^d}\Big\{\bsh\big(y-(0,z')\big)-\widehat{{\sf g}_{yx}}\bsh\big(x-(0,z')\big)\Big\}
\star\widehat{{\sf g}_{y(0,z')}}\bsf(0,z')dz'\\
&=\mcQ_{<\gamma}\sum_{\substack{|k|_\mfs<\gamma-\beta_0 \\ |n|_\mfs\ge\gamma-\beta_0-|k|_\mfs}}
\frac{(y-x)^n}{k!n!}\int_{\bbR^d}\int_{[0,1]^{d+1}}
\partial^{k+n}h\big(x_t-(0,z')\big)m_n(dt)X^k\star\widehat{{\sf g}_{y(0,z')}}\bsf(0,z')dz'.
\end{align*}
Thus it is sufficient to show that
\begin{equation}\label{proof:Taylor:thm Whitney extension}
\bigg\|\int_{\bbR^d}\partial^{k+n}h\big(x_t-(0,z')\big)\widehat{{\sf g}_{y(0,z')}}\bsf(0,z')dz'\bigg\|_{\beta-|k|_\mfs}
\lesssim d(x,y)^{\gamma-\beta-|n|_\mfs}.
\end{equation}
For this purpose, we decompose
\begin{align*}
&\bigg\|\int_{\bbR^d}\partial^{k+n}h\big(x_t-(0,z')\big)\widehat{{\sf g}_{y(0,z')}}\bsf(0,z')dz'\bigg\|_{\beta-|k|_\mfs}\\
&\le\bigg\|\int_{\bbR^d}\partial^{k+n}h\big(x_t-(0,z')\big)
\widehat{{\sf g}_{y(0,x_t')}}\big\{\widehat{{\sf g}_{(0,x_t')(0,z')}}\bsf(0,z')-\bsf(0,x_t')\big\}
dz'\bigg\|_{\beta-|k|_\mfs}\\
&\le\|\widehat{\sf g}\|_\gamma\|\bsf\|_{\mcD_0^\gamma}\sum_{\alpha\in[\beta-|k|_\mfs,\gamma)}\int_{\bbR^d}\big|\partial^{k+n}h\big(x_t-(0,z')\big)\big|
d(x,y)^{\alpha-\beta+|k|_\mfs}|x_t'-z'|^{\gamma-\alpha}dz'\\
&\lesssim\|\widehat{\sf g}\|_\gamma\|\bsf\|_{\mcD_0^\gamma}\sum_{\alpha\in[\beta-|k|_\mfs,\gamma)}
(x_0\wedge y_0)^{\frac{\gamma-\alpha-|k|_\mfs-|n|_\mfs}2}d(x,y)^{\alpha-\beta+|k|_\mfs}.
\end{align*}
Since $\gamma-\alpha-|k|_\mfs-|n|_\mfs\le-\alpha+\beta_0\le0$ under the condition that $|n|_\mfs\ge\gamma-\beta_0-|k|_\mfs$ and $\alpha\in A$, we obtain the desired bound \eqref{proof:Taylor:thm Whitney extension}.
\end{Dem}

\medskip

\begin{cor}\label{cor:thm singular reconstruction}
Let $\sf M=(g,\Pi)$ be a model over $\mathscr{T}$ with a regular product $\star$ satisfying Assumption \REFA.
For any $\bsf\in\mcD^{\gamma,\eta}(T,{\sf g})$ and $a>0$, the restriction $\bsf\vert_{[a,\infty)\times\bbR^d}$ satisfies the bound
$$
\trino{\bsf\vert_{[a,\infty)\times\bbR^d}}_{\mcD_{[a,\infty)}^\gamma}\le Ca^{\frac{\eta\wedge\beta_0-\gamma}2}\trino{\bsf}_{\mcD^{\gamma,\eta}}
$$
for a positive constant $C=C(\|\widehat{\sf g}\|_\gamma)$ depending polynomially on $\|\widehat{\sf g}\|_\gamma$, independent of $a$ and $\bsf$.
Therefore by writing $\textbf{\textsf{R}}_a^{\sf M}\defeq \textbf{\textsf{R}}^{\sf M}\circ E_a$ for the extension map $E_a$ from $\mcD_{[a,\infty)}^\gamma(T,{\sf g})$ to $\mcD^\gamma(T,{\sf g})$, for any Schwartz function $\varphi\in\mcS(\bbR\times\bbR^d)$, $x\in\bbR\times\bbR^d$, and $\lambda\in(0,1]$, one has the bounds
\begin{equation}\label{proof:eq1:cor:thm singular reconstruction}
\left|\big\langle \textbf{\textsf{R}}_a^{\sf M}\bsf-{\sf\Pi}_x^{\sf g}E_a\bsf(x),\varphi_x^\lambda \big\rangle\right|
\le C'\trino{\bsf}_{\mcD^{\gamma,\eta}}a^{\frac{\eta\wedge\beta_0-\gamma}2}\lambda^\gamma
\end{equation}
and
\begin{equation}\label{proof:eq2:cor:thm singular reconstruction}
\left|\big\langle \textbf{\textsf{R}}_a^{\sf M}\bsf,\varphi_x^\lambda \big\rangle\right|
\le C'\trino{\bsf}_{\mcD^{\gamma,\eta}}a^{\frac{\eta\wedge\beta_0-\gamma}2}\lambda^\gamma.
\end{equation}
for a positive constant $C'$ depending polynomially on $\|{\sf g}\|_\gamma+\|{\sf\Pi}^{\sf g}\|_\gamma$, independent of $a$ and $\bsf$.
\end{cor}

\medskip

\begin{Dem}
It is sufficient to show the bound of $\big\|\bsf(y)-\widehat{{\sf g}_{yx}}\bsf(x)\big\|_\beta$ for any $\beta<\gamma$ and $x,y\in[a,\infty)\times\bbR^d$ such that $d(x,y)\ge \omega(x,y)$. We decompose
\begin{align*}
\big\|\bsf(y)-\widehat{{\sf g}_{yx}}\bsf(x)\big\|_\beta
&\le\big\|\bsf(y)\big\|_\beta+\big\|\widehat{{\sf g}_{yx}}\bsf(x)\big\|_\beta\\
&\lesssim\brarb{\bsf}_{\mcD^{\gamma,\eta}}\bigg(a^{\frac{(\eta-\beta)\wedge0}2}+\sum_{\beta\le\alpha<\gamma}d(x,y)^{\alpha-\beta}a^{\frac{(\eta-\alpha)\wedge0}2}\bigg).
\end{align*}
Since $d(x,y)\ge\omega(x,y)\ge a$, we have
\begin{align*}
\big\|\bsf(y)-\widehat{{\sf g}_{yx}}\bsf(x)\big\|_\beta
&\lesssim\brarb{\bsf}_{\mcD^{\gamma,\eta}}d(x,y)^{\gamma-\beta}\sum_{\beta\le\alpha<\gamma}d(x,y)^{\alpha-\gamma}a^{\frac{(\eta-\alpha)\wedge0}2}\\
&\lesssim\brarb{\bsf}_{\mcD^{\gamma,\eta}}d(x,y)^{\gamma-\beta}a^{\frac{\eta\wedge\beta_0-\gamma}2}.
\end{align*}
Thus the modelled distribution $\bsf$ restricted to $[a,\infty)\times\bbR^d$ has the norm of size $a^{\frac{\eta\wedge\beta_0-\gamma}2}$.
\end{Dem}

\medskip

We turn to the proof of the reconstruction theorem for singular modelled distributions.

\medskip

\begin{Dem}[of Theorem \ref{thm singular reconstruction}]
The proof is just an analogue of the proof of Proposition 6.9 in Hairer' seminal work \cite{Hai14}, so we omit a number of details. The only difference is that $p_t(x,\cdot)$ is not compactly supported. For the sake of generality, for any Schwartz function $\varphi\in\mcS(\bbR\times\bbR^d)$, $x\in\bbR\times\bbR^d$, and $\lambda\in(0,1]$, we prove the bound
\begin{equation}\label{proof:eq1:thm singular reconstruction}
\left|\big\langle \textbf{\textsf{R}}^{\sf M}\bsf-{\sf\Pi}_x^{\sf g}\bsf(x),\varphi_x^\lambda \big\rangle\right|
\lesssim \trino{\bsf}_{\mcD^\gamma}\omega(x)^{\eta\wedge\beta_0-\gamma}\lambda^\gamma
\end{equation}
for any $\lambda\in(0,\omega(x)]$, and the bound
\begin{equation}\label{proof:eq2:thm singular reconstruction}
\left|\big\langle \textbf{\textsf{R}}^{\sf M}\bsf,\varphi_x^\lambda \big\rangle\right|
\lesssim \trino{\bsf}_{\mcD^\gamma}\lambda^\eta
\end{equation}
for any $\lambda\in(0,1]$. To lighten notation, we omit the proportional constants depending on the model $\sf M$.
By linearity, we can assume that $\bsf=0$ on $(-\infty,0]\times\bbR^d$. 
For any $a>0$, we consider the distributions $\textbf{\textsf{R}}_a^{\sf M}\bsf\defeq \textbf{\textsf{R}}^{\sf M}(E_a\bsf)$ defined in Corollary \ref{cor:thm singular reconstruction}.
Because of the local property of the reconstruction operator, Corollary \ref{CorLocalDependenceReconstruction}, these distributions are compatible over all $a>0$ in the sense that $\langle \textbf{\textsf{R}}_a^{\sf M}\bsf,\varphi\rangle=\langle \textbf{\textsf{R}}_b^{\sf M}\bsf,\varphi\rangle$ for any $0<a<b$ and Schwartz functions $\varphi$ supported on $[b,\infty)\times\bbR^d$, so the quantity $\big\langle \widetilde{\textbf{\textsf{R}}}{}^{\sf M}\bsf,\varphi \big\rangle$ is defined for any $\varphi$ supported on $(0,\infty)\times\bbR^d$. Since $\bsf$ vanishes on $\bbR\times\bbR^d_0$, one defines $\big\langle \widetilde{\textbf{\textsf{R}}}{}^{\sf M}\bsf,\varphi \big\rangle=0$ if $\varphi$ is supported on $(-\infty,0)\times\bbR^d$. 
To consider arbitrary $\varphi$, fix a family $\{\phi_{n,k}\}_{n\in\bbN,k\in\bbZ^d}$ of functions of the forms
$$
\phi_{n,k}=2^{-n(d+2)}\phi_{x_{n,k}}^{2^{-n}}, \quad x_{n,k}=\big(2^{-2n},2^{-n}k\big)\in\bbR\times\bbR^d
$$
where $\phi$ is a smooth function supported on $\big\{x\in\bbR\times\bbR^d;d(0,x)<1\big\}$, and such that $\sum_{n,k}\phi_{n,k}(x)=1$, if $0<x_0<1/2$. 
To show the bound \eqref{proof:eq2:thm singular reconstruction}, fixing an integer $n_0$ such that $2^{-n_0}\simeq \lambda$, and setting $\widetilde{\phi}_{n_0}=1-\sum_{n\ge n_0,k\in\bbZ^d}\phi_{n,k}$, we decompose
$$
\big\langle \widetilde{\textbf{\textsf{R}}}{}^{\sf M}\bsf,\varphi_x^\lambda \big\rangle
=\sum_{n\ge n_0,k\in\bbZ^d}\big\langle \widetilde{\textbf{\textsf{R}}}{}^{\sf M}\bsf,\varphi_x^\lambda\,\phi_{n,k} \big\rangle
+ \big\langle \widetilde{\textbf{\textsf{R}}}{}^{\sf M}\bsf,\varphi_x^\lambda\,\widetilde{\phi}_{n_0} \big\rangle.
$$
For the second term, since $\widetilde{\phi}_{n_0}$ is supported in $\{y\in\bbR\times\bbR^d;y_0\ge a_{n_0}\simeq2^{-2n_0}\}$ and $\varphi_x^\lambda\,\widetilde{\phi}_{n_0}=\psi_x^\lambda$ for some Schwartz function $\psi$ which is uniform over $\lambda,x,n_0$, by \eqref{proof:eq2:cor:thm singular reconstruction} one has
\begin{align*}
\left|\big\langle \widetilde{\textbf{\textsf{R}}}{}^{\sf M}\bsf,\varphi_x^\lambda\,\widetilde{\phi}_{n_0} \big\rangle\right|
\lesssim2^{-n_0(\eta\wedge\beta_0-\gamma)}\lambda^\gamma|\!|\!|\bsf|\!|\!|_{\mcD^{\gamma,\eta}}
\simeq2^{-n_0(\eta\wedge\beta_0)}|\!|\!|\bsf|\!|\!|_{\mcD^{\gamma,\eta}}.
\end{align*}
For the first term, since $\varphi_x^\lambda\phi_{n,k}$ is supported in $[2^{-2n},\infty)\times\bbR^d$ and roughly equal to 
$$
2^{-n(d+2)}\varphi_x^\lambda(x_{n,k})\phi_{x_{n,k}}^{2^{-n}},
$$ 
one has
\begin{align*}
\left|\big\langle \widetilde{\textbf{\textsf{R}}}{}^{\sf M}\bsf,\varphi_x^\lambda\,\phi_{n,k}\big\rangle\right|
&\lesssim2^{-n(d+2)}|\varphi_x^\lambda(x_{n,k})|2^{-n(\eta\wedge\beta_0-\gamma)}(2^{-n})^\gamma|\!|\!|\bsf|\!|\!|_{\mcD^{\gamma,\eta}}\\
&\lesssim2^{-n(d+2)}|\varphi_x^\lambda(x_{n,k})|2^{-n(\eta\wedge\beta_0)}|\!|\!|\bsf|\!|\!|_{\mcD^{\gamma,\eta}}
\end{align*}
The sum $2^{-nd}\sum_{k\in\bbZ^d}|\varphi_x^\lambda(x_{n,k})|$ is roughly bounded by $\int_{\bbR^d}|\varphi^\lambda(x)|dx'\lesssim\lambda^{-2}$. Hence
\begin{align*}
\sum_{n\ge n_0,\,k\in\bbZ^d}\left|\big\langle \widetilde{\textbf{\textsf{R}}}{}^{\sf M}\bsf,\varphi_x^\lambda\,\phi_{n,k}\big\rangle\right|
&\lesssim|\!|\!|\bsf|\!|\!|_{\mcD^{\gamma,\eta}}\lambda^{-2}\sum_{n\ge n_0}2^{-n(\eta\wedge\beta_0+2)}\\
&\lesssim|\!|\!|\bsf|\!|\!|_{\mcD^{\gamma,\eta}}\lambda^{-2}2^{-n_0(\eta\wedge\beta_0+2)}
\lesssim|\!|\!|\bsf|\!|\!|_{\mcD^{\gamma,\eta}}\lambda^{\eta\wedge\beta_0}.
\end{align*}
Here we use the assumption $\eta\wedge\beta_0>-2$.

The proof of \eqref{proof:eq1:thm singular reconstruction} is almost the same with some modifications.
We fix the same $n_0$ as above and decompose $\varphi_x^\lambda$ into $\sum_{n\ge n_0,\, k\in\bbZ^d}\varphi_x^\lambda\phi_{n,k}+\varphi_x^\lambda\widetilde{\phi}_{n_0}$.
For the dual with the second term, by taking $a\simeq x_0$ in \eqref{proof:eq1:thm singular reconstruction}, one has
\begin{align*}
\left|\big\langle \widetilde{\textbf{\textsf{R}}}{}^{\sf M}\bsf-{\sf\Pi}_x^{\sf g}\bsf(x),\varphi_x^\lambda\,\widetilde{\phi}_{n_0} \big\rangle\right|
\lesssim\omega(x)^{\eta\wedge\beta_0-\gamma}\lambda^\gamma|\!|\!|\bsf|\!|\!|_{\mcD^{\gamma,\eta}}.
\end{align*}
For the dual with the remaining terms, one decomposes

\begin{align*}
\widetilde{\textbf{\textsf{R}}}{}^{\sf M}\bsf-{\sf\Pi}_x^{\sf g}\bsf(x)
&=\left(\widetilde{\textbf{\textsf{R}}}{}^{\sf M}\bsf-{\sf\Pi}_{x_{n,k}}^{\sf g}\bsf(x_{n,k})\right)
+{\sf\Pi}_{x_{n,k}}^{\sf g}\bsf(x_{n,k})
-{\sf\Pi}_{x_{n,k}}^{\sf g}\Big(\widehat{{\sf g}_{x_{n,k}x}}\bsf(x)\Big)  
 \\
&\eqdef (a)+(b)+(c).
\end{align*}
For $(a)$, one has
\begin{align*}
\big\vert\langle(a),\varphi_x^\lambda\phi_{n,k}\rangle\big\vert
\lesssim2^{-n(d+2)}|\varphi_x^\lambda(x_{n,k})|2^{-n(\eta\wedge\beta_0)}\|{\sf\Pi}^{\sf g}\||\!|\!|\bsf|\!|\!|_{\mcD^{\gamma,\eta}}
\end{align*}
as before. The sum $2^{-nd}\sum_{k\in\bbZ^d}|\varphi_x^\lambda(x_{n,k})|$ is roughly bounded by $\int_{\bbR^d}|\varphi^\lambda(x)|dx'\lesssim\lambda^{-2}|\lambda^{-2}x_0|^{-N}$ for any $N>0$. By picking $N=\frac{\gamma-\eta\wedge\beta_0}2$, one has
\begin{align*}
\sum_{n\ge n_0,\,k\in\bbZ^d}\big\vert\langle(a),\varphi_x^\lambda\phi_{n,k}\rangle\big\vert
&\lesssim\lambda^{-2}|\lambda^{-2}x_0|^{\frac{\eta\wedge\beta_0-\gamma}2}
\sum_{n\ge n_0}2^{-n(\eta\wedge\beta_0+2)}|\!|\!|\bsf|\!|\!|_{\mcD^{\gamma,\eta}}\\
&\lesssim\lambda^{-2}|\lambda^{-2}x_0|^{\frac{\eta\wedge\beta_0-\gamma}2}2^{-n_0(\eta\wedge\beta_0+2)}|\!|\!|\bsf|\!|\!|_{\mcD^{\gamma,\eta}}\\
&\lesssim\omega(x)^{\eta\wedge\beta_0-\gamma}\lambda^\gamma|\!|\!|\bsf|\!|\!|_{\mcD^{\gamma,\eta}}.
\end{align*}
Using the bound $\|f(x)\|_\beta\lesssim|x_0|^{(\eta-\beta)/2\wedge0}$, one gets the same bounds as above for $(b)$ and $(c)$.

\ssk

It remains to show the uniqueness of $\textbf{\textsf{R}}^{\sf M}\bsf\in\mcC^{\eta\wedge\beta_0}$ satisfying \eqref{EqSingularReconstructionCondition}.
We start from the identity
$$
\big\vert \big\langle\textbf{\textsf{R}}^{\sf M}\bsf - (\textbf{\textsf{R}}^{\sf M})'\bsf , p_t(x,\cdot)\big\rangle\big\vert \lesssim \omega(x)^{\eta\wedge\beta_0-\gamma}t^{\frac\gamma4}
$$
satisfied uniformly in $x\in\bbR\times\bbR^d$ such that $\omega(x)^4\ge t$ by any other reconstruction 
operator $(\textbf{\textsf{R}}^{\sf M})'$. Since $\textbf{\textsf{R}}^{\sf M}\bsf,(\textbf{\textsf{R}}^{\sf M})'\bsf\in\mcC^{\eta\wedge\beta_0}$, we also have
$$
\big\vert \big\langle\textbf{\textsf{R}}^{\sf M}\bsf - (\textbf{\textsf{R}}^{\sf M})'\bsf , p_t(x,\cdot)\big\rangle\big\vert \lesssim t^{\frac{\eta\wedge\beta_0}4}
$$
uniformly in $x\in\bbR\times\bbR^d$ and $0< t\le1$. From these bounds, for any $\varepsilon\in(0,\gamma)$ we can show that
$$
\big\vert \big\langle\textbf{\textsf{R}}^{\sf M}\bsf - (\textbf{\textsf{R}}^{\sf M})'\bsf , p_t(x,\cdot)\big\rangle\big\vert \lesssim \omega(x)^{\eta\wedge\beta_0-\varepsilon}t^{\frac{\varepsilon}4}
$$
uniformly in $x\in\bbR\times\bbR^d$ and $0< t\le1$. If one chooses small $\varepsilon>0$ such that $\eta\wedge\beta_0-\varepsilon>-2$, since $\int_{\bbR\times\bbR^d}\omega(x)^{\eta\wedge\beta_0-\varepsilon}\varphi(x)dx<\infty$ for any Schwartz function $\varphi\in\mcS(\bbR\times\bbR^d)$, one has
$$
\big\langle \textbf{\textsf{R}}^{\sf M}\bsf - \big(\textbf{\textsf{R}}^{\sf M}\big)'\bsf , \varphi \big\rangle = \underset{t\rightarrow 0}{\lim} 
\int \Big\langle  \textbf{\textsf{R}}^{\sf M}\bsf - (\textbf{\textsf{R}}^{\sf M})'\bsf , p_t(x,\cdot) \Big\rangle\,\varphi(x)dx = \underset{t\rightarrow 0}{\lim} \, O\big(t^{\frac\varepsilon4}\big) = 0.
$$
\end{Dem}

\medskip

In the end of this section, we provide a sketch of the proof of Theorem \ref{thm wellposedness of K}'.

\medskip

\begin{Dem}[of Theorem \ref{thm wellposedness of K}']
The proof is carried out by a method similar to the proof of Theorem \ref{thm wellposedness of K} with modifications on the bounds of $\mcJ^{\sf M}$ and $\mcN^{\sf M}$ terms.
Similarly to the proof of Theorem \ref{thm wellposedness of K}, we use the decomposition $\mcJ^{\sf M}=\int_0^1\mcJ^{\sf M}_rdr$ and $\mcN^{\sf M}=\int_0^1\mcN^{\sf M}_rdr$ and write $(\tau)_{X^n}$ for the $X^n$-component of $\tau\in T$. As for the bound of $\brarb{\mcK^{\sf M}\bsf}_{\mcD^{\gamma,\eta}}$, we focus on the bound of
$$
[\circleddash]_r \defeq  
\Big( \mcJ^{\sf M}_r(x)\bsf(x) + (\mcN^{\sf M}_r\bsf)(x)\Big)_{X^n}
$$
for any $|n|_\mfs<\gamma+2$.
If $r\le\omega(x)^4$, we estimate the $\mcJ^{\sf M}$ and $\mcN^{\sf M}$ terms separately. By using the bound \eqref{EqSingularReconstructionCondition} for $\mcN^{\sf M}$, we have

\begin{align*}
[\circleddash]_r\lesssim\sum_{\beta\in(|n|_\mfs-2,\gamma)}\omega(x)^{(\eta-\beta)\wedge0}r^{\frac{\beta-|n|_\mfs-2}4}+\omega(x)^{\eta\wedge\beta_0-\gamma}r^{\frac{\gamma-|n|_\mfs-2}4}.
\end{align*}
Its integration over $r\in[0,\omega(x)^4]$ has an upper bound $\omega(x)^{\eta\wedge\beta_0+2-|n|_\mfs}$.
If $r>\omega(x)^4$, we decompose
\begin{align*}
[\circleddash]_r =\frac1{n!}\Big\langle\textbf{\textsf{R}}^{\sf M}\bsf,\partial_x^nq_r(x,\cdot)\Big\rangle
-\sum_{\beta<\gamma,\,|n|_\mfs\ge\beta+2}\frac1{n!}\Big\langle{\sf\Pi}_x^{\sf g}\big(\bsf(x)\big)_\beta,\partial_x^nq_r(x,\cdot)\Big\rangle.
\end{align*}
Note that, in the sum over $\beta$, the term associated with $\beta$ such that $|n|_\mfs=\beta+2$ vanishes because of \eqref{eq: integral of K vanishes}.
Since $\textbf{\textsf{R}}^{\sf M}\bsf\in\mcC^{\eta\wedge\beta_0}$, we have
\begin{align*}
[\circleddash]_r \lesssim r^{\frac{\eta\wedge\beta_0-|n|_\mfs-2}4}
+\sum_{\beta<\gamma,\,|n|_\mfs>\beta+2}\omega(x)^{(\eta-\beta)\wedge0}r^{\frac{\beta-|n|_\mfs-2}4}.
\end{align*}
Its integration over $r\in[0,\omega(x)^4]$ has an upper bound $\omega(x)^{(\eta\wedge\beta_0+2-|n|_\mfs)\wedge0}$.

As for the bound of $\|\mcK^{\sf M}\bsf\|_{\mcD^{\gamma,\eta}}$, we focus on the bound of $(\circleddash)_r$ similarly to the proof of Theorem \ref{thm wellposedness of K}.
For the integrations over $r\le d(x,y)^4$ and $d(x,y)^4\le r\le\omega(x,y)^4$, we use the same $(\Asterisk)$-decomposition and the $(\bigstar)$-decomposition respectively. Then we have the same bounds except the existence of the factor $\omega(x,y)^{\eta\wedge\beta_0-\gamma}$. For the integration over $r\ge\omega(x,y)^4$, we further decompose $(\bigstar)_r^1$ into
$$
(\bigstar)_r^1
=\frac1{n!} \Big\langle \textbf{\textsf{R}}^{\sf M}\bsf , (\partial^n q_r)_{y,x}^{\gamma+2-|n|_{\frak{s}}} \Big\rangle  
- \frac1{n!} \Big\langle {\sf\Pi}_x^{\sf g}\bsf(x) , (\partial^n q_r)_{y,x}^{\gamma+2-|n|_{\frak{s}}} \Big\rangle.  
$$
For the first term, since $\textbf{\textsf{R}}^{\sf M}\bsf\in\mcC^{\eta\wedge\beta_0}$, by using the same remainder formula \eqref{lem:anisoTaylorAppqr}, we have the upper bound $\sum_{|\ell|_\mfs>\gamma+2-|n|_\mfs}d(x,y)^{|\ell|_\mfs}r^{\frac{\eta\wedge\beta_0-|n|_\mfs-|\ell|_\mfs-2}4}$.
Its integration over $r\in[\omega(x,y)^4,1]$ has an upper bound
$$
\sum_{|\ell|_\mfs>\gamma+2-|n|_\mfs}d(x,y)^{|\ell|_\mfs}\omega(x,y)^{\eta\wedge\beta_0-|n|_\mfs-|\ell|_\mfs+2}
\lesssim d(x,y)^{\gamma+2-|n|_\mfs}\omega(x,y)^{\eta\wedge\beta_0-\gamma}
$$
since $d(x,y)\le\omega(x,y)$. We also have a similar bound for the remaining terms.
\end{Dem}

\bigskip

\subsection{Proof of Proposition \ref{PropUniversalPreLie} and Proposition \ref{PropUniversalPreLie2}}
\label{Appendix Proof of mpL}

We provide a sketch of the proof Proposition \ref{PropUniversalPreLie} and Proposition \ref{PropUniversalPreLie2}, following \cite{FoissyTypedTrees} and \cite{BM23}. Recall that, the basis $\mcV$ of $V$ is the set of all rooted trees with node types $\frak{T}_{\rm n}$ and edge types $\frak{T}_{\rm e}$, and with 
decorations $\frak{n}:N_\tau\to\bbN^{d+1}$ and $\frak{e}:E_\tau\to\bbN^{d+1}$.
Write ${\sf N}=\frak{T}_{\rm n}\times\bbN^{d+1}$ and ${\sf E}=\frak{T}_{\rm e}\times\bbN^{d+1}$.
Moreover, the basis $\mcB$ of $U$ is a subset satisfying Assumption \textsf{\textbf{(D)}}.
For simplicity, we prove the following proposition for $V$ and $U$, instead of $V^*$ and $U^*$.

\medskip

\begin{prop*}
\begin{enumerate}
\setlength{\itemsep}{0.15cm}
\item The space $V$ (or $U$) with the operators $\{\grafting\}_{\sf e\in E}$ (or $\{\flatgrafting\}_{\sf e\in E}$) is an ${\sf E}$-multi-pre-Lie algebra.
\item Let $(W,\{\rightslice_{\sf {e}}\}_{\sf{e}\in\sf{E}})$ be an $\sf E$-multi-pre-Lie algebra, and let $\{\varphi_{\frak{b}^n}\}_{(\frak{t},n)\in{\sf N}}$ (or $\{\varphi_{\frak{b}^n}\}_{(\frak{t},n)\in{\sf N}\cap\mcB}$) $\subset W$.
Then there exists a unique $\sf E$-multi-pre-Lie morphism $\varphi:V$ (or 
$U$) $\to W$ such that $\varphi(\frak{b}^n)=\varphi_{\frak{b}^n}$ for any $\frak{b}^n\in{\sf N}$ (or ${\sf N}\cap\mcB$).
\end{enumerate}
\end{prop*}

\medskip

The ${\sf E}$-multi-pre-Lie property of $(V,\{\grafting\}_{\sf e\in E})$
$$
(\bstau\overset{\sf e}{\curvearrowright}\bssigma)\overset{{\sf e}'}{\curvearrowright}\bseta
-\bstau\overset{\sf e}{\curvearrowright}(\bssigma\overset{{\sf e}'}{\curvearrowright}\bseta)
=
(\bssigma\overset{{\sf e}'}{\curvearrowright}\bstau)\overset{\sf e}{\curvearrowright}\bseta
-\bssigma\overset{{\sf e}'}{\curvearrowright}(\bstau\overset{\sf e}{\curvearrowright}\bseta)
$$
is proved in a similar way to Proposition 2.2 of \cite{FoissyTypedTrees} and Corollary 2.9 of \cite{BM23}, so we omit the proof. The same property for $(U,\{\flatgrafting\}_{\sf e\in E})$ follows from it. Indeed, by Assumption \textbf{\textsf{(D1)}} for the canonical projection $\pi_U:V\to U$, we have
$$
(\bstau\overset{\sf e}{\curvearrowright}_\flat\bssigma)\overset{{\sf e}'}{\curvearrowright}_\flat\bseta
-\bstau\overset{\sf e}{\curvearrowright}_\flat(\bssigma\overset{{\sf e}'}{\curvearrowright}_\flat\bseta)
=
\pi_U\big(
(\bstau\overset{\sf e}{\curvearrowright}\bssigma)\overset{{\sf e}'}{\curvearrowright}\bseta
-\bstau\overset{\sf e}{\curvearrowright}(\bssigma\overset{{\sf e}'}{\curvearrowright}\bseta)
\big).
$$

In order to prove \textit{\textsf{(b)}}, we introduce the Guin-Oudom extension of the multi-pre-Lie structure. The following is the content of Section 2.2 of \cite{FoissyTypedTrees} and Section 3.2 of \cite{BM23}.

\medskip

\begin{defn*}
Let $(W,\{\rightslice_{\sf {e}}\}_{\sf{e}\in\sf{E}})$ be an $\sf E$-multi-pre-Lie algebra.
Let $(W,{\sf e})=\{(a,{\sf e})\}_{a\in W}$ be a copy of the linear space $W$, and denote
$$
W^{\oplus{\sf E}} \defeq  \bigoplus_{\sf e\in E}(W,{\sf e}).
$$
Moreover, let $S(W^{\oplus{\sf E}})$ be the symmetric algebra of $W^{\oplus{\sf E}}$, with unit $\bf 1$. Then one can define the following linear maps.
\begin{itemize}
\setlength{\itemsep}{0.15cm}
\item Define the linear map $\rightslice_{\sf e}:W\otimes S(W^{\oplus{\sf 
E}})\to S(W^{\oplus{\sf E}})$ inductively as follows.
\begin{align*}
&c\rightslice_{\sf e}{\bf1}=0,   \\
&c\rightslice_{\sf e}\prod_{i=1}^N (b_i,{\sf e}_i)=\sum_{i=1}^N(c\rightslice_{\sf e}b_i,{\sf e}_i)\prod_{j\neq i}(b_j,{\sf e}_j),
\end{align*}
where $c,b_1,\dots,b_N\in W$ and ${\sf e}_1,\dots,{\sf e}_N\in{\sf E}$.
\item Define the linear map $\rightslice:S(W^{\oplus{\sf E}})\otimes W\to 
W$ inductively as follows.
\begin{align*}
&{\bf1}\rightslice b=b,\\
&(c,{\sf e})\rightslice b=c\rightslice_{\sf e}b,\\
&\prod_{i=1}^N (c_i,{\sf e}_i)\rightslice b = c_1\rightslice_{{\sf e}_1}\left(\prod_{i=2}^N(c_i,{\sf e}_i)\rightslice b\right) - \left(c_1\rightslice_{{\sf e}_1}\prod_{i=2}^N(c_i,{\sf e}_i)\right)\rightslice b,
\end{align*}
where $c,c_1,\dots,c_N,b\in W$ and ${\sf e},{\sf e}_1,\dots,{\sf e}_N\in{\sf E}$. (The last quantity is invariant under the permutations of $(c_1,{\sf e}_1),\dots,(c_N,{\sf e}_N)$ because of the multi-pre-Lie property of $W$, so the extension $\rightslice$ is well-defined.)
\end{itemize}
\end{defn*}

\medskip

The above extensions keep the multi-pre-Lie morphism property.
Indeed, if $\varphi:V\to W$ is an $\sf E$-multi-pre-Lie morphism, then defining the extension $\varphi:S(V^{\oplus{\sf E}})\to S(W^{\oplus{\sf E}})$ by $\varphi((c,{\sf e})) \defeq  (\varphi(c),{\sf e})$ for any $(c,{\sf e})\in(V,{\sf E})$, and denoting by $\curvearrowright$ the extension $\curvearrowright:S(V^{\oplus{\sf E}})\otimes V\to V$, one has
$$
\varphi(X\curvearrowright b)=\varphi(X)\rightslice\varphi(b)
$$
for any $X\in S(V^{\oplus{\sf E}})$ and $b\in V$.

\ssk

Since $\frak{T}_{\rm e}=\{\mcI\}$, we identify $\sf E$ with $\bbN^{d+1}$. The following formula can be proved by a similar argument to Lemma 2.6 of \cite{FoissyTypedTrees} by taking the `Taylor deformation' map $\Theta$ introduced in \cite{BM23} (see Theorem 2.7 and Proposition 3.8 of therein) into account. 

\medskip

\begin{lem}
For any $(\frak{t},n)\in{\sf N}\cap\mcB$, $\bstau_1,\dots,\bstau_a\in\mcB$, and $p_1,\dots,p_a\in\bbN^{d+1}$, one has
$$
\pi_U\left(\frak{b}^n\star\Bigstar_{i=1}^a \mcI_{p_i}(\bstau_i)\right)
=\sum_{\substack{q_1\le p_1,\dots,q_a\le p_a \\ q_1+\cdots+q_a\le n}}
(-1)^{|q_1|+\cdots+|q_a|}\binom{n}{q_1,\dots,q_a}
\prod_{i=1}^a (\bstau_i,p_i-q_i)\curvearrowright_\flat \frak{b}^{n-q_1-\cdots-q_a},
$$
where $\curvearrowright_\flat$ denotes the extension $\curvearrowright_\flat : S(U^{\oplus{\sf E}})\otimes U\to U$, and $\binom{n}{q_1,\dots,q_a}$ 
is the multinomial coefficient
$$
\binom{n}{q_1,\dots,q_a}=\frac{n!}{q_1!\cdots q_a!(n-q_1-\cdots-q_a)!}.
$$
\end{lem}

\medskip

Then we can prove the uniqueness part of \textit{\textsf{(b)}} immediately. Indeed, if $\varphi:U\to W$ is an $\sf E$-multi-pre-Lie morphism, then it extends to a multi-pre-Lie morphism from $S(U^{\oplus{\sf E}})$ to $S(W^{\oplus{\sf E}})$ and satisfies
$$
\varphi\left(\frak{b}^n\star\Bigstar_{i=1}^a \mcI_{p_i}(\bstau_i)\right)
=\sum_{\substack{q_1\le p_1,\dots,q_a\le p_a \\ q_1+\cdots+q_a\le n}}
(-1)^{|q_1|+\cdots+|q_a|}\binom{n}{q_1,\dots,q_a}
\prod_{i=1}^a (\varphi(\bstau_i),p_i-q_i)\rightslice \varphi(\frak{b}^{n-q_1-\cdots-q_a})
$$
for any $\frak{b}^n\star\Bigstar_{i=1}^a \mcI_{p_i}(\bstau_i)\in\mcB$.
The right hand side provides the recursive definition of the map $\varphi$, so we can conclude that $\varphi$ is determined by the values $\varphi(\frak{b}^n)$ for any $\frak{b}^n\in{\sf N}$. On the other hand, given $\{\varphi(\frak{b}^n)\}_{(\frak{t},n)\in{\sf N}}$, we can prove that the map $\varphi$ defined by the above formula satisfies indeed the multi-pre-Lie property. See Proposition 2.5 and Corollary 2.10 of \cite{BM23} for details.
We do not provide the details here because only the uniqueness part of \textit{\textsf{(b)}} is used in this paper, especially in Corollary \ref{CorrenormalizedVectorFields}.

\bigskip

\subsection{Proof of Lemma \ref{coassosictive cointeraction Dpm}}
\label{Appendix Proof of coassociativity}

\subsubsection{Reduced coproducts}

First we consider trees with $\frak{n}$ and $\frak{e}$-decorations, without $\frak{o}$-decoration. Recall that $SC$ is a set of strongly conforming trees and $C$ is a set of conforming trees. Set
\begin{align*}
{^\circ}T \defeq  \text{\rm span}(SC),   \qquad
{^\circ}{\sfT}^+ \defeq  \text{\rm span}(C),   \qquad
{^\circ}{\sfU}^- \defeq  \bbR[SC].
\end{align*}


\begin{defn*}
We define the following splitting maps.
\begin{itemize}
\setlength{\itemsep}{0.15cm}
	\item[{\sf 1.}] The linear map ${^\circ}{\sfD} : {^\circ}T\to {^\circ}T\otimes {^\circ}{\sfT}^+$, is defined for $\tau_{\frak{e}}^{\frak{n}}\in SC$ by
\begin{equation*} \begin{split}
{^\circ}\sfD(\tau_{\frak{e}}^{\frak{n}}) \defeq & \sum_{\mu\in ST(\tau)} \sum_{\frak{n}_\mu, \frak{e}'_{\partial\mu}} 
\frac{1}{\frak{e}'_{\partial\mu}!}{{\frak{n}}\choose{\frak{n}_\mu}}
\, \mu^{\frak{n}_\mu+\pi\frak{e}'_{\partial\mu}}_{\frak{e}}
\otimes 
\big(\tau/^{\text{\rm blue}}\mu\big)^{[\frak{n}-\frak{n}_\mu]_\mu}_{\frak{e}+\frak{e}'_{\partial\mu}},
\end{split}  \end{equation*}
where $ST(\tau)$ is the set of all {\sl subtrees} $\mu$ of $\tau$ which contain the root of $\tau$, and the second sum is over functions $\frak{n}:N_\mu\to\bbN\times\bbN^d$ with $\frak{n}_\mu\leq \frak{n}$ and functions $\frak{e}'_{\partial\mu}: \partial\mu\rightarrow\bbN\times\bbN^d$. The algebra morphism 
$$
{^\circ}{\sfD}^+ : {^\circ}{\sfT}^+\to{^\circ}{\sfT}^+\otimes{^\circ}{\sfT}^+
$$ 
is defined by the same formula for $\tau_{\frak{e}}^{\frak{n}}\in C$.

	\item[{\sf 2.}] The algebra morphism 
	$$
	{^\circ}{\sfD}^-: {^\circ}{\sfU}^-\to{^\circ}{\sfU}^-\otimes {^\circ}{\sfU}^-
	$$ 
	is defined by ${^\circ}\sfD^-{\bf1}_-={\bf1}_-\otimes{\bf1}_-$, and for $\tau_{\frak{e}}^{\frak{n}}\in SC$
\begin{equation*} \begin{split}
{^\circ}{\sfD}^-(\tau_{\frak{e}}^{\frak{n}})
&\defeq  \sum_{\varphi\in SF(\tau)}\sum_{\frak{n}_\varphi, \frak{e}'_{\partial\varphi}} \frac{1}{\frak{e}'_{\partial\varphi}!}{{\frak{n}}\choose{\frak{n}_\varphi}}
\varphi^{\frak{n}_\varphi+\pi\frak{e}'_{\partial\varphi}}_{\frak{e}} 
\otimes \big(\tau/^{\text{\rm red}}\varphi\big)^{[\frak{n}-\frak{n_\varphi}]_\varphi}_{\frak{e}+\frak{e}'_{\partial\varphi}},
\end{split}  \end{equation*}
where $SF(\tau)$ is the set of all {\sl subforests} $\varphi$ of $\tau$ {\sl which contain all red nodes of $\tau$}, and the sum over $\frak{n}_\varphi$ and $\frak{e}_{\partial\varphi}'$ is taken as in item $\sf 1$.

	\item[{\sf 3.}] The algebra morphism 
	$$
	{^\circ}\overline{\sfD}^{\,-} : {^\circ}{\sfT}^+\to{^\circ}{\sfU}^{\,-}\otimes {^\circ}{\sfT}^+
	$$ 
	is defined by the same formula as ${^\circ}{\sfD}^-$, but the first sum is restricted to the set $\overline{SF}(\tau)$ of all subforests $\varphi\in SF(\tau)$ {\sl which is disjoint with the root of $\tau$}.
\end{itemize}
\end{defn*}

\ssk

Our aim is to show the coassociativities of ${^\circ}\sfD^\pm$ and ${^\circ}{\sfD}$ and ${^\circ}\overline{\sfD}^{\,-}$. To avoid a confusing calculation, we separate the coproducts into graph part and decoration part. Define simpler coproducts acting on undecorated trees by
\begin{align*}
{^*}{\sf D}^+\tau \defeq  \sum_{\sigma\in ST(\tau)}\sigma\otimes\big(\tau/^{\text{blue}}\sigma\big),\quad
{^*}{\sf D}^-\tau \defeq  \sum_{\varphi\in SF(\tau)}\varphi\otimes\big(\tau/^{\text{red}}\varphi\big),\quad
{^*}\overline{\sf D}^{\,-}\tau \defeq  \sum_{\varphi\in\overline{SF}(\tau)}\varphi\otimes\big(\tau/^{\text{red}}\varphi\big).
\end{align*}
Given an undecorated tree $\tau$, denote by $\bbX_{(n,k)}$ the map {\color{red}adding}
to the node $n\in N_\tau$ the $\frak{n}$-decoration $k\in\bbN\times\bbN^d$, and denote by $\bbI_{(e,\ell)}$ the map giving to the edge $e\in E_\tau$ the $\frak{e}$-decoration $\ell\in\bbN\times\bbN^d$. Then any decorated tree $\tau_{\frak{e}}^{\frak{n}}$ is of the form
\begin{align}\label{decoration map on tau}
\bbF\tau =\bbF_1\cdots\bbF_N\tau,
\end{align}
where $\tau$ is an undecorated tree, and $\bbF_1,\dots,\bbF_N$ are family of $\bbX$-type or $\bbI$-type operators, applying to pairwise different nodes or edges. Moreover, we define the coproducts of such operators by
\begin{align*}
\bbD\bbX_{(n,k)} &= \sum_{k'\le k}\binom{k}{k'}\,\bbX_{(n,k')}\otimes\bbX_{(n,k-k')},\\
\bbD\bbI_{(e,\ell)} &= \bbI_{(e,\ell)}\otimes\bbX_{(e_-,0)} + \sum_{\ell'}\frac1{\ell'!}\,\bbX_{(e_-,\ell')}\otimes\bbI_{(e,\ell+\ell')},
\end{align*}
where $e_-$ denotes the node from where the edge $e$ leaves. For the products of pairwise disjoint such operators, define 
\begin{align}\label{skeletonmultiplicative}
\bbD\bbF\defeq (\bbD\bbF_1)\dots(\bbD\bbF_N).
\end{align}
There are two remarks about this identity. First, the right hand side does not depend on the order of $\bbF_1,\dots,\bbF_N$. To prove it, we have only to show that $\bbD\bbX_{(e_-,k)}\bbD\bbI_{(e,\ell)}=\bbD\bbI_{(e,\ell)}\bbD\bbX_{(e_-,k)}$ and $\bbD\bbI_{(e,k)}\bbD\bbI_{(f,\ell)}=\bbD\bbI_{(f,\ell)}\bbD\bbI_{(e,k)}$ for $e,f$ such that $e_-=f_-$, because for other pairs of operators $\bbF$ and $\bbG$, the operators $\bbD\bbF$ and $\bbD\bbG$ apply to different nodes or edges. These identities can be checked easily by definitions. Second, \eqref{skeletonmultiplicative} holds even if there are distinct $i$ and $j$ such that $\bbF_i$ and $\bbF_j$ apply to the same {\sl node}. Let $\bbF_i=\bbX_{(n,k)}$ and $\bbF_j=\bbX_{(n,\ell)}$. Since $\bbX_{(n,k)}\bbX_{(n,\ell)}=\bbX_{(n,k+\ell)}$, we have
\begin{align*}
\bbD(\bbX_{(n,k)}\bbX_{(n,\ell)})=\sum_{m\le k+\ell}\binom{k+\ell}{m}\bbX_{(n,m)}\otimes\bbX_{(n,k+\ell-m)}.
\end{align*}
On the other hand, we have
\begin{align*}
(\bbD\bbX_{(n,k)})(\bbD\bbX_{(n,\ell)})=\sum_{a\le k,\,b\le\ell}\binom{k}{a}\binom{\ell}{b}\bbX_{(n,a+b)}\otimes\bbX_{(n,k+\ell-a-b)}.
\end{align*}
We can see that the right hand sides coincide by using Chu-Vandermonde identity $\sum_{a+b=m}\binom{k}{a}\binom{\ell}{b}=\binom{k+\ell}{m}$.
Therefore we have $(\bbD\bbX_{(n,k)})(\bbD\bbX_{(n,\ell)})=\bbD\bbX_{(n,k+\ell)}=(\bbD\bbX_{(n,\ell)})(\bbD\bbX_{(n,k)})$.

\ssk

At this stage, we see that the coproducts ${^\circ}{\sf D}^{(\cdot,+,-)}$ 
apply to the decorated tree \eqref{decoration map on tau} by the forms
\begin{align}\label{D commutes with F}
{^\circ}{\sf D}^{(\cdot,+)}(\bbF\tau)=(\bbD\bbF)({^*}{\sf D}^+\tau),\quad
{^\circ}{\sf D}^{-}(\bbF\tau)=(\bbD\bbF)({^*}{\sf D}^{-}\tau).
\end{align}
In the right hand side of \eqref{D commutes with F}, be careful that $\bbD\bbF$ acts on subtrees and contracted trees. 
For an $\bbX$-type operator, if $n\notin N_\sigma$ then set $\bbX_{(n,k)}\sigma={\bf1}_{k=0}\sigma$. On a contracted tree $\tau/\sigma$, the $\bbX$-type operator acts of the form $\bbX_{([n],k)}$, where $[n]$ denotes the equivalence class in the contraction $\tau\to\tau/\sigma$. Hence
\begin{align*}
(\bbD\bbX_{(n,k)})(\sigma\otimes(\tau/\sigma))
=
\begin{cases}
\sum_{k'\le k}\binom{k}{k'}\bbX_{(n,k')}\sigma\otimes\bbX_{([n],k-k')}(\tau/\sigma),&n\in N_\sigma,\\
\sigma\otimes\bbX_{([n],k)}(\tau/\sigma),&n\notin N_\sigma.
\end{cases}
\end{align*}
For an $\bbI$-type operator, if $e\notin E_\sigma$ (resp. $e\notin E_{\tau/\sigma}$) then set $\bbI_{(e,\ell)}\sigma=0$ (resp. $\bbI_{(e,\ell)}(\tau/\sigma)=0$). Combing with the definition of $\bbX_{(n,k)}$, we have
\begin{align*}
(\bbD\bbI_{(e,\ell)})(\sigma\otimes(\tau/\sigma))
=
\begin{cases}
\bbI_{(e,\ell)}\sigma\otimes(\tau/\sigma), &\textrm{for } e\in E_\sigma,\\
\sum_{\ell'}\frac1{\ell'!}\bbX_{(e_-,\ell')}\sigma\otimes\bbI_{(e,\ell+\ell')}(\tau/\sigma), &\textrm{for } e\in\partial\sigma,\\
\sigma\otimes\bbI_{(e,\ell)}(\tau/\sigma), &\textrm{for } e\in E_\tau\setminus(E_\sigma\cup\partial\sigma).
\end{cases}
\end{align*}
These conventions show that the identities \eqref{D commutes with F} hold.


\subsubsection{Coassociativity}

\begin{lem}\label{proof of coassociativity of sfD}
One has the coassociativity formulas
$$
\begin{aligned}
\big({^\circ}\sfD\otimes\text{\rm Id}\big){^\circ}\sfD &= \big(\text{\rm Id}\otimes{^\circ}\sfD^+\big){^\circ}\sfD,   
&\big({^\circ}\sfD^+\otimes\text{\rm Id}\big){^\circ}\sfD^+ &= \big(\text{\rm Id}\otimes{^\circ}\sfD^+\big){^\circ}\sfD^+,   \\
\big({^\circ}\sfD^-\otimes\text{\rm Id}\big){^\circ}\sfD^- &= \big(\text{\rm Id}\otimes {^\circ}\sfD^-\big){^\circ}\sfD^-,   \quad
&\big({^\circ}\sfD^-\otimes\text{\rm Id}\big){^\circ}\overline{\sfD}^{\,-} 
&= \big(\text{\rm Id}\otimes {^\circ}\overline{\sfD}^{\,-}\big){^\circ}\overline{\sfD}^-.
\end{aligned}
$$
\end{lem}

\ssk

\begin{Dem}
We prove the identity
\begin{align} \label{EqLemmaCoassociativity}
({^\circ}\sfD\otimes\text{\rm Id}){^\circ}\sfD = (\text{\rm Id}\otimes{^\circ}\sfD^+){^\circ}\sfD;
\end{align}
the other identities are proved similarly. By the commutation relation \eqref{D commutes with F}, we have
\begin{align*}
\big(\iden\otimes{^\circ}{\sf D}^+\big){^\circ}{\sf D}\, \bbF\tau
= \big(\iden\otimes{^\circ}{\sf D}^+\big)(\bbD\bbF)({^*}{\sf D}^+\tau)
= \big((\iden\otimes\bbD)\bbD\bbF\big)\big(\iden\otimes{^*}{\sf D}^+\big){^*}{\sf D}^+\tau.
\end{align*}
Hence it is sufficient for proving \eqref{EqLemmaCoassociativity} to show 
the two identities
\begin{align}
\label{proof of coassociativity of sfD1}&\big({^*}{\sf D}^+\otimes\text{\rm Id}\big){^*}{\sf D}^+ = \big(\text{\rm Id}\otimes{^*}{\sf D}^+\big){^*}{\sf D}^+,   \\
\label{proof of coassociativity of sfD2}&(\iden\otimes\bbD)\bbD\bbF=(\bbD\otimes\iden)\bbD\bbF.
\end{align}
It is not difficult to show \eqref{proof of coassociativity of sfD1} by the definition of $ST(\tau)$, by noting that
$$
ST(\tau/\sigma) = \big\{\eta/\sigma\,;\,\sigma\subset\eta\subset\tau\big\}
$$
and $(\tau/\sigma)/(\eta/\sigma)=\tau/\eta$. Next we show \eqref{proof of coassociativity of sfD2}. By the multiplicativity \eqref{skeletonmultiplicative}, for any family $\bbF_1,\dots,\bbF_N$ of $\bbX$-type or $\bbI$-type operators, applying to pairwise different nodes or edges, one has
\begin{align*}
(\iden\otimes\bbD)\bbD\bbF
&=(\iden\otimes\bbD)\big((\bbD\bbF_1)\cdots(\bbD\bbF_N)\big)\\
&=\big((\iden\otimes\bbD)(\bbD\bbF_1)\big)\cdots\big((\iden\otimes\bbD)(\bbD\bbF_N)\big).
\end{align*}
In the last equality, we use the fact that \eqref{skeletonmultiplicative} holds even if multiple operators act on the same node.
If \eqref{proof of coassociativity of sfD2} holds for $\bbF_1,\dots,\bbF_N$, then
\begin{align*}
(\iden\otimes\bbD)\bbD\bbF
&=\big((\bbD\otimes\iden)(\bbD\bbF_1)\big)\cdots\big((\bbD\otimes\iden)(\bbD\bbF_N)\big)\\
&=(\bbD\otimes\iden)\big((\bbD\bbF_1)\cdots(\bbD\bbF_N)\big)=(\bbD\otimes\iden)\bbD\bbF.
\end{align*}
Therefore it is sufficient to show \eqref{proof of coassociativity of sfD2} for $\bbF=\bbX_{(n,k)}$ and $\bbI_{(e,\ell)}$. For $\bbF=\bbX_{(n,k)}$, we have
\begin{align*}
(\iden\otimes\bbD)\bbD\bbX_{(n,k)}&=\sum_{k'\le k}\binom{k}{k'}\bbX_{(n,k')}\otimes\bbD\bbX_{(n,k-k')}\\
&=\sum_{k'\le k,\, k''\le k-k'}\binom{k}{k'}\binom{k-k'}{k''}\bbX_{(n,k')}\otimes\bbX_{(n,k'')}\otimes\bbX_{(n,k-k'-k'')}\\
&=\sum_{k',k'';k'+k''\le n}\binom{k}{k',k''}\bbX_{(n,k')}\otimes\bbX_{(n,k'')}\otimes\bbX_{(n,k-k'-k'')}.
\end{align*}
We have the same expansion from $(\bbD\otimes\iden)\bbD\bbX_{(n,k)}$, so \eqref{proof of coassociativity of sfD2} holds for $\bbX$-type operators.
We can prove the same result for $\bbI$-type operators by similar computations.
\end{Dem}

\medskip

Now we consider the extended decoration.


\medskip

\begin{lem}\label{proof of coassociativity bsD}
One has the coassociativity formulas
$$
\begin{aligned}
(\sfD\otimes\text{\rm Id})\sfD &= \big(\text{\rm Id}\otimes\sfD^+\big)\sfD,   
&\big(\sfD^+\otimes\text{\rm Id}\big)\sfD^+ &= \big(\text{\rm Id}\otimes\sfD^+\big)\sfD^+,   \\
(\sfD^-\otimes\text{\rm Id})\sfD^- &= \big(\text{\rm Id}\otimes\sfD^-\big)\sfD^-,   \quad
&\big(\sfD^-\otimes\text{\rm Id}\big)\overline{\sfD}^{\,-} &= \big(\text{\rm Id}\otimes\overline{\sfD}^{\,-}\big)\overline{\sfD}^{\,-}.
\end{aligned}
$$
\end{lem}

\medskip

\begin{Dem}
We consider the first and third identities; the other identities are proved similarly. In this proof, denote by $\bar{\tau}=\tau_{\frak{e}}^{\frak{n}}$ a generic decorated tree without $\frak{o}$-decoration, and write $\bar\tau^{\frak{o}}$ for $\tau_{\frak{e}}^{\frak{n},\frak{o}}$. As in Section \ref{SectionConcreteRS} 
and Section \ref{SectionConcreteRenorS}, we use a shorthand notation
$$
{^\circ}{\sfD}\bar\tau = \sum_{\bar\sigma\le\bar\tau}\bar\sigma\otimes(\bar\tau/\bar\sigma),\quad
{^\circ}{\sfD}^-\bar\tau = \sum_{\bar\varphi\trianglelefteq\bar\tau}\bar\varphi\otimes(\bar\tau/\bar\varphi).
$$
Then we can write
$$
{\sfD}\bar\tau^{\frak{o}} = \sum_{\bar\sigma}\bar\sigma^{\frak{o}}\otimes(\bar\tau/\bar\sigma)^{\frak{o}\vert_{\tau\setminus\sigma}},\quad
{\sfD}^-\bar\tau^{\frak{o}} = \sum_{\bar\varphi}\bar\varphi^{\frak{o}}\otimes(\bar\tau/\bar\varphi)^{\frak{o}+\frak{o}(\bar\varphi)}.
$$
Recall that $\frak{o}(\bar\varphi):N_{\tau/\varphi}\to\bbZ[\beta_0]$, is a function giving the value $\vert\bar\tau_j\vert$, where $\bar\tau_j$ is 
a connected component if $\varphi$, to the node $[\tau_j]\in N_{\tau/\varphi}$. We obtain the coassociativity of $\sfD$ from the coassociativity of ${^\circ}\sfD$, noting that
$$
(\bar\tau/\bar\sigma)/(\bar\eta/\bar\sigma)=\bar\tau/\bar\eta,\quad
\frak{o}\vert_{(\tau\setminus\sigma)\setminus(\eta\setminus\sigma)}=\frak{o}\vert_{\tau\setminus\eta}
$$
for any $\bar\eta\le\bar\sigma\le\bar\tau$. To prove the coassociativity of $\sfD^-$, noting that
\begin{align*}
(\sfD^-\otimes\iden)\sfD^-\bar\tau^{\frak{o}}
= \sum_{\bar\psi\trianglelefteq\bar\varphi\trianglelefteq\bar\tau}
\bar\psi^{\frak{o}}\otimes(\bar\varphi/\bar\psi)^{\frak{o}+\frak{o}(\bar\psi)}
\otimes(\bar\tau/\bar\varphi)^{\frak{o}+\frak{o}(\bar\varphi)}
\end{align*}
and
\begin{align*}
(\iden\otimes\sfD^-)\sfD^-\bar\tau^{\frak{o}}
&= \sum_{\bar\psi\trianglelefteq\bar\tau} \bar\psi^{\frak{o}}\otimes\sfD^-(\bar\tau/\bar\psi)^{\frak{o}+\frak{o}(\bar\psi)}\\
&= \sum_{\bar\psi\trianglelefteq\bar\varphi\trianglelefteq\bar\tau} \bar\psi^{\frak{o}}\otimes(\bar\varphi/\bar\psi)^{\frak{o}+\frak{o}(\bar\psi)} \otimes (\bar\tau/\bar\varphi)^{\frak{o}+\frak{o}(\bar\psi)+\frak{o}(\bar\varphi/\bar\psi)},
\end{align*}
it is sufficient to show that $\frak{o}(\bar\varphi)=\frak{o}(\bar\psi)+\frak{o}(\bar\varphi/\bar\psi)$ as a function on $N_{\tau/\varphi}$. This holds true because $\vert\bar\varphi\vert' = \vert\bar\psi\vert' + \vert\bar\varphi/\bar\psi\vert'$.
\end{Dem}


\subsubsection{Co-interaction}
\label{SubsectionAppendixCointeraction}

\begin{lem}
One has the co-interaction formulas
\begin{equation}\label{proof of cointercation} \begin{split}
\mcM^{(13)}\big({^{\circ}}\sfD^-\otimes{^{\circ}}\overline{\sf D}^{\,-}\big){^{\circ}}\sfD &= (\iden\otimes {^{\circ}}\sfD){^{\circ}}\sfD^-,   \\
\mcM^{(13)}\big({^{\circ}}\overline{\sf D}^{\,-}\otimes{^{\circ}}\overline{\sf D}^{\,-}\big){^{\circ}}\sfD^+ &= \big(\iden\otimes {^{\circ}}\sfD^+\big){^{\circ}}\overline{\sf D}^{\,-},
\end{split} \end{equation}
and
\begin{equation} \label{proof of cointercation bsD} \begin{split}
\mcM^{(13)}\big(\sfD^-\otimes\overline{\sfD}^{\,-}\big)\sfD &= (\iden\otimes \sfD)\sfD^-,   \\
\mcM^{(13)}\big(\overline{\sfD}^{\,-}\otimes\overline{\sfD}^{\,-}\big)\sfD^+ &= (\iden\otimes \sfD^+)\overline{\sfD}^{\,-}.
\end{split} \end{equation}
\end{lem}

\medskip

\begin{Dem}
Consider the first identity of \eqref{proof of cointercation} and the first identity of \eqref{proof of cointercation bsD}; the two other identities are proved similarly. By the commutation relations \eqref{D commutes with F}, identity \eqref{proof of cointercation} rewrites
\begin{align}\label{proof of cointercation 1}
\mcM^{(13)}\big((\bbD\otimes\bbD)\bbD\bbF\big)
\big({^*}{\sf D}^-\otimes{^*}\overline{\sf D}^{\,-}\big){^*}{\sf D}^+\tau
= \big((\text{\rm Id}\otimes\bbD)\bbD\bbF\big)(\text{\rm Id}\otimes{^*}{\sf D}^+){^*}{\sf D}^-.
\end{align}
By the multiplicativity \eqref{skeletonmultiplicative} of $\bbD$, it is sufficient to show \eqref{proof of cointercation 1} for the operators $\bbF=\bbX_{(n,k)}$ and $\bbI_{(e,\ell)}$. By definition,
\begin{align*}
\big({^*}{\sf D}^-\otimes{^*}\overline{\sf D}^{\,-}\big){^*}{\sf D}^+\tau 
= \sum_{\sigma\in ST(\tau)}\sum_{\varphi\in SF(\sigma),\psi\in\overline{SF}(\tau/\sigma)}
\varphi\otimes(\sigma/^{\text{red}}\varphi)\otimes\psi\otimes(\tau/^{\text{blue}}\sigma)/^{\text{red}}\psi.
\end{align*}
Note that $\varphi$ and $\psi$ are disjoint subforests of $\tau$ because of the definition of $\overline{SF}$. Thus we have
\begin{align*}
&\mcM^{(13)}\big((\bbD\otimes\bbD)\bbD\bbX_{(n,k)}\big)
(\varphi\otimes(\sigma/^{\text{red}}\varphi)\otimes\psi\otimes(\tau/^{\text{blue}}\sigma)/^{\text{red}}\psi)\\
&=\mcM^{(13)}\sum_{k=a+b+c+d}\frac{k!}{a!b!c!d!} \, \bbX_{(n,a)}\varphi\otimes\bbX_{(n,b)}(\sigma/^{\text{red}}\varphi)
\otimes\bbX_{(n,c)}\psi\otimes\bbX_{(n,d)}\big((\tau/^{\text{blue}}\sigma)/^{\text{red}}\psi\big)   \\
\end{align*}

\begin{align*}
&=\sum_{k=a+b+d}\frac{k!}{a!b!d!} \, \bbX_{(n,a)}(\varphi\psi)
\otimes\bbX_{(n,b)}(\sigma/^{\text{red}}\varphi)\otimes\bbX_{(n,d)}\big((\tau/^{\text{blue}}\sigma)/^{\text{red}}\psi\big)   \\
&=\big((\text{\rm Id}\otimes\bbD)\bbD\bbX_{(n,k)}\big)
(\varphi\psi\otimes(\sigma/^{\text{red}}\varphi)\otimes\big(\tau/^{\text{blue}}\sigma)/^{\text{red}}\psi\big),
\end{align*}
since either of $a$ and $c$ has to be $0$ in the second line.
It is not difficult to show a similar equality for $\bbF=\bbI_{(e,\ell)}$.
Hence we have
\begin{align*}
\mcM^{(13)}\big((\bbD\otimes\bbD)\bbD\bbF\big)
\big({^*}{\sf D}^-\otimes{^*}\overline{\sf D}^{\,-}\big){^*}{\sf D}^+\tau
= \big((\text{\rm Id}\otimes\bbD)\bbD\bbF\big) \mcM^{(13)}\big({^*}{\sf 
D}^-\otimes{^*}\overline{\sf D}^{\,-}\big){^*}{\sf D}^+\tau.
\end{align*}
Since it is not difficult to show the co-interaction formula
$$
\mcM^{(13)}\big({^*}{\sf D}^-\otimes{^*}\overline{\sf D}^{\,-}\big){^*}{\sf D}^+\tau
= \big(\iden\otimes{^*}{\sf D}^+\big){^*}{\sf D}^-,
$$
identity \eqref{proof of cointercation 1} follows as a consequence.

Next we consider \eqref{proof of cointercation bsD}. By definition,
\begin{align*}
\mcM^{(13)}\big(\sfD^-\otimes\overline{\sfD}^{\,-}\big)\sfD\bar\tau^{\frak{o}}
= \sum_{\bar\sigma\le\bar\tau, \bar\varphi\trianglelefteq\bar\sigma, \bar\psi\trianglelefteq\bar\tau/\bar\sigma}
\bar\varphi^{\frak{o}}\bar\psi^{\frak{o}\vert_{\tau\setminus\sigma}}
\otimes(\bar\sigma/\bar\varphi)^{\frak{o}+\frak{o}(\bar\varphi)}
\otimes\big((\bar\tau/\bar\sigma)/\bar\psi\big)^{\frak{o}\vert_{\tau\setminus\sigma}+\frak{o}(\bar\psi)}
\end{align*}
and
\begin{align*}
(\iden\otimes\sfD)\sfD^-\bar\tau^{\frak{o}}
=\sum_{\bar\zeta\trianglelefteq\bar\tau, \bar\eta\le\bar\tau/\bar\zeta}
\bar\zeta^{\frak{o}}\otimes\bar\eta^{\frak{o}+\frak{o}(\bar\zeta)}\otimes\big((\bar\tau/\bar\zeta)/\bar\eta\big)^{(\frak{o}+\frak{o}(\bar\zeta))\vert_{(\tau/\zeta)\setminus\eta}}.
\end{align*}
The cointeraction between ${^\circ}\sfD$ and ${^\circ}\sfD^-$ implies that the {\it change of variables}
$$
\bar\zeta\leftrightarrow\bar\varphi\,\bar\psi,\quad
\bar\eta\leftrightarrow\bar\sigma/\bar\varphi
$$
is possible. Since $\sigma$ and $\psi$ are disjoint,
\begin{align*}
&(\bar\varphi\bar\psi)^{\frak{o}} = \bar\varphi^{\frak{o}}\,\bar\psi^{\frak{o}\vert_{\tau\setminus\sigma}},\quad
(\bar\sigma/\bar\varphi)^{\frak{o}+\frak{o}(\bar\varphi\bar\psi)}
= (\bar\sigma/\bar\varphi)^{\frak{o}+\frak{o}(\bar\varphi)},   \\
& \big((\bar\tau/\bar\sigma)/\bar\psi\big)^{(\frak{o}+\frak{o}(\bar\varphi\bar\psi))\vert_{(\tau/(\varphi\psi))\setminus(\sigma/\varphi)}}
= \big((\bar\tau/\bar\sigma)/\bar\psi\big)^{\frak{o}\vert_{\tau\setminus\sigma}+\frak{o}(\bar\psi)}.
\end{align*}
Thus \eqref{proof of cointercation bsD} follows.
\end{Dem}


\section{Comments}
\label{SectionAppendixComments}

\textbf{\textsf{Section \ref{SectionIntro} -- }} Regularity structures theory has its roots in T. Lyons' theory of rough paths and rough differential equations \cite{Lyons98}. This theory deals with controlled ordinary differential equations
$$
dx_t = V(x_t)dh_t,
$$
with controls $h$ of low regularity, say $\alpha$-H\"older. For $\alpha>1/2$, Young integration theory allows to make sense of the equation as a fixed point problem for an integral equation. As one expects a solution path to be $\alpha$-H\"older, the product $V(x_t)dh_t$ makes sense as a distribution on $\bbR_+$ iff $\alpha+(\alpha-1)>0$, that is $\alpha>1/2$. One of Lyons' deep insights was to realize that what really governs the dynamics is not the $\bbR^\ell$-valued control $h$, say, but rather a finite 
collection of its iterated integrals. The latter are ill-defined when $\alpha\leq 1/2$, and a rough path is the a priori datum of quantities playing their role. Natural algebraic and size constraints on these objects are then sufficient to set the entire theory. These constraints are similar 
to the constraints that define the $\sf g$-part of a model. Several reformulations of rough paths theory were given after Lyons' seminal work: Davie's numerical scheme approach \cite{Davie}, Gubinelli's controlled paths 
approach \cite{GubinelliControlled, GubinelliBranched}, Friz \& Victoir's 
limit ODE picture \cite{FrizVictoir}, and Bailleul's approximate flow-to-flow approach \cite{BailleulRMI}, amongst others. Gubinelli's versatile notion of controlled paths was a direct source of inspiration for the construction of regularity structures.   

\medskip

Other tools than regularity structures have been developed for the study of singular stochastic PDEs. None of them offers presently a complete alternative to regularity structures. 

\ssk

{\begin{itemize}
   \item Gubinelli, Imkeller and Perkowski laid in \cite{GIP} the foundations of paracontrolled calculus, that was developed by Bailleul \& Bernicot \cite{BB1, BB2, BB3}. While the fundamental notions of regularity structures involve pointwise expansions, paracontrolled calculus uses paraproducts as a mean for making sense of what it means to look like a reference quantity. See \cite{LNPCGP} for lecture notes on the subject and \cite{GPRoeckner} for an overview on the subject, both by Gubinelli \& Perkowski. 
   
   In a nutshell, the starting point of the paracontrolled approach to the study of singular stochastic PDEs is the decomposition of a product of two distributions $f,g$ into
   $$
   fg = \big({\sf P}_fg+{\sf P}_gf\big) + {\sf \Pi}(f,g).
   $$
   This decomposition is obtained in a Fourier picture of the product by splitting the convolution into what happens far from the diagonal $\big({\sf P}_fg+{\sf P}_gf\big)$ from what happens near the diagonal ${\sf \Pi}(f,g)$. This decomposition isolates in the resonant term ${\sf \Pi}(f,g)$ 
what does not make sense in a general product, the paraproduct terms ${\sf P}_fg, {\sf P}_gf$ being always well-defined. The definition of ${\sf P}_fg$ allows to think of it as a modulation of $g$ by $f$ and give meaning to what it means for a distribution/function $u$ to look like another distribution/function $g$
   $$
   u= {\sf P}_{u'}g + u^\sharp,
   $$
   for a function $u'$ and a distribution/function $u^\sharp$ that is more 
regular than $g$. The role of modelled distributions is played in a paracontrolled setting by systems $(u_a)_{a\in\mathscr{A}}$ of paracontrolled distributions/functions
   \begin{equation} \label{EqPCSystemForm}
   u_a = \sum_{\tau\in\mathscr{T};\vert a\tau\vert\leq n\alpha}{\sf P}_{u_{a\tau}}[\tau] + u_a^\sharp
   \end{equation}
   indexed by the set $\mathscr{A}$ of words over an alphabet $\mathscr{T}$, with remainders $u_a^\sharp$ sufficiently regular. The reference distributions/functions $[\tau]$ somehow play the role of ${\sf \Pi}\tau$ and 
the $(u_a)_{a\in\mathscr{A}}$ the role of the $(u_\tau)_{\tau\in T,\vert\tau\vert<\gamma}$. (Bailleul \& Hoshino's work \cite{BH1, BH2} on the relations between paracontrolled calculus and regularity structures make that link clear.) Identity \eqref{EqPCSystemForm} is an analogue of the notion of modelled distribution. While the definition of the latter involves pointwise comparisons, here the comparison is somehow done in `momentum space', although not in a pointwise sense. The core point of the paracontrolled analysis of a (system of) singular stochastic PDE(s) is that we end up dealing with ill-defined terms of the form
   \begin{equation} \label{EqIllDefinedTermPC}
   {\sf \Pi}'\big(u,\zeta^{(i)}\big)
   \end{equation}
   for operators ${\sf \Pi}'$ that have similar properties as the resonant operator $\sf \Pi$, and possibly multi-dimensional functionals $\zeta^{(i)}$ of the noise $\zeta$. It turns out that while an expression like \eqref{EqIllDefinedTermPC} does not make sense for a generic $u$, it makes sense on a restricted class of $u$ of the form \eqref{EqPCSystemForm} {\it provided} one can make sense of the terms ${\sf\Pi}'\big(\tau,\zeta^{(i)}\big)$. The analysis of a given (system of) singular stochastic PDE(s) gives an inductive definition of the $\zeta^{(i)}$ and the $\tau$'s. Compared to the regularity structures setting, the datum of all the ${\sf\Pi}'\big(\tau,\zeta^{(i)}\big)$ plays the role of the datum of a model. The inductive/tree structure of the elements of a regularity structure takes here the form of the inductive definition of the $\tau$'s and $\zeta^{(i)}$'s. A systematic treatment of renormalization operations within paracontrolled calculus has not been invented yet. The links between the regularity structure and paracontrolled settings detailed in Bailleul \& Hoshino's works \cite{BH1, BH2} allow however to transport the renormalization machinery of regularity structures into the setting paracontrolled calculus. What is missing presently is an independent, purely paracontrolled, approach of the renormalization problem.      \vspace{0.1cm}
   
   \item Otto \& Weber \cite{OttoWeber} developed jointly with Sauer and Smith \cite{OSSW} a variant of regularity structures that is more in the flavour of rough paths theory. See in also their most recent joint works \cite{LOT,LOTT} with Linares, Tempelmayr and Tsatsoulis to see how far they were able to go. Most concepts and objects from regularity structures have counterparts in their setting. It was specifically designed and used for the analysis of a number of quasilinear singular stochastic PDEs. Some of these equations can be approached using the original first order paracontrolled calculus as in Bailleul, Debussche and Hofmanov\'a's work \cite{BDH} or a variant of it using paracomposition operators, as in Furlan \& Gubinelli's work \cite{FG}. See also \cite{BailleulMouzard, GerencserHairer} for extensions of paracontrolled calculus and regularity structures designed for the study of a whole class of quasilinear singular stochastic PDEs.\vspace{0.1cm}
   
   \item Kupiainen \& Marcozzi managed in \cite{Kupiainen1, Kupiainen2} to implement a renormalization group approach to the (KPZ) and $\Phi^4_3$ equations. The starting point of their strategy consists in decomposing the resolution operator ${\bf K}^{-1}$ involved in the Picard formulation of the equation as a sum of operators ${\bf K}^{-1}_n$ turning distributions into smooth functions that vary essentially only up to scale $2^{-n}$. The approximate renormalized dynamics will take the form
   $$
   u_N = \bigg(\sum_{n=0}^N {\bf K}^{-1}_n\bigg)\Big(F_N(u_N,\partial 
u_N\,;\,\zeta)\Big)
   $$
   in a simplified problem where the initial condition was taken to be null. The noise is left untouched, with no problem for defining the nonlinearity in the right hand side since $u_N$ is smooth. So one has
   $$
   u_N = u_N^0 + \cdots + u_N^N,
   $$
   where each term $u_N^n$ is morally varying only up to scale $2^{-n}$. The point is now to see that one can choose the nonlinearity $F_N$ in such a way that each $u_N^n$ is the solution of an equation of the form
   $$
   u_N^n = \Big(\sum_{m=n}^N {\bf K}^{-1}_n\Big)\Big(F_N^n(u_N^n,\partial u_N^n\,;\,\zeta)\Big)
   $$
   for a nonlinearity $F_N^n$, and is converging in a proper space as $N$ goes to infinity. This is done via the use of rescaling operators, taking profit from the exact scaling property of the heat kernel, by turning the problem of convergence of each $u_N^n$ into the problem of the convergence of the family of rescaled versions of the functions $F_N^n$ -- taking profit of the fact that the former is a continuous function of the later. The overall convergence of $u_N$ as $N$ goes to $\infty$ is somehow similar to the well-known fact that a sum of functions $\sum_n u^n$ converges in an $\alpha$-H\"older space if $u^n$ is localized in Fourier space on a 
ball of size $2^n$ and has uniform norm of order $2^{-\alpha n}$.

This approach was improved a lot by some recent works of P. Duch \cite{Duch1, Duch2, Duch3}, who traded the above discrete scale decomposition for a continuous scale decomposition and uncovered a certain structure on the cumulants of some functionals of the noise that pave the way to the development of a robust approach to some classes of singular equations. One can look at Chandra \& Ferdinand's work \cite{CF24} for an application of this approach to the generalized (KPZ) equation.   
\end{itemize} }

\medskip

\textbf{\textsf{Section \ref{SectionBasicsRS} -- }} The functional setting adopted here draws inspiration from \cite{BB1,BB2, BB3} and \cite{OttoWeber}. The main result of this section is the reconstruction theorem.

\ssk

Several proofs of the reconstruction theorem are available now, in addition to Hairer's original proof. Gubinelli, Imkeller and Perkowski gave in \cite{GIP} an alternative construction of the reconstruction map using a paraproduct-like operator. Singh \& Teichmann showed in \cite{SinghTeichmann} how it can be understood as the continuous extension of an elementary reconstruction operator defined on a set of smooth modelled distributions. Otto and Weber have an analogue of the reconstruction map in their rough paths-like setting \cite{OttoWeber, OSSW}. Caravenna \& Zambotti's recent work \cite{CaravennaZambotti} provide a robust version of the reconstruction theorem in a setting free of any reference to regularity structures. The notion of coherent germ turns it into a particularly versatile tool. See \cite{HairerLabbe, HenselRosati, LiuPromelTeichmann} for versions of the reconstruction theorem in functional settings different from H\"older spaces. See also our previous work \cite{BH1} for a paracontrolled representation of the reconstruction operator that refines over a similar flavoured representation given in Theorem 6.10 of \cite{GIP}.

\ssk

The reconstruction theorem takes its place in the history of a family of statements producing `transcendantal' objects, i.e. objects constructed by limiting procedures, from families of objects satisfying constraints involving no limiting procedures. The one-step Euler scheme for solving ordinary differential equations characterizes for instance uniquely their flows under sufficient regularity conditions on the vector fields. In its simplest form, for the equation $\dot y=-y$, in $\bbR$, it yields the elementary identity $(1-t/n)^n\rightarrow e^{-t}$, as $n$ goes to $\infty$. 
It takes a more elaborate form in Hille's approximation $\big((\textrm{Id}+\frac{t}{n}\,G)^{-1}\big)^n$ of the semigroup $(e^{-tG})_{t>0}$ generated by an unbounded operator $G$ under well-known conditions. Chernov's theorem \cite{Chernov} on families of strongly continuous perturbations of the identity used for constructing $(e^{-tG})_{t>0}$ has a similar flavour. So is the $C^1$-approximate flow-to-flow machinery of \cite{BailleulRMI}, that provides a far reaching generalization of Lyons' extension theorem in rough paths theory and Gubinelli and Feyel \& de la Pradelle' sewing lemma \cite{GubinelliControlled, FdlP, FdlPM}. All these statements characterize uniquely a transcendantal object as the unique object close to a family of objects satisfying a `$O(1)$ condition', involving no limiting procedure. The characterizing identity \eqref{EqReconstructionCondition} for the reconstruction is of that form when the reconstruction operator 
is unique. This kind of situation allows to build a calculus for the transcendantal objects from an elementary calculus on their generators.   

\ssk

The space of models over a given regularity structure is nonlinear. Bailleul \& Hoshino showed in \cite{BH1, BH2} how to parametrize this space by 
a linear space using the tools of paracontrolled calculus. The set of $\bf K$-admissible models on a given regularity structure turns out in particular to be parametrized by the data for each $\tau$ of negative homogeneity of a $\vert\tau\vert$-H\"older distribution, describing somehow the most regular part of the distribution ${\sf \Pi}\tau$. This has a number of consequences, such as an extension theorem similar to Lyons' extension theorem in rough paths theory.   

\ssk

Proposition \ref{PropFModelled}, giving a definition of the image of a modelled distribution by a nonlinear map, has a counterpart in paracontrolled calculus, generalizing Bony's paralinearisation formula to an arbitrary order -- see Section 2 of Bailleul \& Bernicot's work \cite{BB3}.   

\medskip

\textbf{\textsf{Section \ref{SectionIntegration} -- }} The proof of the continuity result for $\mcK^{\sf M}$ is an adaptation of the material from 
Hairer's groundbreaking work \cite{Hai14} to the functional setting adopted here. It is called by Hairer the multilevel Schauder estimates.  The construction of $\bf K$-admissible models from Section \ref{SubsectionAdmissibleModels} is adapted from Bailleul \& Hoshino's work \cite{BH1}, which gives amongst others a parametrization of the set of all admissible models on any reasonable concrete regularity structure. See \cite{BH2} for more results on the structure of the space of models and modelled distributions. Note that in the different components $u_\tau$ of a modelled distribution also appear in the paracontrolled approach, in which they are involved in the global description of a possible solution, as opposed to their local meaning in the regularity structure setting. 

\medskip

\textbf{\textsf{Section \ref{SectionSolvingPDEs} -- }} This section essentially follows the line of the corresponding results in \cite{Hai14}, Section 7 therein. Note that the setting presented here does not allow to take as initial condition a Dirac mass for instance. One needs for that purpose to set a Besov counterpart of the theory, as opposed to the H\"older flavoured version presented here. See Hairer \& Labb\'e's work \cite{HairerLabbe}, Hensel \& Rosati's work \cite{HenselRosati} or Singh \& Teichmann's work \cite{SinghTeichmann}. 

In a different direction, an number of works have been done on quasilinear singular stochastic PDEs \cite{BDH, BailleulMouzard, FG, GerencserHairer, LOT, LOTT, OttoWeber, OSSW, OSSW2}.

\medskip

\textbf{\textsf{Section \ref{SectionConcreteRenormStructure} -- }} The notion of renormalization structures and compatible regularity and renormalization structures introduced in this section is new. It encodes in a simple way the mechanics at work in Bruned, Hairer and Zambotti's work \cite{BHZ}.   

\medskip

\textbf{\textsf{Section \ref{SectionMultiAndrenormalizedEquations} -- }} This section contains the core insights of Bruned, Chandra, Chevyrev and Hairer's work \cite{BCCH18}, implemented here on the example of the generalized (KPZ) equation. The relevance of the notion of pre-Lie algebra was 
first noticed in the work \cite{BrunedChevyrevFriz} of Bruned, Chevyrev and Friz on rough paths. The article \cite{Cayley} would provide a pre-history of the pre-Lie algebra. The comodule-bialgebra structure of the Butcher-Connes-Kreimer Hopf algebra was first investigated in the work \cite{CalaqueEbrahimiFardManchon} of Calaque, Ebrahimi-Fard and Manchon; it played a key motivating role in the work of Bruned, Chevyrev and Friz. The description of the free pre-Lie algebra in this setting is due to Chapoton and Livernet \cite{ChapotonLivernet}. The notion of multi-pre-Lie algebra was introduced in the work \cite{BCCH18} of Bruned, Chandra, Chevyrev and 
Hairer, where the free multi-pre-Lie was first described.

The $\frak{o}$-decorations introduced here under the form of ${\color{red} \bullet}^{n,\alpha}$ is forced by our construction of compatible regularity and renormalization structures for the generalized (KPZ) equation, given in Section \ref{SectionBuildingRS}. It has no dynamical meaning. Bailleul \& Bruned showed in \cite{BailleulBruned} how to obtain the renormalized equation without using extended decorations for a large class of renormalization procedures including the BPHZ scheme. This work is based on the recursive renormalization scheme introduced by Bruned \cite{BrunedRecursive}.

\ssk

The setting described here is robust enough to deal with equations driven by multiple noises, or systems of equations driven by multiple noises. We 
take Funaki's example \cite{Funaki} of the random motion of a rubber on a manifold as an archetype -- see also \cite{HairerString, BrunedGabrielHairerZambotti}. The unknown $u$ is a spacetime function with values in $\bbR^d$, solution of the system
$$
(\partial_t-\partial_x^2) u = \Gamma(u)\big(\partial_xu,\partial_xu\big) + \Sigma(u)\xi,
$$
where $\Gamma(z)$ is a symmetric matrix on $\bbR^d$, and $\Sigma(z)$ a linear map from $\bbR^k$ to $\bbR^d$, for any $z\in\bbR^d$, and $\xi = (\xi^1,\dots, \xi^k)$ is an $k$-dimensional tuple of identically distributed independent one-dimensional spacetime white noises. We still have only one operator $(\partial_t-\partial_x^2)$ in this example, so the edge type set is here the same as in the study of the generalised (KPZ) equation. 
The node set is changed from $\{\circ,\bullet\}\times\bbN^{d+1}$ to $\{\circ^1,\dots, \circ^k, \bullet\}\times\bbN^{d+1}$ to account for the fact that we have $k$ noises $\xi^1,\dots, \xi^k$ in the system. Things get a bit messier if the system involves different operators, with different regularising properties, and noises with different regularities. The fundamental ideas involved in the analysis remain the same, while the notations 
needed to take care of this richer setting become heavier. All this is explained in full details in \cite{BHZ}. 

\medskip

\textbf{\textsf{Section \ref{SectionBHZCharacter} -- }} This section gives what seems to us to be one of the two core results of \cite{BHZ}, Theorem \ref{ThmPropertyBPHZRenorm} here. More general renormalization schemes were introduced by Bruned in \cite{BrunedRecursive} and Bailleul \& Bruned have shown in \cite{BailleulBruned} how to get back the renormalized equation in a very simple way in a setting with no extended decorations.

\medskip

\textbf{\textsf{Section \ref{SectionManifold} -- }} The fact that the family of solutions to the generalized (KPZ) equation forms a finite dimensional manifold of some function space had not been noticed so far. More generally, this is true for any singular stochastic PDE that can be treated by the methods of regularity structures This short section emphasizes that fact.

\medskip

\textbf{\textsf{Section \ref{SectionBuildingRS} -- }} This section builds 
on the fundamental work \cite{BHZ}, with a number of simplifications. The 
notion of subcritical equation is subtle to check in the general case of a system of equations, as one needs to keep track of how a given symbol of a regularity structure `flows' in the different pieces of a system, involving possibly operators with different regularizing properties. The meaning of subcriticality remains, though.

\medskip

{\small 
\noindent \textcolor{gray}{$\bullet$} {\sf I. Bailleul} -- Univ. Brest, LMBA - UMR 6205, Brest, France.   \\
\noindent {\it E-mail}: ismael.bailleul@univ-brest.fr    

\medskip

\noindent \textcolor{gray}{$\bullet$} {\sf M. Hoshino} --  Graduate School of Engineering Science, Osaka University, Japan   \\
{\it E-mail}: hoshino@sigmath.es.osaka-u.ac.jp   }

\end{document}